%% file: phDThesis.tex
\documentclass[12pt, oneside]{book}
\usepackage{setspace}

\usepackage[utf8]{inputenc} 
\usepackage{lmodern} 
\usepackage[T1]{fontenc} 
\usepackage[spanish,english]{babel} 
\usepackage{textgreek} 
\usepackage{amsmath}
\usepackage{amsthm, amssymb, amsfonts}
\usepackage{graphicx, color}
\usepackage{caption, subcaption} 
\usepackage{multicol} 
\usepackage{makeidx} 
\usepackage{stmaryrd}

\makeindex
\usepackage{array, tabularx, multirow} 
\usepackage[width=155mm,top=25mm,bottom=25mm, headheight=15pt,
]{geometry}
\usepackage[pagebackref=true]{hyperref} 
\usepackage{enumitem} 
\hypersetup{
	hidelinks,
  colorlinks  = true,    
  urlcolor    = red,    
  linkcolor   = blue,    
  citecolor   = blue,      
  pdfauthor   = {Jorge Antonio Cruz Chapital},%
  pdfsubject  = {Ph. D. Thesis},%
  pdftitle    = {Between omega and the hill of treasures},  %
,}

\usepackage{tikz-cd} 

\usepackage{lipsum}  

\usepackage[capitalise,noabbrev, nameinlink]{cleveref} 

\Crefname{part}{Part}{Parts}
\Crefname{step}{Step}{Steps}
\Crefname{prop}{Proposition}{Propositions}
\Crefname{prob}{Problem}{Problems}

\usepackage{url}

\usepackage{fancyhdr} 

\usepackage{titling} 
\usepackage{datetime}

\graphicspath{ {img/} } 

\makeatletter
  \def\cleardoublepage{\clearpage\if@twoside \ifodd\c@page\else
  \hbox{}
  \vspace*{\fill}
  \begin{center}
  \end{center}
  \vspace{\fill}
  \thispagestyle{empty}
  \newpage
  \if@twocolumn\hbox{}\newpage\fi\fi\fi}
\makeatother

\numberwithin{equation}{section}

\theoremstyle{plain}
\newtheorem{theorem}[equation]{Theorem}
\newtheorem{lemma}[equation]{Lemma}
\newtheorem{proposition}[equation]{Proposition}
\newtheorem{corollary}[equation]{Corollary}

\theoremstyle{definition}
\newtheorem{definition}[equation]{Definition}

\newtheorem{problem}[equation]{Problem}

\theoremstyle{remark}

\newtheorem{rem}[equation]{Remark}

\usepackage[bb=boondox]{mathalfa}

\renewcommand{\leq}{\leqslant}
\renewcommand{\geq}{\geqslant}

\renewcommand{\epsilon}{\varepsilon}

\usepackage[shadow, colorinlistoftodos,disable]{todonotes}

\newcolumntype{b}{X}
\newcolumntype{s}{>{\hsize=.5\hsize}X}

\usepackage{mathtools}
\usepackage{amsfonts}
\usepackage{amssymb}
\usepackage{amsthm}
\usepackage{amsmath}
\usepackage{dirtytalk}
\usepackage{mathrsfs}  
\usepackage{lineno}
\newtheorem{remark}{Remark}
 \usepackage{multicol}  
 \usepackage{bbm}
 \usepackage{tabularx}
\newenvironment{claimproof}[1][\proofname]
{\begin{proof}[#1]}
{\end{proof}}

\begin{document}

\author{Jorge Antonio Cruz Chapital}
\title{On construction schemes: Building the uncountable from finite pieces}
\date{Agosto, 2019}

\frontmatter
\input{chapters/portada}

\cleardoublepage
\include{chapters/dedi}

\cleardoublepage
\include{chapters/acknowlmts}

\cleardoublepage
\include{chapters/dinero}

\cleardoublepage
\include{chapters/abstract}

\renewcommand{\chaptermark}[1]{\markboth{#1}{}}
\renewcommand{\sectionmark}[1]{\markright{#1}}
\selectlanguage{english}
\tableofcontents
\addcontentsline{toc}{chapter}{Resumen/Abstract}
\cleardoublepage
\selectlanguage{spanish}
\selectlanguage{english}
\include{chapters/intro}
%
\mainmatter
\input{chapters/Preliminaries}
\input{chapters/metricsandschemes}

\input{chapters/newgaps}

\input{chapters/Trees_and_lines}
\input{chapters/Ramsey}
\input{chapters/Oscillation}
\input{chapters/Onsigmamonotone}
\input{chapters/Fragments_of_Martins_axiom_and_ultrafilters}

\input{chapters/A_deeper_analysis_of_construction_schemes}
\newpage
\input{chapters/problems}

\backmatter

\cleardoublepage
\phantomsection
%

\cleardoublepage
\phantomsection
\addcontentsline{toc}{chapter}{Bibliography}

\bibliographystyle{plainurl}
\bibliography{chapters/phDThesis}

\end{document}

%% file: chapters/portada.tex
\selectlanguage{spanish}

\begin{titlepage}
\setlength{\parindent}{0pt} \setlength{\parskip}{0pt}

\begin{center}
 \vfill 
 
 \begin{minipage}{\textwidth}
  
  \newcolumntype{V}{>{\centering\arraybackslash} m{.17\textwidth} }
  \newcolumntype{C}{>{\centering\arraybackslash} m{.56\textwidth} }
  
   \begin{tabular}{ V C V }
  \includegraphics[width=.15\textwidth]{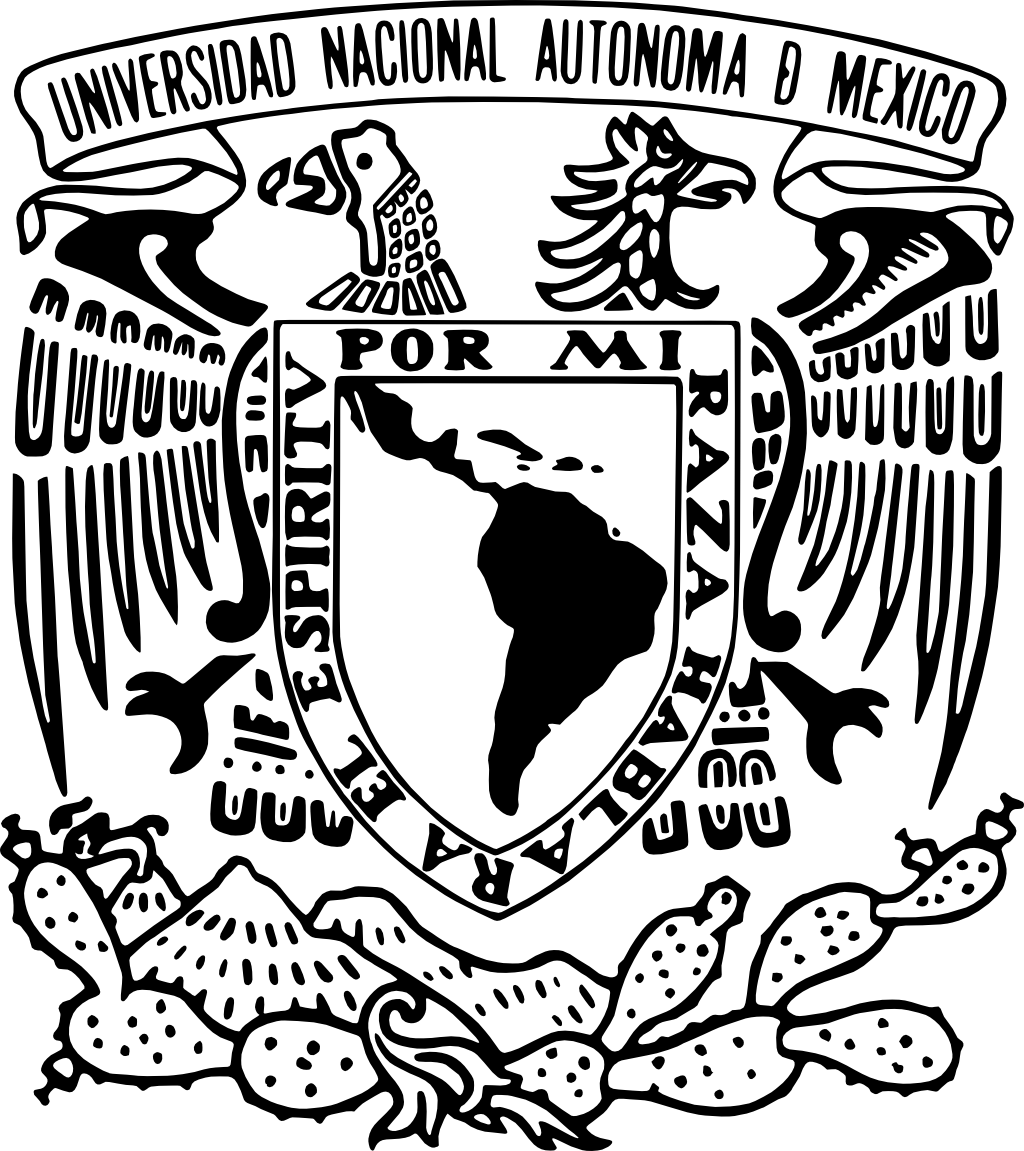}  %
  &\Large Universidad Nacional Aut\'onoma de M\'exico y Universidad Michoacana de San Nicol\'as de Hidalgo &%
  \includegraphics[width=.15\textwidth]{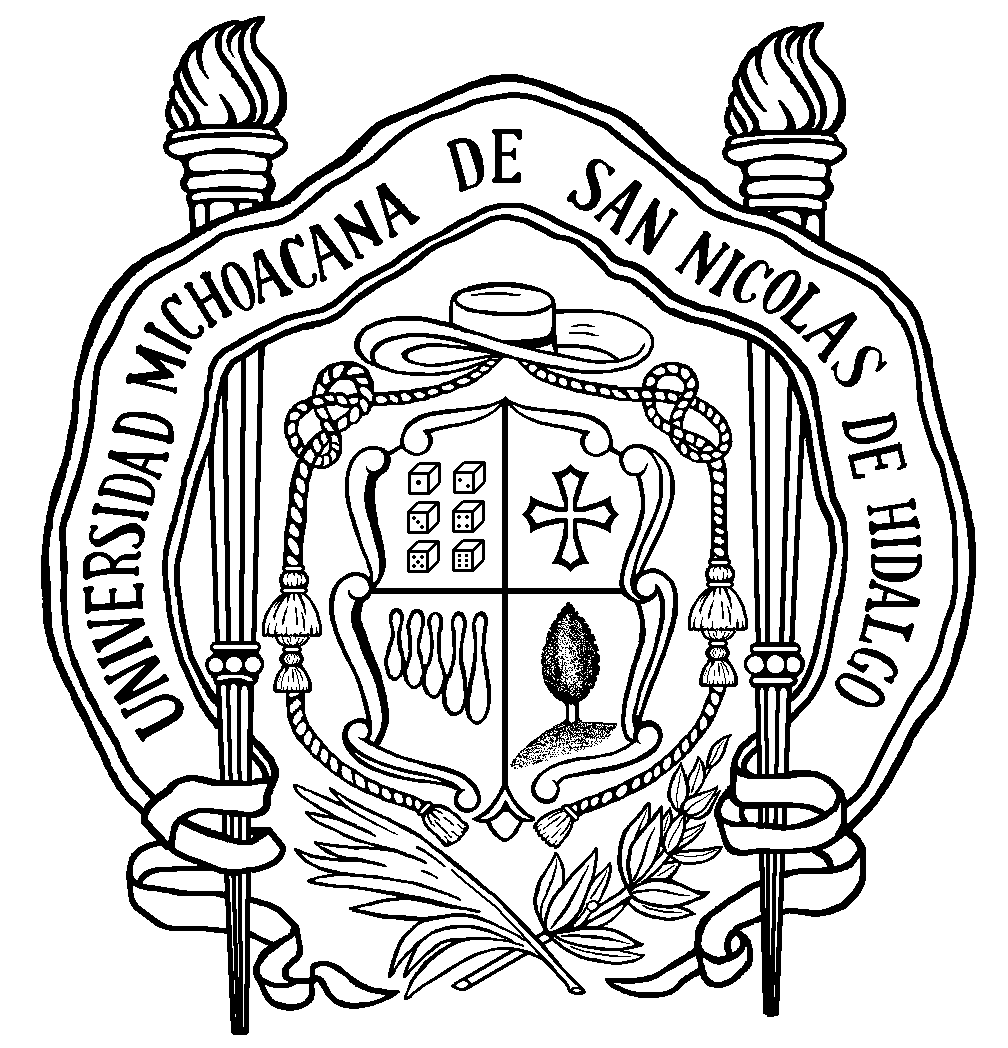}
 \end{tabular}

 \end{minipage}


 \vfill

\begin{minipage}{0.7\textwidth}
\begin{center}
\large Posgrado Conjunto en Ciencias Matem\'aticas UMSNH-UNAM
\end{center}
\end{minipage}
\begin{center}

\vfill
\medskip \rule{.9\textwidth}{2pt}
\bigskip

{\Large \bfseries  \thetitle   }

\medskip \rule{.9\textwidth}{2pt}
\end{center}

\begin{center}
\vspace*{1cm}
{\huge T E S I S}\\
\vspace*{1cm}
que para obtener el grado de \\ \bigskip
{\large \emph{Doctor en Ciencias Matemáticas}} \\ \bigskip

presenta\\ \medskip
\theauthor \\
\texttt{jorgeacruzchapital@hotmail.com}

\bigskip

Asesores:\\
Osvaldo Guzm\'an Gonz\'alez \\

Michael Hru\v{s}\'ak \\

\end{center}

\vfill

\begin{center}
{Morelia, Michoacán, México}\\
{Junio, 2024}
\end{center}
\cleardoublepage
\end{center}
\end{titlepage}

%% file: chapters/dedi.tex
\newpage
\thispagestyle{empty}
\vspace*{2cm}
\begin{flushright}

\end{flushright}
\vfill

%% file: chapters/acknowlmts.tex
\chapter*{Agradecimientos}
Puedo decir sin lugar a dudas que estoy donde estoy, y soy quien soy, gracias a toda la gente que me ha amado y apoyado a lo largo de los años. Entre dicha gente, el lugar mas especial lo ocupa mi  familia en constante crecimiento .  Les reitero mi gratitud, ya que sin ustedes nada de esto hubiera sido posible.

Fue por culpa Manuel Lara  y Jonas Martínez que conocí y me enamoré de la Teoría de Conjuntos. Gracias a Manuel, tambien conocí, al menos de nombre, a mis tres mayores heroes en el mundo de la Teoría de Conjuntos: Jonathan Cancino, Osvaldo Guzmán y Michael Hru\v{s}\'ak. A mitad de mi licenciatura tuve la oportunidad de tomar clase con Osvaldo y Jonathan, y desde ese momento reafirme mi deseo, no solo de hacer un posgrado en Teor\'ia de conjuntos, sino de hacerlo bajo la tutela de Michael. Afortunadamente sucedio. Durante estos ultimos 4 años, Michael y Osvaldo fueron mi asesores en el Doctorado. Les agradezco profundamente por su paciencia, comprension, empatía y por todo lo que pude aprender de ustedes a lo largo de este tiempo.

Un Agradecimiento especial merecen tambien mis sinodales: David Fernández, Fernando Hernández y Reynaldo Rojas. Sus valiosas aportaciones mejoraron esta tesis sustancialmente. De igual manera, siento la necesidad de agradecer a aquellos profesoras y profesores que más impacto tuvieron en mi formación matemática: Natalia Jonard, Jorge Marcos Martínez, Gabriela Campero, Ángel Tamariz, Diana Avella, David Meza, Luis Paredes, Pierre Bayard, Javier Páez, Luis Turcio, Rodrigo Hernández, Ariet Ramos, Leonardo Salmerón y Luis Nava.  

Creo firmemente que el One Piece si son los amigos que hicimos en el camino. Le agradezco a Richi por su incondicional amistad a lo largo de todos estos años, al igual que por la interminable ayuda que me ha brindado. Él, sin duda, ocupa un lugar muy especial dentro de estos agradecimientos. En igual medida le agradezco a Marcos y Uhthoff. Me llena de felicidad saber que las amistades pueden perdurar tanto tiempo. Durante mi estancía en el $CCM$ conoci a grandes amigos. David Valencía, por ejemplo, llegó al punto de llevarme atole a mi casa una vez que estaba triste. Por todo esto y más, le estoy muy agradecido. Definitivamente, no hubiera divertido tanto como lo hice si no hubiera sido por los amigos que conoci durante esta etapa de mi vida: Norberto, Miguel, Mario, Yhon, Carlos, Emmanuel, Kevin, Yulu, Sonia, Angel, Francisco, Tero, Jesus, Daniel, Tristan, Tatsuya, Julia,  y muchos más.  Otro agradecimiento especial se lo debo a Cesar, pues además de darme su amistad, fue mi guia y ejemplo a seguir durante gran parte de mi posgrado. Él, Thelma, Mati y Dony me adoptaron cuando era un vagabundo en las frias calle de Toronto. Por un mes, me hicieron parte de su familia, me dieron comida y techo, y jamás me pidieron nada a cambio. A ustedes cuatro les estaré eternamente agradecido.

Agradezco de nuevo a Tero, ya que sin su ayuda este trabajo no tendría la calidad que tiene. El me pasó la plantilla con la que hizo su tesis, y durante estos dos últimos años no dudó en resolver al instante cualquier duda que llegue a tener relacionada a Latex.

Por último, me gustaría agradecer a Stevo Todor\v{c}evi\'c. Primero, por las bellas matemáticas que ha desarrollado a lo largo de los años, sin la cuales, esta tesis jamás habría existido. Tambien le agradezco por la oportunidad que me dio de poder trabajar junto con él y Osvaldo en esquemas de construcción.  Durante el desarrollo de este proyecto, Stevo no solo ha sido coautor, también ha sido un gran mentor y maestro. Muchas gracias Stevo.

%% file: chapters/dinero.tex
\newpage
\thispagestyle{empty}
\vspace*{8cm}
\noindent La investigación de este trabajo se realizó en parte gracias al apoyo económico brindado por el gobierno de M\'exico, por parte del Consejo Nacional de Humanidades, Ciencias y Tecnolog\'ias (\textbf{CONAHCyT}). Por esto y más, les estoy eternamente agradecido.

%% file: chapters/abstract.tex
\chapter*{Abstract}
\thispagestyle{plain}
\begin{center}
    \Large
    \textbf{\thetitle}
 
    \large
 
    \vspace{0.4cm}
    {\theauthor}
 
    \vspace{0.9cm}
\end{center}

\selectlanguage{spanish}
\section*{\centering \begin{normalsize}Resumen\end{normalsize}}
\begin{quotation}
\noindent
En esta tesis se desarrollar\'a un an\'alisis estructural de los esquemas de construci\'on introducidos  en \cite{schemenonseparablestructures}. La importancia de este estudio se ver\'a reflejada al construir una gran cantidad de objetos combinatorios distintos que han sido de gran inter\'es en las matem\'aticas. Tambien se continuar\'a con el estudio de los axiomas de captura asociados a esquemas de construcci\'on. De dichos axiomas se construir\'an varios objetos cuya existencia se sabe independiente de los axiomas usuales de la teor\'ia de conjuntos.
\end{quotation}
\begin{quotation}
 \noindent \textit{Palabras Clave: uncountable, construction scheme, diamond principle, morasses, capturing, trees, gap, coloring.}
\end{quotation}

\vfill
\selectlanguage{english}

\section*{\centering \begin{normalsize}Abstract\end{normalsize}}
\begin{quotation}
\noindent In this thesis, a structural analysis of construction schemes (as introduced in \cite{schemenonseparablestructures}) is developed. The importance of this study will be justified  by constructing several distinct combinatorial objects which have been of great interest in mathematics. We then continue the study of capturing axioms associated to construction schemes. From them, we construct several uncountable structures whose existence is known to be independent from  the usual axioms of Set Theory.
\end{quotation}
\begin{quotation}
 \noindent \textit{Keywords: uncountable, construction scheme, diamond principle, morasses, capturing, trees, gap, coloring.}
\end{quotation}
\thispagestyle{empty}

%% file: chapters/intro.tex
\chapter{Introduction}

Throughout history, infinity has captivated and amazed humanity, sparking curiosity and contemplation among generations. A fabulous trip lasting more than 2000 years lead, in the 19th century, to beautiful results regarding this concept. Georg Cantor delved deeper into the nature of infinity by proving, among other things, that infinite sets may have different \say{sizes}. Formally, he showed that there can not be a bijective function from the natural to the real numbers. Since then, a great part of effort has been put into understanding the mysteries lying in the realm surrounding the two aforementioned infinities. In this thesis, we will explore $\omega_1$, the first uncountable cardinal lying in \textit{between $\omega$ and the hill of treasures} that are the real numbers.\\

The purpose of this thesis is to study $\omega_1$ by further developing the theory of construction schemes as introduced by Stevo Todor\v{c}evi\'c in \cite{schemenonseparablestructures}. By means of descriptive set-theoretic results, it is often impossible to deduce the existence of interesting objects of size $\omega_1$ by presenting them throughout \say{nice and simple} definitions. For that reason, there has been a lot of interest in presenting and studying methods for constructing such objects:
\begin{itemize}

\item Probably the most straightforward approach is to build them by transfinite recursion in $\omega_1$ many steps. Typical examples of this method are the construction of an Aronszajn tree from  \cite{jechsettheory} or \cite{kunensettheory} and the construction of a Hausdorff gap from \cite{ScheepersGaps} or \cite{integersvan}.

\item A second approach would be to show the consistency of an uncountable object with the desired properties, analyze the complexity of the sentences defining such properties, and then appeal to Jerome H. Keisler's completeness theorem for $L^\omega_\omega(\mathcal{Q})$ (see \cite{Keisleruncountablymany}) to conclude that such an object in fact exists. An example of this method is the coherent family of functions from \cite{coherentfamilyoffunctions}.
A variation of this method involving the $\Diamond$-principle was presented by Menachem Magidor and Jerome Malitz in \cite{magidorlQdiamond}.

\item A third approach is by using the method of walks on ordinals and ordinal metrics as introduced by Tordor\v{c}evi\'c (see \cite{Walksonordinals}). Historically, this method has been proved to be one of the most useful when dealing with constructions of uncountable objects. Examples of this are the solution of the famous $L$-space problem from \cite{solutionlspace} and the construction of a \say{rainbow coloring} from \cite{partitioningpairs}.

\item A fourth approach, which is arguably the most related to the present work, is to construct the desirable uncountable object by finite approximations through the use of simplified morasses as introduced by Ronald Jensen (see \cite{aspectsofconstructility}) and later simplified by Daniel Velleman in \cite{simplifiedmorasses}.

\end{itemize}
Just as simplified morasses, construction schemes are objects that let us build uncountable objects through finite approximations. Roughly speaking, a construction scheme $\mathcal{F}$ is a collection of finite subsets of $\omega_1$ with strong coherent properties (see Definition \ref{constructionschemedef}). These coherent properties allow us to build directed families of finite structures indexed by the elements of the scheme. Such structures are built recursively by amalgamation processes dictated by the construction scheme. In the end, the uncountable object that we construct is in some sense the direct limit of the family that we defined. The main difference between construction schemes and morasses is that in morasses we only amalgamate two structures at a given time, whereas in construction schemes, the number of amalgamations may vary.\\
That same attribute which sets appart construction schemes and morasses, allow construction schemes to \textit{consistently} have further properties, powerful enough to imply the existence of a large amount of objects which are known to exist under extra set-theoretic assumptions. These properties, which we call \textit{$n$-capturing}, \textit{capturing} and \textit{fully capturing}, can be described as \say{finitizations} (relative to a construction scheme) of the well-known $\Diamond$-prinicple (see \ref{capturedfamiliesdef}). The capturing axioms $CA_n$, $CA$ and $FCA$ assert the existence of construction schemes which are $n$-capturing, capturing and fully capturing respectively. All of these axioms are implied by the $\Diamond$-principle and also hold in any forcing extension of the universe obtained by adding at least $\omega_1$ many Cohen reals.\\
Given an $n$-capturing construction scheme $\mathcal{F}$ we also can define the $n$-parameterized Martin's number associated to it as $\mathfrak{m}^n_\mathcal{F}$. This cardinal invariant is defined as the Martin's number corresponding to the family of $ccc$-forcing notions which force our construction scheme to be $n$-capturing. 
Even though this cardinal is strictly related to a combinatorial structure over $\omega_1$ (the construction scheme), the assertion \say{$\mathfrak{m}^n_\mathcal{F}>\omega_1$} is strong enough to imply the existence of objects such as Ramsey ultrafilters over $\omega$. Furthermore, as we will show, there are  meaningful statements about uncountable objects which Martin's axiom can not decide, that are true under \say{$\mathfrak{m}^n_\mathcal{F}>\omega_1$}, and which are false under axioms such as $PFA$, and $CH$ (so in particular, $\Diamond$-principle).\\

\section*{Structure of the thesis}

In this thesis we  will develop to a great extent the theory of construction schemes and their structural properties, and we will obtain important results about set-theoretic assumptions associated with them. The importance of this study will be justified by constructing several distinct combinatorial objects which have been of a great interest to set theorists and topologists. A large amount of these objects will be constructed in a rather different form from the original constructions, and some of them will make their first appearance here. \\

\noindent
In {\bf Chapter \ref{prelimaries}} we will present some of the preliminaries, and fix the notation that will be used through out the thesis.\\

\noindent
 {\bf Chapter \ref{metricandschemeschapter}} is intended to serve as an introduction to the theory of ordinal metrics and construction schemes over $\omega_1$. In this chapter, we will define these two notions and analyze their general behavior. In particular, we will prove that construction schemes are in a one-to-one correspondence between certain kinds of ordinal metrics. We will also define the notions of $n$-capturing, capturing, and fully capturing construction schemes, as well as the capturing axioms $CA_n$, $CA$, and $FCA$ associated to them. Finally, we will introduce the parametrized Martin's axioms $\mathfrak{m}^n_\mathcal{F}$ associated to $n$-capturing construction schemes. Along the way, we will present the main theorems regarding the existence of distinct types of construction schemes. Nevertheless, such theorems will be proved and discussed in further detail until Chapters \ref{martinschapter} and \ref{deeperanalysis}. \\

\noindent
In {\bf Chapter \ref{applicationschapter}} we will provide a handful of applications of construction schemes both in $ZFC$ and under extra set-theoretic assumptions. These applications are of interest in set theory, topology, infinite combinatorics, algebra, and analysis. The constructions are mostly independent of each other, so that the reader can start from the ones that she or he finds more interesting. For the convenience of the reader, we list here some of the constructions that appear in this chapter. References for further study and historic remarks will be provided as we encounter them.
\begin{enumerate}
\item \textbf{Hausdorff gaps} An interesting feature of the Boolean algebra $\mathscr{P}(\omega)/\text{FIN}$ is that it is not complete. The easiest way to see this is as follows: Let $\{A_\alpha\,:\,\alpha\in 2^\omega\}\subseteq [\omega]^{\omega}$ be an $AD$ family. For every $H\subseteq 2^\omega$, consider the set $\mathcal{A}_H=\{[A_\alpha]\,:\,\alpha\in H\}$ (where $[A_\alpha]$ is the class of $A_\alpha$ in $\mathscr{P}(\omega)/\text{FIN}$). It follows form a simple counting argument that there must be an $H\subseteq 2^\omega$ such that $\mathcal{A}_H$ has no supremum. Although this is a very simple argument, we are often interested in more concrete examples of the incompleteness of $\mathscr{P}(\omega)/\text{FIN}$. The nicest examples are provided by gaps. Moreover, gaps are important because they represent obstructions that we may encounter when embedding structures in $\mathscr{P}(\omega)/\text{FIN}$ (a very illustrative example of this situation is Theorem 8.8 of \cite{PartitionProblems}). The classic construction of a Hausdorff gap requires a very clever argument to take care of $2^\omega$ many tasks in only $\omega_1$ many steps. With construction schemes, we will be able to provide a very simple construction of a Hausdorff gap, which under extra set-theoretic assumptions satisfies some interesting and (until now) unstudied combinatorial properties.

    \item \textbf{Luzin-Jones almost disjoint families.} This topic is about almost disjoint families (AD) and their separation properties. A \textit{Luzin family} is an $AD$ family of size $\omega_1$ in which no two uncountable subfamilies can be separated. On the other hand, a \text{Jones} family is an $AD$ family with the property that every countable subfamily of it can be separated from its complement. It is easy to prove that both of this kind of families exist. However, building a family that is both Luzin and Jones at the same time is much more complicated, since there is a tension between this two properties. A very difficult and highly complex construction of a Luzin-Jones family appears in \cite{guzman2019mathbb}. With construction schemes, we will be able to build an almost disjoint family satisfying strong  properties related to a concept that we will call \textit{the Luzin representation of a partial order}. From these properties we will not only be able to show that the constructed family is in fact Luzin-Jones, but from it, we will derive the existence of a handful of objects related to gaps. In particular, we will generalize, in strong way, the results of Hausdorff regarding the existence of gaps in $\mathscr{P}(\omega)/\text{FIN}$. It is worth pointing out that Luzin-Jones families can also be used to build interesting examples in functional analysis (see \cite{anonstable}).

\item \textbf{Donut-inseparable gaps.} Naively speaking, we say that gap over $\mathscr{P}(\omega)\backslash \text{FIN}$ is donut-inseparable if it can be reconstructed from two disjoint inseparable subfamilies of an $AD$ family. We will show, without assuming any extra axioms, that there are some gaps which are donut-inseparable and some which are not. This result suggest the following weakining of the previously defined notion: A gap is strongly donut-separable if each cofinal subgap is donut-separable. With the help of construction schemes and their parametrized Martin's axioms, we will show that the statement \say{There are no strongly donut-separable gaps} is independent from $ZFC$.
\item\textbf{Gap cohomology groups.} In \cite{cohomologytalayco}, Talayco defined what it means for two $(\omega_1,\omega_1)$-gaps to be cohomologically different. This notion allowed him to define \textit{the gap cohomology group} of an $\omega_1$-tower. We will generalize the notion of the gap cohomology group to a greater class of substructures of $\mathscr{P}(\omega)/\text{FIN}$ and show that all these groups can be as big as possible.
    
    \item \textbf{Luzin coherent family of functions.} We now look at a generalization of the Hausdorff gaps discussed previously. A Luzin coherent family of functions is a coherent system of functions supported by a pretower, in which we impose a strong non triviality condition. The importance of these families is that they provide many cohomologically different gaps. They were first studied by Talayco in \cite{cohomologytalayco}. Later Farah proved that such families exist (see \cite{coherentfamilyoffunctions}). The proof of Farah is highly non-constructive and indirect, since it appeals to Keisler's completeness Theorem. We will build such families using a construction scheme. No previous direct construction was known.

\item \textbf{Independent coherent family of functions}.
 We now return to the study of gaps that are obtained from a coherent family of functions. However, this time we want our family of gaps to be \say{independent}. This means that we can either fill or freeze any subfamily without filling or freezing any of the remaining gaps in the family. A similar result was obtained by Yorioka assuming the $\Diamond$-principle in \cite{yoriokadestructiblegaps} (the analogue result for Suslin trees was proved by Abraham and Shelah in \cite{Shelahabrahamincompactness}).

 \item \textbf{Countryman lines}. Let $(X,<)$ be a total linear order. Except for the trivial cases $X^2$ is not a linear order. In this way, it makes sense to ask how many chains we need to cover it. A Countryman line is an uncountable linear order whose square can be covered with only countably many chains. These orders seem so paradoxical at first glance that Countryman conjectured they do not exist. However, it was first proved by Shelah that Countryman lines do exist (see \cite{decomposingsquares}).

    \item \textbf{Aronszajn trees.} Aronszajn trees are the most well-known examples of the incompactness of $\omega_1$. An Aronszajn tree is a tree of height $\omega_1$, its levels are countable, yet it has no cofinal branches. A simple way in which we can guarantee that a tree (of height $\omega_1$) has no cofinal branches is to make it \say{special}, which means that it can be covered with countably many antichains (equivalently, they can be embedded in the rational numbers). Aronszajn proved this kind of trees exists.

    \item\textbf{Suslin trees.} A Suslin tree is an Aronszajn tree in which every antichain is countable. The Suslin Hypothesis ($SH$) is the statement that there are no Suslin trees. We now know that $SH$ is independent from $ZFC$. A related concept, the Suslin lines, were introduced by Suslin while studying the ordering of the real numbers. Kurepa was the one to realize that there is a Suslin line if and only if there is a Suslin tree. Applications and constructions from Suslin trees are abundant in the literature. We will use capturing schemes to build two types of these trees: Coherent Suslin and full Suslin trees. These two families of trees are diametrically opposed. While forcing with a Coherent tree, it completely destroys the $ccc$-ness of it, while with a full Suslin tree, many subtrees of it remain $ccc$.
    \item \textbf{Suslin lower semi-lattices}. If in the definition of a Suslin tree, we relax the condition of being a tree to just a being a lower semi-lattice, we get the notion of a Suslin lower semi-lattice. They were introduced by Dilworth, Odell and Sari (see \cite{dilworth2007lattice}) in the context of Banach spaces. They were then studied by Raghavan and Yorioka (see \cite{raghavan2014suslin}). Among other things they proved that the $\Diamond$-principle implies that $\mathscr{P}(\omega)$ contains a Suslin lower semi-lattice. We were able to obtain the same result from a $2$-capturing scheme.
    \item \textbf{Entangled sets}. A well-known theorem of Cantor is that any two countable dense linear orders with no end points are isomorphic. The straightforward generalization to linear orders of size $\omega_1$ is false. For this reason we want to restrict to suborders of the real numbers. We say that  an uncountable $B\subseteq \mathbb{R}$ is $\omega_1$-dense if $|U\cap \mathbb{B}|=\omega_1$ whenever $U$ is an open interval whose intersection with $U$ is nonempty. A remarkable theorem of Baumgartner is that $PFA$ implies that any two $\omega_1$-dense sets of reals are isomorphic (see \cite{ApplicationsofPFA}). This statement is now known as Baumgartner axiom ($BA(\omega_1)$). An entangled set is a subset of $\mathbb{R}$ with very strong combinatorial properties. The existence of an entangled set implies the failure of the Baumgartner axiom. Entangled sets were introduced by Abraham and Shelah in order to show that $BA(\omega_1)$ does not follow from Martin's axiom (see \cite{MAdoesnotImplyBA}). We will use a capturing scheme to build an entangled set. In this way, the existence of certain capturing schemes contradict the Baumgartner axiom. 
 \item \textbf{$ccc$ destructible $2$-bounded coloring without injective sets}. A coloring $c:[\omega_1]^2\longrightarrow\omega_1$ is called $2$-bounded if every color appears at most $2$ times. A set $A\subseteq \omega_1$ is $c$-injective if no color appears twice in $[A]^2$. Galvin was the first to wonder if there is a $2$-bounded coloring without an uncountable injective set. He proved that such coloring exists assuming the Countinuum Hypothesis. On the other hand, the third author proved that no such coloring exist under $PFA$.  years later, Abraham, Cummings and Smyth proved that $MA$ is consistent with the existence of a $2$-bounded coloring without uncountable injective sets (see \cite{abrahampolychromatic}). After hearing this result, Friedman asked for a concrete example of a $2$-bounded coloring without an uncountable injective set, but that such set can be added with a $ccc$ partial order. In \cite{abrahampolychromatic} such example is constructed assuming $CH$ and the failure of the Suslin hypothesis. We will  find an example with a $3$-capturing scheme.

\item \textbf{Oscillation theory of $2$-capturing construction schemes} In the book \cite{PartitionProblems}, Todor\v{c}evi\'{c} developed an oscillation theory which is based on an unbounded family of functions. A plethora of applications of this theory have been found throughout the years. For example, the oscillation theory is key in the proof that $PFA$ implies that the continuum is $\omega_2$ (see \cite{topicssettheory}). We will develop a similar theory using a $2$-capturing scheme. An important difference between the classic oscillation theory and the one from capturing schemes, is that this new one is based on a bounded family of functions. Using this new oscillation theory, we can prove the existence of the following objects:
\begin{enumerate}
    \item \textbf{Sixth Tukey type}. The Tukey ordering is a useful tool to compare directed partial orders. Its purpose is to study how a directed partial order behave cofinally. It was introduced by Tukey in \cite{Tukey} in order to study convergence in topology. The Tukey classification of countable directed partial orders is very simple: Every countable directed partial order is Tukey equivalent to $1$ or to $\omega$. The Tukey classification of directed sets of size $\omega_1$ becomes much more interesting. We now have at leas five Tukey types: $1$, $\omega$, $\omega_1$, $\omega\times\omega_1$ and $[\omega_1]^{<\omega}$. We may wonder if there is a directed partial order of size $\omega_1$ that is not Tukey equivalent to one if this five. In \cite{CategorycofinaltypesII}, Isbell proved that $CH$ entails the existence of a sixth Tukey type. This was greatly improved by Todor\v{c}evi\'c in \cite{directedsetscofinaltypes}, where he proved that $CH$ implies that there are $2^{\omega_1}$ distinct Tukey types. On the other hand, in the same paper he showed that $PFA$ implies that there are no sixth Tukey types. Here, we found a sixth Tukey type from a capturing scheme.
    \item $\omega_1\not\rightarrow (\omega_1,\omega+2)^2_2$. Given $\alpha<\omega_1$, the partition $\omega_1\rightarrow (\omega_1,\alpha)^2_2$ means that for every $c:[\omega_1]^2\longrightarrow 2$, either there is an uncountable $0$-monochromatic set, or there is a $1$-monochromatic set of order type $\alpha$. A celebrated result of Erd\"os and Rado (extending a theorem by Dushnik, Miller and Erd\"os) is that $\omega_1\rightarrow (\omega_1,\omega+1)^2_2$. We may wonder if this theorem can be improved. This turns out to be independent from $ZFC$. The Proper Forcing Axiom implies that $\omega_1\rightarrow (\omega_1,\alpha)^2_2$ for every $\alpha<\omega_1$ (see \cite{notesforcingaxioms}), while $\mathfrak{b}=\omega_1$ implies $\omega_1\not\rightarrow (\omega_1,\omega+2)^2_2$. We will prove a similar result from our oscillation theory. 
    \item \textbf{Non productivity of $ccc$ partial orders}. When is the product of two $ccc$ partial orders again $ccc$? This question has been of interest to set theorists for a long time. On one hand, Martin's axiom implies that the product of $ccc$ partial orders is $ccc$. On the other hand, the failure of the Suslin Hypothesis implies the opposite. Consistent examples of two $ccc$ partial orders whose product is not $ccc$ have been constructed by Galvin under $CH$ (see \cite{kunensettheory}) and by Todor\v{c}evi\'c under $\mathfrak{b}=\omega_1$ in  \cite{PartitionProblems}. As an applications of the oscillation theory that we developed, we encounter new situations in which there are $ccc$ partial orders whose product is not $ccc$. 
    \item \textbf{Suslin towers}. In \cite{oscillationsintegers}, Todor\v{c}evi\'c developed an analogue of his oscillation theory, now this time based on non-meager towers. While studying this oscillation, Borodulin-Nadzieja and Chodounsk\'y introduced the notion of a Suslin tower. A tower $\mathcal{T}=\{T_\alpha\,:\,\alpha\in\omega_1\}$ is Suslin if every uncountable subset of $\mathcal{T}$ contains two pairwise $\subseteq$-incomparable elements. The existence of a Suslin tower is independent from $ZFC$. It can be proved that the Open Graph Axiom forbids the existence of such families, while a Suslin tower can be constructed assuming $\mathfrak{b}=\omega_1$ (see \cite{GapsandTowers}). We will build a Suslin tower as an application of the oscillation theory obtained from a $2$-capturing scheme.
    \item \textbf{S-spaces}. Hereditarily separable and hereditarily Lindel\"of are two properties that in some sense are dual to each other. We may wonder if they always coincide. The question is only of interest in the realm of regular spaces. An $S$-space is a regular hereditarily separable, non Lindel\"of space. The study of $S$-spaces (and the dual notion, $L$-spaces) used to be one of the most active areas of set-theoretic topology. $S$-spaces can be constructed under several set-theoretic hypothesis (like $CH$ or the negation of the Suslin Hypothesis). It was Todor\v{c}evi\'c who for the first time succeeded in proving that it is consistent that there may be no $S$-spaces (see \cite{PartitionProblems}). Here, we will apply the oscillation theory obtained from a $2$-capturing scheme to construct several $S$-spaces.
    \end{enumerate} 
\item \textbf{$\sigma$-monotone spaces.} Monotone and $\sigma$-monotone spaces are a particular kind of metric spaces which were defined in \cite{monotonemetricspaces} by Ale\v{s} Nekvinda and Ond\v{r}ej Zindulka. In \cite{universalzindulka}, Zindulka used such spaces to prove the existence of universal measure zero sets of large Hausdorff dimension. Here, we will show that the capturing axiom $CA$ implies the existence of a metric space of cardinality $\omega_1$ which has no uncountable monotone subspaces. 

 \end{enumerate} 




\noindent 
In {\bf Chapter \ref{martinschapter}} we will study the effect that forcing with $ccc$ posets has on the capturing properties related to construction schemes. This is done for two purposes: The first one is to prove that the notions of \textit{capturing} and \textit{capturing with partitions} are different, and the second one is to show that the inequality $\mathfrak{m}_\mathcal{F}>\omega_1$ is consistent with $ZFC$. Lastly, we will prove that under this later assumption, we can deduce the existence of Ramsey ultafilters over $\omega$. A particularity of these results is that such ultrafilters can be explicitly defined from a combinatorial structure over $\omega_1$ (i.e., a construction scheme).\\

\noindent
In {\bf Chapter \ref{deeperanalysis}} we will define the forcing $\mathbb{P}(\mathcal{F})$ and use it to construct construction schemes without assuming any extra axioms. Later we will prove that Jensen's $\Diamond$-principle implies the capturing axiom $FCA(part)$. This was first claimed in \cite{schemenonseparablestructures}. However, the sketch of the proof given there is incomplete. Finally, we will show that $PID$ implies the nonexistence of $2$-capturing construction schemes.\\

\noindent 
Finally, in \textbf{Chapter }\ref{problemschapter}, we briefly discuss some open problems and future lines of research.\\

\noindent
The following chart summarizes many of the constructions which can be carried out with a construction scheme, and that will appear in this thesis. All of the results stated in the chart are new contributions.\\\\

\noindent
\begin{tabularx}{1\textwidth}{
  | >{\small\raggedright\arraybackslash}X 
  || >{\small\centering\arraybackslash\hsize=0.5\hsize}X
  | >{\small\centering\arraybackslash\hsize=0.5\hsize}X |}
 \hline
  \textbf{Object} & \textbf{Extra set-theoretic assumption } & \textbf{Result number}\\
  \hline \hline
  Donut-separable gap & - & \ref{firsthausdorffdonutseparable} \\
  \hline
  Donut-inseparable gap & - &\ref{donutinseparablegap}\\
  \hline
 Almost disjoint family coding any $\omega_1$-like order & - & \ref{luzinjonestheorem} and  \ref{luzinreptheorem}\\
 \hline
$(X,Y)$-gap for  all $\omega_1$-like orders & - & \ref{allgapsthm}\\
 \hline
 $*$-Lower semilattices with large gap cohomology groups & - & \ref{hausdorffcoherenttheorem}\\
\hline

 Not strongly donut-separable Hausdorff gap & $\mathfrak{m}_\mathcal{F}>\omega_1$ & \ref{mfstronglydonuttheorem}\\ 
 \hline
 Ramsey ultrafilter  & $\mathfrak{m}_\mathcal{F}>\omega_1$ & \ref{UnFramseytheorem} and \ref{pivframseytheorem}\\
 \hline
 Rainbow coloring $o^*$ & $CA_2$ &\ref{coloringo*}\\
 \hline
 Sixth Tukey type  & $CA_2$ & \ref{sixthtukeytypethm}  \\
\hline
$\omega_1\not\rightarrow(\omega_1,\omega+2)^2_2$ & $CA_2$ & \ref{notdusnikmillertheorem}\\
\hline
Non productivity of $ccc$ & $CA_2$ & \ref{cccnotproductive} \\
\hline
Suslin towers & $CA_2$ & \ref{suslintowercoro} \\
\hline 
$S$-spaces & $CA_2$ & \ref{sspace1}, \ref{sspace2} and \ref{sspace3}\\
\hline
Suslin lower semi-lattices & $CA_2$ & \ref{suslinlatticescheme} \\
\hline
Failure of $BA(\omega_1)$ & $CA_2$ & \ref{failurebaomega1scheme} \\
\hline
$ccc$ destructible $2$-bounded coloring & $CA_3$& \ref{policromaticscheme}  \\
\hline
Uncountable metric space without uncountable monotone subspaces & $CA$& \ref{nonmonotonescheme}\\
\hline
Coherent Suslin Tree & $FCA(part)$ & \ref{coherentsuslinscheme}\\
\hline

Full Suslin Tree & $FCA$ & \ref{fullsuslinscheme}\\
\hline

Entangled sets & $FCA$ & \ref{entangledscheme}\\
\hline 
Independent coherent family of functions & $FCA$& \ref{independentcoherentschemefca}\\
\hline
\end{tabularx}\\\\

 In addition to applications, the main contributions of this work are listed below. All of the undefined notions will be defined throught the thesis.\\\\

 \noindent
\begin{tabularx}{1\textwidth}{
  | >{\small\raggedright\arraybackslash}X 
  || >{\small\centering\arraybackslash\hsize=0.5\hsize}X
  | >{\small\centering\arraybackslash\hsize=0.5\hsize}X |}
 \hline
  \textbf{Result} & \textbf{Extra set-theoretic assumption } & \textbf{Result number}\\
  \hline \hline
  There are no strongly donut-separable gaps Either &$PFA$ or $CH$& \ref{weaklydonutch} and \ref{stronglyseppfa}\\
  \hline MA is independent from the statement \say{There is a strongly donut-separable gap}& - &\ref{independentmastronglysep}\\
  \hline
  There are no destructible gaps& $\mathfrak{m}_\mathcal{F}>\omega_1$ & \ref{nodestructiblemf}\\
  \hline
  There are no Suslin trees & $\mathfrak{m}_\mathcal{F}>\omega_1$ &\ref{nosuslinmf}\\
  \hline
  $\mathfrak{c}=\mathfrak{m}_\mathcal{F}^n> \omega_1$ is consistent with an arbitrarily large continuum & - & \ref{consistencymnfcontinuum}\\
  \hline
  Consistency of the existence of a construction scheme $\mathcal{F}$ which is $n+1$-capturing, $\mathcal{P}$-capturing, but not $\mathcal{P}'$-capturing &-&\ref{partitionncapturingcorolary}\\
  \hline
  There are no $2$-capturing construction schemes& Either $\mathfrak{m}>\omega_1$ or $PID$& \ref{corollarymanotca2} and \ref{coropidnotca2}\\
  \hline
  Full capturing with partitions axiom ($FCA(part)$)& $\Diamond $& \ref{diamondfcapartheoremproooof}\\
  \hline
  \end{tabularx}\\\\
  \quad

  Most of the results presented here were obtained in a joint work with Osvaldo Guzm\'an and Stevo Todor\v{c}evi\'c. Part of them  appear in \cite{schemescruz} which has already been submitted.

%% file: chapters/Preliminaries.tex
\chapter{Preliminaries and Notation}\label{prelimaries}
We devote this chapter to fix most of the preliminaries, notation and terminology that will be used throughout the rest of this work. Broadly speaking, the results appearing in this text can be separated into two categories. The ones regarding the development of the theory of \textit{construction schemes} and the ones regarding applications of this theory to set theory and topology. Of course, familiarity with the basics of these two branches of mathematics is assumed. Nothing more is asked for the reader solely interested in applications. However, familiarity with the forcing method is mandatory for the reader interested in learning more about the development of the theory. For set theory, the texts Halbeisen \cite{halbeisensettheory}, Jech \cite{jechsettheory} and Kunen \cite{kunensettheory} provide, in particular, all the prerequisites needed. The texts Engelking \cite{engelkingtopology}, Nagata \cite{nagatatopology} and Willard \cite{willardtopology} do the same for the topological background.

\section{Set-Theoretic Notation}

Ordinals are usually denoted by the first lower case Greek letters $\alpha,\beta,\gamma,\delta,\dots$,  whereas the middle letters $\kappa,\lambda,\mu,\dots$ are reserved for cardinals (commonly infinite). The set of natural numbers is denoted by $\omega$ and its elements are mainly denoted by the lowercase letters in the English alphabet. $\omega_1$ stands for the first uncountable ordinal. The set of all limit ordinals smaller than $\omega_1$ is denoted by $\text{Lim}$. 

We denote the power set of $X$ by $\mathscr{P}(X)$. The cardinality of $X$ is denoted by $|X|$. Given a cardinal $\kappa$, we let $[X]^\kappa$ be the set of all subsets of $X$ of size $\kappa$ and $[X]^{<\kappa}$ denotes $\bigcup_{\alpha<\kappa}[X]^\alpha$. The set $[X]^{\leq \kappa}$ is defined in a similar way. By $\text{FIN}(X)$ we mean $[X]^{<\omega}\backslash \{\emptyset\}$. $\text{FIN}(\omega)$ is just denoted as $\text{FIN}$.  We write $X\subseteq^* Y$ if $X\backslash Y$ is finite and $X=^*Y$ if $X\subseteq^*Y$ and $Y\subseteq^*X$. Following this notation, we write $X\subsetneq^*Y$ whenever $X\subseteq^*Y$ but $X\not=^*Y$.

\begin{definition}[$\Delta$-system] A family $\mathcal{D}$ is called a \textit{$\Delta$-system} with root $R$ if $|\mathcal{D}|\geq 2$ and $X\cap Y= R$ whenever $X,Y\in \mathcal{D}$ are different.
\end{definition}

Let $f:X\longrightarrow Y$ be a function. As usual, for $x\in X$ we denote the image of $x$ under $f$ as $f(x)$. For $A\subseteq X$ and $B\subseteq Y$ the \textit{direct image} of $A$ under $f$, that is, $\{f(x)\,:\,x\in A\}$, is denoted as $f[A]$ and the \textit{inverse image} of $B$, that is, $\{x\in X\,:\,f(x)\in B\}$ is denoted by $f^{-1}[B]$. Sometimes we write $dom(f)$ instead of $X$ $im(f)$ instead of $f[X]$. Whenever we write $g;X\longrightarrow Y$ we mean that $g$ is a function with $dom(g)\subseteq X$ and $im(g)\subseteq Y$. The set of all functions from $X$ to $Y$ is denoted by $X^Y$.  For an ordinal $\beta$ we define $X^{<\beta}$ and $X^{\leq \beta}$ as $\bigcup_{\alpha<\kappa}X^\alpha$ and $X^{<\beta}\cup X^\beta$ respectively.

\section{Order-Theoretic Notation}\label{settheoreticnotation}
Let $(X,<)$ be a partial order and $x,y\in X$. We say $x$ and $y$ are \textit{comparable} if either $x\leq y $ or $y\leq x$.  We say $x$ and $y $ are \textit{compatible} and write it as $x\parallel y$ if there is $r\in X$ with both $r\leq x$ and $r\leq y$. The \textit{incompatibility} relation between $x$ and $y$ is written as $x\perp y$. We define the following sets:\begin{multicols}{2} 
\begin{align*} (-\infty, x)_X&=\{z\in X\,:\,z<x\}\\
 (-\infty, x]_X&=\{z\in X:\,z\leq x\} \\
 (x,\infty)_X&=\{z\in X\,:\,x<z\} \\
 [x,\infty)_X&=\{z\in X\,:\, x\leq z\}
 \end{align*}
 \columnbreak
\begin{align*}\\
 (x,y)_X&=(x,\infty)\cap (-\infty, y) \\
 (x,y]_X&=(x,\infty)\cap (-\infty,y] \\
 [x,y)_X&=[x,\infty)\cap (-\infty, y) \\
 [x,y]_X&=[x,\infty)\cap (-\infty, y] 
\end{align*}
\end{multicols}
We say that $A\subseteq X$ is an \textit{interval} in $X$ whenever $(x,y)\subseteq A$ for all $x,y\in A$. $A$ is an \textit{initial segment of }$X$ if $(-\infty, x)\subseteq A$ for every $x\in A$. $A$ is a \textit{final segment of }$X$ if $(x,\infty)\subseteq A$ whenever $x\in A$. $A$ is cofinal in $X$ if for any $[x,\infty)\cap A$ is non-empty for all $x\in X$. $A$ is a \textit{chain} if any two elements of $A$ are comparable. On the other hand, if any two elements of $A$ are incomparable, we call $A$ a \textit{pie}. Finally, $A$ is an \textit{antichain} in $X$ if any two elements in $A$ are incompatible.

Let $B\subseteq X$ and $x\in B$. $x$ is said to be a \textit{minimal} element of $B$ if there is no $y\in B$ with $y<x$. Furthermore, if $x\leq y$ for any $y\in B$ we say that $x$  is the \textit{minimum of }$B$ and denote it as $\min(B)$. Analogously, $x$ is said to be a maximal element of $B$ if there is no $y\in B$ with $y>x$. Moreover, $x\geq y $ for any $y\in  B$ we call $x$ the \textit{maximum} of $B$ and denote it as $\max(B)$. An element $z\in X$ is an \textit{upper bound} (resp. \textit{lower bound}) of $B$ if $z=\max(B\cup\{z\})$ (resp. $z=\min(B\cup\{z\})$. In case the set of upper bounds of $B$ has a minimum, it is denoted as $\sup B$. In a similar way, if the set of lower bounds of $B$ has a maximum, it is denoted as $\inf B$. Lastly, for an element $x\in X$ we denote by $$suc(x)=\{y\in (x,\infty)\,:\,y\textit{ is a minimal element of }(x,\infty)\,\},$$
$$pred(x)=\{y\in (-\infty,x)\,:\,y\textit{ is a maximal element of }(-\infty,x)\,\}.$$
Note that $y\in suc(x)$ if $x<y$ and there is no $z\in X$ with $x<z<y$. Analogously $y\in pred(x)$ if $y<x$ and there is no $z\in X$ such that $y<z<x$.
Let $A,B\subseteq X$ be arbitrary subsets of $X$. We say $A<B$ if $B$ is non-empty and $a<b$ for all $a\in A$ and $b\in B$. We say that $B$ \textit{end-extends} $A$ (or $A$ is \textit{end-extended} by $B$) and write it as $A\sqsubseteq B$ if $A$ is an initial segment of $B$.

\begin{definition}[Root-tail-tail $\Delta$-system] A $\Delta$-system $\mathcal{D}\subseteq \mathscr{P}$ with root $R$ is said to be \textit{root-tail-tail} if $R\sqsubseteq A$ and either $A\backslash R<B\backslash R$ or $B\backslash R<A\backslash R$ for any two distinct $A,B\in \mathcal{D}$. In this same situation, if $\mathcal{D}$ is indexed by some ordinal $\gamma$ as, say, $\langle D_\alpha\rangle_{\alpha<\gamma}$ we will make the implicit assumption that $D_\alpha\backslash R<D_\beta\backslash R$ whenever $\alpha<\beta.$
\end{definition}

$X$ is said to be \textit{well-founded} if any non-empty subset of $X$ has a minimal element. In this case, there is a unique function $rank:X\longrightarrow Ord$ satisfying:$$rank(x)=\sup(\,rank(y)+1\,:\,y<x\,)$$
Note that this function is order preserving and its image is an ordinal. For each ordinal $\alpha$, we define the \textit{level $\alpha$ of $X$} as $X_\alpha=\{x\in X\,:\,rank(x)=\alpha\}$ and the \textit{height of $X$} as $Ht(X)=\min(\,\alpha\in Ord\,:\, X_\alpha=\emptyset\,)=rank[X]$. 
$X$ is said to be a \textit{well-order} if  any non-empty subset of $X$ has a minimum. In this case, there is a unique ordinal $\beta$ and a unique order preserving function $\phi:\alpha\longrightarrow X$ called \textit{the enumeration of $X$}. We call $\beta$ the order type of $X$ and denote it as $ot(X)$.\\\\
\fbox{\parbox{6in}{{\bf Very important notational remark} In the previous situation, following what some people may call abuse of notation, we identify $X$ with $\phi$. In this way, $X(\alpha)$ denotes $\phi(\alpha)$  and $X^{-1}(x)$ denotes $\phi^{-1}(x)$ whenever $\alpha\in \beta$ and $x\in X$. In the same way, $X[S]$ denotes $\phi[S]$ and $X^{-1}[Y]$ denotes $\phi^{-1}[Y]$ for all $S\subseteq \alpha$ and $Y\subseteq X.$ }}

\begin{definition}[Product order and lexicographical order]
Whenever $(X,\leq_X)$ and $(Y,\leq_Y)$ are two partial orders, we define the following two partial orders over the product $X\times Y$;
the \textit{order product }$\leq$ given by $(x,y)\leq (w,z)$ if $x\leq_X w$ and $y\leq_Y z$ and the \textit{lexicographical product} $\leq_{lex}$ given by $(x,y)\leq_{lex} (w,z)$ if either $x<_X w$ or $x=w$ and $y\leq_Y z$.
\end{definition}
\begin{definition}[Lower semi-lattice]$(X,\leq,\wedge)$ is said to be a \textit{lower semi-lattice} if $(X,\leq)$ is a partial order and $\wedge: X^2\longrightarrow X$ is such that for all $x,y\in X$, $\inf(x,y)$ exists and it is equal to $x\wedge y$.
\end{definition}

\section{Boolean algebras and filters}

Boolean algebras are one of the central objects of study in modern set theory and set-theoretic topology. Strictly speaking, all the results appearing in this thesis can be presented by avoiding any explicit mention to these objects. However, I think that, at least in this case, they provide a perfect setting for the motivation and analysis of some of the main topics which will be discussed eventually. We will briefly review some of the basic definitions and theorems. A detailed treatment of the subject can be found in Koppelberg \cite{handbookbooleanalgebras}.
\begin{definition}A Boolean algebra is a structure $(\mathbb{A},+,\cdot,-,\mathbb{0},\mathbb{1})$ with two binary operations\footnote{Here we decided to adopt the \say{algebraic} notation for the operations in Boolean algebras. Many authors prefer the usage of the symbols $\vee$, $\wedge$, and $\neg$ to refer to these operations (see Givant and Halmos \cite{introductionbooleanalgebras}).} $+$ and $\cdot$, a unary operation $-$, and two (distinct) distinguished elements $\mathbb{0}$ and $\mathbb{1}$ such that for all $x,y,z\in A$ the following conditions hold:
\begin{multicols}{2}
\begin{itemize}
\item $x+(y+z)=(x+y)+z,$
\item $x+y=y+x,$
\item $x\cdot( y+z)=x\cdot y+x\cdot z,$
\item $x+y\cdot x=x,$
\item $x+(-x)=\mathbb{1}$,
\columnbreak
\item $x\cdot(y\cdot z)=(x\cdot y)\cdot z,$
\item $x\cdot y=y\cdot x,$
\item $x+(y\cdot z)=(x+y)\cdot(x+z),$
\item $x\cdot(y+x)=x,$
\item $x\cdot(-x)=\mathbb{0}.$
\end{itemize}
\end{multicols}
\end{definition}
The most tipical example of a Boolean algebra is $(\mathscr{P}(X),\cup,\cap,\,^c,\emptyset,X)$ where $^c$ represents the complementation with respect to $X$.

Any Boolean algebra has an associated partial order given by $$x\leq y \text{ if }x\cdot y=x \text{(equivalently }x+y=y\text{)}.$$ In this order, $\mathbb{0}$ and $\mathbb{1}$ are the minimum and maximum of $\mathbb{A}$, respectively. In addition, Boolean operations can be completely recovered by the order in the sense that $\sup\{x,y\}=x+y$ and $\inf\{x,y\}=x\cdot y$. In the case of $\mathscr{P}(X)$ the order induced by the boolean operations coincide with \say{ $\subseteq$}.
\begin{definition}Let $\mathbb{A}$ be a Boolean algebra and let $\mathcal{F},\mathcal{I}\subseteq \mathbb{A}$. We say that
\begin{multicols}{2}$\mathcal{F}$ is a \textit{filter} if:
\begin{itemize}
\item $\mathbb{1}\in \mathcal{F}$ and $\mathbb{0}\notin\mathcal{F}$,
\item $x\cdot y\in \mathcal{F}$ whenever $x,y\in \mathcal{F}$,
\item If $x\in \mathcal{F}$ and $y\geq x$ then $y\in \mathcal{F}$.
\end{itemize}
\columnbreak 

$\mathcal{I}$ is an \textit{ideal} if:
\begin{itemize}
\item $\mathbb{0}\in \mathcal{I}$ and $\mathbb{1}\notin\mathcal{I}$,
\item $x+y\in \mathcal{I}$ whenever $x,y\in \mathcal{I}$,
\item If $x\in \mathcal{I}$ and $y\leq x$ then $y\in \mathcal{I}$.
\end{itemize}
\end{multicols}
A maximal filter is called an ultrafilter, and a maximal ideal is called a prime ideal.
\end{definition}
In this text, we will work mainly with filters over Boolean algebras of the form $\mathscr{P}(X)$. In this context, a filter $\mathcal{F}$ over $\mathscr{P}(X)$ (or simply a filter over $X$) is a non-empty family of non-empty subsets of $X$ which is closed upwards and under finite intersections. We call such $\mathcal{F}$ \textit{principal} if there is $A\subseteq X$ for which $\mathcal{F}=\{B\subseteq X\,:\,A\subseteq B\}.$ If there is no such $A$, $\mathcal{F}$ is called non-principal. Analogously, an ideal $\mathcal{I}$ over $\mathscr{P}(X)$ (or simply an ideal over $X$) is a family of proper subsets of $X$  that is closed downards and under finite unions. \\
Let $\mathcal{F}$ be a filter over some set $X$. We define the \textit{family of positive sets with respect to $\mathcal{F}$}, that is, $\mathcal{F}^+$, as the set of all $A\subseteq X$ which intersect each member of $\mathcal{F}$. In the same way, if $\mathcal{I}$ is an ideal in some set $X$, we define $\mathcal{I}^+$ as $\mathscr{P}(X)\backslash \mathcal{I}.$

\begin{definition}[Quotient algebras] Let $\mathbb{A}$ be a Boolean algebra and let $\mathcal{I}\subseteq \mathbb{A}$ be an ideal. $\mathcal{I}$ naturally induces an equivalence relation \say{$\sim$} on $\mathbb{A}$ given by$$x\sim y\text{ if and only if }x\cdot(-y)+y\cdot(-x)\in \mathcal{I}.$$ 
We denote the equivalence class of $x$ under such relation as $[x]$. The set of equivalence classes is denoted as $\mathbb{A}/\mathcal{I}$. The operations on $\mathbb{A}/\mathcal{I}$ given by $-[x]=[x]$, $[x]+[y]=[x+y]$, and  $[x]\cdot[y]=[x\cdot y]$ are well-defined. Furthermore, $(\mathbb{A}/\mathcal{I},+,\cdot,-,[\mathbb{0}],[\mathbb{1}])$ is a Boolean algebra.
\end{definition}
The Boolean algebra in which we will focus most of our attention is $\mathscr{P}(\omega)/\text{FIN}$. Note that in this case, if $X,Y\in \mathscr{P}(\omega)$ then $[X]\leq [Y]$ if and only if $X\subseteq^* Y$. 

\section{Forcing}
Throughout this thesis, we will force downwards. Thus, a \textit{forcing notion} (or simply a \textit{ forcing}) is triplet $(\mathbb{P},\leq, \mathbb{1}_\mathbb{P})$ so that $(\mathbb{P},\leq)$ is a partial order with $\mathbb{1}_\mathbb{P}$ as its maximum. Usually, we will denote $(\mathbb{P},\leq,\mathbb{1}_\mathbb{P})$ simply as $\mathbb{P}$ unless there is any risk of confusion. A subset $D$ of $\mathbb{P}$ is said to be \textit{dense} if for any $q\in \mathbb{P}$ there is $p\in D$ so that $p\leq q$. Lastly, we call a subset $G$ of $\mathbb{P}$ a filter if $G$ is non-empty, it is closed upwards and for any $p,q\in G$ there is $r\in G$ so that $r\leq p,q$.
\begin{definition}[Cohen forcing]Given $\kappa$ an infinite cardinal, we define the \textit{$\kappa$-Cohen forcing} as $$\mathbb{C}_\kappa=\{\,p;\kappa\longrightarrow 2\,:\,p\textit{ is finite}\}$$
and order with the reverse inclusion. The forcing $\mathbbb{C}_\omega$ is simply called \textit{Cohen forcing} and denoted as $\mathbb{C}$.
\end{definition}

Let $V$ be a transitive model of $ZFC$ and $\mathbb{P}\in V$. We say that a filter $G$ is \textit{$\mathbb{P}$-generic over $V$} if $G$ intersects any dense set of $\mathbb{P}$ lying in $V$. In this situation, there is a minimal transitive model of $ZFC$ which extends $V$ and has $G$ as an element. This model is denoted as $V[G].$\\
The forcing relation is denoted as $\Vdash$. Given a formula $\phi$, we write $\mathbb{P}\Vdash\text{\say{$\phi$}}$ instead of $\mathbb{1}_\mathbb{P}\Vdash\text{\say{$\phi$}}$. We adopt the convention of denoting $\mathbb{P}$-names with dots above them. However, if $X$ is an element of the universe $V$, we will make an abuse of notation and denote the canonical $\mathbb{P}$-name of $X$ simply as $X$. Whenever $\dot{Y}$ is a $\mathbb{P}$-name and $G$ is a $\mathbb{P}$-generic filter over $V$, we denote the interpretation of $\dot{Y}$ with respect to $G$ as $\dot{Y}^G$.
Given an ordinal $\gamma$,  \textit{finite support iterations} of length $\gamma$ are denoted as $$(\langle \mathbb{P}_\xi\rangle_{\xi\leq \gamma},\langle \dot{\mathbb{Q}}_\xi\rangle_{\xi<\gamma}).$$ For us, elements $p\in \mathbb{P}_\gamma$ are finite partial functions with domain contained in $\gamma$.

We now define some classes of forcing notions which will be important in further discussions.

\begin{definition}[$ccc$ forcings]A forcing $\mathbb{P}$ is said to be $ccc$ if any antichain in $\mathbb{P}$ is at most countable.
\end{definition}
\begin{definition}[Knaster forcings]Let $2\leq n\in \omega$. Given a forcing $\mathbb{P}$, we say that $\mathcal{B}\subseteq \mathbb{P}$ is \textit{$n$-linked} if for any $p_0,\dots, p_{n-1}\in \mathcal{B}$ there is $p\in \mathbb{P}$ (not necessarily in $\mathcal{B}$) so that $p\leq p_i$ for each $i<n$. $\mathbb{P}$ is said to be \textit{$n$-Knaster} if for any $\mathcal{A}\in [\mathbb{P}]^{\omega_1}$ there is $\mathcal{B}\in [\mathcal{A}]^{\omega_1}$ which is $n$-linked. $2$-Knaster forcings are simply called \textit{Knaster}.
\end{definition}

\begin{definition}[$\sigma$-centered forcings.] Given a forcing $\mathbb{P}$, a subset $\mathcal{C}$ of $\mathbb{P}$ is said to be \textit{centered} if it is $n$-linked for each $2\leq n\in\omega$. We say that $\mathbb{P}$ is $\textit{$\sigma$-centered}$ if $\mathbb{P}=\cup_{l\in\omega}\mathcal{C}_l$ where $\mathcal{C}_l$ is centered for each $l$.
    
\end{definition}
\begin{definition}[Precaliber]  Let $\kappa$ be an infinite cardinal. We say that a forcing $\mathbb{P}$ has \textit{precaliber $\kappa$} if for each sequence $\langle p_\alpha\rangle_{\alpha\in \mathbb{P}}$ there is $X\in[\kappa]^\kappa$ for which $\{p_\alpha\,:\,\alpha\in X\}$ is centered.
    
\end{definition}
\begin{rem}Any $\sigma$-centered forcing has precaliber $\omega_1$, and each forcing with precaliber $\omega_1$ is $n$-Knaster for each $2\leq n\in\omega$. Furthermore, any Knaster forcing is $ccc.$
\end{rem}
The following theorem will be frequently used. 
\begin{theorem}Let $\mathbb{P}$ be an uncountable $ccc$ forcing. Then there is $p\in \mathbb{P}$ which forces the generic filter to be uncountable.
    
\end{theorem}
Given a forcing $\mathbb{P}$, we define the \textit{Martin's number of $\mathbb{P}$}, namely $\mathfrak{m}(\mathbb{P})$, as the minimal cardinal $\kappa$ for which there is a family $\mathcal{D}$ of size $\kappa$ of dense sets of $\mathbb{P}$ so that there is no filter $G$ intersecting each member of $\mathcal{D}$.  To learn more about Martin's numbers, see \cite{topicssettheory} and \cite{notesforcingaxioms}. Regarding this notion, we define the following cardinal invariants:
$$\mathfrak{m}=\min(\,\mathfrak{m}(\mathbb{P})\,:\,\mathbb{P}\textit{ is }ccc \,).$$
$$\mathfrak{m}_{K_n}=\min(\,\mathfrak{m}(\mathbb{P})\,:\,\mathbb{P}\textit{ is }n\text{-Knaster}\,).$$
\begin{theorem} Suppose that $\mathbb{P}$ is an uncountable $ccc$ forcing so that $\mathfrak{m}(\mathbb{P})>\omega_1$. Then $\mathbb{P}$ has an uncountable filter.
\end{theorem}
Two well-known axioms regarding the Martin's numbers of certain classes of forcing notions will appear in some discussions, namely $MA$ and $PFA$. We will state them bellow for the sake of completeness. The notion of \textit{property} will not be defined here. However, the reader may find a full treatment of this subject in \cite{propershelah}. \\

\noindent
{\textbf{Martin's Axiom [MA]:}} $\mathfrak{m}=\mathbb{c}>\omega_1$.\\\\
\noindent
{\textbf{Martin's Axiom for $n$-Knaster forcings}:} $\mathfrak{m}_{K_n}>\omega_1.$\\\\
\noindent
{\textbf{Proper Forcing Axiom [PFA]:}} $\mathfrak{m}(\mathbb{P})>\omega_1$ for any proper forcing notion $\mathbb{P}.$
\section{P-ideal Dichotomy}\label{pidealsection}
Let $\mathcal{I}$ be an ideal over a set $X$. We say that $\mathcal{I}$ is a \textit{$P$-ideal} if for each countable $\mathcal{A}\subseteq \mathcal{I}$ there is $B\in \mathcal{I}$ such that $A\subseteq^*B$ for any $A\in \mathcal{A}$.
Any such set $B$ is called a \textit{pseudo-union} of $\mathcal{A}$. The $P$-ideal dichotomy ($PID$) is a well-known consequence of $PFA$ regarding $P$-ideals. Stevo Todor\v{c}evi\'c introduced it in its current form in \cite{dichotomypidfirst}. There he showed that $PID$ is consistent with $CH$, although $PFA$ is not. A less general version of $PID$ was previously formulated by Uri Abraham and Todor\v{c}evi\'c in \cite{partitionpropertiesch}.

   \begin{definition}Let $\mathcal{I}$ be an ideal over a set $X$. We define the \textit{orthogonal ideal of $\mathcal{I}$}, namely  $\mathcal{I}^\perp$,  as the set of all $A\subseteq X$ such that $A\cap B=^*\emptyset$ for any $B\in \mathcal{I}.$  
\end{definition}

\noindent
{\textbf{$P$-ideal Dichotomy [PID]:}} Let $X$ be a set and $\mathcal{I}$ be a $P$-ideal over $X$ such that $[X]^{<\omega}\subseteq \mathcal{I}\subseteq [X]^{\leq \omega}$. Then one of the two following conditions hold:\begin{itemize}
    \item There is $Y\in [X]^{\omega_1}$ so that $[Y]^\omega\subseteq \mathcal{I}$.
    \item There is $\{ Z_n\,:\,n\in\omega\}\subseteq\mathcal{I}^\perp$ for which $X=\bigcup\limits_{n\in\omega} Z_n.$
\end{itemize}

\section{Diamond principle}

Undoubtedly, the most important axiom regarding this work is Jensen's $\Diamond$-principle. This principle was introduced by Ronald Jensen in \cite{jensendiamond}. In there, he showed that this principle holds in the constructible universe $L$. We will recall the version of the $\Diamond$-principle that we will be using.

A $\Diamond$-sequence (over $\omega_1$) is a sequence $\langle D_\alpha\rangle_{\alpha\in\text{Lim}} $ of countable subsets of $\omega_1$ for which the following properties hold:\begin{enumerate}[label=$(\alph*)$]
    \item For each $\alpha\in\text{LIM}$, $D_\alpha\subseteq \alpha.$
    \item Given $S\subseteq \omega_1$, the set $\{\alpha\in \text{Lim}\,:\,S\cap \alpha=D_\alpha\,\}$
    is stationary.
\end{enumerate}
\noindent
{\textbf{Jensen's $\Diamond$-principle [$\Diamond$]:}} There is a $\Diamond$-sequence.

%% file: chapters/metricsandschemes.tex
\chapter{Ordinal metrics and construction schemes}\label{metricandschemeschapter}
\begin{center}
\includegraphics[scale=0.9]{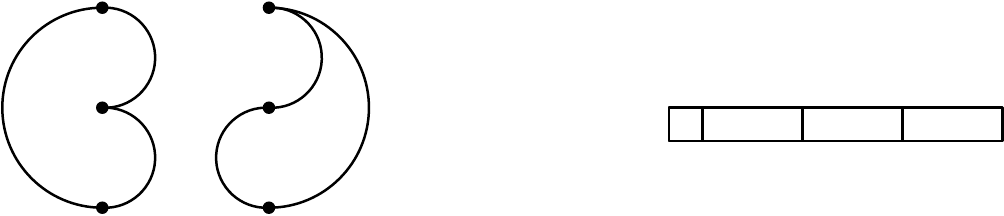}
\end{center}\vspace{0.1\linewidth}
In this chapter we introduce ordinal metrics and construction schemes as well as study the relationship between them.

Stevo Todor\v{c}evi\'{c} implicitly introduced ordinal metrics in \cite{partitioningpairs}\footnote{In Section 2 of such paper, it is proved  in (2.3) that the function $\rho$ defined immediately before satisfies the properties of what we now call an ordinal metric.}, as a tool for showing the existence of a coloring of the pairs of countable ordinals with uncountably many colors so that every uncountable set contains pairs of every color. In \cite{schemenonseparablestructures}, he introduced construction schemes as a tool for building uncountable structures by finite approximations. Broadly speaking, construction schemes provide blueprints for the domain of such approximations, as well as a procedure to recursively improve them through the use of finite amalgamations. This concept can be viewed as a generalisation of the so-called $(\omega,1)$-morasses as defined by Ronald Jensen (see \cite{aspectsofconstructility}) and simplified by Daniel Velleman in \cite{simplifiedmorasses}. The main difference between these two notions is that with morasses we always amalgamate \textit{two} objects when improving approximations, whereas in construction schemes the number of amalgamations allowed may vary.  

The relationship between ordinal metrics and construction schemes has already been highlighted in \cite{schemenonseparablestructures}. It is worth pointing out that Charles Morgan presented in \cite{morassessquareforcingaxioms} a detailed analysis of the relation between ordinal metrics and gap-1 morasses. Readers who want to 
know more about morasses and their applications may look at \cite{aspectsofconstructility}, \cite{somebanachspacesaddedbyacohenreal}, \cite{morassesincombinatorialsettheory}, \cite{ablackboxtheoremformorasses}, \cite{widescatteredspacesandmorasses}, \cite{simplifiedmorasses},  \cite{souslintreesconstructedfrommorasses}, and \cite{morassesdiamondandforcing}.

\section{Ordinal metrics}
For the rest of this chapter, $X$ stands for a non-empty set of ordinals.
\begin{definition}[Closure operation] Let $\rho:X^2\longrightarrow \omega$ be an arbitrary function. For each $F\in \text{FIN}(X)$ and $k\in \omega$ we define the \textit{$k$-closure} of $F$ as:
\begin{align*}(F)_k & =\{\, \alpha\in X\, :\, \exists \beta\in F\backslash \alpha \,(\rho(\alpha,\beta)\leq k )\,\}\\
(F)^-_k & =(F)_k\cap \max (F)=(F)_k\backslash \{\max(F)\}
\end{align*}
The \textit{diameter} of $F$ is defined as $$\rho^F=\max(\, \rho(\alpha,\beta)\,:\,\alpha,\beta\in F\,).$$
 $F$  is said to be \textit{$k$-closed} whenever $F=(F)_k$.  Moreover, if $F$ is $\rho^F$-closed we just say that $F$ is \textit{closed}.  For $\alpha\in X$, we will write $(\alpha)_k$ and $(\alpha)^-_k$ instead of $(\{\alpha\})_k$ and $(\{\alpha\})^-_k$ respectively.
\end{definition}
\begin{rem}Both the $k$-closure and the diameter are monotone operations. That is, if $F\subseteq G$ then $(F)_k\subseteq (G)_k$ and $\rho^F\leq \rho^G$.
\end{rem}

Following the convention established in the preliminaries, if $(\alpha)_k=\{\beta_0,\dots,\beta_n\}$ with $\beta_0<\dots<\beta_n$, we write $(\alpha)_k(i)=\beta_i$ for each $i\leq n$. This notation will be used throughout this thesis.
\begin{definition}[$k$-cardinality] Let $\rho$ be as in the previous definition. Given $\alpha\in X$ and $k\in\omega$ we define the \textit{$k$-cardinality} of $\alpha$ as $$\lVert\alpha \rVert_k=|(\alpha)^-_k|.$$
\end{definition}
Now we define ordinal metrics following the same approach (although slightly different notation) as in \cite{Walksonordinals}. In practice, we will only be interested in ordinal metrics whose domain is $\omega_1^2$. However, such metrics are constructed via recursion. For that reason, it is useful to present the definition in this generality.
\begin{definition}[Ordinal metric]\footnote{Some authors prefer to call ordinal metrics \say{$T$-functions} or \say{$\rho$-functions}.}\label{ordinalmetricdef} We say that $\rho:X^2\longrightarrow \omega$ is an \textit{ordinal metric} (over $X$) if:
\begin{enumerate}[label=$(\alph*)$]
\item $\forall\alpha,\beta\in X\,(\,\rho(\alpha,\beta)=0\textit{ if and only if }\alpha=\beta\,),$
\item $\forall \alpha,\beta\in X\,(\,\rho(\alpha,\beta)=\rho(\beta,\alpha)\,),$
  \item $\forall\alpha,\beta,\gamma\in X\,(\;\alpha\leq \beta,\gamma \rightarrow \rho(\alpha,\beta)\leq \max(\rho(\alpha,\gamma),\rho(\beta,\gamma))\;),$
\item $\forall \beta\in X\,\forall k\in\omega\,(\;\lVert\beta\rVert_k<\omega\;)$.
\end{enumerate}
\end{definition}

As the reader may note, the previous definition resembles that of an ultra-metric. However, one of the triangle inequalities is missing. In this way, one may interpret $(\beta)_k$ as the ball of radius $k$ centered on $\beta$ intersected with $\beta+1$. In the theory of metric spaces, we use the metric to approximate or locate points in a given space. To do that, we study the behaviour of balls centered on a given point as they become smaller. Ordinal metrics tend to work in a different way. Here, we want to construct a structure whose elements are parameterized by the points of the domain of our ordinal metric, say $X$. In order to do that, we make approximations of such structure by analyzing the interaction between the elements $\beta\in X$ and the points in $(\beta)_k$ as $k$ grows larger.  For that reason it may not be always suitable to carry the intuition of metric spaces into the context of ordinal metrics, but rather use the word \say{metric} to help us remember the main aspects of the definition. With the following proposition we show the difference between these two notions in a formal sense.  

\begin{proposition}No ordinal metric over $\omega_1$ is a metric.
\begin{proof}  
 Suppose $\rho$ is as previously mentioned and let $M$ be a countable elementary submodel of $H(\omega_2)$ having $\rho$ as an element. We will show that $\rho$ does not satisfy the triangle inequality of metric spaces. Take $\gamma\in \omega_1\backslash M $, consider an arbitrary $\alpha\in M\cap \omega_1$ and define $k$ as $\rho(\alpha,\gamma)$. Using elementarity we can find $\beta \in M\cap \omega_1$ above $\max((\gamma)_{2k}\cap M)$ for which $k=\rho(\alpha,\beta)$. By definition of the $k$-closure, we have that $\rho(\beta,\gamma)>2k=\rho(\alpha,\beta)+\rho(\alpha,\gamma)$. This finishes the proof.
 \end{proof}
 \end{proposition}

If $\rho$ is an ordinal metric over $X$ and $\alpha\in X$ then $(-\infty,\alpha]_X=\bigcup_{k\in\omega}(\alpha)_k$. Due to the property (d) of Definition \ref{ordinalmetricdef} this set must be countable. Hence, the existence of such $\rho$ implies that $ot(X)$ is at most $\omega_1$.

\begin{rem}If $\rho$ is an ordinal metric over $X$ and $F\in \text{FIN}(X)$ then $F\subseteq (F)_k$ for each $k\in\omega$. Furthermore, as the $k$-closure operation does not add points above the maximum of the set to which it is applied, it follows that $\max((F)_k)=\max(F)$ for any such $F$.
\end{rem}
The following proposition summarizes the most basic properties regarding the closure of sets in the context of ordinal metrics.
\begin{proposition}\label{closureprop1} Let $\rho:X^2\longrightarrow \omega$ be and ordinal metric, $F\in \text{FIN}(X)$ and $k\geq \rho^F$. Then:
\begin{enumerate}[label=$(\arabic*)$]
\item $\rho^F=\max(\,\rho(\alpha,\max(F))\,:\,\alpha\in F).$
\item $(F)_k=\{\,\alpha\leq \max(F) \,:\,\rho(\alpha,\max(F))\leq k\,\}=(\max(F))_k.$
\item If $\beta\in F$, then $(F)_k\cap (\beta+1)=(F\cap (\beta+1))_k.$
\item $(F)_k$ is $k$-closed and $\rho^{(F)_k}\leq k$, so in particular $(F)_k$ is closed. Furthermore, if $k=\rho^F$, then $\rho^{(F)_k}=k$.
\item  $F$ is closed if and only if $F=(\max(F))_{\rho^F}$.
\item If $G\in \text{FIN}(X)$ is $k$-closed, then $F\cap G\sqsubseteq F$.
\end{enumerate}
\begin{proof}The points (1) and (2) are direct consequences of the property (c) of Definition \ref{ordinalmetricdef} and the points (4) and (5) are  direct consequences of the point (1) and (2). Thus, we will only prove the points (3) and (6).

\begin{claimproof}[Proof of $(3)$.] It is clear that $(F\cap (\beta+1))_k\subseteq (F)_k\cap (\beta+1)$ so we will only prove that $(F)_k\cap (\beta+1)\subseteq (F\cap (\beta+1))_k$. For this take an arbitrary $\alpha\in (F)_k\cap (\beta+1)$. By the point (2) we
have that $\rho(\alpha,\max(F))\leq k$.  Using the property (c) of Definition 
\ref{ordinalmetricdef} we can conclude that $$\rho(\alpha,\beta)\leq \max(\rho(\alpha,\max(F)),\,\rho(\beta,\max(F)))\leq k.$$ Hence $\alpha\in (F\cap(\beta+1))_k$ and we are done.
\end{claimproof}
\begin{claimproof}[Proof of $(6)$]Fix $\beta\in F\cap G$ and let $\alpha\in F$ with $\alpha\leq \beta$ be arbitrary. Then $\rho(\alpha,\beta)\leq \rho^F\leq k$. As $G$ is $k$-closed, $\alpha\in G$. In this way $\alpha\in F\cap G$. This finishes the proof.
\end{claimproof}
 \end{proof}
\end{proposition}
In this way, in order to compute $\rho^F$ and $(F)_k$, we only need to use the maximum of $F.$
\begin{corollary}\label{initialsegmenteclosedisclosed}Let $\rho$ be an ordinal metric over $X$ and $F,G\in \text{FIN}(X)$ with $F\sqsubseteq G$. If $G$ is closed then so is $F$.
\begin{proof}As $G$ is closed then $G=(G)_{\rho^G}$. Let $\beta=\max F$. Since $F\sqsubseteq G$ then $\beta\in G$ and $F=G\cap (\beta+1)$. By the point (3) of Proposition \ref{closureprop1} we have that $$(F)_{\rho^G}=(G\cap(\beta+1))_{\rho^G}=(G)_{\rho^G}\cap (\beta+1)=G\cap (\beta+1)=F.$$
This means that $F$ is $\rho^G$-closed. But $F\subseteq G$ so $\rho^F\leq \rho^G$. Consequently $F$ is closed.
\end{proof}
\end{corollary}
We are interested in ordinal metrics which belong to three particular classes. In the rest of this section we will study them. But first, we need the following definition.
\begin{definition}Let $F\in \text{FIN}(X)$. We say that $F$ is \textit{maximally closed} if it is closed and there is no closed $G\in \text{FIN}(X)$ with $\rho^G=\rho^F$ such that $F\subsetneq G$.
\end{definition}

\begin{definition}[Locally finite metric]\label{locallyfinitedef}Let $\rho:X^2\longrightarrow \omega$ be an ordinal metric. We say that $\rho$ is \textit{locally finite} if $$\sup(\,|F|\,:\,F\in \text{FIN}(X), F\textit{ is closed and }\rho^F\leq k\,)$$
is finite for every $k\in \omega$.
\end{definition}
\begin{rem}\label{cofinalremarkregular} If $\rho$ is locally finite and $F\in \text{FIN}(X)$ then there is a maximally closed set $G\in \text{FIN}(X)$ such that $F\subseteq G$ an $\rho^F= \rho^G$. In fact, due to the point (6) of Proposition \ref{closureprop1} we can conclude that $(F)_{\rho^F}\sqsubseteq G$ for any such $G$. In particular, any closed set can be end-extended to a maximally closed set of the same diameter.
\end{rem}
The following result is an immediate consequence of the previous remark.
\begin{proposition}Let $\rho:X^2\longrightarrow \omega$ be a locally finite ordinal metric. For each $F\in \text{FIN}(X)$ we have that \begin{align*}\rho^F&=\min( \,\rho^G\,:\,G\textit{ is maximally closed and }F\subseteq G\,)\\
&=\min(\, \rho^G\,:\,G\textit{ is maximally closed and }F\sqsubseteq G\,)
\end{align*}
\end{proposition}

\begin{proposition}\label{regularimpliesunbounded}Let $\rho:X^2\longrightarrow \omega$ be a locally finite ordinal metric. For all $k\in\omega$ and every $A\in [X]^\omega$ there exist $\alpha,\beta\in A$ such that $\rho(\alpha,\beta)>k.$
\begin{proof} Let $k$ and $A$ be as in the hypotheses. Suppose towards a contradiction that $\rho(\alpha,\beta)\leq k$ for all $\alpha,\beta\in A$. Then $\rho^{A[n]}\leq k$ for each $n\in \omega$. Due to the point (6) of Proposition \ref{closureprop1} we conclude that $\langle (A[n])_k\rangle_{k\in\omega}$ is a strictly increasing sequence of closed sets all of whose diameter is at most $k$. This contradicts the fact that $\rho$ is locally finite. So we are done.
\end{proof}
\end{proposition}
Ordinal metrics satisfying the conclusion of the previous proposition deserve a special name.
\begin{definition}[unbounded metric]Let $\rho:X^2\longrightarrow \omega$ be an ordinal metric. We say that $\rho$ is \textit{unbounded} if for all $k\in\omega$ and every $A\in[X]^{\omega}$ there are $\alpha,\beta\in A$ with $\rho(\alpha,\beta)>k$.
    \end{definition}
    The following proposition can be found in \cite{Walksonordinals} (Lemma 3.4.11, page 73). The proof  that we present here avoids the use of ultrafilter quantifiers. 
\begin{proposition}Let $\rho:X^2\longrightarrow \omega$ be an unbounded ordinal metric and $k\in\omega$. For all $l\in\omega$ and each infinite $A\subseteq [X]^l$ family of pairwise disjoint sets there is $B\in [A]^\omega$ such that for all distinct $a,b\in B$ and every $\alpha\in a$ and $\beta \in b$, we have that $\rho(\alpha,\beta)>k$.
\begin{proof}Let $l$ and $A$ be as in the hypotheses.  Enumerate $A$ as $\langle a_n\rangle_{n\in\omega}$. For any $i<l$ there are two particular colorings that we may consider. The first one is the coloring $c_i:[\langle a_n(i)\rangle_{n\in\omega}]^2\longrightarrow 2$ given by $c_i(a_m(i),a_n(i))=0$ if and only if  $m<n$ and $a_m(i)<a_n(i)$. The second one is the coloring $o_i:[\langle a_n(i)\rangle_{n\in\omega}]^2\longrightarrow 2$  given by $o_i(a_m(i),a_n(i))=0$ if and only if $\rho(a_m(i),a_n(i))>k$. Since there are no infinite decreasing sequences of ordinals and $\rho$ is an unbounded metric then all infinite monochromatic sets with respect to these two coloring choose color $0$. In this way, by applying Ramsey's Theorem multiple (but finitely many) times, we can get a set $B\in [A]^\omega$ enumerated as $\langle b_n\rangle_{n\in\omega}$ satisfying the following conditions for each $i<l$ and all $m<n\in\omega$:
\begin{itemize}
\item $b_m(i)<b_n(i)$, 
\item $\rho(b_m(i),b_n(i))>k$.
\end{itemize}

and such that for all distinct $i,j<l$ either:\\
\begin{center}\begin{minipage}{3.5cm} \begin{center} \textbf{(A)}\end{center} For all $m<n\in \omega$, $ \rho(b_m(i),b_n(j))\leq k.$
\end{minipage}\hspace{3cm} \begin{minipage}{3.5cm}\begin{center} \textbf{(B)}\end{center} For all $m<n\in \omega,$ $\rho(b_m(i),b_n(j))> k.$
\end{minipage}
\end{center}
\bigskip
The proof wil be finished if we are able show that (B) holds for all distinct $i,j<l$. We divide the proof of this fact into three cases. In each of them, it is enough to show that there are $m<n$ for which $\rho(b_m(i),b_n(j))> k$. Indeed, this implies that (A) can not hold for $i,j,$ or equivalently, $(B)$ holds. \\\\
\underline{Case 1}: If $i<j$.

\begin{claimproof}[Proof of case.]Note that $b_0(i)<b_1(i)<b_2(i)<b_2(j)$. Then, $$k<\rho(b_0(i),b_1(i))\leq \max(\,\rho(b_0(i),b_2(j)),\,\rho(b_1(i),b_2(j))\,).$$
In other words, either $\rho(b_0(i),b_2(j))>k$ or $\rho(b_1(i),b_2(j))>k$. 
\end{claimproof}
\noindent
\underline{Case 2}: If $i>j$ and $b_2(j)<b_0(i).$

\begin{claimproof}[Proof of case.]
Here we also have that $b_1(j)<b_2(j)$. Hence, $$k<\rho(b_1(j),b_2(j))\leq \max(\,\rho(b_0(i),b_1(j)),\,\rho(b_0(i),b_2(j))\,).$$
As in the previous case, since either $\rho(b_0(i),b_1(j))>k$ or $\rho(b_0(i),b_1(j))>k$ then $B$ must hold.
\end{claimproof}
\noindent
\underline{Case 3}: If $i>j$ and $b_0(i)<b_2(j).$
\begin{claimproof}[Proof of case.]
 Here, as $b_0(i)<b_1(i)$ the inequality associated to the first case also holds. The rest of the argument is equal to the first two cases.
 \end{claimproof}
\end{proof}
\end{proposition}

\begin{definition}Let $\rho:X^2\longrightarrow \omega$ be an ordinal metric and let $F,G\in \text{FIN}(X)$. We say that $F$ and $G$ are \textit{$\rho$-isomorphic} if  $|F|=|G|$ and for all $i,j<|F|$, $$\rho(F(i),F(j))=\rho(G(i),G(j)).$$
Equivalently, if $h:F\longrightarrow G$ is the only increasing bijection\footnote{Note that if $|F|=|G|$ then $h$ is  defined by the formula $h(F(i))=G(i)$.} then $$\rho(\alpha,\beta)=\rho(h(\alpha),h(\beta))$$
for all $\alpha,\beta\in F$.
\end{definition}
\begin{definition}[Homogeneous metric]Let $\rho:X^2\longrightarrow \omega$ be an ordinal metric. We say that $\rho$ is \textit{homogeneous} if all $F,G\in \text{FIN}(X)$ maximally closed sets with the same diameter are $\rho$-isomorphic.
\end{definition}
\begin{lemma}\label{homoegeneitylemma}Suppose that $\rho:X^2\longrightarrow \omega$ is a locally finite and homogeneous ordinal metric.  If $F,G\in\text{FIN}(X)$ are two closed sets (not necessarily maximal) with $|F|=|G|$  and $\rho^F=\rho^G$ then $F$ and $G$ are $\rho$-isomorphic. 

\begin{proof}Let $k=\rho^F=\rho^G$. As $\rho$ is locally finite, there are maximally closed sets $\overline{F},\overline{G}\in \text{FIN}(X)$ end-extending $F$ and $G$ respectively and such that $\rho^{\overline{F}}=\rho^{\overline{G}}=k$. In particular, this means that $F(i)=\overline{F}(i)$ and $G(i)=\overline{G}(i)$ for each $i<|F|$. Now, since $\rho$ is homogeneous then $\rho(F(i),F(j))=\rho(\overline{F}(i),\overline{F}(j))=\rho(\overline{F}(i),\overline{F}(j))=\rho(G(i),G(j))$ for any $i,j<|F|\leq |\overline{F}|$. This finishes the proof. 
\end{proof}
\end{lemma}
\begin{rem}If $\rho$, $F$ and $G$ are as in the previous lemma and $\phi:F\longrightarrow G$ is the only increasing bijection then $\rho^H=\rho^{\phi[H]}$ for each $H\subseteq F$. As $F$ and $G$ are closed it is easy to see that $\phi[H]$ is closed whenever $H$ is. This is because $\rho^H\leq \rho^F$. Furthermore, if $H$ is maximally closed then $\phi[H]$ is also maximally closed. This is because $\phi[H]$ can be extended to a maximally closed set $\overline{\phi[H]}$ of diameter $\rho^H$. By homogeneity, such set has the same cardinality as $H$. But then $\phi[H]=\overline{\phi[H]}$.
\end{rem}
\begin{proposition}\label{propnoname}Let $\rho:X^2\longrightarrow \omega$ be an ordinal metric and let $\mathcal{A}\subseteq \text{FIN}(X)$ be a family of closed sets of diameter $k$. If $F\in \text{FIN}(X)$  is such that  $\rho^F\leq k$ and $F\subseteq \bigcup \mathcal{A}$ then there is $G\in \mathcal{A}$ for which $F\subseteq G$.
\begin{proof}By hypothesis, there is $G\in \mathcal{A}$ with $\max(F)\in G$. As $\rho^F\leq k=\rho^G$ and $G$ is closed, we know that $F\cap G\sqsubseteq F$ due to the point (6) of Proposition \ref{closureprop1}. It follows directly that $F\subseteq G$.
\end{proof}
\end{proposition}
As a corollary of the previous proposition we get that there are no leftovers in the unions of maximally closed sets of the same diameter. 
\begin{corollary}\label{uniquenessdescomosition}Let $\rho:X^2\longrightarrow \omega$ be an ordinal metric and let $F\in \text{FIN}(X)$. If $i<\rho^F$ and $\mathcal{A}\subseteq \text{FIN}(X)$ is a family of maximally closed sets of diameter $i$ such that $F=\bigcup\mathcal{A}$, then $$\mathcal{A}=\{\,G\in \text{FIN}(X)\,:\,G\subseteq F,\, G\textit{ is maximally closed and }\rho^G=i\,\}.$$
\begin{proof}We only need to prove the inclusion from right to left. For this, let $G\in \text{FIN}(X)$ be such that $G\subseteq F$, $G$ is maximally closed and $\rho^G=i$. By Proposition \ref{propnoname} we know there is $H\in \mathcal{A}$ with $G\subseteq H$. As $G $ is maximally closed and $\rho^G=i=\rho^H$, then $G=H$. This finishes the proof.
\end{proof}
\end{corollary}

\begin{definition}[Regular metric]\label{regularmetricdef}Let $\rho:X^2\longrightarrow \omega$ be an ordinal metric. We say that $\rho$ is \textit{regular} if for each $k\in \omega\backslash 1$ and each maximally closed set $F\in \text{FIN}(X)$ of diameter $k$ there are $j_F\in \omega\backslash 2$ and $F_0,\dots, F_{j_F-1}\in \text{FIN}(X)$ such that:
\begin{itemize}
\item $F=\bigcup\limits_{i<j_F} F_i,$
\item For each $i<j_F$, $F_i$ is a maximal closed set with $\rho^{F_i}=k-1$,
\item $\langle F_i\rangle_{i<j_F}$ forms a root-tail-tail $\Delta$-system with root $R(F)$. That is, $$R(F)<F_0\backslash R(F)<\dots<F_{j_F-1}\backslash R(F).$$
\end{itemize}
\end{definition}
\begin{rem} Observe that both the number $j_F$ and sequence $\langle F_i\rangle_{i<j_F}$ are unique for each $F$ due to Corollary \ref{uniquenessdescomosition}. So from now on, we will call $j_F$ the decomposition number of $F$ and $\langle F_i\rangle_{i<j_F}$ its canonical decomposition.

\end{rem}
\begin{lemma}\label{maximallyclosedregularlemma}Let $\rho:X^2\longrightarrow \omega$ be a regular, locally finite ordinal metric. If $F\in FIN(X)$ is a closed set with positive diameter and $i<\rho^F$, then there is a maximally closed set  $G\in\text{FIN}(X)$  with $\rho^G=i$ and $G\sqsubseteq F$.
\begin{proof}
The proof is carried by induction over $k=\rho^F$. For this, fix $k\in\omega$ and suppose that we have proved the lemma for each closed set $H\in \text{FIN}(X)$ with $\rho^H<k$.  Let $F\in \text{FIN}(X)$ be a closed set of diameter $k$ and let $i<\rho^F$. Furthermore, $R(F)$ is also uniquely determined so we will call it the root of $F$.

As $\rho$ is locally finite, there is a maximally closed set $\overline{F}$ with $F\sqsubseteq \overline{F}$ and $\rho^{\overline{F}}=k$. Since $\rho$ is regular and $k>0$ we can consider $j_{\overline{F}}\in \omega\backslash 2$ and $ \overline{F}_0,\dots,\overline{F}_{j_{\overline{F}}-1}\in \text{FIN}(X)$ as in Definition \ref{regularmetricdef}. Now, $F\not\subseteq \overline{F}_0$ as $\rho^{\overline{F}_0}=k-1$. As both $F$ and $\overline{F}_0$ are initial segments of $\overline{F}$ it follows that $\overline{F}_0\sqsubseteq F$. Observe that if $i=k-1$ then $\overline{F}_0$ testifies  the conclusion of the Lemma. Otherwise, if $i<k-1$ we can use the induction hypotheses to get a maximally closed set $G$ for which $\rho^G=i$ and $G\sqsubseteq \overline{F}_0\sqsubseteq F$. This finishes the proof.
\end{proof}
\end{lemma}
In this way, for a regular, locally finite ordinal metric, we can not \say{skip a distance} in a closed set.
\begin{corollary}Let $\rho:X^2\longrightarrow \omega$ be a regular, locally finite ordinal metric. If $F\in FIN(X)$ is closed and $i\leq\rho^F$, then there are $\alpha,\beta\in F$ for which $\rho(\alpha,\beta)=i.$
\end{corollary}

 \begin{theorem}\label{closedsizetheorem}Let $\rho:X^2\longrightarrow \omega$ be a regular, locally finite and homogeneous ordinal metric. If $F,G\in FIN(X)$ are closed sets of the same cardinality then $\rho^F=\rho^G.$    
\begin{proof}Suppose towards a contradiction that $\rho^F\not=\rho^G$. Without loss of generality we can assume that $\rho^F<\rho^G$. As $\rho$ is locally finite there is maximally closed set $\overline{F}$ with $F\sqsubseteq \overline{F}$ and $\rho^{\overline{F}}=\rho^F$. Due to Lemma \ref{maximallyclosedregularlemma} there is also a maximally closed set $H$ such that $H\sqsubseteq G$ and $\rho^H=\rho^F$. Observe that $|H|<|G|$. On the other hand, by the homogeneity of $\rho$ we have that $|\overline{F}|=|H|$. Thus, $|F|\leq|\overline{F}|=|H|<|G|$ which is a contradiction. This finishes the proof.
\end{proof}
\end{theorem}
\begin{corollary}\label{ballhomogeneitylemma}Let $\rho:X^2\longrightarrow \omega$ be a regular, locally finite and homogeneous ordinal metric. If $\alpha,\beta\in X$ and $k\in\omega$ is such that $\lVert \alpha\rVert_k=\lVert \beta\rVert_k$ then $(\alpha)_k$ and $(\beta)_k$ are $\rho$-isomorphic.
\begin{proof}By Theorem \ref{closedsizetheorem} it follows that $\rho^{(\alpha)_k}=\rho^{(\beta)_k}$ so by Lemma \ref{homoegeneitylemma} we are done.
\end{proof}
\end{corollary}
\section{Construction schemes}\label{constructionschemessection}

By means of Lemma \ref{homoegeneitylemma} and Theorem \ref{closedsizetheorem}, closed sets in regular locally finite homogeneous ordinal metrics are completely determined by their size. Given such a metric $\rho$, let us call $m^\rho_k$ the common cardinality of maximally closed sets of diameter $k$. The sequence $\langle m^\rho_k\rangle_{k\in\omega}$ is strictly increasing and the first point of such sequence, that is $m^\rho_0$, is $1$ due the point (a) of Definition \ref{ordinalmetricdef}. A closed set will have diameter $k+1$ if and only if its cardinality lies in the interval $(m^\rho_k,m^\rho_{k+1}]$, and will be maximal only when its size is exactly $m^\rho_{k+1}$.
\begin{center}
\begin{minipage}[b]{0.8\linewidth}
\centering
\includegraphics[scale=0.8]{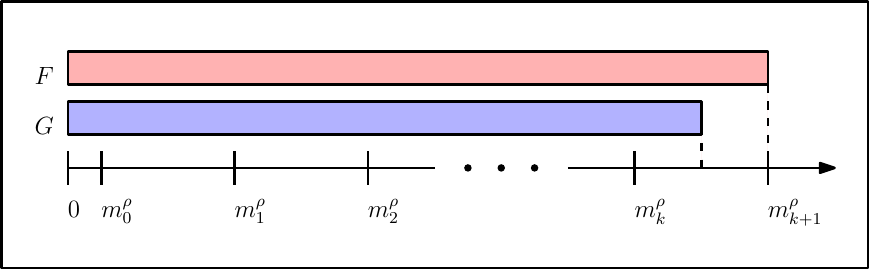}
\textit{ \small Here, $F$ represents a maximally closed set of diameter $k+1$, while $G$ represents a closed (not maximal) set of diameter $k+1.$}
\end{minipage}
\end{center}

By homogeneity we also know that for each maximally closed set $F$ of positive diameter $k$, both the decomposition number $j_F$ as well as the cardinality of the root $|R(F)|$ are uniquely determined by $k$. Hence, we can rename these two numbers as $n^\rho_k$ and $r^\rho_k$ respectively. 
\begin{center}
\begin{minipage}[b]{0.8\linewidth}
\centering
\includegraphics[width=11cm, height=6cm]{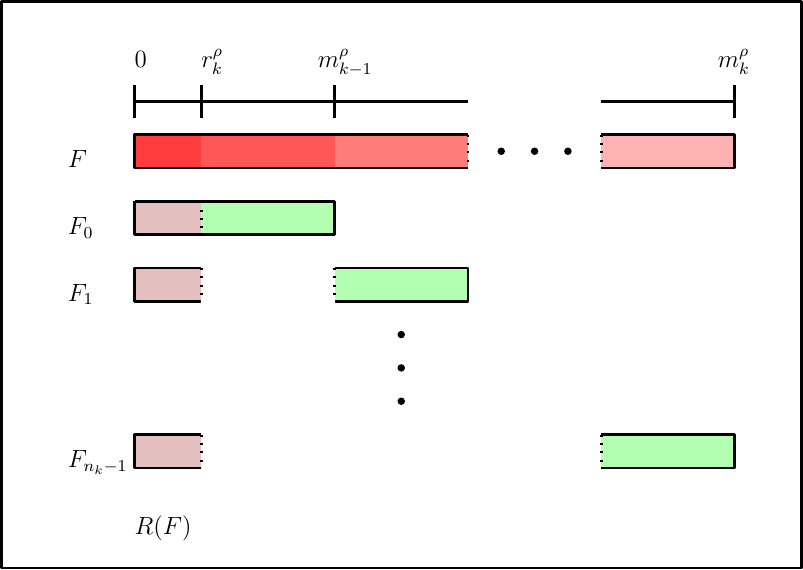}

\textit{\small Here, $F$ represents a maximally closed of diameter $k$. Each piece of its canonical decomposition is painted into two colors. The color brown represents the root of $F$ and the green one represents the remaining part. }
\end{minipage}
\end{center}
From the previous picture it is clear that the numbers $m^\rho_k$, $n^\rho_k$ and $r^\rho_k$ are related through the equation: $$m^\rho_{k+1}=r^\rho_{k+1}+(m^\rho_k-r^\rho_{k+1}) n^\rho_{k+1}.$$
An abstraction of the properties satisfied by the sequence $\langle m^\rho_k,n^\rho_{k+1},r^\rho_{k+1}\rangle_{k\in\omega}$ lead to the next definition.
\begin{definition}[Type]\label{definitiontype}We call a sequence $\tau=\langle m_k,n_{k+1},r_{k+1}\rangle_{k\in\omega}\subseteq \omega^3$ a \textit{type} if:\begin{enumerate}[label=$(\alph*)$]
\item $m_0=1,$ 
\item $\forall k\in \omega\backslash 1 \big(\;n_k\geq 2\;\big),$
\item $\forall k\in \omega \big( \; m_k>r_{k+1}\; \big), $
\item $\forall k\in\omega\backslash 1 \big( \; m_{k+1}=r_{k+1}+(m_k-r_{k+1})n_{k+1}\;\big).$
\end{enumerate}
We say that $\tau$ is \textit{good} if:
\begin{enumerate}
\item[($e$)] $\forall r\in \omega\; \exists^\infty k\in \omega\big(\; r_k=r\;\big), $
\end{enumerate}
Additionally, we say that a partition of $\omega$, namely $\mathcal{P}$, is compatible with $\tau$ if: 
\begin{enumerate}
    \item[($e$')]$\forall P\in \mathcal{P}\,\forall r\in\omega\,\exists^\infty k\in P (\;r_k=r\;).$
\end{enumerate}
\end{definition}

Now, let $\mathcal{F}^\rho$ be the collection of all maximally closed sets and consider it as ordered set with respect to $\subseteq$. First observe that $\mathcal{F}^\rho$ is cofinal in $\text{FIN}(X)$ by Remark \ref{cofinalremarkregular}. Note also that $\mathcal{F}^\rho$ is well-founded since it consists of finite sets. Using Lemma \ref{maximallyclosedregularlemma} one can show recursively that $rank(F)=\rho^F$ for each $F\in \mathcal{F}^\rho$.

We now define a construction scheme as a family of finite sets of $X$ satisfying the same basic properties as $\mathcal{F}^\rho$.\footnote{This approach is not historically accurate. In \cite{schemenonseparablestructures}, construction schemes were defined in the second section, applied in the next three sections, and finally related to ordinal metrics in the sixth one.} But first we make a brief comment regarding the notation.  Each family $\mathcal{F}\subseteq \text{FIN}(X)$ is well-founded with respect to $\subseteq$ (moreover, $Ht(\mathcal{F})\leq \omega)$. In this way, we can define $\mathcal{F}_k$ as in the preliminaries for each $k\in\omega$. That is, $\mathcal{F}_k=\{F\in \mathcal{F}\,:\,rank(F)=k\}$.
\begin{definition}[Construction scheme] 
\label{constructionschemedef}Let $\tau=\langle m_k,n_{k+1},r_{k+1}\rangle_{k\in\omega}$ be a type,
and $X$ be a set of ordinals.
We say that  $\mathcal{F}\subseteq \text{FIN}(X)$ is a \textit{construction scheme} (or simply \textit{a scheme}) over
$X$ of type $\tau$  if:\begin{enumerate}[label=$(\alph*)$]
\item $\mathcal{F}$ is cofinal in $\text{FIN}(X),$
\item $\forall k\in\omega\;\forall F\in \mathcal{F}_k\big(\;|F|=m_k\;\big)$,
\item $\forall k\in\omega\;\forall F,E\in \mathcal{F}_k\big(\; E\cap F\sqsubseteq E,F\;\big)$,
\item $\forall k\in\omega\;\forall F\in \mathcal{F}_{k+1}\;\exists F_0,\dots,F_{n_{k+1}-1}\in \mathcal{F}_k$ such that $$F=\bigcup\limits_{i<n_{k+1}}F_i.$$
Moreover, $\langle F_i\rangle_{i<n_{k+1}}$ forms a $\Delta$-system with root $R(F)$ such that  $|R(F)|=r_{k+1}$ and $$R(F)<F_0\backslash R(F)<\dots < F_{n_{k+1}-1}\backslash R(F).$$
\end{enumerate}
If $n\in \omega$ and $n_{k+1}=n$ for each $k\in\omega$, we will call  $\mathcal{F}$ an \textit{n-construction scheme} (or simply an \textit{n-scheme}).\footnote{2-construction schemes are the objects which we also know as $(\omega,1)$-gap morasses.}
\end{definition}

\begin{rem}\label{remarkcanonicaldecomposition}If $F\in\mathcal{F}_{k+1}$ and $\langle F_i\rangle_{i<n_{k+1}}$ is as in the point (d) of the previous definition, then for each $i<n_{k+1}$,   $$F_i=F[r_{k+1}]\cup F[\,[a_i, a_{i+1})\,]$$
where $a_i=r_{k+1}+i\cdot (m_k-r_{k+1})$. In particular this means that the sequence $\langle F_i\rangle_{i<n_{k+1}}$ is uniquely determined. For this reason, we will call it \textit{the canonical decomposition of $F$.}
    
\end{rem}
By the previous discussion and results in the first section of this chapter we conclude:
\begin{proposition}\label{ordinalmetricimplieschemes}Let $\rho:X^2\longrightarrow \omega$ be a regular, locally finite and homogeneous ordinal metric. The set $$\mathcal{F}^\rho=\{\,F\in \text{FIN}(X)\,:\,F\textit{ is maximally closed}\,\} $$is a construction scheme over $X$ of type $\langle m^\rho_k, n^\rho_{k+1}, r^\rho_{k+1}\rangle_{k\in\omega}$. 
\end{proposition}

The following Theorem was proved in \cite{schemenonseparablestructures}. Another proof can be found in \cite{lopezschemethesis}. We will provide a complete proof of it in the Subsection \ref{subsectionschemeszfc} for the convenience of the reader.
 
\begin{theorem}\label{theoremschemesinzfcstevo} For any good type there is a construction scheme over $\omega_1$ of that type.
\end{theorem}

Our next goal is to prove the converse of Proposition \ref{ordinalmetricimplieschemes}. That is, we will prove that every construction scheme naturally defines a regular, locally finite and homogeneous ordinal metric. Moreover, we will show the construction scheme coincides with the family of maximally closed sets of such metric. In order to do that, we will need some previous results.

The next lemma is an easy consequence of the condition (d) in Definition \ref{constructionschemedef}.
\begin{lemma}\label{firstlemmaconstruction}Let $\mathcal{F}$ be a construction scheme over $X$, $l\in \omega$ and $F\in \mathcal{F}_l$. For each $k\leq l$, it happens that $F=\bigcup\{H\in \mathcal{F}_k\,:\,H\subseteq F\}$.
\end{lemma}
 There is a version of condition condition (c) in Definition \ref{constructionschemedef} in the case where the sets have different rank.
\begin{lemma}
\label{initialsegmentlemma}Let $\mathcal{F}$ be a construcion scheme over $X$ and let $k\leq l\in \omega$. If $G\in \mathcal{F}_k$ and $F\in \mathcal{F}_l$ then $G\cap F\sqsubseteq G$.
\begin{proof}As a consequence of Lemma \ref{firstlemmaconstruction} we know that \begin{align*}G\cap F&=G\cap\big(\bigcup\{\, H\in \mathcal{F}_k\,:\,H\subseteq F\}\big)\\
&=\bigcup\limits_{H\in\mathcal{F}_k\cap \mathscr{P}(F)} G\cap H.
\end{align*}
As each of the sets forming the union of the last equality is an initial segment of $G$ by condition (c) of Definition \ref{constructionschemedef}, such union is also an initial segment of $G$. 
\end{proof}
\end{lemma}
\begin{rem}In general it is not true that if $k<l$, $G\in \mathcal{F}_k$ and $F\in \mathcal{F}_l$ then $G\cap F$ is an initial segment of $F.$ An example of this ocurrs when $G$ is the second piece of the canonical decomposition of $F$, namely $F_1$.
\end{rem}
\begin{proposition}\label{unionlevelscheme}Let $\mathcal{F}$ be a construction scheme over $X$ and let $k\in\omega$ be such that $m_k\leq |X|$. Then $X=\bigcup \mathcal{F}_k$.
\begin{proof}Let $\alpha\in X$ and let $Y$ be a subset of $X$ of cardinality $m_k$. Due to condition (a) of Definition \ref{constructionschemedef} there is $l\in\omega$ and $F\in \mathcal{F}_l$ for which $\{\alpha\}\cup Y\subseteq F$. Note that $l\geq k$ as $|Y|=m_k$. Therefore, by Lemma \ref{firstlemmaconstruction} we conclude  that $F=\bigcup\{G\in \mathcal{F}_k\,:\,G\subseteq F\}$. In particular $\alpha\in G$ for some $G\in \mathcal{F}_k$ with $G\subseteq F$. This finishes the proof. 
\end{proof}
\end{proposition}
\begin{rem}\label{remarksingletons}By applying the previous proposition to the case where $k=0$ we conclude that $[X]^1=\mathcal{F}_0$.
\end{rem}

\begin{definition}\label{defmetricscheme}Let $\mathcal{F}$ be a construction scheme over $X$. We define $\rho_\mathcal{F}:X^2\longrightarrow \omega$ as:
$$\rho_\mathcal{F}(\alpha,\beta)=\min(\,k\in\omega\;:\;\exists F\in\mathcal{F}_k(\;\{\alpha,\beta\}\subseteq F\;)\,).$$
If $\mathcal{F}$ is clear from context, we will write $\rho_\mathcal{F}$ simply as $\rho.$
\end{definition}
\begin{rem}Note that $\rho(\alpha,\beta)$ is well defined since $\mathcal{F}$ is cofinal in $\text{FIN}(X)$.
    
\end{rem}
\begin{lemma}\label{closureschemelemma} Let $\mathcal{F}$ be a construction scheme over $X$. If $\beta\in X$ and $k\in\omega$, then $$(\beta)_k=F\cap (\beta+1)$$
for each $F\in \mathcal{F}_k$ with $\beta\in F.$
\begin{proof}
By Definition \ref{defmetricscheme} it is clear that $F\cap (\beta+1)\subseteq(\beta)_k$. It only remains to prove the other inclusion. For this, let $\alpha\in (\beta)_k$ and $G\in \mathcal{F}_{\rho(\alpha,\beta)}$ be such that $\{\alpha,\beta\}\subseteq G$.  By the definition of the $k$-closure we know that $\rho(\alpha,\beta)\leq k$. Hence, $G\cap F\sqsubseteq G$ due to Lemma \ref{initialsegmentlemma}. As $\beta\in G\cap F$ and $\alpha\leq \beta$ it must happen that $\alpha\in G\cap F$. In particular $\alpha\in F$, so we are done.
\end{proof}
 \end{lemma}
\begin{proposition}Let $\mathcal{F}$ be a construction scheme over $X$. Then $\rho=\rho_\mathcal{F}$ is an ordinal metric.
\begin{proof}The condition (b) of Definition \ref{ordinalmetricdef} is trivially satisfied, the condition (a) is a consequence of the Remark \ref{remarksingletons} and the condition (d) follows from Lemma \ref{closureschemelemma}. Thus, we only need to prove that the condition (c) of such definition also holds. For this, let $\alpha,\beta,\gamma\in X$ be such that $\alpha\leq \beta, \gamma. $ Consider $F\in \mathcal{F}_{\rho(\alpha,\gamma)}$ and $G\in \mathcal{F}_{\rho(\beta,\gamma)}$ for which $\{\alpha,\gamma\}\subseteq F$ and $\{\beta,\gamma\}\subseteq G$. We need to consider two cases.\\\\
\underline{Case 1:} If $\rho(\alpha,\gamma)\leq \rho(\beta,\gamma)$.

\begin{claimproof}[Proof of case.]
Here $F\cap G\sqsubseteq F$ by Lemma \ref{initialsegmentlemma}. As $\gamma\in F\cap G$ and $\alpha\leq \gamma$ we have that $\alpha\in F\cap G$. In particular $\alpha\in G$ so $\{\alpha,\beta\}\subseteq G$. Therefore $\rho(\alpha,\beta)\leq \rho(\beta,\gamma)$.
\end{claimproof}

\noindent \underline{Case 2:} If $\rho(\beta,\gamma)\leq \rho(\alpha,\gamma)$.
\begin{claimproof}[Proof of case.]In this case $F\cap G\sqsubseteq G$. If $\beta\leq \gamma$ then $\beta\in F\cap G$ because $\gamma\in F\cap G$. Therefore $\{\alpha,\beta\}\subseteq G$ which means $\rho(\alpha,\beta)\leq \rho(\beta,\gamma)$. On the other hand, if $\beta\geq \gamma$ we can use Proposition \ref{unionlevelscheme} to pick $F'\in \mathcal{F}_{\rho(\alpha,\gamma)}$ for which $\beta\in F'$. Since $F'\cap (\beta+1)=(\beta)_{\rho(\alpha,\gamma)}$ by Lemma \ref{closureschemelemma} then $\gamma\in F'$. Again, by the same lemma,  $F'\cap (\gamma+1)=(\gamma)_{\rho(\alpha,\gamma)}$ which means that $\alpha\in F'$. Consequently $\{\alpha,\beta\}\subseteq F'$ so $\rho(\alpha,\beta)\leq \rho(\alpha,\gamma)$.
\end{claimproof}
\end{proof}
\end{proposition}
\begin{proposition}\label{rhosdecompositionprop}Let $\mathcal{F}$ be a construction scheme over $X$, $\rho=\rho_\mathcal{F}$ and let $F\in\mathcal{F}_k$ for some $k>0$. Consider $\langle F_i\rangle_{i<n_k}$ the canonical decomposition of $F$ as described in point (d) of Definition \ref{constructionschemedef}. If $\alpha\in F_i\backslash R(F)$ and $\beta\in F_j\backslash R(F)$ for $i<j<n_k$ then $\rho(\alpha,\beta)=k.$
\begin{proof}
First observe that $\rho(\alpha,\beta)\leq k$ as testified by $F.$ Now, $(\beta)_{k-1}=F_j\cap (\beta+1)$ due to the Lemma \ref{closureschemelemma}. In this way, $\alpha\not\in(\beta)_{k-1}$, which means that $\rho(\alpha,\beta)>k-1$. This finishes the proof.
\end{proof}
\end{proposition}
\begin{lemma}\label{diameterscheme}Let $\mathcal{F}$ be a construction scheme over $X$ and $\rho=\rho_\mathcal{F}$. For every $k\in\omega$ and $F\in \mathcal{F}_k$, we have that $\rho^F=k.$ 
\begin{proof}For $k=0$ the result follows from Remark \ref{remarksingletons}, so we will only consider the case where $k>0$. Let $F\in \mathcal{F}_k$.  By Definition \ref{defmetricscheme} it should be clear that $\rho^F\leq k$. To prove the other inequality let $\langle F_i\rangle_{i<n_k}$ be the canonical decomposition of $F$ as described in point (d) of Definition \ref{constructionschemedef}.  Now, consider $\alpha\in F_0\backslash R(F)$ and $\beta\in F_1\backslash R(F)$. By Proposition \ref{rhosdecompositionprop}, $\rho(\alpha,\beta)=k$. Hence  $k\leq \rho(\alpha,\beta)\leq \rho^F$.
\end{proof}
\end{lemma}
\begin{lemma}\label{schememaximallemma}Let $\mathcal{F}$ be a construction scheme over $X$ and $\rho=\rho_\mathcal{F}$. Each element of $\mathcal{F}$ is maximally closed with respect to $\rho$.
\begin{proof}
Let $k\in\omega$ and $F \in \mathcal{F}_k$. Consider $\alpha=\max(F)$. By Lemma \ref{closureschemelemma}, $F=F\cap (\alpha+1)=(\alpha)_k$. In this way $F$ is closed thanks to the point (4) of Proposition \ref{closureprop1}. To see that $F$ is maximally closed let $G\in \text{FIN}(X)$ be a closed set with $\rho^G=k$ and $F\subseteq G$. Consider $\beta=\max(G)$. Through the use of Proposition \ref{unionlevelscheme} we can take $H\in \mathcal{F}_k$ for which $\beta\in H$. Then $(\beta)_k=H\cap (\beta+1)$ due to Lemma \ref{closureschemelemma}. Furthermore, by the point (1) of Proposition \ref{closureprop1} we know that $\rho(\delta,\beta)\leq k$ for each $\delta\in G$. In particular this means that $F\subseteq G\subseteq H$. Since both $F$ and $H$ have size $m_k$, we conclude that $F=H=G$.  This finishes the proof.
\end{proof}
\end{lemma}

\begin{proposition}\label{maximallyclosedschemeprop}Let $\mathcal{F}$ be a construction scheme over $X$ and $\rho=\rho_\mathcal{F}$. Then $$\mathcal{F}_k=\{\,H\in\text{FIN}(X)\,:\,H\textit{ is maximally closed and }\rho^H=k\,\}$$
for each $k\in\omega$.

\begin{proof} The inclusion from left to right is a direct consequence of Lemmas \ref{diameterscheme} and \ref{schememaximallemma}. To prove the one from right to left, let $H\in \text{FIN}(X)$ be a maximally closed set of diameter $k$. Consider $\alpha=\max(H)$ and take $F\in \mathcal{F}_k$ for which $\alpha\in F.$ Then $(\alpha)_k=F\cap (\alpha+1)=F$. This means that $H\subseteq F$ due to the point (1) of Proposition \ref{closureprop1}. Since $\rho^F=k$ by Lemma \ref{diameterscheme} and $H$ is maximally closed, it follows that $H=F$. Thus, the proof is over.
\end{proof}
\end{proposition}
\begin{theorem}Let $\mathcal{F}$ be a construction scheme over $X$. Then $\rho=\rho_\mathcal{F}$ is locally finite, regular and homogeneous.

\begin{proof}The function $\rho$ is locally finite because each closed set $H$ is contained in $F\cap(\max(H)+1)$ where $F\in \mathcal{F}_{\rho^H}$ is such that $\max(H)\in F$. In particular, this means that $|H|$ is bounded by $m_{\rho^H}$. Now, $\rho$ is regular due to Proposition \ref{maximallyclosedschemeprop} and the point (d) of Definition \ref{constructionschemedef}. The only thing left to do is to prove that $\rho$ is homogeneous. This proof is carried by induction over the diameter $k$ of the maximally closed sets involved. The case where $k=0$ is trivial as maximally closed sets of diameter $0$ are singletons. So suppose that $k\in \omega$, and we have proved homogeneity for all maximally closed sets of diameter less than $k$. Let $F,G$ be two maximally closed sets with $\rho^F=\rho^G=k$.  By Proposition \ref{maximallyclosedschemeprop}, $F,G\in \mathcal{F}_k$. Consequently $|F|=|G|=m_k$. Thus, we can let $h:F\longrightarrow G$ be the only increasing bijection. We claim that $\rho(\alpha,\beta)=\rho(h(\alpha),h(\beta))$ for each $\alpha,\beta\in F$. For this, take $\langle F_i\rangle_{i<n_k}$ and $\langle G_i \rangle_{i<n_k}$ the canonical decompositions of $F$ and $G$ respectively. It is easy to see that $h[F_i]=G_i$ for each $i<n_k$ and that $h[R(F)]=R(G)$. Consequently, if $\alpha,\beta\in F_i$ for some $i<n_k$ we can use the induction hypothesis over $F_i$ and $G_i$ to conclude that $\rho(\alpha,\beta)=\rho(h(\alpha),h(\beta))$. On the other hand, if there is no $i$ for which both $\alpha$ and $\beta$ belong to $F_i$, there are distinct $i,j<n_k$ for which $\alpha\in F_i\backslash R(F)$ and $\beta\in F_j\backslash R(F)$. Observe that $h(\alpha)\in G_i\backslash R(G)$ and $h(\beta)\in G_j\backslash R(F)$. Therefore, by Proposition \ref{rhosdecompositionprop} we have that $\rho(\alpha,\beta)=k=\rho(h(\alpha),h(\beta))$. This finishes the proof. 
\end{proof}
\end{theorem}

The following corollary follows from the proof of the previous theorem as well as of Proposition 
\ref{ordinalmetricimplieschemes}.

\begin{corollary} If $\mathcal{F}$ is a construction scheme over $X$ and $\rho:X^2\longrightarrow \omega$ is a locally finite, homogeneous and regular ordinal metric, then:
\begin{itemize}
    \item $\mathcal{F}^{\rho_\mathcal{F}}=\mathcal{F}$,
    \item $\rho_{\mathcal{F}^\rho}=\rho.$
\end{itemize}
\end{corollary}
\section{Two more canonical functions}
Apart from the ordinal metric associated to a construction scheme, there are two important functions that we need to analyze before we enter the world of applications. The first one being the $\Delta$ function and the second one being the $\Xi$ function.

For the rest of this section, fix a construction scheme $\mathcal{F}$ over $X$ and let $\rho=\rho_\mathcal{F}$ be its associated ordinal metric.\\

The $\Delta$ function just measures the least moment in which the $k$-cardinality of two ordinals differs. 

\begin{definition}[The $\Delta$ function] We define $\Delta:X^2\longrightarrow \omega+1$ as:
$$\Delta(\alpha,\beta)=\begin{cases} \min(\,k\in\omega\,:\,\lVert \alpha\rVert_k\not=\lVert \beta\rVert_k\,)&\textit{ if }\alpha\not=\beta\\
\omega&\textit{ if }\alpha=\beta
\end{cases}$$
\end{definition}
\begin{rem}$\Delta$ is well defined since $\lVert\alpha\rVert_{\rho(\alpha,\beta)}\not=\lVert\beta\rVert_{\rho(\alpha,\beta)}$ whenever $\alpha\not=\beta.$ Moreover, $\Delta(\alpha,\beta)\leq \rho(\alpha,\beta)$.
\end{rem}
\begin{lemma}\label{lemmadelta1}Let $\alpha,\beta\in X$ be distinct ordinals and $k\in\omega$. If $\lVert\alpha\rVert_k=\lVert\beta\rVert_k$ then $k<\Delta(\alpha,\beta)$.\begin{proof}By Corollary \ref{ballhomogeneitylemma} $(\alpha)_k$ is $\rho$-isomorphic to $(\beta)_k$. Let $h:(\alpha)_k\longrightarrow (\beta)_k$ be the only increasing bijection. Note that $h(\alpha)=\beta$. Now consider an arbitrary $i\leq k$. Then $(\alpha)_i\subseteq (\alpha)_k$. In this way, $h[(\alpha)_i]=(h(\alpha))_i=(\beta)_i$. This means that $\lVert \alpha\rVert_i=\lVert\beta\rVert_i$. Therefore $k<\Delta(\alpha,\beta)$.
 \end{proof}
\end{lemma}

\begin{lemma}\label{countrymanlemma3} Let $\alpha,\beta,\delta\in X$ be distinct ordinals such that $\Delta(\alpha,\beta)< \Delta(\beta,\delta)$. Then $\Delta(\alpha,\delta)=\Delta(\alpha,\beta)$.
\begin{proof}$\Delta(\alpha,\delta)\leq \Delta(\alpha,\beta)$ since $\lVert \alpha\rVert_{\Delta(\alpha,\beta)}\not=\lVert \beta\rVert_{\Delta(\alpha,\beta)}= \lVert \delta \rVert_{\Delta(\alpha,\beta)}$. Now, to see that $\Delta(\alpha,\delta)\geq \Delta(\alpha,\beta)$ take an arbitrary $i<\Delta(\alpha,\beta)$.  Then $\lVert \alpha\rVert_i=\lVert \beta \rVert_i=\lVert \delta\rVert_i$.  This finishes the proof.
\end{proof}
\end{lemma}

\begin{lemma}\label{lemmainequalitiesdelta}Let $\alpha,\beta\in X$ be distinct ordinals, $k<\Delta(\alpha,\beta)$ and $h:(\alpha)_k\longrightarrow (\beta)_k$ be the only increasing bijection.
\begin{multicols}{2}
If $\delta\in (\alpha)_k$ then:
\begin{enumerate}[label=$(\alph*)$]
\item $\Delta(\delta,h(\delta))\geq \Delta(\alpha,\beta)$,
\item $\rho(\alpha,\beta)\geq \rho(\delta,h(\delta))$.
\end{enumerate}
\columnbreak
In other words, if $i<\lVert\alpha\rVert_k$ then:
\begin{enumerate}[label=$(\alph*)$]
\item $\Delta((\alpha)_k(i),(\beta)_k(i))\geq \Delta(\alpha,\beta)$,
\item $\rho(\alpha,\beta)\geq \rho((\alpha)_k(i),(\beta)_k(i))$.
\end{enumerate}
\end{multicols}
\begin{proof}\begin{claimproof}[Proof of $(a)$] Let $l=\Delta(\alpha,\beta)-1$ and $h':(\alpha)_l\longrightarrow (\beta)_l$ be the only increasing bijection. As $k<\Delta(\alpha,\beta)$ then $(\alpha)_k\subseteq (\alpha)_l$ and $h'[(\alpha)_k]=(\beta)_k$. Hence $h'|_{(\alpha)_k}=h$ so $h'(\delta)=h(\delta)$. Therefore $h'[(\delta)_l]=(h(\delta))_l$. In this way, $\lVert \delta\rVert_l=\lVert h(\delta)\rVert_l$. Thus, by Lemma \ref{lemmadelta1} we conclude that $\Delta(\alpha,\beta)\leq \Delta(\delta,\phi(\delta))$. 
\end{claimproof}
\begin{claimproof}[Proof of $(b)$] Let $F\in \mathcal{F}_{\rho(\alpha,\beta)}$ be such that $\alpha,\beta\in F$. Since $k<\Delta(\alpha,\beta)\leq \rho(\alpha,\beta)$ then $(\alpha)_k\cup (\beta)_k\subseteq F$. In particular $\delta,h(\delta)\in F$ so $\rho(\alpha,\beta)\geq \rho(\delta,h(\delta)).$ 
\end{claimproof}
\end{proof}
\end{lemma}

\begin{rem}\label{remarkbijectionidentityballs}Recall that if $\alpha$ and $\beta$ are so that $\lVert \alpha\rVert_k=\lVert \beta\rVert_k$ and $h:(\alpha)_k\longrightarrow (\beta)_k$ is the only increasing bijection, then the following happens:
\begin{multicols}{2}
\begin{itemize}
\item $h(\alpha)=h(\beta),$
\item $(\alpha)_k\cap (\beta)_k\sqsubseteq (\alpha)_k,(\beta)_k,$
\item $h|_{(\alpha)_k\cap (\beta)_k}$ is the identity function,
\item $h[(\alpha)_k\backslash (\beta)_k]=(\beta)_k\backslash (\alpha)_k.$
\end{itemize}
\columnbreak
\hspace{1cm}
\begin{minipage}[t]{0.4\linewidth}
\centering
\includegraphics[width=6cm, height=3cm]
{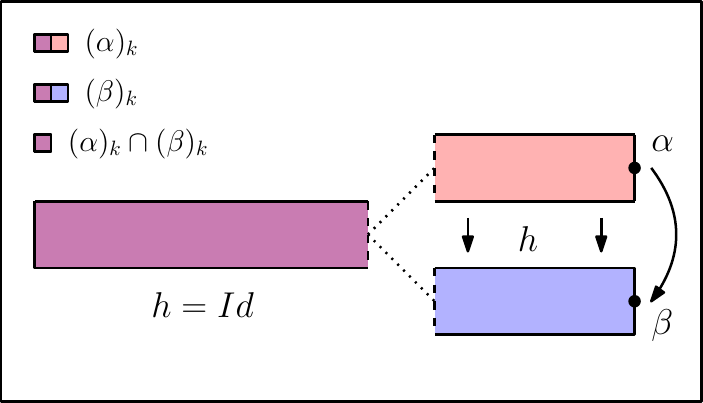}

\end{minipage}
\end{multicols}
\end{rem}
\begin{lemma}\label{lemmahdeltarhoinequalities}Let $\alpha,\beta\in X$ be distinct ordinals and $k<\Delta(\alpha,\beta)$. If $h:(\alpha)_k\longrightarrow (\beta)_k$ is the only increasing bijection and $\delta,\gamma\in (\alpha)_k$ are such that $\delta\leq \gamma$ and $h(\delta)\not=\delta$ then the following happens:
\begin{enumerate}[label=$(\alph*)$]
\item$h(\gamma)\not=\gamma,$
\item$\rho(\alpha,\beta)\geq \rho(\gamma,h(\gamma))\geq \rho(\delta,h(\delta))\geq\Delta(\delta,h(\delta)) \geq \Delta(\gamma,h(\gamma))\geq \Delta(\alpha,\beta).$
\end{enumerate}
\begin{proof} Remember that $C=(\alpha)_k\cap (\beta)_k$ is an initial segment of both $(\alpha)_k$ and $(\beta)_k$. In this way, it is easy to see that $h(\xi)=\xi$ if and only if $\xi \in C.$ Note that since $\delta\not\in C$ and $\gamma >\delta$, then $\gamma\notin C$. Therefore $h(\gamma)\not=\gamma$. This proves the point (a).

The inequalities from (b) follow directly from Lemma \ref{lemmainequalitiesdelta} and the fact that $h|_{(\gamma)_k}$ is the only increasing bijection from $(\gamma)_k$ to $(h(\gamma))_k$.

\end{proof}
\end{lemma}
\begin{definition}[The $\Xi$ function]\label{Xifunction} Let $\alpha\in X$.
$\Xi_\alpha:\omega\longrightarrow \omega\cup\{-1\}$ is the function defined as:

$$\Xi_\alpha(k)=\begin{cases}0 &\textit{if }k=0\\
-1 &\textit{if }k>0\textit{ and }\lVert\alpha\rVert_k\leq r_k\\
\frac{\lVert\alpha\rVert_k-\lVert\alpha\rVert_{k-1}}{m_{k-1}-r_k} &\textit{otherwise}
\end{cases}$$
It is not hard to check that if $k\in \omega\backslash 1$ and $F\in \mathcal{F}_k$ is such that $\alpha\in F$, then:

$$\Xi_\alpha(k)=\begin{cases}-1&\textit{ if }\alpha\in R(F)\\
i&\textit{ if }\alpha\in F_i\backslash R(F)
    
\end{cases}$$
\end{definition}
\begin{multicols}{2}

The main reason for defining the $\Xi$ function is that without appealing to any extra axioms it is really hard to give useful properties regarding the behavior of the $k$-cardinality of pairs of ordinals. The $\Xi$ function reveals two important and natural critical points \say{leading the dance} of the $k$-cardinality of ordinals $\alpha$ and $\beta$ as $k$ grows larger. As the next lemma shows, below $\Delta(\alpha,\beta)$ and above $\rho(\alpha,\beta)$, this dance is smooth and pleasant whereas between $\Delta(\alpha,\beta)$ and $\rho(\alpha,\beta)$ it seems that almost anything is possible.
\columnbreak
\vspace{1cm}
\begin{center}
\begin{minipage}[b]{0.8\linewidth}
\centering
\includegraphics[width=6.5cm, height=5.5cm]{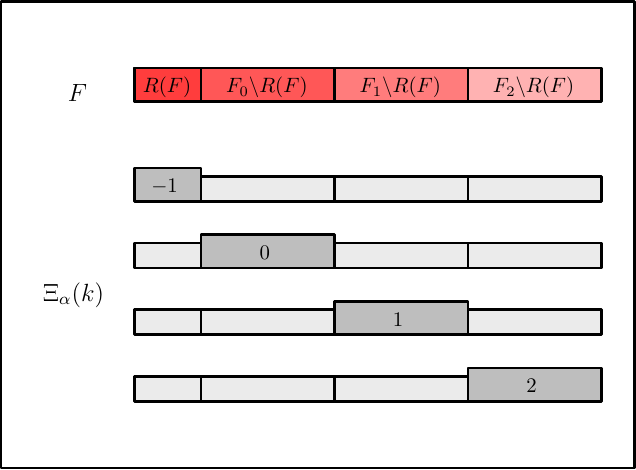}\\
\textit{\small In here, $F$ represents a member of $F_k$ with $n_k=3$.}
\end{minipage}
\end{center}
\end{multicols}

\begin{lemma}\label{lemmaxi}Let  $\alpha<\beta\in X$ and $k\in \omega\backslash 1$. Then:
\begin{enumerate}[label=$(\alph*)$]
\item If $k<\Delta(\alpha,\beta)$, then  $\Xi_\alpha(k)=\Xi_\beta(k).$
\item If $k=\rho(\alpha,\beta)$, then $0\leq \Xi_\alpha(k)<\Xi_\beta(k).$
\item If $k>\rho(\alpha,\beta)$, then either $\Xi_\alpha(k)=-1$ or $\Xi_\alpha(k)=\Xi_\beta(k).$ 
\item If $k=\Delta(\alpha,\beta)$ then $0\leq \Xi_\alpha(k)\not=\Xi_\beta(k)\geq 0.$ 
\end{enumerate}
\begin{proof}
(a) is direct from the definition of $\Delta(\alpha,\beta)$ so we will only prove the remaining points.

\begin{claimproof}[Proof of $(b)$.] Fix $F\in \mathcal{F}_{\rho(\alpha,\beta)}$ for which $\alpha,\beta\in F$. We know there are $i,j<n_k$ such that $\alpha\in F_i$ and $\beta\in F_j$. Since $F_i,F_j\in \mathcal{F}_{\rho(\alpha,\beta)-1}$, minimality of $\rho$ implies that $\alpha\not\in F_j$ and $\beta\not\in F_i$. In particular, this means $i\not=j $ and $\alpha,\beta\not\in R(F)$. Thus, by Definition \ref{Xifunction} we have that $\Xi_\alpha(k)=i$ and $\Xi_\beta(k)=j$. Moreover, since $\alpha<\beta$ then $i<j$.
\end{claimproof}
\begin{claimproof}[Proof of $(c)$.]For this, suppose that $\Xi_\alpha(k)\not=-1$ and let $F\in \mathcal{F}_k$ be such that $\alpha,\beta\in F$. Since $\alpha<\beta$, we also have $\Xi_\beta(k)\not=-1$. Thus, there is $j<n_k$ for which $\beta\in F_j\backslash R(F) $. Recall that $F_j\in \mathcal{F}_{k-1}$. In this way, $\alpha\in (\beta)_{k-1}=(\beta+1)\cap F_j$. As $\alpha\not\in R(F)$, this means  $\Xi_\alpha(k)=j=\Xi_\beta(k)$. So we are done.
\end{claimproof}
\begin{claimproof}[Proof of $(d)$.] Let $F,G\in \mathcal{F}_k$ be such that $\alpha\in F$ and $\beta\in G$. Also, let $h:F\longrightarrow G$ be the increasing bijection. Then $h[(\alpha)_i]=(h(\alpha))_i$ for any $i\leq k$. In this way, $$k=\Delta(\alpha,\beta)=\Delta(h(\alpha),\beta)\leq \rho(h(\alpha),\beta)\leq k.$$
That is, $k=\rho(h(\alpha),\beta)$. So by the part (b) of this lemma we conclude that $\Xi_{h(\alpha)}(k)$ and $\Xi_\beta(k)$ are both distinct and greater or equal to $0$. To finish just note that $\Xi_\alpha(k)=\Xi_{h(\alpha)}(k)$.
\end{claimproof}
\end{proof}
\end{lemma}
As an easy consequence of the previous lemma we have the following corollary. 
\begin{corollary} Suppose that $X=\omega_1$ and let $\alpha\in X$. Then:
\begin{enumerate}[label=$(\arabic*)$]
\item $\Xi_\alpha(k)=-1$ for infinitely many $k\in \omega$.
\item If $\alpha$ is infinite, then $\Xi_\alpha(k)\geq 1$ for infinitely many $k\in\omega$.
\item If $\tau$ (the type of $\mathcal{F}$) is a good type, then $\Xi_\alpha(k)\geq 0$ for infinitely many $k\in\omega$.
\end{enumerate}
\begin{proof}In order to prove (3) just observe that since $\tau$ is a good type then $r_k=0$ for infinitely many $k$'s. For any such $k$ it is necessarily true that $\Xi_\alpha(k)\geq 0$. We now prove the remaining points.
\begin{claimproof}[Proof of $(1)$] Let $l\in\omega$. We will find $k>l$ for which $\Xi_\alpha(k)=-1$.  Since $X=\omega_1$ there is an uncountable $S\subseteq \omega_1\backslash \alpha$ and $s\in \omega$ such that $\rho(\alpha,\beta)=\rho(\alpha,\delta)=s$ for any $\beta,\delta\in S$. In virtue of Lemma \ref{regularimpliesunbounded} there are $\beta<\delta\in S$ for which $\rho(\beta,\delta)>\max( l,s)$. Let $k=\rho(\beta,\delta)$. By the point (b) in Lemma \ref{lemmaxi} we know that $0\leq\Xi_\beta(k)<\Xi_\delta(k)$. If $\Xi_\alpha(k)\not=-1$ then we would have that $\Xi_\beta(k)=\Xi_\alpha(k)=\Xi_\delta(k)$ due to part (c) of that same lemma, which is impossible. Hence, we are done.
\end{claimproof}
\begin{claimproof}[Proof of $(2)$] Fix $l\in\omega$. We shall find $k>l$ for which $\Xi_\alpha(k)\geq 1$. Since  $(\alpha)_l$ is finite and $\alpha$ is infinite there is $\xi\in \alpha\backslash (\alpha)_l$. By the definition of the $l$-closure we have that $k=\rho(\xi,\alpha)>l$. Then $\Xi_\alpha(k)\geq 1$ do to the point (b) of Lemma \ref{lemmaxi}.
\end{claimproof}

\end{proof}
\end{corollary}

\begin{lemma}\label{lemmainequalitiesdeltatwo}Let $\alpha,\beta\in X$ be distinct ordinals, $k<\Delta(\alpha,\beta)$ and $h:(\alpha)_k\longrightarrow (\beta)_k$ be the only increasing bijection. If $\gamma\in (\alpha)_k$ is such that $\gamma\not=h(\gamma)$ then the following are equivalent:
\begin{enumerate}[label=$(\alph*)$]
\item $\Delta(\gamma,h(\gamma))>\Delta(\alpha,\beta)$.
\item $\Xi_\gamma(\Delta(\alpha,\beta))=-1$.
\end{enumerate}
Furthermore, if $\Xi_\gamma(\Delta(\alpha,\beta)\geq 0$ then $\Xi_\gamma(\Delta(\alpha,\beta))=\Xi_\alpha(\Delta(\alpha,\beta))$ and $\Xi_{h(\gamma)}(\Delta(\alpha,\beta))=\Xi_\beta(\Delta(\alpha,\beta)).$
\begin{proof}
\begin{claimproof}[Proof of $(a)\Rightarrow (b)$.] We argue by contradiction. Suppose that $\Xi_\gamma(\Delta(\alpha,\beta))\geq 0$. Since $\rho(\alpha,\gamma)\leq k<\Delta(\alpha,\beta)$ then $\Xi_\gamma(\Delta(\alpha,\beta))=\Xi_\alpha(\Delta(\alpha,\beta))\geq 0$ by the part (c) of Lemma \ref{lemmaxi}. Now, since we are assuming that $\Delta(\gamma,h(\gamma))>\Delta(\alpha,\beta)$, then $0\leq \Xi_\gamma(\Delta(\alpha,\beta))=\Xi_{h(\gamma)}(\Delta(\alpha,\beta))$ due to the point (a) of Lemma \ref{lemmaxi}. Thus, we can argue in the same way as before to conclude that $\Xi_{h(\gamma)}(\Delta(\alpha,\beta))=\Xi_\beta(\Delta(\alpha,\beta))$. Therefore, according the part (d) of the same lemma, $$\Xi_\gamma(\Delta(\alpha,\beta))=\Xi_\alpha(\Delta(\alpha,\beta))\not=\Xi_\beta(\Delta(\alpha,\beta))=\Xi_{h(\gamma)}(\Delta(\alpha,\beta)).$$
We conclude using the part (a) of Lemma \ref{lemmaxi} that $\Delta(\gamma,h(\gamma))=\Delta(\alpha,\beta)$.
\end{claimproof}
\begin{claimproof}[Proof of $(b)\Rightarrow (a)$.] Suppose that $\Xi_\gamma(\Delta(\alpha,\beta))=-1$. Then the conclusion of the part (d) of Lemma \ref{lemmaxi} can not hold when applied to $\gamma$, $h(\gamma)$ and $\Delta(\alpha,\beta).$ In virtue of this, $\Delta(\alpha,\beta)$ must be distinct from $\Delta(\gamma,h(\gamma)).$
\end{claimproof}
\end{proof}
\end{lemma}
\begin{corollary}Let $\xi<\alpha<\beta\in X$. If $\rho(\xi,\beta)<\rho(\alpha,\beta)$ then $\rho(\xi,\alpha)<\rho(\alpha,\beta)$.
\begin{proof}Since $\rho$ is an ordinal metric, $\rho(\xi,\alpha)\leq \max(\rho(\xi,\beta),\rho(\alpha,\beta))=\rho(\alpha,\beta)$. Suppose towards a contradiction that $\rho(\alpha,\beta)=\rho(\xi,\alpha)$. Then by the part (b) of Lemma \ref{lemmaxi} we have that $$0\leq \Xi_\xi(\rho(\alpha,\beta))<\Xi_\alpha(\rho(\alpha,\beta))<\Xi_\beta(\rho(\alpha,\beta)).$$
On the other hand, since $\rho(\alpha,\beta)>\rho(\beta,\xi)$ then $\Xi_\xi(\rho(\alpha,\beta))=-1$ or $\Xi_\xi(\rho(\alpha,\beta))=\Xi_\beta(\rho(\alpha,\beta))$ due to the part (c) of Lemma \ref{lemmaxi}. Both cases contradict the previous inequality so we are done.
\end{proof}
\end{corollary}
The definition of the $\Xi$ function can be extended (with some restrictions) to arbitrary elements of $\text{FIN}(X)$. More precisely, by virtue of Proposition \ref{propnoname} we know that:
\begin{lemma}Let $A\in \text{FIN}(X)$, $k>\rho^A$ and $F\in \mathcal{F}_k$ such that $A\subseteq F_i$. Then there is $i<n_k$ for which $A\subseteq F$. Furthermore, it is easy to see that if $A\not\subseteq R(F)$ then this $i$ is unique and does not depend on the choice of $F$. 
\end{lemma}
This leads to the following definition.
\begin{definition}[set-valued $\Xi$ function]Let $A\in \text{FIN}(X)$ for each $k>\rho^A$ we define $$\Xi_A(k)=\begin{cases}-1&\textit{ if }A\subseteq R(F)\\
i&\textit{ if }A\not\subseteq R(F)\textit{ and }A\subseteq F_i
\end{cases}$$
Here, $F\in \mathcal{F}_k$ is such that $A\subseteq F.$ 
\end{definition}
\begin{rem}\label{remxisetvalued}It is not hard to check that actually $\Xi_A(k)=\Xi_{\max(A)}(k)$  for each $k> \rho^A$.    
\end{rem}

The following lemma is easy to prove and it is left to the reader as it is just a variation of  Lemma 
\ref{lemmaxi}.
\begin{lemma}\label{xilemmaset}For any distinct $A,B\in \text{FIN}(X)$ and $k\geq\rho^{A\cup B}$ the following happens:
\begin{enumerate}[label=$(\alph*)$]

\item If $0\leq \Xi_A(k)\not=\Xi_B(k)\geq 0$ then $k=\rho^{A\cup B}.$
\item If $k>\rho^{A\cup B}$ then $\Xi_A(k)=\Xi_B(k)$ or $\min(\Xi_A(k),\Xi_B(k))=-1$.
\end{enumerate}
\end{lemma}
\section{Capturing Axioms}
Suppose that $A$ and $B$ are two different elements of $\text{FIN}(X)$ and $k$ is a natural number bigger (or equal) than the diameter of both $A$ and $B$. Then $(A)_k$ and $(B)_k$ are closed sets by the part (4) of Proposition \ref{closureprop1}. By Lemma \ref{homoegeneitylemma} and Theorem \ref{closedsizetheorem} we know that these two sets have the same size if and only if they are $\rho$-isomorphic. However, if this situation occurs, nothing assures that if we take $h:(A)_k\longrightarrow (B)_k$ to be the increasing bijection then $h[A]=B$. In particular, if $k\geq\rho^{A\cup B}$, then either $|(A)_k|\not=|(B)_k|$ (so these closures are not even $\rho$-isomorphic) or $(A)_k=(B)_k$\footnote{Since $k\geq \rho^{A\cup B}\geq \max(\rho^A,\rho^B)$, then $(A)_k=(\max(A))_k$ and $(B)_k=(\max(B))_k$.}. Note that in the latter case, the increasing bijection $h$ between $(A)_k$ and $(B)_k$ is the identity function which means that $h[A]\not= B$.  In virtue of this observation, the biggest $k$ for which we can expect $(A)_k$ and $(B)_k$ to satisfy this stronger isomorphism condition is $\rho^{A\cup B}-1$. 

\begin{definition}Let $A,B\in FIN(X)$ be such that $\rho^A,\rho^B<\rho^{A\cup B}$. We say that $A$ and $B$ are \textit{strongly $\rho$-isomorphic} if for $l=\rho^{A\cup B}$ it happens that:
\begin{itemize}
    \item $|(A)_{l-1}|=|(B)_{l-1}|$. That is, $(A)_{l-1}$ is $\rho$-isomorphic to $(B)_{l-1}$,
    \item $h[A]=B$ where $h$ is the increasing bijection from $(A)_{l-1}$ to $(B)_{l-1}$.
\end{itemize}   
\end{definition}
The proof of the following lemma is easy.
\begin{lemma}If $A$ and $B$ are strongly $\rho$-isomorphic, then $\rho^A=\rho^B$ and $(A)_k$ is isomorphic to $(B)_k$ for each $k<\rho^{A\cup B}.$
\end{lemma}
The following proposition gives us a better picture of how two strongly $\rho$-isomorphic sets look inside the closure of their union. 
\begin{proposition}\label{equivalencestronglyisomorphic} Let $A,B\in \text{FIN}(X)$ be such that $\max(\rho^A,\rho^B)<\rho^{A\cup B}=l$. Then $A$ and $B$ are strongly $\rho$-isomorphic if and only if  for any $F\in \mathcal{F}_l$ with $A\cup B\subseteq F$ the following happens:
\begin{itemize}
\item $\Xi_A(l),\Xi_B(l)\not=-1$.

\item $h[A]=B$ where $h$ is the increasing bijection from $F_{\Xi_A(l)}$ to $F_{\Xi_{A}(l)}$.

\end{itemize}
\begin{proof}\begin{claimproof}[Proof of $\Leftarrow$.] Let $F\in \mathcal{F}_l$ be such that $A\cup B\subseteq F$. Note that $(A)_{l-1}\sqsubseteq F_{\Xi_A(l)}$ and $(B)_{l-1}\sqsubseteq F_{\Xi_B(l)} $. By this and since $|(A)_{l-1}|=|(B)_{l-1}|$ we have that $h[(A)_{l-1}]=(B)_{l-1}$ where $h:F_{\Xi_A(l)}\longrightarrow F_{\Xi_B(l)}$ is the increasing bijection. Hence, the restriction of $h|_{(A)_{l-1}}$ is the increasing bijection from $(A)_{l-1}$ to $(B)_{l-1}$. In this way we conclude that $h[A]=B $. Now suppose towards a contradiction that $\min(\Xi_A(l),\Xi_B(l))=-1$. Without loss of generality we can assume that $\Xi_A(l)=-1$. That is, $A\subseteq R(F)$. Since $h|_{R(F)}$ is the identity function this means that $A=B$. This is a contradiction to the fact that $\max(\rho^A,\rho^B)<\rho^{A\cup B}$. We conclude that $\min(\Xi_A(l),\Xi_B(l))\geq 0$ so we are done.
\end{claimproof}

\begin{claimproof}[Proof of $\Rightarrow$.] For this, let $F\in \mathcal{F}_l$ be such that $A\cup B\subseteq F$ and $h':F_{\Xi_A(l)}\longrightarrow F_{\Xi_B(l)}$ be the increasing bijection. Since $h'[A]=B$ it is easy to see that $h'[(A)_{l-1}]=(B)_{l-1}$.  In this way, $|(A)_{l-1}|=|(B)_{l-1}|$. To finish just note that $h=h'|_{(A)_{l-1}}$ is the increasing bijection to $(B)_{l-1}$ so $h[A]=B$.
\end{claimproof}

\end{proof}
\end{proposition}
\begin{rem}\label{remarkcapturing}The two conditions imposed to $F$ in the previous proposition are equivalent to the existence of distinct $i,j<n_l$ (namely, $\Xi_A(l)$ and $\Xi_B l$) and a unique $S\subseteq m_{l-1}$ such that $F_i[S]=A$ and $F_j[S]=B$ (such $S$ being $F^{-1}_{\Xi_A(l)}[A]$). 
\end{rem}

\begin{lemma}\label{intersectionrhoisomorphiclemma}Let $A,B\in \text{FIN}(X)$ be strongly $\rho$-isomorphic and let $l=\rho^{A\cup B}$. If $F\in \mathcal{F}_l$ is such that $A\cup B\subseteq F$, then $A\cap B=A\cap R(F)$. In particular, $A\cap B\sqsubseteq A$ and $\Xi_{A\cap B}(l)=-1$. 
\begin{proof}
Let $h:F_{\Xi_A(l)}\longrightarrow F_{\Xi_B(l)}$ be the increasing bijection. The inclusion from left to right is clear since $A\cap B=A \cap (A\cap B)\subseteq A\cap (F_{\Xi_A(l)}\cap F_{\Xi_B(l)})=A\cap R(F)$. To prove the other one just note that $h|_{R(F)}$ is the identity function. Hence, if $\alpha\in A\cap R(F)$ then $\alpha=h(\alpha)\in B$ due to the second point of Proposition \ref{equivalencestronglyisomorphic}.

\end{proof}
\end{lemma}
As a corollary we have:
\begin{corollary}\label{corollarydisjointrhoisomorphic}Let $A,B\in \text{FIN}(X)$ be strongly $\rho$-isomorphic and let $l=\rho^{A\cup B}$. If $F\in \mathcal{F}_l$ is such that $A\cup B\subseteq F$ and $A\cap B=\emptyset$ then $A\subseteq F_{\Xi_A(l)}\backslash R(F)$ and $B\subseteq F_{\Xi_A(l)}\backslash R(F)$. In particular, $\Xi_\alpha(l)=\Xi_A(l)$ and $\Xi_\beta(l)=\Xi_B(l)$ for any $\alpha\in A$ and $\beta \in B$.
\end{corollary}
\begin{rem}\label{remarkcapturing2}In virtue of Proposition \ref{equivalencestronglyisomorphic}, if $A=\{\alpha\}$ and $B=\{\beta\}$ for $\alpha$ and $\beta$ distinct ordinals, then $A$ and $B$ are strongly $\rho$-isomorphic if and only if $\Delta(\alpha,\beta)=\rho(\alpha,\beta)$. 
\end{rem}

 \begin{definition}[Captured families]\label{capturedfamiliesdef} Let $\mathcal{D}$ be a finite subset of $\text{FIN}(X)$ and $l\in\omega$. We say that that $\mathcal{D}$ is \textit{captured} at level $l$ if $2\leq|\mathcal{D}|\leq n_l$ and:
 \begin{enumerate}[label=$(\arabic*)$]
 \item $l=\rho^{\cup \mathcal{D}}>\rho^D$ for each $D\in \mathcal{D}$,
 \item $\{\Xi_D(l)\,:\, D\in \mathcal{D}\}=\{0,\dots,|\mathcal{D}|-1\}=|\mathcal{D}|$.\footnote{by Lemma \ref{xilemmaset}, this condition implies that $\rho^{E\cup D}=l$ for any distinct $D,E\in \mathcal{D}$.}
 \item For any distinct $D, E\in \mathcal{D}$, $D$ is strongly $\rho$-isomorphic to $E$.
 \end{enumerate}

 Equivalently, for any $F\in \mathcal{F}_l$ with $\bigcup \mathcal{D}\subseteq F$ there is $S\subseteq m_{l-1}$ such that $$\mathcal{D}=\{F_i[S]\,:\,i<|\mathcal{D}|\}.$$
Additionally, if $|\mathcal{D}|=n_l$ we  say that $\mathcal{D}$ is \textit{fully captured }at level $l$. Whenever $D\in FIN(X)$, we say that $D$ is captured (resp. fully captured) in case $\{\{\alpha\}\,:\,\alpha\in D\}$ is captured (resp. fully captured). Finally, if $\mathcal{D}$ is  captured at level $l$ and we write $\mathcal{D}$ as a list of elements, say $\mathcal{D}=\{D_0,D_1,\dots, D_{n-1}\}$,  we always assume that $\Xi_{D_i}(l)=i$ for each $i<n$.\end{definition}
\begin{center}
\begin{minipage}[b]{0.8\linewidth}
\centering
\includegraphics[width=12 cm, height=5.5cm]{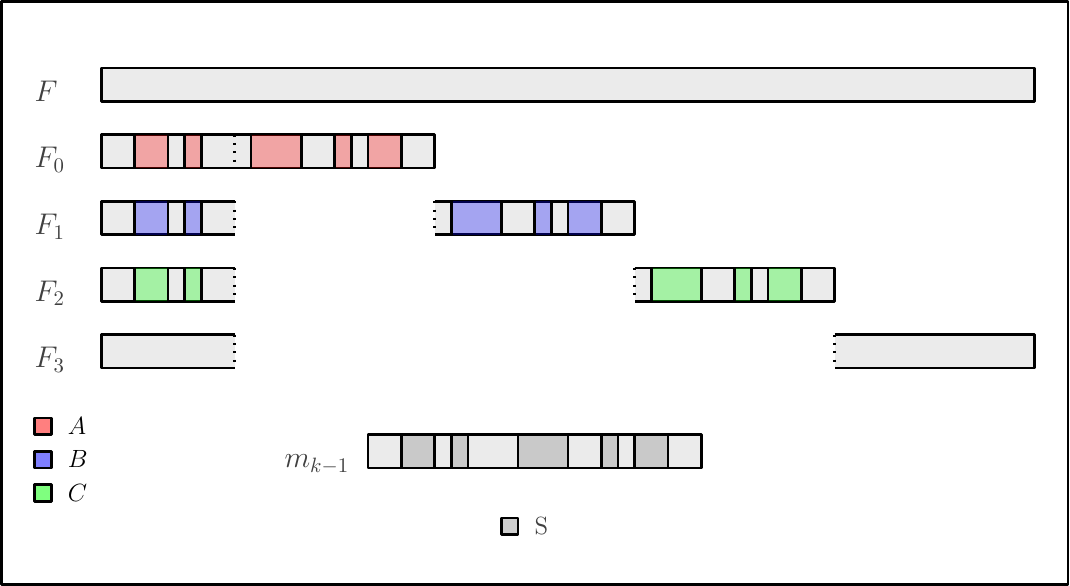}

\textit{\small In here, the set $\mathcal{D}=\{A,B,C\}$ is being captured as testified by $F$ and $S$. Note that $\mathcal{D}$ forms a $\Delta$-system whose root is contained in $R(F)$. } 
\end{minipage}
\end{center}

We will frequently use the following lemma without any explicit mention to it. Its proof is left to the reader. 
\begin{lemma}Let $\mathcal{D}$ be a finite subset of $\text{FIN}(X)$ which is captured at some level $l\in \omega$. If $n$ is the common cardinality of the members of $\mathcal{D}$ and $S$ is a non-empty finite subset of $n$, then $$\{D[S]\,:\,D\in \mathcal{D}\}$$
is also captured at level $l$ as long as this set has at least two elements.
\end{lemma}
The following proposition is a direct consequence of the Remark \ref{remarkcapturing2}.
\begin{proposition}\label{deltarhoequalityprop}Let $l\in\omega$ and $D\in \text{FIN}(X)$ be such that $|D|\geq 2$.  Then $D$ is captured at level $l$ if and only if:
\begin{enumerate}[label=$(\arabic*)$]
\item For any distinct $\alpha,\beta\in D$, $\Delta(\alpha,\beta)=l=\rho(\alpha,\beta)$,
\item For each $i<|D|$, $\Xi_{D(i)}(l)=i$ . In other words, $\{ \, \Xi_\alpha(l)\,:\,\alpha\in D\,\}=|D|$.
\end{enumerate}
\end{proposition}
\begin{lemma}\label{capturedfamiliestosetslemma}Let $\mathcal{D}$ be a finite subset of $\text{FIN}(X)$ and $l\in \omega$.  For any $A\in \mathcal{D}$, let $\alpha_A=\max(A)$. Suppose that there are $j\in \omega$, $a<m_j$ and $C\subseteq a+1$ so that the following conditions hold for any $A\in \mathcal{D}$:
\begin{enumerate}[label=$(\alph*)$]
    \item $\rho^A\leq j,$
    \item $\lVert \alpha_A\rVert_j=a$,
    \item $(\alpha_A)_j[C]=A.$ 
\end{enumerate}
If  $D=\{\alpha_A\,:\,A\in \mathcal{D}\}$ is captured at some level $l$, then $\mathcal{D}$ is also captured at level $l$.
\begin{proof}First observe that that $l>j$. In order to prove that $\mathcal{D}$ is captured at level $l$, we will prove that the three points in Definition \ref{capturedfamiliesdef} are satisfied.\\

\noindent
\underline{Claim 1}: $l=\rho^{\cup\mathcal{D}}$ and $l>\rho^A$ for each $A\in \mathcal{D}$.
\begin{claimproof}[Proof of claim] We know that $l>\rho^A$ for each $A\in \mathcal{D}$ due to the condition (a) written above. $l\leq \rho^{\cup \mathcal{D}}$ because if $A,B\in \mathcal{D}$ are distinct then $\alpha_A,\alpha_B\in \bigcup \mathcal{D}$ and $\rho(\alpha_A,\alpha_B)=l$. Lastly, if $F\in \mathcal{F}_l$ is such that $D\subseteq F$ then $A\subseteq (\alpha_A)_{j}\subseteq (\alpha_A)_l\subseteq F$ for any $A\in \mathcal{D}$. In this way, $\bigcup \mathcal{D}\subseteq F$. Thus, $\rho^{\mathcal{D}}\leq l$.
\end{claimproof}
\noindent
\underline{Claim 2}: $\{\Xi_A(l)\,:\,A\in \mathcal{D}\,\}=n_l$.
\begin{claimproof}[Proof of claim] As $l>\rho^A$ for each $A\in \mathcal{D}$, then $\Xi_A(l)=\Xi_{\alpha_A}(l)$ for any $A\in \mathcal{D}$ due to the Remark \ref{remxisetvalued}. The claim follows from this fact.
\end{claimproof}
\noindent
\underline{Claim 3}: For any distinct $A,B\in \mathcal{D}$, $A$ is strongly $\rho$-isomorphic to $B$.
\begin{claimproof}[Proof of claim] Let $h:(A)_{l-1}\Longrightarrow (B)_{l-1}$ be the increasing bijection. As $D$ is fully captured, we know that $\{\alpha_A\}$ is strongly isomorphic to $\{\alpha_B\}$. In this way, $h(\alpha_A)=\alpha_B$. Since $k\leq l-1$, we conclude that $h[(\alpha_A)_k]=(\alpha_B)_k$. Thus, $$h[A]=h[\,(\alpha_A)_k[C]\,]=(\alpha_B)_k[C]=B.$$
This finishes the proof.
\end{claimproof}
    
\end{proof}
    
\end{lemma}

\begin{definition}[$n$-capturing schemes]Let $\mathcal{P}$ be a partition of $\omega$ and $n\in\omega$. We say that $\mathcal{F}$ is \textit{$n$-$\mathcal{P}$-capturing} if for each uncountable $\mathcal{S}\subseteq \text{FIN}(X)$ and $P\in \mathcal{P}$ there are infinitely many  $l\in P$ with some $\mathcal{D}\in [\mathcal{S}]^n$ which is captured at level $l$. Whenever $\mathcal{P}=\{\omega\}$, we simply say that $\mathcal{F}$ is \textit{$n$-capturing.}
\end{definition}
\begin{definition}[capturing schemes]Let $\mathcal{P}$ be a partition of $\omega$. We say that $\mathcal{F}$ is \textit{$\mathcal{P}$-capturing} (resp. \textit{capturing}) if $\mathcal{F}$ is $n$-$\mathcal{P}$-capturing (resp. \textit{$n$-capturing}) for each $n\in\omega$.
\end{definition}

\begin{definition}[fully capturing schemes] Let $\mathcal{P}$ be a partition of $\omega$. We say  that $\mathcal{F}$ is \textit{$\mathcal{P}$-fully capturing} if for each uncountable $\mathcal{S}\subseteq \text{FIN}(X)$ and $P\in \mathcal{P}$  there are infinitely many $l\in P$ with some $\mathcal{D}\in \text{FIN}(\mathcal{S})$ which is fully captured at level $l$. Whenever $\mathcal{P}=\{\omega\}$, we simply say that $\mathcal{F}$ is \textit{fully capturing.}
\end{definition}

The following Lemma was first proved in \cite{schemenonseparablestructures} (Lemma 7.1) and it presents useful equivalences of the previous definitions.
\begin{lemma}\label{equivalencecapturing}Let $\mathcal{F}$ be a construction scheme and $\mathcal{P}$ be a partition of $\omega$ compatible with $\tau$. Then:
\begin{itemize}
\item For each $n\in\omega$, $\mathcal{F}$ is $n$-$\mathcal{P}$-capturing if and only if for each $S\in[X]^{\omega_1}$ and $P\in \mathcal{P}$ there is $D\in [S]^n$ which is captured at some level $l\in P$.
\item $\mathcal{F}$ is $\mathcal{P}$-fully capturing if and only if for each $S\in [X]^{\omega_1}$ and $P\in \mathcal{P}$ there is $D\in \text{FIN}(S)$ which is fully captured at some level $l\in P$.
\end{itemize}

\begin{proof} We will only prove the second point as the first one is proved in a similar way. In order to prove the nontrivial direction of this statement, let $\mathcal{S}$ be an uncountable subset of $\text{FIN}(X)$ and fix $k\in\omega$. For any $A\in S$, let $\alpha_A=\max(A)$. Due to the pigeonhole principle, we can find an uncountable $\mathcal{S}'\subseteq S$, $k<j\in \omega$, $a<m_j$ and $C\subseteq a$ so that the following conditions hold for any $A\in \mathcal{S}'$:
\begin{enumerate}[label=$(\alph*)$]
    \item $\rho^A\leq j.$
    \item $\lVert \alpha_A\rVert_j=a$,
    \item $(\alpha_A)_j[C]=A.$
\end{enumerate}
By the point (b), it follows that $\rho(\alpha_A,\alpha_B)>j$ for any two distinct $A,B\in \mathcal{S}'$. According to our assumptions, there is $\mathcal{D}\in \text{FIN}(\mathcal{S})$ for which $D=\{\alpha_A\,:\,A\in \mathcal{D}\}$ is fully captured at some level $l$. Then $\mathcal{D}$ is captured in level $l$ by virtue of the Lemma \ref{capturedfamiliestosetslemma}. Thus, the proof is over.

\end{proof}

\end{lemma}
We are now ready to state the Capturing axioms. The axiom $FCA(part)$ was introduced in \cite{schemenonseparablestructures}. The other axioms were later studied in 
\cite{irredundantsetsoperator},
\cite{forcingandconstructionschemes}, \cite{banachspacescheme} and \cite{lopezschemethesis}.\\\\
{\bf Fully Capturing Axiom} {[\bf FCA]}: There is a fully capturing construction scheme over $\omega_1$ of every possible good type.\\\\
{ \bf Fully Capturing Axiom with Partitions [FCA(part)]}: There is a $\mathcal{P}$-fully capturing construction scheme over $\omega_1$ for every good type $\tau$  and each partition $\mathcal{P}$ compatible with $\tau$.\\\\
{\bf $n$-Capturing Axiom [CA$_n$]}: There is an $n$-capturing construction scheme over $\omega_1$ of every possible good type satisfying that $n\leq n_k$ for each $k\in \omega\backslash1$.\\\\
{\bf $n$-Capturing Axiom with Partitions [CA$_n$(part)]}: There is a $\mathcal{P}$-$n$-capturing construction scheme over $\omega_1$ for every good type $\tau$ satisfying that $n\leq n_k$ for each $k\in \omega\backslash 1$  and each partition $\mathcal{P}$ compatible with $\tau$.\\\\
{\bf Capturing Axiom [CA]}: $CA_n$ holds for each $n\in\omega$  and there is a capturing construction scheme over $\omega_1$ for every good type satisfying that the sequence $\langle n_{k+1}\rangle_{k\in\omega}$ is non-decreasing and unbounded.\\\\
{\bf Capturing Axiom with partitions [CA(part)]}: $CA_n(part)$ holds for each $n\in\omega$ and there is a capturing construction scheme over $\omega_1$ for every good type $\tau$ satisfying that the sequence $\langle n_{k+1}\rangle_{k\in\omega}$ is non-decreasing and unbounded and each partition $\mathcal{P}$-compatible with $\tau$.\\\\

Note that the following relations hold between the previously defined capturing axioms where the arrows denote the implication relation.
\begin{center}
\begin{tikzcd}
CA_2(part) \arrow[d, Rightarrow] &CA_3(part) \arrow[l, Rightarrow] \arrow[d, Rightarrow]&\dots \arrow[l, Rightarrow] \arrow[d, Rightarrow] & CA(part)\arrow[l, Rightarrow]\arrow[d, Rightarrow] & FCA(part)\arrow[d, Rightarrow] \arrow[l, Rightarrow]\\
CA_2 & CA_3 \arrow[l, Rightarrow] &\dots\arrow[l, Rightarrow] & CA\arrow[l, Rightarrow] & FCA\arrow[l, Rightarrow]
\end{tikzcd}
\end{center}
 
In \cite{forcingandconstructionschemes}, Damian Kalajdzievski and Fulgencio Lopez proved that none of the arrows going from $CA_{n+1}$ to $CA_n$ (resp. from $CA_{n+1}(part)$ to $CA_n(part)$) can be reversed. They also proved the following theorem.
\begin{theorem}\label{Cohenforcingfcatheorem}Let $\kappa$ be an uncountable cardinal. Then $\mathbb{1}_{\mathbb{C}_\kappa}\vDash\text{\say{ $ \text{FCA}(part)$ }}$.
\end{theorem}

For the sake of completeness, we will prove such results later on. Furthermore, we will show that it is consistent that there is an $n$-capturing construction scheme which is not $\mathcal{P}$-$n$-capturing for any partition of $\omega$ in at least two infinite pieces. For this,  we need to introduce a new cardinal invariant.

\begin{definition}[Parametrized Martin's numbers]\label{parametrizedmartinsdef}Let $\mathcal{F}$ be a construction scheme and $2\leq n\in\omega$. We define $\mathfrak{m}^n_\mathcal{F}$ as follows:
$$\mathfrak{m}^n_\mathcal{F}=\begin{cases}
\omega &\textit{if }\mathcal{F}\textit{ is not }n\textit{-capturing}\\
\min(\mathfrak{m}(\mathbb{P})\,:\,\mathbb{P}\textit{ is }ccc\textit{ and }\mathbb{P}\Vdash \text{\say{$\mathcal{F}\textit{ is }n-capturing$}})&\textit{if }\mathcal{F}\textit{ is }n\textit{-capturing}
\end{cases}$$
$\mathfrak{m}^2_\mathcal{F}$ is denoted simply as $\mathfrak{m}_\mathcal{F}.$ 
\end{definition}
Of course, if $\mathcal{F}$ is $n$-capturing then $\mathfrak{m}^n_\mathcal{F}\geq \omega_1.$ The following theorem will be proved in Chapter \ref{deeperanalysis}.
\begin{theorem}\label{theoremofmF}For any $n$-capturing construction scheme $\mathcal{F}$ there is a $ccc$-forcing $\mathbb{P}$ for which $$\mathbb{P}\Vdash\text{\say{ $\mathfrak{m}_{\mathcal{F}}^n>\omega_1$ }}.$$
\end{theorem}
We end this chapter by announcing the main theorem of the thesis.
\begin{theorem}\label{diamondfcaparttheorem} The $\Diamond$-principle implies $\text{FCA}(part)$.
\end{theorem}

%% file: chapters/newgaps.tex
\chapter{Applications}\label{applicationschapter}

\section{Gaps, towers and almost disjoint families} \label{sectiongaps}
This section is dedicated to the study of gaps and towers. In a general framework, one of the main reasons for studying gaps in Boolean algebras is because of their close relation with the existence of extensions of both homomorphisms and embeddings from one Boolean algebra to another one. More concretely, in \cite{Sikorskiextension} Roman Sikorski implicitly gave a criterion for the extension of homomorphisms in Boolean algebras. Such criterion tells us that the objects the we now know as gaps are the only thing  which can prevent an homomorphism to be extended.  As we will center our attention solely on gaps over $\mathscr{P}(\omega)/\text{FIN}$, all the definitions included here are stated only for such algebra.

The reader interested in knowing more about towers and gaps is referred to \cite{GapsandTowers}, \cite{rothbergergapsinfragmentedideals}, \cite{gapsandlimits}, \cite{Agapcohomologygroup},  \cite{ScheepersGaps}, \cite{SomeResultsonGaps}, \cite{Walksonordinals}, \cite{StevoIlias}, \cite{PartitionProblems}, \cite{AnalyticGaps} and \cite{combinatorialprinciplesonomega1}.

We start by recalling the main definitions regarding this chapter.
\begin{definition}[Towers]Let $X$ be a countable set and $\mathcal{T}\subseteq[X]^\omega$. We say that $\mathcal{T}$ is a \textit{tower}\footnote{It is also common in the literature to call these objects pre-towers, and towers are pre-towers which are maximal with respect to the end-extension. As we won't be dealing with maximal pre-towers, we've opted to simply refer to them as "towers" to simplify the notation.} if it is well-ordered with respect to $\subseteq^*$. Furthermore, for an ordinal $\kappa$, we say that $\mathcal{T}$ is a $\kappa$-tower whenever it is a tower and $ot(\mathcal{T})=\kappa.$
\end{definition}

\begin{definition}[Gaps and Pregaps]
Let $X$ be a countable set and $\mathcal{L},\mathcal{R}\subseteq [X]^\omega$. We say that $(\mathcal{L},\mathcal{R})$ is \textit{pregap}, and write it as $\mathcal{L}\perp \mathcal{R}$, if $L\cap R=^*\emptyset$ for all $L\in \mathcal{L}$ and $R\in \mathcal{R}$. An element $C\in [X]^\omega$ is said to \textit{separate} $(\mathcal{L},\mathcal{R})$ if $L\subseteq^*C$ and $C\cap R=^*\emptyset$ for each $L\in\mathcal{L}$ and $R\in \mathcal{R}$. Finally, we say that $(\mathcal{L},\mathcal{R})$ is a \textit{gap} if it is a pregap and there is no $C\in [X]^\omega$ separating it.
\end{definition}
\begin{rem}It is well known that there are no countable gaps. That is, there is no gap $(\mathcal{L},\mathcal{R})$ with both $\mathcal{L}$ and $\mathcal{R}$ countable.
    
\end{rem}

\begin{rem}Whenever the domain of a pregap is not specified (the set $X$ of the previous definition) we assume that such domain is $\omega$.
\end{rem}

\begin{definition}[Type of a gap]Let $(X,\leq_X)$ and $(Y,\leq_Y)$ be two partial orders. We say that a pregap $(\mathcal{L},\mathcal{R})$ is an \textit{$(X,Y)$-pregap} if $(\mathcal{L},\subseteq^*)$ and $(\mathcal{R},\subseteq^*)$ are isomorphic to $(X,\leq_X)$ and $(Y,\leq_Y)$ respectively. Furthermore, if $(\mathcal{L},\mathcal{R})$ is a gap we will refer to it as an $(X,Y)$-gap. Whenever both sides of a pregap $(\mathcal{L},\mathcal{R})$ are index by some set $I$, say $\mathcal{L}=\langle A_i\rangle_{i\in I}$ and $\mathcal{R}=\langle B_i\rangle_{i\in I}$ we may refer to $(\mathcal{L},\mathcal{R})$ as $(A_i,B_i)_{i\in I}$. Of course, if $I$ is a  partial order and we specify that the pregap is an $(I,I)$-gap it is understood such an indexing is order preserving.
\end{definition}

It is easy to construct an $(\omega_1,\omega_1)$-gap using $CH$. However, in 1909, Felix Hausdorff gave a clever recursive construction of an $(\omega_1,\omega_1)$-gap without appealing to any extra axioms. This is quite surprising, as it is required to overcome $\mathfrak{c}$ obstacles in only $\omega_1$-many steps. In Theorem \ref{allgapsthm}, we will prove an analogous result but for large collection of partial orders of size $\omega_1$.

 The gap that Hausdorff constructed satisfied the following property.

\begin{definition}[Hausdorff condition]Let $(\mathcal{L},\mathcal{R})$ be an $(\omega_1,\omega_1)$-pregap on $\omega$. We say that $(\mathcal{L},\mathcal{R})$ is  Hausdorff if $\{R\in \mathcal{R}\,:\,rank(R)<rank(L)\textit{ and }L\cap R\subseteq k\}$ is finite for each $L\in\mathcal{L}$ and $k\in\omega.$
\end{definition}

\begin{definition}[Luzin condition]Let $(\mathcal{L},\mathcal{R})$ be an $(\omega_1,\omega_1)$-pregap on an infinite set $X$. We say that $(\mathcal{L},\mathcal{R})$ is Luzin if $\{R\in \mathcal{R}\,:\,rank(R)<rank(L)\textit{ and }|L\cap R|\leq k\}$ is finite for each $L\in\mathcal{L}$ and $k\in\omega.$
\end{definition}

\begin{rem}\label{remarkhausdorffgaps1}Note that the Luzin condition is stronger than the Hausdroff condition. It is a standard excercise to prove that any $(\omega_1,\omega_1)$-pregap satisfying the Hausdorff condition is in fact a gap. Under the $P$-ideal Dichotomy ($PID$), if $(\mathcal{L},\mathcal{R})$ is an $(\omega_1,\omega_1)$-gap then there are cofinal $\mathcal{L}'\subseteq \mathcal{L}$ and $\mathcal{R}'\subseteq \mathcal{R}$ so that $(\mathcal{L}',\mathcal{R}')$ is a Hausdorff gap. This result was proved by Uri Abraham and Stevo Todor\v{c}evi\'{c} in \cite{partitionpropertiesch}. The same result follows from $\mathfrak{m}>\omega_1.$
    
\end{rem}

In the following theorem we construct a Hausdorff gap using a $2$-construction scheme (that is, a construction scheme of type $\langle m_k, 2, r_{k+1}\rangle_{k\in\omega}$). The first construction of a Hausdorff gap using morasses was performed by Daniel Velleman in \cite{omegamorasses}. In \cite{lopezschemethesis} and \cite{Walksonordinals}  the reader may find a construction of such a gap using construction schemes and ordinal metrics respectively. First, a quick remark.
\begin{rem}\label{remhaus1}If $\langle L_\alpha\rangle_{\alpha\in \omega_1}$ and $\langle R_\alpha\rangle_{\alpha\in \omega_1}$ are $\omega_1$-towers and $L_\alpha\cap R_\alpha=^*\emptyset$ for each $\alpha\in \omega_1$ then $(L_\alpha,R_\alpha)_{\alpha\in\omega_1}$ is a pregap. This is because for each $\alpha<\beta\in\omega_1$, $L_\alpha\cap R_\beta\subseteq^*L_\beta\cap R_\beta.$
\end{rem}
\begin{theorem}\label{hausdorffgapconstruction}Let $\mathcal{F}$ be a $2$-construction scheme. For each $\alpha\in \omega_1$, define $$L_\alpha=\{2k+\Xi_\alpha(k)\,:\,k\in \omega\backslash 1,\textit{ and }\Xi_\alpha(k)\geq 0\},$$
$$R_\alpha=\{2k+(1-\Xi_\alpha(k))\,:\,k\in\omega\backslash 1,\textit{ and }\Xi_\alpha(k)\geq 0\}.$$
Then, $(L_\alpha, R_\alpha)_{\alpha\in\omega_1}$ is a Hausdorff gap.
\begin{proof} Since $r_k=0$ for infinitely many $k's$, it should be clear that each $L_\alpha$ and $R_\alpha$ are infinite.  By virtue of the part (c) in Lemma \ref{lemmaxi}, we have that if $\alpha<\beta$ then  $$L_\alpha\backslash L_\beta\subseteq \{2k+\Xi_\alpha(k)\,:\,k\leq \rho(\alpha,\beta)\textit{ and }\Xi_\alpha(k)\not=\Xi_\beta(k)\,\},$$
$$R_\alpha\backslash R_\beta\subseteq \{2k+(1-\Xi_\alpha(k))\,:\,k\leq \rho(\alpha,\beta)\textit{ and }\Xi_\alpha(k)\not=\Xi_\beta(k)\,\}.$$ 
As the sets on the right are finite, this means  $\langle L_\alpha\rangle_{\alpha\in \omega_1}$ and $\langle R_\alpha\rangle_{\alpha\in \omega_1}$ are both $\omega_1$-towers. Furthermore by definition we have that $L_\alpha\cap R_\alpha=\emptyset$ for each $\alpha\in\omega_1$. This  implies that $(L_\alpha,R_\alpha)_{\alpha\in\omega_1}$ is a pregap due to the Remark \ref{remhaus1}.

The only thing left to prove is that the Hausdorff condition is satisfied. For this purpose take $\beta\in \omega_1$ and $k\in\omega$. We claim  $\{\alpha<\beta\,:\,L_\beta\cap R_\alpha\subseteq k\}\subseteq (\beta)_k$. For this, take an arbitrary $\alpha<\beta$ satisfying $\rho(\alpha,\beta)\geq k.$ By the part (b) of Lemma \ref{lemmaxi},  $\Xi_\alpha(\rho(\alpha,\beta))=0$ and $\Xi_\beta(\rho(\alpha,\beta))=1$. This means $2\rho(\alpha,\beta)+1\in L_\beta\cap R_\alpha$, so we are done.
\end{proof}
\end{theorem}
An interesting feature of the proof given above is that not only we prove that there is a gap using  a $2$-construction scheme, but such gap can be explicitly defined from it.

For a countably infinite set $X$, we say that a family $\mathcal{A}\subseteq [X]^\omega$ is an almost disjoint family if $A\cap B=^*\emptyset$ whenever $A,B\in \mathcal{A}$ are different. Almost disjoint families are one of the central objects of study in modern combinatorial set theory.  Constructing  almost disjoint families with special properties is usually difficult and in most cases had lead to the development of powerful tools (see  \cite{InvariancePropertiesofAlmostDisjointFamilies}, \cite{Madfamilieswithstrongcombinatorialproperties}, \cite{madnessandnormality},
\cite{frechetlike},  \cite{nonpartitionable}, \cite{cechfunction},\cite{ThereisavanDouwenMADfamily} and \cite{SANEPlayer}). Almost disjoint families have also played a central roll in the solution of many problems of Topology and Analysis. An example of this is the solution of the selection problem posed by Jan van Mill and Evert Wattel in \cite{van1981selections} and solved by Michael Hru\v{s}\'ak and Iv\'an Mart\'inez-Ru\'iz in \cite{hruvsak2009selections}. The reader interested in learning more about almost disjoint families is refered to \cite{TopologyofMrowkaIsbellSpaces}, \cite{Combinatoricsoffiltersandideals} and \cite{AlmostDisjointFamiliesandTopology}.

Our next goal is to use almost disjoint families with the aim of proving that, basically, there are gaps of any possible type for which the cofinality of the two partial ordered sets involved is $\omega_1$.

\begin{definition}Let $\mathcal{A}$ be an almost disjoint family of size $\omega_1$ over a set $X$. We say that:\begin{itemize}
\item $\mathcal{A}$ is \textit{inseparable} if for any two disjoint $\mathcal{L},\mathcal{R}\in [\mathcal{A}]^{\omega_1}$, the pair $(\mathcal{L},\mathcal{R})$ forms a gap. 
    \item $\mathcal{A}$ is \textit{Luzin} if we can enumerate it as $\langle A_\alpha\rangle_{\alpha\in\omega_1}$ in such way that $\{\alpha\in\beta\,:\,|A_\alpha\cap A_\beta|< n\}$ is finite for each $\beta\in\omega_1$ and $n\in\omega$.
    \item $\mathcal{A}$ is \textit{Jones} if for any two disjoint  $\mathcal{L}\in [\mathcal{A}]^{\omega}$ and $\mathcal{R}\in [\mathcal{A}]^{\omega_1}$, the pregap $(\mathcal{L},\mathcal{R})$ can be separated. If $\mathcal{A}$ is indexed as $\langle A_\alpha\rangle_{\alpha\in \omega_1}$ then $\mathcal{A}$ is Jones if and only if for any $\beta\in \omega_1$ the pregap $(\langle A_\alpha\rangle_{\alpha\leq \beta},\langle A_{\delta}\rangle_{\delta>\beta})$ can be separated.
    \item $\mathcal{A}$ is \textit{Luzin-Jones} if it is both Luzin and Jones.
\end{itemize} 
\end{definition}
\begin{rem}It is easy to check that any Luzin family is in fact inseparable. In \cite{luzin1947subsets}, Judith Roitman and Lajos Soukup showed that under $\mathfrak{m}>\omega_1$ any uncountable $AD$ family either contains a Luzin family or contains no inseparable family.
\end{rem}
   The first construction of a Luzin family was done in \cite{luzin1947subsets} by Nikol\'ai Nikol\'ayevich Luzin. A Jones family was implicitly constructed by F. B. Jones in \cite{jones1937concerning}. Although, at first glance, Luzin and Jones properties seem to be incompatible, a construction of a Luzin-Jones family was obtained in \cite{guzman2019mathbb} by Osvaldo Guzm\'an
, Michael Hru\v{s}\'ak and Piotr Koszmider (building from work by Koszmider in \cite{onconstructionswith2cardinals}). We would like to point out that the highly complex construction of the mentioned family is carried out through the use of simplified morasses. 
Here we give an elementary construction of such object.

\begin{theorem}\label{luzinjonestheorem}There is a Luzin-Jones family.
\begin{proof}Let $\mathcal{F}$ be a $2$-construction scheme. For every $k\in \omega\backslash 1$, let $$N_k= \{k\}\times (m_{k-1}\backslash r_k)\times [k\cdot m_{k-1},\,k\cdot m_k).$$
We aim to construct a  Luzin-Jones family $\langle A_\alpha\rangle_{\alpha\in \omega_1}$ over the union of the $N_k$'s, namely $N$. For this, take an arbitrary $\alpha\in\omega_1$. Given $k\in \omega\backslash 1$, we define $A_\alpha^k\subseteq N_k$ as follows:
\begin{enumerate}[label=$(\alph*)$]
\item If $\Xi_\alpha(k)=-1$, let $A^k_\alpha=\emptyset$.
     \item If $\Xi_\alpha(k)=0$, let $A^k_\alpha=\{k\}\times \{\lVert \alpha\rVert_k\}\times [k\cdot m_{k-1},\,k\cdot m_k)$.
     \item If $\Xi_\alpha(k)=1$, let $A^k_\alpha=\{k\}\times [r_k,m_{k-1})\times \big[k\cdot\lVert\alpha \rVert_k,\, k (\lVert \alpha\rVert_k+1)\big).$
\end{enumerate} 
Now, we define $A_\alpha$ as $\bigcup\limits_{k\in\omega\backslash 1}A^k_\alpha$ and $\mathcal{A}$ as $\langle A_\alpha\rangle_{\alpha\in\omega_1}$.\\

Our first task is to show that $\mathcal{A}$ is an almost disjoint family. First observe that since $r_k=0$ infinitely often it follows that each $A_\alpha$ is in fact infinite. Now let $\alpha<\beta\in\omega_1$ and take an arbitrary $k>\rho(\alpha,\beta)$. In virtue of the part (c) by Lemma \ref{lemmaxi} we have that $\Xi_\alpha(k)=-1$ or $\Xi_\alpha(k)=\Xi_\beta(k)$. In either case it follows that $A^{k}_\alpha\cap A^{k}_\beta =\emptyset.$ Therefore we conclude that $$A_\alpha\cap A_\beta=\bigcup_{k\in\omega\backslash 1}\big(A^k_\alpha\cap A^k_\beta\big)\subseteq \bigcup_{k\leq \rho(\alpha,\beta)}N_{k}.$$

\noindent
\begin{minipage}[b]{0.57\linewidth}
Now we prove that $\mathcal{A}$ is Luzin. We claim that whenever $\alpha<\beta$ and $k=\rho(\alpha,\beta)$ then $|A_\alpha\cap A_\beta|\geq k$. If this happens then $\{\alpha<\beta\,:\,|A_\alpha\cap A_\beta|<k\}\subseteq (\beta)_k$ for any $\beta\in\omega_1$. Let $\alpha, \beta$ and $k$ be as previously stated. The part (b) of Lemma \ref{lemmaxi} assures  $\Xi_\alpha(k)=0<\Xi_\beta(k)=1$. In this way $$A^k_\alpha\cap A^k_\beta=\{k\}\times \{|(\alpha)^-_k|\}\times \big[k\cdot\lVert\alpha \rVert_k,\, k (\lVert \alpha\rVert_k+1)\big).$$
The cardinality of this set is $k$ and it is contained in $A_\alpha\cap A_\beta$. Thus, we are done. 
\end{minipage}\hspace{0.03\linewidth}\begin{minipage}[t]{0.4\linewidth}
\centering
\includegraphics[width=6cm, height=5.4cm]
{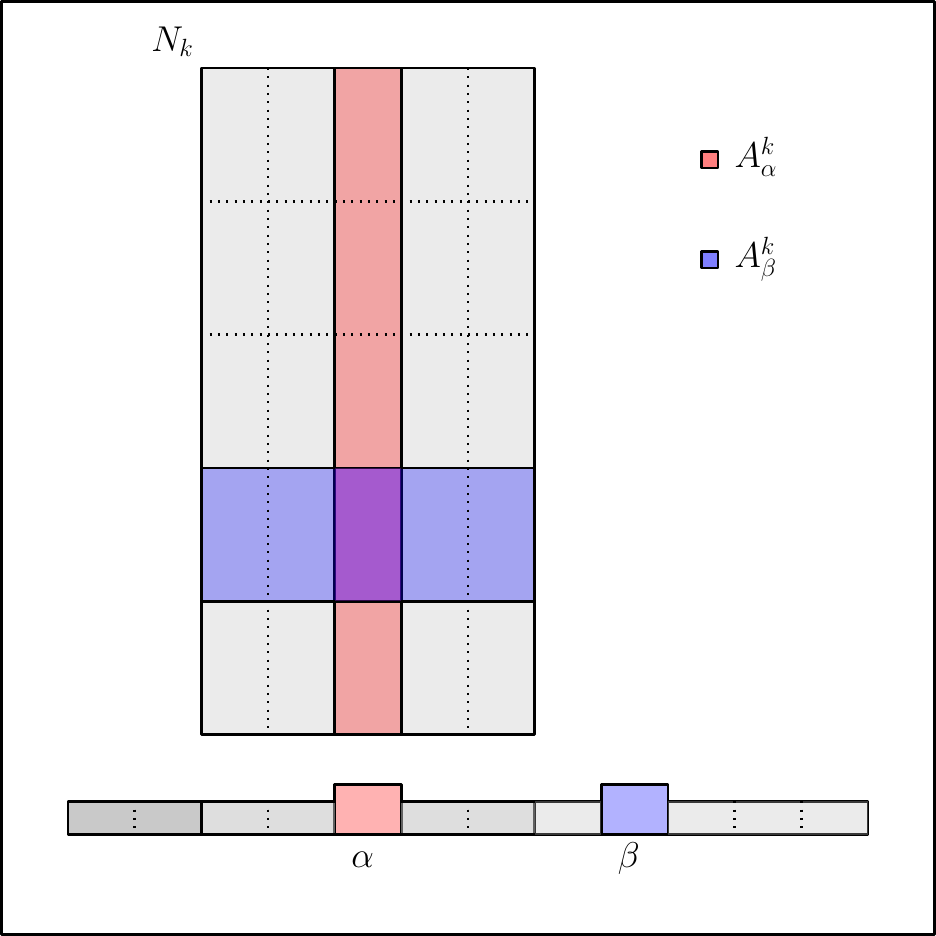}
\textit{\small In here, $\Xi_\alpha(k)=0$ and $\Xi_\beta(k)=1$.}

\end{minipage}
\vspace{0.5cm}\\
We will end the proof by showing that $\mathcal{A}$ is a Jones family. For this purpose define  $C_\beta$ as $$\bigcup\limits_{k\in\omega}\big(\bigcup\limits_{\alpha\in (\beta)_k } A^{k+1}_\alpha\big)$$ for each $\beta\in\omega_1$. We will show that $C_\beta$ separates $(\langle A_\alpha\rangle_{\alpha\leq \beta},\langle A_\alpha\rangle_{\delta> \beta})$. Indeed, if $\alpha\leq\beta$ then $A^{k+1}_\alpha\subseteq C_\beta$ for 
each $k\geq \rho(\alpha,\beta)$, so in particular $A_\alpha\subseteq^*C_\beta$. On the other hand, if $\delta>\beta$ then $A^{k+1}_\delta\cap A^{k+1}_\alpha=\emptyset$ for any $k\geq \rho(\delta,\beta)$ and each $\alpha\in (\beta)_k$. This is because for any such $\alpha$ we have that $\rho(\alpha,\delta)\leq k$. Consequently $A_\delta\cap C_\beta=^*\emptyset$. This finishes the proof.
\end{proof}
\end{theorem}
The previous construction yields an almost disjoint family which, in principle, satisfies a lot more properties than just being Luzin-Jones. In the rest this section, we will analyse those properties. The following definitions are essential for that task.
\begin{definition}[Succesor-like elements]Let $(X,\leq)$ be a partial order. We say that $x\in X$ is \textit{succesor-like} if:
\begin{itemize}
    \item $pred(x)$ is finite,
    \item for any $y<x$ there is $z\in pred(x)$ for which $y\leq z$.
\end{itemize}
 
\end{definition}

\begin{definition}[Luzin representation]\label{Luzinrepdef} Let $(X,\leq)$ be a partial order of size $\omega_1$. A \textit{Luzin representation of $X$} is an ordered pair $(\mathcal{T},\mathcal{A})$ consisting of two families of infinite subsets of $\omega$ indexed as $\langle T_x\rangle_{x\in X}$ and $\langle A_x\rangle_{x\in X}$ respectively. Furthermore, $\mathcal{A}$ is a Luzin family and for all $x,y\in X$ the following properties hold:
\begin{enumerate}[label=$(\alph*)$]
\item $A_x\subseteq T_x.$

\item If $y\not\leq x$, then $A_y\subseteq^*T_y\backslash T_x$.
\item If $\inf(x,y)$ exists, then $T_x\cap T_y=^*T_{\inf(x,y)}.$

\item If $x$ is succesor-like, then $$T_x\backslash\big(\bigcup\limits_{z\in pred(x)}T_z\big)=^*A_x.$$
\item If $x$ and $y$ are incompatible, then $T_x\cap T_y=^*\emptyset.$

\end{enumerate}
In particular $(\mathcal{T},\subseteq^*)$ is order isomorphic to $X$ by virtue of the points $(b)$ and $(c)$. In the case there is $(\mathcal{T},\mathcal{A})$ a Luzin representation of $X$ we say that $\mathcal{A}$ \textit{codes} $X$.
\end{definition}

\begin{definition}[$\omega_1$-like order]\label{omega1likedef} We say that a partial order $(X,\leq)$ is $\omega_1$-like if it is well-founded, $|X|=\omega_1$ and $|(-\infty,x)|\leq \omega$ for all $x\in X$.
\end{definition}

\begin{proposition}\label{propositionwellfounded}Let $(Y,\leq)$ be a partial order of cardinality $\omega_1$. There is a well-founded cofinal $X\subseteq Y$ with $|(-\infty, x)_X|\leq \omega$ for each $x\in X$.
\begin{proof} Enumerate $Y$ as $\langle y_\alpha\rangle_{\alpha\in \omega_1}$.  Define $X$ as the set of all $y_\beta$'s such that $y_\beta\not\leq y_\alpha$ for each $\alpha<\beta$. Of course $X$ is well-founded as for each $x,z\in X$ with $x=y_\alpha$ and $z=y_\beta$ it happens that if $x<z$ then $\alpha<\beta$. Because of this we also have that $(-\infty,\, y_\beta)_X\subseteq \{y_\alpha\,:\,\alpha<\beta\}$ for each $y_\beta\in X$. Lastly $X$ is cofinal in $Y$ because for each $\beta\in \omega_1$ the element $y_{\xi_\beta}\in X$ where  $\xi_\beta=\min(\,\alpha\leq \beta\,:\, y_\alpha \geq y_\beta\,)$.
\end{proof}
\end{proposition}
The following lemma is easy.
\begin{lemma}\label{bijectionlemma}Let $(X,\leq)$ be $\omega_1$-like. Then there is a bijection $\phi:X\longrightarrow \omega_1$ with $\phi(x)<\phi(y)$ for all $x<y\in X.$
\end{lemma}

Let $(X,\leq)$ $\omega_1$-like and $\phi:X\longrightarrow \omega_1$ as above. In the next lemma, in some sense, we will pull back a construction scheme from $\omega_1$ to some sort of \say{construction scheme} on $X$. This suggests that the theory of construction schemes could be generalized to other partial orders of size $\omega_1$. The extend and utility of this idea remains mainly unexplored.

\begin{lemma}\label{Mklemma}Let $\mathcal{F}$ be a $2$-construction scheme and let $X$ and $\phi$ be as in Lemma \ref{bijectionlemma}. For each $x\in X$ and $k\in\omega$ define $M^k_x=\{\, z\leq x\,:\, \phi(z)\in (\phi(x))_k\,\}$. The following properties hold for each $x,y\in X$ and $k\in\omega$:
\begin{enumerate}[label=$(\arabic*)$]
\item  If $\inf(x,y)$ exists and $k>\rho^{\phi[\{x,y,\inf(x,y)\}]}$, then $M^k_x\cap M^k_y=M^k_{\inf(x,y)}.$
\item[$(\frac{3}{2})$] If $y\leq x$ and $k>\rho^{\phi[\{x,y\}]}$, then $M^k_y\subseteq M^k_x$.
\item If $y\not\leq x$, then $y\in M^k_y\backslash M^k_x.$
\item If $x$ is succesor-like and $k>\rho^{\phi[pred(x)\cup\{x\}]}$, then $M^k_x\backslash\big( \bigcup\limits_{z\in pred(x)} M^k_z \big)=\{x\}$.
\item If $x$ and $y$ are incompatible,  then $M^k_x\cap M^k_y=\emptyset.$
\end{enumerate}
\begin{proof}
The points (2) and (4) are trivial and the point ($\frac{3}{2}$) follows directly from (1). Therefore,  we will only prove (1) and (3).

\begin{claimproof}[Proof of  $(1)$] Since $\phi$ preserves the order then $\phi(\inf(x,y))\leq \min(\phi(x),\phi(y))$ which implies that $(\phi(\inf(x,y))_k\sqsubseteq (\phi(x))_k\cap (\phi(y))_k$. In this way, if $z\in M^k_{\inf(x,y)}$ then $z\leq x,y$ by definition of the infimum, and $z$ belongs to both $(\phi(x))_k$ and $(\phi(y))_k$ because $z\in (\phi(\inf(x,y)))_k$. Therefore $z\in M^k_x\cap M^k_y$. On the other hand, if $z\in M^k_x\cap M^k_y$ then $z\leq x$ and $z\leq y$, so $z\leq \inf(x,y)$. Furthermore, $\phi(z)\in (\phi(x))_k\cap (\phi(y))_k$ and $\phi(z)\leq \phi(\inf(x,y))$. As $(\phi(\inf(x,y)))_k$ is an initial segment of such intersection then $\phi(z)$ belongs to it. Consequently $z\in M^k_{\inf(x,y)}$. 
\end{claimproof}
\begin{claimproof}[Proof of $(3)$] Note that since $x\not\leq z$ for any $z\in pred(x)$ then $x\in M^k_x\backslash M^k_z$ for any such $z$. This proves the inclusion from right to left. To show that the one from left to right also holds, let $$w\in M^k_x\backslash\big( \bigcup\limits_{z\in pred(x)} M^k_z \big).$$
Suppose towards a contradiction that $w<x$ and consider $z\in pred(x)$ with $w\leq z$. Since $k>\rho(\phi(x),\phi(z))$ we conclude that $(\phi(z))_k\sqsubseteq (\phi(x))_k$. But then  $\phi(w)\in (\phi(z))_k$ because $\phi(w)\leq \phi(z)$ and $\phi(w)\in (\phi(x))_k$. Consequently $w\in M^k_z$ which is a contradiction. This finishes the proof.
  \end{claimproof}  

\end{proof}
\end{lemma}
In the next theorem we show that not only a great variety of partial orders posses a Luzin representation, but that the  particular Luzin family that we constructed in the Theorem \ref{luzinjonestheorem} codes all of them. 
For that reason, it will be convenient to recall some basic facts regarding such family.

We have a countable set $N$ and $\langle N_k\rangle_{k\in\omega}$ a partition of $N$ into finite sets. For each $\alpha\in \omega_1$ we have $A_\alpha^k\subseteq N_k$ and we defined $A_\alpha=\bigcup\limits_{k\in \omega} A^{k+1}_\alpha$ and $\mathcal{A}=\langle A_\alpha\rangle_{\alpha\in \omega_1}$. One key feature of $\mathcal{A}$ is that: \begin{center}If $\alpha,\beta\in\omega_1$ are distinct and $k>\rho(\alpha,\beta)$ then $A^k_\alpha\cap A^k_\beta=\emptyset.$
\end{center}

\begin{theorem}\label{luzinreptheorem}Let $\mathcal{A}$ be the Luzin family constructed in Theorem \ref{luzinjonestheorem} and let $(X,\leq)$ an $\omega_1$-like order. Then there is a family $\mathcal{T}=\langle T_x\rangle_{x\in X}\subseteq [\omega]^{\omega}$ and a re-indexing of $\mathcal{A}$ as $\langle\hat{A}_x\rangle_{x\in X}$ so that $(\mathcal{T},\mathcal{A})$ is a Luzin representation of $X$.
\begin{proof}Let $\phi:X\longrightarrow \omega_1$ be as in Lemma \ref{bijectionlemma} and
$ M^k_x$ be as in Lemma 
\ref{Mklemma} for any $x\in X$ and $k\in \omega$. Fix $x\in X$ and let $\hat{A}_x=A_{\phi(x)}.$
This defines the re-indexing of $\mathcal{A}$. Now, for each $k\in\omega$ and $x\in X$, let 
$$T^k_x=\bigcup \{\,A_{\phi(z)}^{k+1}\,:\,z\in M^k_x\,\}=\bigcup\{\,A^{k+1}_\xi\,:\,\xi\in \phi[M^k_x]\,\}$$
$$T_x=\bigcup\limits_{k\in\omega}T^k_x.$$
We claim that $(\mathcal{T},\mathcal{A})$ is a Luzin representation of $X$. In the following paragraphs we will show that the points $(a)$, $(b)$, $(c)$, $(d)$ and $(e)$ of Definition \ref{Luzinrepdef} are satisfied for such pair. 
\begin{claimproof}[Proof of $(a)$]Given $k>0$ we have that $A^k_{\phi(x)}\subseteq T^{k-1}_x$ because $x$ is always an element of $M^k_x$. Therefore $\hat{A}_x=A_{\phi(x)}\subseteq T_x$.   
\end{claimproof}
 Before proving the remaining points let us fix some notation. Given $P\in[\omega_1]^{<\omega}$ and $k\in\omega$, let $A^k_P=\bigcup\{A^k_\alpha\,:\,\alpha\in P\}$. Note that if $P,Q\in [\omega_1]^{<\omega}$ and $k>\rho^{P\cup Q}$ then, by virtue of the key feature of $\mathcal{A}$ highlighted just before this theorem, we have the following:

 $$A^k_P\cap A^k_Q=A^k_{P\cap Q},$$
$$A^k_P\cup A^k_Q=A^k_{P\cup Q},$$
$$A^k_{P}\backslash A^k_Q=A^k_{P\backslash Q}.$$

We will apply these three equalities for finite sets of the form $\phi[M^k_x]$. More precisely, observe that if $x,y\in X$ and $k>\rho(\phi(x),\phi(y))$ then $\rho^{\phi[M^k_x]\cup \phi[M^k_y]}<k+1$. This is because if we take  $F\in \mathcal{F}_k$ with $\phi(x),\phi(y)\in F$, then $\phi[M^k_x]\cup \phi[M^k_y]=\phi[M^k_x\cup M^k_y]\subseteq (\phi(x))_k\cup (\phi(y))_k\subseteq F.$ 

\begin{claimproof}[Proof of $(c)$]Let $x,y\in X$ and suppose that $\inf(x,y)$ exists.  Take an arbitrary $k>\rho^{\phi[\{x,y,\inf(x,y)\}]}$. By means of the point $(1)$ of Lemma  \ref{Mklemma}, $M^k_x\cap M^k_y=M^k_{\inf(x,y)}$. In this way: \begin{align*}
T^k_x\cap T^k_y =A^{k+1}_{M^k_x}\cap A^{k+1}_{M^k_y}
=A^{k+1}_{M^k_x\cap M^k_y}
=A^{k+1}_{M^k_{\inf(x,y)}}=T^k_{\inf(x,y)}.
\end{align*}
Therefore, $T_x\cap T_y=\bigcup\limits_{k\in\omega} \big( T^k_x\cap T^k_y\big)=^*\bigcup \{T^k_{\inf(x,y)}\,:\, k>\rho^{\phi[x,y,\inf(x,y)]}\}=^*T_{\inf(x,y)}.$
\end{claimproof}
The remaining parts of the theorem are proved in a completely similar way. Because of this, we leave the calculations to the reader. 
\end{proof}
\end{theorem}
The following result greatly extends Hausdorff's theorem about the existence of an $(\omega_1,\omega_1)$-gap.
\begin{theorem}\label{allgapsthm}Let $(X,<_X)$ and $(Y,<_Y)$ be two partial orders with $cof(X)=cof(Y)=\omega_1$. Then there are cofinal $X'\subseteq X$ and $Y'\subseteq Y$ for which there is an $(X',Y')$-gap.
\begin{proof} Without loss of generality we may assume that $X\cap Y=\emptyset$. By applying Proposition \ref{propositionwellfounded} we get $X'\subseteq X$ and $Y'\subseteq Y$ two $\omega_1$-like cofinal subsets of $X$ and $Y'$ respectively. Let $Z=X'\cup Y'$ 
and $<_Z$ be the partial order over $Z$ induced by $X'$ and $Y'$. That is, $x<_Z y$ if and only if $x,y\in X'$ and $x<_X y$ or $x,y\in Y'$ and $x<_Y y$. It is straightforward that $(Z,<_Z)$ is an $\omega_1$-like order. 
Furthermore, if $x\in X'$ and $y\in Y'$ then there is no $z\in Z$ with $z\leq x, y.$ In virtue of Theorem \ref{luzinreptheorem}, there is a Luzin representation of $Z$, namely $(\mathcal{T},\mathcal{A})$. Due to the point (e) of Definition \ref{Luzinrepdef} it follows that $(\langle T_x\rangle_{x\in X'},\langle T_y\rangle_{y\in Y'})$ forms an $(X',Y')$-pregap. 
We claim that this is in fact a gap. For this, just note that any set separating $(\langle T_x\rangle_{x\in X'},\langle T_y\rangle_{y\in Y'})$ would also separate $(\langle A_x\rangle_{x\in X'}, \langle A_y\rangle_{y\in Y'})$ by the point (a) of Definition \ref{Luzinrepdef}. But this is impossible since $\mathcal{A}$ is Luzin (hence inseparable) and both $X'$ and $Y'$ are uncountable.
\end{proof}
\end{theorem}

\subsection{Donut-separability}
In this subsection we will explore in more depth the relation between almost disjoint families and gaps.
 When applying the argument of Theorem \ref{allgapsthm} to the orders $X=Y=\omega_1$ we get the following result.
\begin{theorem}\label{donutinseparablegap}There is an $(\omega_1,\omega_1)$-pregap $(L_\alpha,R_\alpha)_{\alpha\in\omega_1}$ such that the set $$\{L_{\alpha+1}\backslash L_\alpha\,:\,\alpha\in\omega_1\}\cup \{R_{\alpha+1}\backslash R_\alpha\,:\,\alpha\in\omega_1\}$$
forms a Luzin family. In particular $(L_\alpha,R_\alpha)_{\alpha\in\omega_1}$ is a gap.
\begin{proof}Let $Z=\omega_1\times 2$  and define the order over $Z$ given by $(\alpha,i)<_Z (\beta,j)$ if and only if $\alpha<\beta$ and $i=j$. $(Z,<_z)$ is an $\omega_1$-like order so that for any $\alpha,\beta\in \omega_1$, there is no $z\in Z$ with $z\leq (\alpha,0),(\beta,1)$. In particular, by Theorem \ref{luzinjonestheorem} there is a Luzin representation $(\mathcal{T},\mathcal{A})$ of $Z$. For any $\alpha\in\omega_1$, let $L_\alpha= T_{(\alpha,0)}$ 
and $R_{\alpha}=T_{(\alpha,1)}$. Then $(L_\alpha,R_\alpha)_{\alpha\in\omega_1}$ is a pregap due to the same arguments given in Theorem \ref{allgapsthm}. In this particular case, note that for any $\alpha\in\omega_1$ and $i\in 2$,  the set $pred((\alpha+1,i))$ satisfies the conditions of point (d) in Definition \ref{Luzinrepdef}. Even more, $pred((\alpha+1,i))=\{(\alpha,i)\}$. In this way, $L_{\alpha+1}\backslash L_{\alpha}=^*A_{(\alpha,0)}$ and $R_{\alpha+1}\backslash R_{\alpha}=^*A_{(\alpha,1)}$. We conclude that $$\{L_{\alpha+1}\backslash L_\alpha\,:\,\alpha\in\omega_1\}\cup \{R_{\alpha+1}\backslash R_\alpha\,:\,\alpha\in\omega_1\}$$ is an uncountable subset of $\mathcal{A}$. Therefore, it must also  be a Luzin family.
\end{proof}

\end{theorem}
The core idea of Theorems \ref{allgapsthm} and \ref{donutinseparablegap} is the following: If  $(\mathcal{L},\mathcal{R})$ is a pregap and for each $L\in \mathcal{L}$ and $R\in \mathcal{R}$ we can find infinite sets $A_L\subseteq L$ and $A_R\subseteq R$ such that the family $\{A_L\}_{L\in \mathcal{L}}\cup \{A_R\}_{R\in \mathcal{R}}$ forms an inseparable family, then the original pregap is in fact a gap. Such phenomenom was studied in \cite{hausdorffgapsreconstructed} by  Piotr Kalemba and Szymon Plewik for the case of $(\omega_1,\omega_1)$-gaps. The main difference between the gaps they studied and the one we constructed in the previous theorem, is that the associated Luzin family can be explicitly defined from such gap. This leads to the following definition.

\begin{definition}[donut-separable gaps]We call an $(\omega_1,\omega_1)$-gap $(L_\alpha,R_\alpha)_{\alpha\in\omega_1}$ \textit{donut-separable} if there is a set separating  $(L_{\alpha+1}\backslash L_\alpha, R_{\alpha+1}\backslash R_\alpha)_{\alpha\in\omega_1}$. If no such set exists we call the gap \textit{donut-inseparable}.
\end{definition}
\begin{corollary}There is a donut-inseparable Hausdorff gap.
\end{corollary}
In contrast, there is also (in ZFC) a Hausdorff gap which is donut-separable.

\begin{theorem}\label{firsthausdorffdonutseparable}The Hausdorff gap $(\mathcal{A},\mathcal{B})$ constructed in Theorem \ref{hausdorffgapconstruction} is donut-separable.
\begin{proof}Let $\alpha\in\omega_1$. Now consider an arbitrary $k> \rho(\alpha,\alpha+1)$. Note that $\lVert \alpha\rVert_k +1=\lVert \alpha+1 \rVert_k$. By the point (c) of Lemma \ref{lemmaxi} we know that either $\Xi_{\alpha+1}(k)=\Xi_{\alpha}(k)$ or these two numbers are different and $\Xi_{\alpha}(k)=-1$. It follows that in the latter case $r_k-1=\lVert \alpha\rVert_k<\lVert \alpha+1\rVert_k=r_k<m_k$, which means that $\Xi_{\alpha+1}(k)=0$.  From this fact we conclude that \begin{align*}L_{\alpha+1}\backslash L_{\alpha}&=^*\{2k+\Xi_{\alpha+1}(k)\,:k>\rho(\alpha,\alpha+1)\textit{ and }\Xi_{\alpha}(k)\not=\Xi_{\alpha+1}(k)\}\\
&\subseteq \{ 2k\,: k\in\omega\textit{ and }\Xi_{\alpha+1}(k)=0\}\subseteq \{2k\,: k\in\omega\}.
\end{align*}
In the same way, we have that $R_{\alpha+1}\backslash R_\alpha\subseteq ^* \{2k+1\,:k\in\omega\}.$
Hence, the set of even numbers separates $(L_{\alpha+1}\backslash L_\alpha,R_{\alpha+1}\backslash R_\alpha)_{\alpha\in\omega_1}$.
\end{proof}
\end{theorem}
\begin{lemma}Let $(\mathcal{L},\mathcal{R})$ be an $(\omega_1,\omega_1)$-gap and $\mathcal{L}'$, $\mathcal{R}'$ be cofinal subsets of $\mathcal{L}$ and $\mathcal{R}$ respectively. 
 If $(\mathcal{L},\mathcal{R})$ is donut-inseparable then so is $(\mathcal{L}',\mathcal{R}').$
 \end{lemma}
 The previous lemma motivates the following definition.
\begin{definition}[Strongly donut-separable gaps] Let $(\mathcal{L},\mathcal{R})$ be an $(\omega_1,\omega_1)$-gap. We say that $(\mathcal{L},\mathcal{R})$ is strongly \textit{donut-separable} if for any two cofinal subsets of $\mathcal{L}$ and $\mathcal{R}$, say $\mathcal{L}'$ and $\mathcal{R}'$, the gap $(\mathcal{L}',\mathcal{R}')$ is donut-separable. 
\end{definition}
It is natural to wonder if there is a strongly donut-separable gap. We will prove that this statement is independent from $ZFC.$
 \begin{lemma}\label{separabilityequivalence} Let $(L_\alpha,R_\alpha)_{\alpha\in\omega_1}$ be an $(\omega_1,\omega_1)$-gap. The following statements are equivalent:
 \begin{enumerate}[label=$(\arabic*)$]
 \item $(L_\alpha,R_\alpha)_{\alpha\in\omega_1}$ is strongly donut-separable. 
 \item For any $S\in[\omega_1]^{\omega_1}$, the gap $(L_\alpha,R_\alpha)_{\alpha\in S}$ is donut-separable.
 \item For any club $S\subseteq \omega_1$, the gap $(L_\alpha,R_\alpha)_{\alpha\in S}$ is donut-separable.
 \end{enumerate}
 
\begin{proof}The only nontrivial part of the proof is to show that $(3)$ implies $(1)$. For this, let $X$ and $Y$ be uncountable subsets of $\omega_1$. Our goal is to prove that there is a $C\in [\omega]^{\omega}$ which separates $(L_{X(\alpha+1)}\backslash L_{X(\alpha)},R_{Y(\alpha+1)}\backslash R_{Y(\alpha)})_{\alpha\in\omega_1}$. Let $\langle M_\alpha\rangle_{\alpha\in\omega_1}$ be a continuous chain of elementary submodels of a largely enough $H(\lambda)$ so that $X,Y\in M_0$. Now define $S$ as $\{0\}\cup \{ M_\alpha\cap \omega_1\,:\,\alpha\in \omega_1\}$. Then $S$ is a club. Furthermore, by elementarity it is straightforward that for any $\alpha\in\omega_1$ there is $\beta \in S$ with:
\begin{multicols}{2}\begin{itemize}
\item $S(\beta)\leq X(\alpha)<X(\alpha+1)\leq S(\beta+1)$.
\item $S(\beta)\leq Y(\alpha)<Y(\alpha+1)\leq S(\beta+1).$
\end{itemize}
\end{multicols}
By the hypotheses, there is $C\in [\omega]^{\omega}$ separating $(L_{S(\alpha+1)}\backslash L_{S(\alpha)}, R_{S(\alpha+1)}\backslash R_{S(\alpha)})_{\alpha\in\omega_1}$. Note that if $\alpha\in\omega_1$ and $\beta\in S$ are as previously stated then $L_{X(\alpha+1)}\backslash L_{X(\alpha)}\subseteq^* L_{S(\beta+1)}\backslash L_{S(\beta)}\subseteq^* C$ and $R_{Y(\alpha+1)}\backslash R_{Y(\alpha)}\subseteq^* R_{Y(\beta+1)}\backslash R_{Y(\beta)}\subseteq^* \omega\backslash C$. Therefore $C$ separates the pregap  $(L_{X(\alpha+1)}\backslash L_{X(\alpha)},R_{Y(\alpha+1)}\backslash R_{Y(\alpha)})_{\alpha\in\omega_1}$. This finishes the proof.
\end{proof}
\end{lemma}
\begin{corollary}\label{corocccpreservingdonut}Let $(\mathcal{L},\mathcal{R})=(L_\alpha,R_\alpha)_{\alpha\in\omega_1}$ be an $(\omega_1,\omega_1)$-gap which is strongly donut-separable. If $\mathbb{P}$ is a $ccc$ forcing, then $$\mathbb{P}\Vdash\text{\say{ $(\mathcal{L},\mathcal{R})$ is strongly donut-separable}}.$$
\begin{proof}Suppose towards a contradiction that there is $G$ a $\mathbb{P}$-generic filter over $V$ such that, in $V[G]$,  $(\mathcal{L},\mathcal{R})$ is not strongly donut-separable.  According to the point (3) of Lemma \ref{separabilityequivalence}, there is a club $S\in V[G]$ such that the gap $(L_\alpha,R_\alpha)_{\alpha\in S}$ is donut-inseparable. Since $\mathbb{P}$ is a $ccc$ forcing, we can find a club $C\in V$ for which $C\subseteq S$. is donut-inseparable and belongs to $V$. This contradiction ends the proof. 
\end{proof}
    
\end{corollary}
Suppose that $(\mathcal{L},\mathcal{R})$ is a pregap over $\omega$ and $C$ is an infinite subset separating it. Define $s:\mathcal{L}\cup\mathcal{R}\longrightarrow \omega$ as: $$s_C(X)=\begin{cases}
\min(n\,: X\backslash n\subseteq C)&\textit{ if }X\in \mathcal{L}\\
\min(n\,:X\backslash n \cap C=\emptyset)&\textit{ if }X\in \mathcal{R}
\end{cases}$$
Then $L\cap R\subseteq \max(s_C(L),s_C(R))$ for any $L\in \mathcal{L}$ and $R\in \mathcal{R}.$ This motivates the following definition.
\begin{definition}[Separating functions]Let $(\mathcal{L},\mathcal{R})$ be a pregap. We say that $s:\mathcal{L}\cup\mathcal{R}\longrightarrow \omega$ is \textit{separating} if $L\cap R\subseteq  \max(s(L),s(R))$ for any $L\in \mathcal{L}$ and $R\in \mathcal{R}$. 
\end{definition}
\begin{lemma}Let $(\mathcal{L},\mathcal{R})$ be a pregap. If there is a separating $s:\mathcal{L}\cup\mathcal{R}\longrightarrow \omega$ then there is $C\in [\omega]^{\omega}$ which separates $(\mathcal{L},\mathcal{R})$.
\begin{proof}Let $s$ be as in the hypotheses. We define $C$ as $$\{n\in\omega\,:\, \exists L\in \mathcal{L}\,(n\in L\textit{ and }s(L)<n)\}=\bigcup\limits_{L\in \mathcal{L}}L\backslash (s(L)+1).$$
Note that $L\subseteq^* C$ for any $L\in \mathcal{L}$. For $R\in \mathcal{R}$ we claim that $R\backslash (s(R)+1)\cap C=\emptyset$ which in particular implies that $R\cap C=^* \emptyset$. Suppose towards a contradiction that this is not the case and let $n$ be an element in the intersection of both sets. On one hand, since $n\in C$ there is $L\in \mathcal{L}$ with $n\in L$ and $s(L)<n$. In particular, $n\in L\cap R$.  On the other hand, since $s$ is separating then $L\cap R\subseteq \max(s(R),s(L))\subseteq n$. This is a contradiction. Therefore the claim is true, which means that $C$ separates $(\mathcal{L},\mathcal{R})$. 
\end{proof}
\end{lemma}

\begin{definition}\label{separabilityforcingdonuts}Let $(\mathcal{D},\mathcal{E})=(D_\alpha,E_\alpha)_{\alpha\in\omega_1}$ be a pregap (not necessarily of type $(\omega_1,\omega_1)$) so that $D_\alpha\cap E_\alpha=\emptyset$ for each $\alpha$. We define the forcing $\mathbb{P}(\mathcal{D},\mathcal{E})=\mathbb{P}(D_\alpha,E_\alpha)_{\alpha\in\omega_1}$ as the set of all functions $p;\omega_1\longrightarrow \omega$ with finite domain and such that $D_\alpha\cap E_\beta\subseteq  \max(p(\alpha),p(\beta))$ for all $\alpha,\beta\in dom(p)$. The order is given by $$p\leq q\text{ if and only if }q\subseteq p.$$
    
\end{definition}

\begin{rem}\label{remdenseseparating}If $(\mathcal{D},\mathcal{E})=(D_\alpha,E_\alpha)_{\alpha\in \omega_1}$ is as in the previous definition and $\beta\in \omega_1$ then $\mathcal{M}_\beta=\{p\in \mathbb{P}(\mathcal{D},\mathcal{E})\,:\,\beta\in dom(p)\,\}$ is dense in 
$\mathbb{P}(\mathcal{E},\mathcal{D})$. If $\mathcal{G}$ is a filter intersecting each $\mathcal{M}_\beta$ then the function $s:\mathcal{D}\cup \mathcal{E}\longrightarrow\omega_1$ given by: $$s(D_\alpha)=s(E_\alpha)=\bigcup\mathcal{G}(\alpha)$$
is well defined and separating.
\end{rem}
The following proposition generalizes a well-known result of Kenneth Kunen (we reiterate that our pregaps do not need to be linearly ordered).

\begin{proposition}\label{separabilitygapsprop} Let $(E_\alpha,D_\alpha)_{\alpha\in\omega_1}$ be a pregap  with $E_\alpha\cap D_\alpha=\emptyset$ for each $\alpha$. The following statements are equivalent:
\begin{enumerate}[label=$(\alph*)$]
\item $\mathbb{P}=\mathbb{P}(E_\alpha,D_\alpha)_{\alpha\in\omega_1}$ is $ccc$.
\item There is $\mathbb{Q}$ $ccc$ with $\mathbb{Q}\Vdash \text{\say{ $(E_\alpha,D_\alpha)_{\alpha\in\omega_1}\textit{ can be separated }$}}$.
\item There is $W$ a transitive model of $ZFC$ extending $V$ with $\omega_1^W=\omega_1^V$ where $$W\models (E_\alpha,D_\alpha)_{\alpha\in\omega_1}\textit{ can be separated.}$$
\end{enumerate}
\begin{proof}Trivially $(a)$ implies $(b)$ and $(b)$ implies $(c)$. In order to prove that $(c)$ implies $(a)$ take an arbitrary uncountable subset $\mathbb{P}$ in $V$, say $\mathcal{A}$. For  any $p\in \mathcal{A}$ let $d_p=dom(p)$. Without loss of generality we can suppose that $\mathcal{A}$ is a $\Delta$-system with root $R$, there is $n\in\omega$ so that $|p|=n$ for every $p\in\mathcal{A}$ and $p(d_p(i))=q(d_q(i))$ for all $p,q\in \mathcal{A}$ and each $i<n$.

Let $W$ be as in the hypotheses of $(c)$ and let $C\in W$ which separates $(E_\alpha,D_\alpha)_{\alpha\in\omega_1}$. Since $\omega_1^V=\omega_1^W$ and $V$ models that $\mathcal{A}$ is uncountable, we can find an uncountable $\mathcal{A}'\subseteq \mathcal{A}$ in $W$, so that $E_{d_p(i)}\backslash C=E_{d_q(i)}\backslash C$ and $D_{d_p(i)}\cap C= D_{d_q(i)}\cap C$ for all $p,q\in \mathcal{A}'$ and $i<n$. Fix two distinct $p,q\in\mathcal{A}'$. We claim that $p\cup q\in \mathbb{P}$. Indeed, let $i,j<n$. Then \begin{align*}E_{d_p(i)}\cap D_{d_q(j)}&=((E_{d_p(i)}\backslash C)\cap D_{d_q(j)})\cup (E_{d_p(i)}\cap (D_{d_q(j)}\cap C))\\
&=((E_{d_q(i)}\backslash C)\cap D_{d_q(j)})\cup (E_{d_p(i)}\cap (D_{d_p(j)}\cap C))\\
&\subseteq \max(q(d_q(i)),q(d_q(j))\cup \max(p(d_p(i)),p(d_p(j)))\\
&= \max(p(d_p(i)),p(d_q(j))).
\end{align*}
This finishes the proof.
\end{proof}

\end{proposition}
Consider the gap $(L_\alpha, R_\alpha)_{\alpha\in\omega}$ constructed in Theorem \ref{hausdorffgapconstruction} and for any $k\in\omega$ let $N_k=2(k+1)=\{0,\dots,2k+1\}$. The following properties follow directly from the proof of that theorem and the definition of the $\Delta$-function.
\begin{proposition}\label{propositioninterhausdorf}Let $\alpha<\beta\in \omega_1$. Then:
\begin{enumerate}[label=$(\arabic*)$]
\item[$(0)$]$L_\alpha\cap R_\alpha=\emptyset.$
\item For each $k\in\omega$, both $\{2k,2k+1\}\cap L_\alpha$ and $\{2k,2k+1\}\cap R_\alpha$ have at most one point.  
\item All of the sets $L_\alpha\backslash L_\beta$,  $R_\alpha\backslash R_\beta$,  $L_\alpha\cap R_\beta$  and $L_\beta\cap R_\alpha$ are subsets of $N_{\rho(\alpha,\beta)}$.
\item $L_\alpha \cap N_{\Delta(\alpha,\beta)-1}=L_\beta\cap N_{\Delta(\alpha,\beta)-1}$ and  $R_\alpha \cap N_{\Delta(\alpha,\beta)-1}=R_\beta\cap N_{\Delta(\alpha,\beta)-1}$.

\end{enumerate}
\end{proposition}
Let us briefly recall Definition \ref{parametrizedmartinsdef}:  If $\mathcal{F}$ is a  $2$-capturing construction scheme over $\omega_1$, the parametrized Martin's number $\mathfrak{m}_\mathcal{F}$ associated to $\mathcal{F}$, is the minimun of the Martin's numbers of $ccc$ forcing notions which force $\mathcal{F}$ to be $2$-capturing in the generic extension.\\
In the following theorem, we will show that the existence of $(\omega_1,\omega_1)$-gaps which are strongly donut-separable is consistent with $ZFC$.  The proof we provide can be simplified by the use of Lemma \ref{lemmaequivalencefunctioncapturingpreserving}. We decided to keep it this way in order to motivate such lemma. 
\begin{theorem}\label{mfstronglydonuttheorem} Let $\mathcal{F}$ be a $2$-capturing $2$-construction scheme for which $\mathfrak{m}_\mathcal{F}>\omega_1$. The Hausdorff gap constructed from $\mathcal{F}$ in Theorem \ref{hausdorffgapconstruction}  is strongly donut-separable.
\begin{proof}Let $S\in [\omega_1]^{\omega_1}$. We will show that the gap $(L_{S(\alpha)},R_{S(\alpha)})_{\alpha\in\omega_1}$ is donut-separable. In other words, we will show that pregap  $(D_\alpha, E_\alpha)_{\alpha\in \omega_1}$ can be separated where $D_\alpha=L_{S(\alpha+1)}\backslash L_{S(\alpha)} $ and  $ E_\alpha=R_{S(\alpha+1)}\backslash R_{S(\alpha)}$ for each $\alpha\in \omega_1$. This is enough due to the Lemma \ref{separabilityequivalence}. In view of Remark \ref{remdenseseparating} and since we are assuming that $\mathfrak{m}_\mathcal{F}>\omega_1$, it is sufficient to prove that for $\mathbb{P}=\mathbb{P}(D_\alpha,E_\alpha)_{\alpha\in \omega_1}$ we have:
\begin{enumerate}[label=$(\arabic*)$]
\item $\mathbb{P}$ is a $ccc$ forcing,
\item $\mathbb{P}\Vdash\text{\say{ $ \mathcal{F}\textit{ is }2\textit{-capturing}$ }}.$
\end{enumerate}

\begin{claimproof}[Proof of $(1)$] For this consider an arbitrary $\mathcal{A}\in[\mathbb{P}]^{\omega_1}$. Given $p\in \mathcal{A}$ define $Z_p=\{ S(\alpha),S(\alpha+1)\,:\,\alpha\in dom(p)\}$. By refining $\mathcal{A}$ we can assume without loss of generality that for any $p,q\in \mathcal{A}$ the following conditions hold:
\begin{enumerate}[label=$(\alph*)$]
    \item $|p|=|q|$,
    \item If $f:dom(p)\longrightarrow dom(q)$ is the increasing bijection then $p(\alpha)=q(f(\alpha))$ for each $\alpha\in dom(p).$ In particular, there is $k\in\omega$ with $k>\max(im(p))$ for each $p\in \mathcal{A}$.
\end{enumerate}
Furthermore, we can suppose that the set $\{dom(p)\,:\,p\in \mathcal{A}\}$ forms a root-tail-tail $\Delta$-system with root $R$ satisfying the following properties for any two distinct $p,q\in \mathcal{A}$:
\begin{enumerate}
    \item[(c)] If $R\not=\emptyset$, then $\max(R)+1<\min(dom(p)\backslash R)$.
    \item [(d)] $\max(dom(p))+1\not\in dom(q)$.
\end{enumerate}
 As the set $\{Z_p\,:\,p\in \mathcal{A}\}$ is uncountable and $\mathcal{F}$ is assumed to be $2$-capturing then there are distinct $p,q\in \mathcal{A}$ for the set $\{Z_p,Z_q\}$ is captured at some level $l>k$. We affirm that $p$ and $q$ are compatible. This will follow from the next claim.\\\\
 \underline{Claim 1:} $p\cup q$ is a condition of $\mathbb{P}$.
 
 \begin{claimproof}[Proof of claim] First note that if $F\in \mathcal{F}_l$ is such that $Z_p\cup Z_q\subseteq F$ then $Z_p\cap Z_q=Z_p\cap R(F)=Z_q\cap R(F)$ by means of Lemma \ref{intersectionrhoisomorphiclemma}. 

Now, let us consider $f:dom(p)\longrightarrow dom(q)$  and $h:F_0\longrightarrow F_1$ the increasing bijections. In order to prove that $p\cup q$ is a condition it is enough to take $\alpha\in dom(p)$ and $\beta\in dom(q)$ and show that both $E_\alpha\cap D_\beta$ and $D_\alpha\cap E_\beta$ are contained in $\max(p(\alpha),q(\beta))$. If either $\alpha$ or $\beta$ belong to the intersection of $dom(p)\cap dom(q)$ there is nothing to do. So we can assume that $\alpha\in dom(p)\backslash dom(q)$ and $\beta\in dom(q)\backslash dom(p)$. 
Then $\{S(\alpha),S(\alpha+1)\} \subseteq Z_p\backslash Z_q\subseteq F_0\backslash R(F)$ and $\{S(\beta),S(\beta+1)\}\subseteq Z_q
\backslash Z_p\subseteq F_1\backslash R(F)$ due to the points (c) and (d).  Thus $\rho(S(\alpha+1),S(\beta+1))=l$ and consequently $D_\alpha\cap E_\beta\subseteq L_{S(\alpha+1)}\cap R_{S(\beta+1)}\subseteq N_l$ due to the point (2) of Proposition \ref{propositioninterhausdorf}.  But $\Xi_{S(\alpha)}(l)=\Xi_{S(\alpha+1)}(l)=0$. 
Therefore $2l\in L_{S(\alpha)}\cap L_{S(\alpha+1)}$. In particular $2l\in L_{S(\alpha+1)}$,  so $2l+1\not\in L_{S(\alpha+1)}$ by the point (1) of Proposition \ref{propositioninterhausdorf}. In this way  $\{2l,2l+1\}\cap D_\alpha=\emptyset$. Thus, $$D_\alpha\cap E_\beta\subseteq N_{l-1}.$$  
The next thing to note is that $h(S(\alpha))=S(f(\alpha))$ and $h(S(\alpha+1))=S(f(\alpha)+1).$ This means that $\Delta(S(\alpha),S(f(\alpha)))=l=\Delta(S(\alpha+1),S(f(\alpha)+1))$. In virtue of the point (3) of proposition \ref{propositioninterhausdorf} we have $$L_{S(\alpha)}\cap N_{l-1}=L_{S(f(\alpha))}\cap N_{l-1},$$
$$L_{S(\alpha+1)}\cap N_{l-1}=L_{S(f(\alpha)+1)}\cap N_{l-1}.$$
Hence, $D_\alpha\cap N_{l-1}=D_{f(\alpha)}\cap N_{l-1}$. From all the equations we have so far we deduce that $D_\alpha\cap E_\beta=D_{f(\alpha)}\cap E_\beta\subseteq \max(q(f(\alpha)),q(\beta))=\max(p(\alpha),q(\beta)).$ In a completely similar way we can show that $E_\alpha\cap D_\beta\subseteq \max(p(\alpha),q(\beta))$. We conclude that $p\cup q\in \mathbb{P}$ so we are done.
\end{claimproof}
\end{claimproof}

\begin{claimproof}[Proof of $(2)$]For this we will show that the equivalence of $2$-capturing stated in Lemma \ref{equivalencecapturing} is forced by $\mathbb{P}$. Let $\dot{X}$ be a name so that $\mathbb{1}_\mathbb{P}\Vdash \text{\say{ $\dot{X}\in [\omega_1]^{\omega_1}$}}$ and let $q\in \mathbb{P}$. We shall find $p\leq q$ so that $$p\Vdash\text{\say{ $\exists M\in[\dot{X}]^2\,(\, M\textit{ is captured })$ }} .$$ As $\dot{X}$ is forced 
to be  uncountable we can find for each $\alpha\in \omega_1$ a condition $p_\alpha\leq q$ and $\xi_\alpha>\alpha$ with $p_\alpha\Vdash \text{\say{ $\xi_\alpha\in \dot{X}$ }}.$ We can suppose without loss of generality that the family $\mathcal{A}=\{p_\alpha\,:\,\alpha\in \omega_1\}$ satisfies the conditions (a), (b), (c) and (d) of the previous paragraphs. We may also assume that for any distinct $\alpha,\beta\in \omega_1:$
\begin{enumerate}
\item[$(d)$]$\xi_\alpha\not=\xi_\beta$ and $|Z_{p_\alpha}\cup \{\xi_\alpha\}|=|Z_{p_\beta}\cup \{\xi_\beta\}|$,
\item[$(e)$] If $g: Z_{p_\alpha}\cup \{\xi_\alpha\}\longrightarrow Z_{p_\beta}\cup \{\xi_\beta\}$ is the increasing bijection then $g[Z_{p_\alpha}]=Z_{p_\beta}$ and $g(\xi_\alpha)=\xi_\beta.$
\end{enumerate}
Since $\mathcal{F}$ is $2$-capturing there are $\alpha<\beta\in \omega_1$ for which $\{Z_{p_\alpha}\cup \{\xi_\alpha\},Z_{p_\beta}\cup \{\xi_\beta\}\}$ is captured at some level $l$ greater than $k$. In virtue of the conditions (d) and (e) it is easy to see that $\{Z_{p_\alpha},Z_{p_\beta}\}$ and $\{\xi_\alpha,\xi_\beta\}$ are also captured at level $l$. As in the first part of the proof, this implies that $p=p_\alpha\cup p_\beta$ is a condition of $\mathbb{P}$. To finish, just note that $p\leq q$ and  $p\Vdash\text{\say{ $\{\xi_\alpha,\xi_\beta\}\in [\dot{X}]^2\textit{ and }\{\xi_\alpha,\xi_\beta\}\textit{ is captured }$}}.$ This finishes the proof.
\end{claimproof}
\end{proof}
\end{theorem}
\begin{definition}[Adequate sets]\label{adequatesetsgap}Let $(\mathcal{L},\mathcal{R})=(L_\alpha,R_\alpha)_{\alpha\in\omega_1}$ be an $(\omega_1,\omega_1)$-gap. We say that $X\in[\omega_1]^{\omega_1}$ is \textit{adequate} for $(\mathcal{L},\mathcal{R})$ if for any $Y\in[\omega_1]^{\omega_1}$, the pregap $$(L_{X(\alpha+1)}\backslash L_{X(\alpha)},R_{X(\alpha+1)}\backslash R_{X(\alpha)})_{\alpha\in Y}$$ is a gap.

\end{definition}
\begin{rem}If an $(\omega_1,\omega_1)$-pregap admits an adequate set, then such pregap is not strongly donut-separable.
\end{rem}
\begin{theorem}[Under $CH$]\label{weaklydonutch} Let $(\mathcal{L},\mathcal{R})=(L_\alpha,R_\alpha)_{\alpha\in\omega_1}$ be an $(\omega_1,\omega_1)$-gap. Then there is a club $X\in [\omega_1]^{\omega_1}$ adequate for $(\mathcal{L},\mathcal{R}).$ In particular, there are no strongly-donut separable gaps.
\begin{proof}Since we are assuming $CH$ we can enumerate $[\omega]^{\omega}$ as $\langle C_\alpha\rangle_{\alpha\in \omega_1}$. We will build $X$ be recursion in such way that for any $\beta\in \omega_1$ and each $\alpha\leq \beta$ one of the following conditions occur:
\begin{center}\begin{minipage}{5cm} \begin{center} \textbf{(A)}\end{center} $$L_{X(\beta+1)}\backslash L_{X(\beta)}\not\subseteq^* C_\alpha.$$
\end{minipage}\hspace{2cm} \begin{minipage}{5cm}\begin{center} \textbf{(B)}\end{center} $$(R_{X(\beta+1)}\backslash R_{X(\beta)})\cap C_\alpha\not=^*\emptyset.$$
\end{minipage}
\end{center}
If $\gamma$ is limit and we have constructed $X(\alpha)$ for any $\alpha<\gamma$ we just define $X(\gamma)$ as $\sup(X(\alpha)\,:\,\alpha<\gamma)$. The interesting case happens when we have constucted $X(\beta)$ for some $\beta\in \omega_1$ and we want to define $X(\beta+1).$ For this case, we have the following claim:\\\\
\underline{Claim}: The pregap $(L_\delta\backslash L_{X(\beta)}, R_\delta\backslash R_{X(\beta)})_{\delta>X(\beta)}$ is a gap.
\begin{claimproof}[Proof of claim] Suppose towards a contradiction that this is not the case and let $C$ be a set separating $(L_\delta\backslash L_{X(\beta)}, R_\delta\backslash R_{X(\beta)})_{\delta>X(\beta)}$. Then $C\cup L_{X(\beta)}$ separates $(\mathcal{L},\mathcal{R})$ which is a contradiction. 
\end{claimproof}
Fix $\alpha\leq \beta$. By the previous claim, $C_\alpha$ does not separate   $(L_\delta\backslash L_{X(\beta)}, R_\delta\backslash R_{X(\beta)})_{\delta>X(\beta)}$. In this way, there is $\delta_\alpha >X(\beta)$ so that either $L_{\delta_\alpha}\backslash L_{X(\beta)}\not\subseteq^*C_\alpha$ or  $R_{\delta_\alpha}\cap C_\alpha\not=^*\emptyset.$ Let us define $X(\beta+1)$ as $\sup(\delta_\alpha\,:\alpha\leq \beta)$. Then $L_{\delta_\alpha}\backslash L_{X(\beta)}\subseteq L_{X(\beta+1)}\backslash L_{X(\beta)}$ and $R_{\delta_\alpha}\backslash R_{X(\beta)}\subseteq R_{X(\beta+1)}\backslash R_{X(\beta)}$ for each $\alpha\leq \beta$. In this way we guarantee that either contition (A) or condition (B) will ocurr for any such $\alpha.$ This finishes the recursion.

We will now prove that $X$ is as desired. Let $Y\in [\omega_1]^{\omega_1}$. We need to show that $(L_{X(\alpha+1)}\backslash L_{X(\alpha)},R_{X(\alpha+1)}\backslash R_{X(\alpha)})_{\alpha\in Y}$ is a gap.
For this let $C$ be an infinite subset of $\omega$. By the assumptions we know there is $\alpha\in \omega_1$ so that $C=C_\alpha$. Let $\beta\in Y$ be such that $\beta>\alpha$. Therefore, either $L_{X(\beta+1)}\backslash L_{X(\beta)}\not\subseteq^*C_\alpha$ or $(R_{X(\beta+1)}\backslash R_{X(\beta)})\cap C_\alpha\not=^*\emptyset.$ In any case, $X(\beta)$ testifies that $C_\alpha$ does not separate the gap that we are considering. This finishes the proof.
\end{proof}
\end{theorem}
By combining the previous result and Theorem \ref{mfstronglydonuttheorem}, we conclude the following:
\begin{theorem}The statement \say{ There is a strongly donut-separable gap} is independent from $ZFC.$ 
\end{theorem}
Now we will analyze the relationship between donut-separable gaps and the Cohen forcing.

\begin{lemma}\label{adequatecohenpreserving} Let $(\mathcal{L},\mathcal{R})=(L_\alpha,R_\alpha)_{\alpha\in\omega_1}$ be an $(\omega_1,\omega_1)$-gap and let $X\in [\omega_1]^{\omega_1}$. If $X$ is adequate for $(\mathcal{L},\mathcal{R})$ then  $\mathbb{C}\Vdash \text{\say{ $X\textit{ is adequate for }(\mathcal{L},\mathcal{R})$}}.$
\begin{proof}Let $G$ be a $\mathbb{C}$-generic filter over $V$ and let $Y\in [\omega_1]^{\omega_1}\cap V[G]$. Suppose towards a contradiction that there is $C\in V[G]$ separating the pregap $$(L_{X(\alpha+1)}\backslash L_{X(\alpha)},R_{X(\alpha+1)}\backslash R_{X(\alpha)})_{\alpha\in Y}.$$ Since $Y$ is uncountable we can find $a,b\in [\omega]^{<\omega}$  and $Y'\in [Y]^{\omega_1}\cap V[G]$ so that for any $\alpha\in Y'$, $(L_{X(\alpha+1)}\backslash L_{X(\alpha)})\backslash C= a.$
Due to well-known facts concerning the Cohen forcing $\mathbb{C}$, we know there is $Z\in [\omega_1]^{\omega_1}\cap V$ with $Z\subseteq Y$. Let us define $$C'=\bigcup\limits_{\alpha\in Z}L_{X(\alpha+1)}\backslash L_{X(\alpha)}.$$
Then $C'\in V$ and trivially $L_{X(\alpha+1)}\backslash L_{X(\alpha)}\subseteq C'$ for any $\alpha\in Z$. Furthermore, $C'\subseteq C\cup a$. Therefore, $(R_{X(\alpha+1)}\backslash R_{X(\alpha)})\cap C'$ is finite for any $\alpha\in Z.$ We conclude that, in $V$, there is a set separating the gap $(L_{X(\alpha+1)}\backslash L_{X(\alpha)},R_{X(\alpha+1)}\backslash R_{X(\alpha)})_{\alpha\in Z}.$ This is a contradiction to the hypotheses. Hence, the proof is over.
\end{proof}

\end{lemma}
\begin{theorem}[Under $CH$] Let $\kappa$ be an uncountable cardinal. Then $$\mathbb{C}_\kappa\Vdash\text{\say{ There are no  strongly donut-separable $(\omega_1,\omega_1)$-gaps }}.$$
In particular, the statement \say{There are no strongly  donut-separable $(\omega_1,\omega_1)$-gaps} is consistent with an arbitrarily large continuum.
\begin{proof}If $\kappa=\omega_1$ the argument is clear, since in this case $\mathbb{C}_\kappa$ forces $CH$. So suppose that $\kappa>\omega_1$ and let $G$ be a $\mathbb{C}_\kappa$-generic filter over $V$. Finally, let $(\mathcal{L},\mathcal{R})=(L_\alpha,R_\alpha)_{\alpha\in \omega_1}\in V[G]$ be an arbitrary $(\omega_1,\omega_1)$-gap. Since $|\mathcal{L}|=|\mathcal{R}|=\omega_1$, then there is a $\mathbb{C}_{\omega_1}$-generic filter over $V$, namely $H$, and  a $\mathbb{C}_{\kappa}$-generic filter over $V[H]$, namely $K$, so that $V[H][K]=V[G]$ and $(\mathcal{L},\mathcal{R})\in V[H].$ By the hypotheses, $V$ is a model of $CH$. Therefore $V[H]$ models $CH $ too. According to Theorem \ref{weaklydonutch} there is $X\in V[H]$ which is adequate for $(\mathcal{L},\mathcal{R}).$ We finish by proving the following claim.\\\\
\underline{Claim}: $X$ testifies that $(\mathcal{L},\mathcal{R})$ is not strongly donut-separable.
\begin{claimproof}[Proof of claim.] Suppose towards a contradiction that there is $C\in V[G]$ separating the pregap $$(L_{X(\alpha+1)}\backslash L_{X(\alpha)}, R_{X(\alpha+1)}\backslash R_{X(\alpha)})_{\alpha\in \omega_1}.$$
Again, since $|C|=\omega$ we can find $\mathbb{C}$-generic filter over $V[H]$, say $H'$, for which $C\in V[H][H']$. By virtue of Lemma \ref{adequatecohenpreserving}, $V[H][H']\models X\textit{ is adequate for }(\mathcal{L},\mathcal{R}).$ In particular $C$ can not separate the pregap that we are considering. This contradiction finishes the proof.
\end{claimproof}
\end{proof}
\end{theorem}
As $\mathfrak{m}_\mathcal{F}>\omega_1$ is consistent with $\mathfrak{c}$ arbitrarily large, we conclude that both the statement \say{there is a strongly donut-separable gap} and its negation are consistent with an arbitrarily large continuum.\\
Our next goal is to show that $PFA$ also implies that there are no strongly donut-separable gaps. For this sake, we need to consider the following definitions.
\begin{definition}[Almost Luzin gaps] Let $(\mathcal{D},\mathcal{E})=(D_\alpha,E_\alpha)_{\alpha\in\omega_1}$ be pregap (not necessarily of type $(\omega_1,\omega_1))$. We say that $(\mathcal{D},\mathcal{E})$ is \textit{almost Luzin} if $(\mathcal{D}',\mathcal{E}')$ is a gap for any two uncountable $\mathcal{D}'\subseteq \mathcal{D}$ and $\mathcal{E}'\subseteq \mathcal{E}$.
\end{definition}
\begin{definition}[Highly adequate sets]Let $(\mathcal{L},\mathcal{R})=(L_\alpha,R_\alpha)_{\alpha\in\omega_1}$ be an $(\omega_1,\omega_1)$-gap. We say that $X\in [\omega_1]^{\omega_1}$ is \textit{highly adequate} for $(\mathcal{L},\mathcal{R})$ if the pregap $$(L_{X(\alpha+1)}\backslash L_{X(\alpha)},R_{X(\alpha+1)}\backslash R_{(X(\alpha)})_{\alpha\in X}$$
is almost Luzin.
\end{definition}
\begin{remark}If $X$ is highly adequate for $(\mathcal{L},\mathcal{R})$ then it is also adequate.  
\end{remark}
It is very easy check that the dichotomy satisfied by club $X$ that we constructed in Theorem \ref{weaklydonutch} already turns $X$ into a highly adequate set for the gap $(\mathcal{L},\mathcal{R})$. Hence, we have the following proposition.
\begin{proposition}[Under $CH$]\label{highlyadequateprop} Let $(\mathcal{L},\mathcal{R})=(L_\alpha,R_\alpha)_{\alpha\in\omega_1}$ be an $(\omega_1,\omega_1)$-gap. Then there is a club $X\in [\omega_1]^{\omega_1}$ highly adequate for $(\mathcal{L},\mathcal{R})$.
\end{proposition}

\begin{definition}[Normal pregaps] Suppose that $
(D_\alpha,E_\alpha)_{\alpha\in\omega_1}$ is a pregap. We say that it is \textit{normal} if $D_\alpha\cap E_\alpha=\emptyset$ for each $\alpha\in\omega_1.$
    
\end{definition}
\begin{rem}\label{normalgaprem}Suppose that $
(D_\alpha,E_\alpha)_{\alpha\in\omega_1}$ is a pregap and consider, for each $\alpha$, $D'_\alpha=D_\alpha\backslash E_\alpha$ and $E'_\alpha=E_\alpha\backslash D_\alpha$. Then $(D'_\alpha,E'_\alpha)_{\alpha\in\omega_1}$ is a normal pregap. Furthermore, both gaps induce the same sets in $\mathscr{P}(\omega)/\text{FIN}$. In particular, $
(D_\alpha,E_\alpha)_{\alpha\in\omega_1}$ is strongly donut-separable if and only if the same holds for $
(D'_\alpha,E'_\alpha)_{\alpha\in\omega_1}$.
\end{rem}

\begin{definition}[Biorthogonal gaps] Let $(\mathcal{D},\mathcal{E})=(D_\alpha,E_\alpha)_{\alpha\in \omega_1}$ be a normal pregap. We say that $(\mathcal{D},\mathcal{E})$ is \textit{biorthogonal} if $(D_\alpha\cap E_\beta)\cup (D_\beta\cap E_\alpha)\not=\emptyset$ for all $\alpha\not=\beta\in \omega_1$.
\end{definition}
\begin{lemma}\label{biorthogonalgaplemma}If $(\mathcal{D},\mathcal{E})=(E_\alpha,D_\alpha)_{\alpha\in\omega_1}$ is a normal biorthogonal pregap, then it is a gap.
\begin{proof}Suppose towards a contradiction that there is $C\in [\omega]^{\omega}$ which separates $(\mathcal{D},\mathcal{E})$. According to the pigheonhole principle, there is $X\in [\omega_1]^{\omega_1}$ for which $D_\alpha\backslash C=D_\beta\backslash C$ and $E_\alpha\cap C=E_\beta\cap C$ for all $\alpha,\beta\in X.$ Fix two distinct $\alpha,\beta\in X$. Then \begin{align*}D_\alpha\cap E_\beta&= (D_\alpha\cap (E_\beta\cap C))\cup (E_\beta\cap (D_\alpha\backslash C))\\&=(D_\alpha\cap (E_\alpha\cap C))\cup (E_\beta\cap (D_\beta\backslash C))=\emptyset
\end{align*}
By symmetry, we also have that $D_\beta\cap E_\alpha=\emptyset$, but this is a contradiction to the normality of $(\mathcal{D},\mathcal{E})$. Thus, the proof is over.
\end{proof}
    
\end{lemma}

\begin{definition}Let $(\mathcal{D},\mathcal{E})=(D_\alpha,E_\alpha)_{\alpha\in \omega_1}$ be a normal gap. We define the forcing $\mathbb{L}(\mathcal{D},\mathcal{E})$ as the set of all $p\in [\omega_1]^{<\omega}$ so that $(D_\alpha\cap E_\beta)\cup (D_\beta\cap E_\alpha)\not=\emptyset$
for all $\alpha\not=\beta\in dom(p).$ The order is given by  $$p\leq q\text{ if and only if }q\subseteq p.$$
    
\end{definition}

\begin{proposition}Let $(\mathcal{D},\mathcal{E})=(D_\alpha,E_\alpha)_{\alpha\in \omega_1}$ be a normal gap.  If $(\mathcal{D},\mathcal{E})$ is almost Luzin, then $\mathbb{L}(\mathcal{D},\mathcal{E})$ is $ccc$. 
\begin{proof}
Let $\mathcal{A}$ be an uncountable subset of $\mathbb{L}(\mathcal{D},\mathcal{E})$. We will show that $\mathcal{A}$ is not an antichain.  By refining $\mathcal{A}$, we may assume without loss of generality that there are $n,m\in \omega$ for which the following properties hold for each $p,q\in \mathcal{A}$:
\begin{enumerate}[label=$(\arabic*)$]
\item $|p|=n$.
\item For all $i,j<n$, $(D_{p(i)}\cap E_{p(j)})\cup (D_{p(j)}\cap E_{p(i)})\subseteq m$.
\item For all $i<n$, $D_{p(i)}\cap m=D_{q(i)}\cap m$ and $D_{p(i)}\cap m=D_{q(i)}\cap m$.\end{enumerate}
Note that for any two given conditions in $p,q\in \mathbb{L}(\mathcal{D},\mathcal{E})$, we have that $p$ and $q$ are compatible if and only if $p\backslash q$ and $q\backslash p$ are compatible conditions. Because of this and due to the $\Delta$-system Lemma, we may also assume that the elements of $\mathcal{A}$ are pairwise disjoint. Now, we proceed to find  distinct $p,q\in \mathcal{A}$ for which $p\cup q$ is a condition. This will be done after proving the following claims. \\

\noindent
\underline{Claim 1}: Let $p,q\in \mathcal{A}$ and $i\not=j<n$, then $(D_{p(i)}\cap E_{q(j)})\cup (D_{q(j)}\cap E_{p(i)})\not=\emptyset.$
\begin{claimproof}[Proof of claim]Just note that by the conditions (2) and (3), we have the following chain of equalities: \begin{align*}
   ( (D_{p(i)}\cap E_{q(j)})\cup (D_{q(j)}\cap E_{p(i)}))\cap m&=(D_{p(i)}\cap (E_{q(j)}\cap m))\cup ((D_{q(j)}\cap m)\cap E_{p(i)})\\&
    =(D_{p(i)}\cap (E_{p(j)}\cap m))\cup ((D_{p(j)}\cap m)\cap E_{p(i)})\\& 
    =(D_{p(i)}\cap E_{p(j)})\cup (D_{p(j)}\cap E_{p(i)})\not=\emptyset.&    \end{align*}
    This proves the claim.    
\end{claimproof}

\noindent
\underline{Claim 2}: Let $\mathcal{A}'\in [\mathcal{A}]^{\omega_1}$ and $i<n$. Then there is $\mathcal{X}\in [\mathcal{A}']^{\omega_1}$ so that for any $k\in \omega$, if  $\{p\in \mathcal{X}\,:\,k\in D_{p(i)}\}$  is non-empty, then it is uncountable. Analogously with the set $\{p\in \mathcal{X}\,:\,k\in E_{p(i)}\}$.
\begin{claimproof}[Proof of claim] Let $M$ be a countable elementary submodel of $H(\omega_2)$ such that $\mathcal{A}'\in M$. We put $\mathcal{X}= \mathcal{A}'\backslash M$. Note that if $k\in \omega$ is such that $k\in D_{q(i)}$ for some $q\in \mathcal{X}$, then the set $\{p\in \mathcal{A}'\,:\,k\in D_{p(i)}\}$ is an element of $M$ which is not contained in it. By elementarity, it follows that this set is uncountable. Therefore  $\{p\in \mathcal{X}\,:\,k\in D_{p(i)}\}=\{p\in \mathcal{A}'\,:\,k\in D_{p(i)}\}\backslash M$ is uncountable as well. The same argument holds for the set $\{p\in \mathcal{X}\,:\,k\in E_{p(i)}\}$. Hence, the proof is over.
\end{claimproof}
\noindent
\underline{Claim 3}: Let $i<n$ and $\mathcal{X}$, $\mathcal{Y}$ be uncountable disjoint subsets of $\mathcal{A}$ satisfying the conclusions of the Claim 2 when applied to $i$. Then there are $\mathcal{X}'\in [\mathcal{X}]^{\omega_1}$ and $\mathcal{Y}'\in [\mathcal{X}]^{\omega_1}$ so that $$(D_{p(i)}\cap E_{q(i)})\cup(D_{q(i)}\cap E_{q(i)})\not=\emptyset$$
for all $p\in \mathcal{X}'$ and $q\in \mathcal{Y}'$.
\begin{claimproof}[Proof of claim.]By the hypotheses of the proposition, we have that $(\langle D_{p(i)}\rangle_{p\in \mathcal{X}}, \langle E_{q(i)}\rangle_{q\in \mathcal{Y}})$ is a gap. In particular, there is $k\in \omega$, $p'\in \mathcal{X}$ and $q'\in \mathcal{Y}$ for which $k\in D_{p'(i)}\cap E_{q'(i)})\cup(D_{q'(i)}\cap E_{q'(i)})$. Without loss of generality we may assume that $k\in D_{p'(i)}\cap E_{q'(i)}$. We define $$\mathcal{X}'=\{p\in \mathcal{X}\,:,k\in D_{p(i)}\},$$
$$\mathcal{Y}'=\{q\in \mathcal{Y}\,:\,k\in E_{q(i)}\}.$$
It is straightforward that $\mathcal{X}'$ and $\mathcal{Y}'$ are the sets that we are looking for.
\end{claimproof}
\noindent
By applying multiple times the claims  2 and 3, we may build two sequences $\mathcal{X}_{n-1}\subseteq\dots \subseteq \mathcal{X}_0\subseteq \mathcal{A}$ and $\mathcal{Y}_{n-1}\subseteq \dots \subseteq \mathcal{Y}_0\subseteq \mathcal{A}$ of uncountable sets for which $\mathcal{X}\cap \mathcal{Y}=\emptyset$ and such that $(D_{p(i)}\cap E_{q(i)})\cup(D_{q(i)}\cap E_{q(i)})\not=\emptyset$
for all $i<n$, $p\in \mathcal{X}_i$ and $q\in \mathcal{Y}_i$. Using Claim 1, it should be clear that if $p\in \mathcal{X}_{n-1}$ and $q\in \mathcal{Y}_{n-1}$, then $p\cup q$ is a condition of $\mathbb{L}(\mathcal{D},\mathcal{E})$. Thus, we are done.

\end{proof}
\end{proposition}
\begin{corollary}\label{corollarycccbiorthogonal}Let $(\mathcal{D},\mathcal{E})$ be an almost Luzin gap.  Then there is a $ccc$ forcing $\mathbb{Q}$ so that $\mathbb{Q}\Vdash\text{\say{ $(\mathcal{D},\mathcal{E})$ has a biorthogonal subgap.}}$
\end{corollary}
\begin{theorem}[Under $PFA$]\label{stronglyseppfa}There are no strongly donut-separable $(\omega_1,\omega_1)$-gaps
\begin{proof} Let $(\mathcal{L},\mathcal{R})=(L_\alpha,R_\alpha)_{\alpha\in\omega_1}$ be an $(\omega_1,\omega_1)$-gap and consider the forcing $\mathbb{P}=2^{<\omega_1}$. According to the Remark \ref{normalgaprem}, we may assume that $(\mathcal{L},\mathcal{R})$ is normal. Now, let $G\subseteq \mathbb{P}$ be a $\mathbb{P}$-generic filter over $V$.
Then $V[G]$ models $CH$. Furthermore, $[\omega]^{\omega}\cap V=[\omega]^{\omega}\cap V[G]$. In this way, $(\mathcal{L},\mathcal{R})$ is still a gap in $V[G]$. According to the Proposition \ref{highlyadequateprop}, we can find a club $X\subseteq \omega_1$ in $V[G]$ which is highly adequate for $(\mathcal{L},\mathcal{R})$. That is, in $V[G]$, the gap $(L_{X(\alpha+1)}\backslash L_{X(\alpha)},R_{X(\alpha+1)}\backslash R_{X(\alpha)})_{\alpha\in \omega_1}$ is almost Luzin. Therefore, by virtue of the Corollary \ref{corollarycccbiorthogonal} there is a $ccc$ forcing $\mathbb{Q}$ which forces this gap to have a biorthogonal subgap. Let $H$ be a $\mathbb{Q}$-generic filter over $V[G]$. Then, in $V[G][H]$, there is $Y\in [\omega_1]^{\omega_1}$ so that the pregap $$(\langle L_{X(\alpha+1)}\backslash L_{X(\alpha)}\rangle, \langle R_{X(\alpha+1)}\backslash R_{X(\alpha)}\rangle)_{\alpha\in Y}$$
is biorthogonal.\\
Returning to $V$, let $\dot{\mathbb{Q}}$ be a $\mathbb{P}$ name for $\mathbb{Q}$ in $V[G]$. Now, let $\dot{X}$ and $\dot{Y}$ be $\mathbb{P}*\dot{Q}$-names for $X$ and $Y$ respectively. Finally, let $p\in \mathbb{P}*\dot{\mathbb{Q}}$ which forces $\dot{X}$ and $\dot{Y}$ to have the properties discussed in the previous paragraph. Given $\xi\in \omega_1$, let $$D_\xi=\{q\leq p\:\,:\exists \xi<\alpha,\beta\in \omega_1\,(q\Vdash\text{\say{ $\alpha\in\dot{X}$, and $\beta\in \dot{Y}$)}}\,\},$$
$$E_\xi=\{q\leq p\,:\,p\Vdash\text{\say{$\xi\not\in\dot{X}$}}\textit{ or there is }\alpha\in \omega_1\textit{ such that }p\Vdash\text{\say{$\dot{X}(\alpha)=\xi$}}\}.$$
It is straightforward that $D_\xi$ and $E_\xi$ are dense below $p$ for any $\xi\in\omega_1$. Since $\mathbb{P}$ is $\sigma$-closed and $\dot{\mathbb{Q}}$ is forced to be $ccc$, then $\mathbb{P}*\dot{\mathbb{Q}}$ is proper. In this way, there is a filter $F$ which intersects  each $D_\xi$ and $E_\xi$. Let $$X'=\{\alpha\,:\,\exists q\in F\,(q\Vdash\text{\say{$\alpha\in \dot{X}$}}\},$$
$$Y'=\{\alpha\,:\,\exists q\in F\,(q\Vdash\text{\say{$\alpha\in \dot{Y}$}}\}.$$
Note that the statement \begin{center}
    \say{$(L_{X'(\alpha+1)}\backslash L_{X'(\alpha)}\cap R_{X'(\beta+1)}\backslash R_{X'(\beta)})\cup  (R_{X'(\alpha+1)}\backslash R_{X'(\alpha)}\cap L_{X'(\beta+1)}\backslash L_{X'(\beta)})=\emptyset$}
\end{center}
is absolute. From this, it easily follows that that the pregap $$(L_{X'(\alpha+1)}\backslash L_{X'(\alpha)}, R_{X'(\alpha+1)}\backslash R_{X'(\alpha)})_{\alpha\in Y}$$ is biorthogonal. In particular, it is a gap by Lemma \ref{biorthogonalgaplemma}. Form this, it follows directly that $(\mathcal{L},\mathcal{R})$ is not strongly donut-separable.

\end{proof}
\end{theorem}

So far we have seen that both $CH$ and $PFA$ imply that the non-existence of strongly donut-separable gaps. On the other hand, the existence of a $2$-capturing construction scheme $\mathcal{F}$ for which $\mathfrak{m}_\mathcal{F}>\omega_1$ implies the existence of a strongly donut-separable gap. We will end this subsection by proving that $MA$ is not sufficient to decide this problem. This result suggests that the cardinal invariant $\mathfrak{m}_\mathcal{F}$ should not be treated just as a sort of weakening of $MA$, but rather as an interesting principle on its own. 
\begin{theorem}\label{independentmastronglysep}$MA$ is independent from the statement \say{there are strongly donut-inseparable gaps}.
\begin{proof}Since $PFA$ implies that there are no strongly donut-inseparable gaps, then $MA$ is consistent with that same statement.
In order to show the consistency of the negation, we start with a model $V$ in which there is a strongly donut-inseparable gap. We now consider a $ccc$ forcing $\mathbb{P}$ which forces  $MA$ and $G$ a $\mathbb{P}$-generic filter over $V$. Then $V[G]$ is a model $MA$. Furthermore, by means of the Corollary \ref{corocccpreservingdonut}, $V[G]$ also models the existence of a strongly donut-inseparable gap. Thus, the proof is over.
    
\end{proof}
    
\end{theorem}

\subsection{The gap cohomology group}
In the subsection we naturally extend the work done by Daniel Talayco in \cite{cohomologytalayco}. In particular, we study the size of groups defined in a similar context as Talayco and prove that, under certain assumptions, these groups are as big as possible.
\begin{definition}[$*$-lower semi-lattice]Let $X$ be an infinite set.
We say that a family $\mathcal{T}\subseteq [X]^{\omega}$ is a \textit{$*$-lower semi-lattice} if  for any $A,B\in \mathcal{T}$ there is $C\in \mathcal{T}$ so that $A\cap B=^*C.$ In other words, $\{[A]\,:\,A\in \mathcal{T}\}$ is a lower semi-lattice in $\mathscr{P}(X)/\text{FIN}$.
\end{definition}
\begin{remark}Towers are particular cases of $*$-lower semi-lattices.
\end{remark}
\begin{definition}[Coherent Subsystems]\label{cohsystemsdef} Let $\mathcal{T}$ be a $*$-lower semi-lattice over a countable set $X$. A function $g:\mathcal{T}\longrightarrow \mathscr{P}(X)$ is said to be a \textit{coherent subsystem of $\mathcal{T}$} if for any $A,B\in \mathcal{T}$, the following happens:
\begin{enumerate}[label=$(\arabic*)$]
\item $g(A)\subseteq A$,
\item if $B\subseteq^* A$ then $g(A)\cap B=^*g(B)$.
\end{enumerate}
 We say that $g$ is \textit{trivial} if there is $C\in \mathscr{P}(X)$ so that $C\cap A=^*g(A)$ for any $A\in \mathcal{T}$. In this case, we say that $C$ \textit{trivializes} $g$. The set of all coherent subsystems of $\mathcal{T}$ is denoted as $C(\mathcal{T})$ and the set of all trivial coherent subsystem is denoted as $Tr(\mathcal{T}).$
\end{definition}
\begin{definition}Let $\mathcal{T}$ be a $*$-lower semilattice over a countable set $X$ and $g$ be a coherent subsystem of $\mathcal{T}$. We denote the sequence $(g(A),A\backslash g(A))_{A\in\mathcal{T}}$ as $\mathcal{C}(\mathcal{T},g).$
    
\end{definition}
The following lemma is easy. We prove it just to emphasise the reason of defining coherent subsystems just for $*$-lower semi-lattices.
\begin{lemma}\label{lemmacoherentsubsystempregap}Let $\mathcal{T}$ be a $*$-lower semi-lattice and let $g$ be a coherent subsystem of $\mathcal{T}$. Then $\mathcal{C}(\mathcal{T},g)$ is a pregap. Furthermore, $g$ is non-trivial if and only if $\mathcal{C}(\mathcal{T},g)$ forms a gap.

\begin{proof}Let $A,B\in \mathcal{T}$. We need to prove that $g(A)\cap (B\backslash g(B))=^*\emptyset.$ For this let $C\in \mathcal{T}$ be so that $C=^*A\cap B$. Then \begin{align*}g(A)\cap B=^* (g(A)\cap A)\cap B=^*g(A)\cap C
=^*g(C)=^*g(B)\cap C\subseteq ^*g(B).
\end{align*}
This finishes the argument. Now we will prove the proposed equivalence.
\begin{claimproof}[Proof of $\Rightarrow$] We will prove this by contra-positive. Suppose that there is a set $C$ separating the pregap $\mathcal{C}(\mathcal{T},g)$. For any $A\in \mathcal{T}$ we have that $g(A)\subseteq^* C$ and $C\cap (A\backslash g(A))=^*\emptyset.$ In this way, $C\cap A=^*g(A)$. We conclude that $C$ trivializes $g$. 
\end{claimproof}
\begin{claimproof}[Proof of $\Leftarrow$] Again, by contra-positive. Suppose that there is $C$ which trivializes $g$. Then $C\cap A=^*g(A)$ for any $A\in \mathcal{T}$. In particular $g(A)\subseteq^* C$ and $C\cap (A\backslash g(A))=\emptyset$. We conclude that $C$ separates $\mathcal{C}(\mathcal{T},g)$.
\
\end{claimproof}
\end{proof}
\end{lemma}

\begin{definition}Let $\mathcal{T}$ be a $*$-lower semi-lattice over a set $X$. For $f,g\in C(\mathcal{T})$ we can define $f\cdot g\in C(\mathcal{T})$ given as: $$(f\cdot g)(A)= f(A)\Delta g(A).$$
\end{definition}
The following lemma is easy.
\begin{lemma}
$C(\mathcal{T})$ is an abelian (Boolean)\footnote{A group $(G,\cdot, e)$ is Boolean if $g\cdot g=e$ for any $g\in G$. Equivalently, $G$ is a vector space over $\mathbb{Z}_2$. In particular, the isomorphism type of $G$ is completely determined by its cardinality.} group with respect to this operation and $Tr(\mathcal{T})$ is a subgroup of it. 
\end{lemma}
Since $C(\mathcal{T})$ is abelian, then $Tr(\mathcal{T})$ is a normal subgroup. Therefore, we can consider the quotient induced by it.
\begin{definition}[Gap cohomology group] We define the \textit{gap cohomology group of $\mathcal{T}$} as the quotient of these two groups. That is,  $$G(\mathcal{T})=C(\mathcal{T})/Tr(\mathcal{T}).$$
We  say that $f,g\in C(\mathcal{T})$ are \textit{cohomologous} if the class of $f$ is equal to the class of $g$ inside $G(\mathcal{T})$, or equivalently, if $f\cdot g\in Tr(\mathcal{T})$. 
\end{definition}
\begin{rem}If $\mathcal{T}$ has size $\omega_1$, then $|G(\mathcal{T})|\leq 2^{\omega_1}.$
    
\end{rem}
In \cite{cohomologytalayco}, Daniel E. Talayco defined the gap cohomology group for the particular case of $\omega_1$-towers. In there he showed that the $2^\omega\leq |G(\mathcal{T})|$ for any $\omega_1$-tower $\mathcal{T}$. In order to do that, he proved the following Theorem.
\begin{theorem}[The $\aleph_0$ gap theorem] let $\mathcal{T}=\langle T_\alpha\rangle_{\alpha\in \omega_1}$ be an $\omega_1$-tower. There is  a sequence of functions $\langle f_\alpha\rangle_{\alpha\in\omega_1}$ with the following properties:
\begin{itemize}
\item $\forall \alpha\in \omega_1\,(\,f_\alpha:T_\alpha\longrightarrow \omega\,)$,
\item $\forall \alpha<\beta\in \omega_1\, (\,f_\alpha|_{T_\alpha\cap T_\beta}=^*f_\beta|_{T_\alpha\cap T_\beta}\,),$
\item $\forall m<n\in\omega\,(\, (f_\alpha^{-1}[\{m\}],f_{\alpha}^{-1}[\{n\}]\,)_{\alpha\in\omega_1}\textit{ is a Hausdorff gap }).$ 

\end{itemize}
\end{theorem}
He also proved that the $\Diamond$-principle implies that the size of this group is always $2^{\omega_1}$. Later, Stevo Todor\v{c}evi\'{c} noted that the $\aleph_0$ gap theorem already gives you the previous conclusion without assuming any extra axioms (see Theorem 32 in \cite{cohomologytalaycotrees}). Namely, if $2^{\omega}=2^{\omega_1}$ the result is clear due to Talayco's calculations. On the other hand, if $2^{\omega}<2^{\omega_1}$ then the result follows from an easy counting argument involving the quotient structure of $G(\mathcal{T}).$
In \cite{Agapcohomologygroup}, Charles Morgan constructed an $\omega_1$-tower whose gap cohomology group can be explicitly calculated. This construction was carried through the use of morasses. Finally, in \cite{coherentfamilyoffunctions}, Ilijas Farah improved the $\aleph_0$ gap theorem. This was achieved by considering the following concept.
\begin{definition}[Coherent families of functions]\label{coherentfamiliesdef} A \textit{coherent family of functions} supported by an $\omega_1$-tower $\langle T_\alpha\rangle_{\alpha\in\omega_1\backslash \omega}$ is a family of functions $\langle f_\alpha \rangle_{\alpha\in\omega_1}$ such that:
\begin{enumerate}[label=$(\arabic*)$]
    \item $\forall \alpha\in \omega_1\backslash \omega\,(\,f_\alpha:T_\alpha\longrightarrow \alpha\,),$
\item $\forall \alpha,\beta\in \omega_1\backslash \omega\,(\,f_\alpha|_{T_\alpha\cap T_\beta}=^*f_\beta|_{T_\alpha\cap T_\beta}\,).$
\end{enumerate}
 Given such family, we define $L^\xi_\alpha$ as $f^{-1}_\alpha[\{\xi\}]$ for all $\xi\in \omega_1$ and $\alpha> \xi$. Additionally, we let $\mathcal{L}^{\xi}=\langle L^\xi_\alpha\rangle_{\alpha>\xi}$
 \end{definition}
 \begin{lemma} If  $\langle f_\alpha\rangle_{\omega_1\backslash \omega}$ is as in the previous definition and $\xi<\mu\in \omega_1$, the pair $(L^\xi_\alpha,L^\mu_\alpha)_{\alpha\in \omega_1\backslash \mu}$ is a pregap. We will denote it as $(\mathcal{L}^\xi,\mathcal{L}^\mu)$ although, in principle, $\mathcal{L}^\xi=\langle L^\xi_\alpha\rangle_{\alpha>\xi}$ instead of  $\langle L^\xi_\alpha\rangle_{\alpha>\mu}.$
 
 \end{lemma}
\begin{definition}Let $\mathfrak{F}=\langle f_\alpha\rangle_{\alpha\in \omega_1\backslash \omega}$ be a coherent family of functions supported by an $\omega_1$-tower. We say that $\mathfrak{F}$ is Luzin (respectively Hausdorff) if $(\mathcal{L}^\xi,\mathcal{L}^\mu)$ is Luzin (respectively Hausdorff) for each $\xi<\mu\in \omega_1$. 

\end{definition}

Ilijas Farah constructed a Haudorff coherent family of functions supported by an $\omega_1$-tower by forcing it and then appealing to Keisler's completeness Theorem for $L^\omega(Q)$. Unlike Talayco's result, the proof idea behind the construction of Farah does not work for any $\omega_1$-tower. This is due to the nature of Keisler's Theorem. 

Here we give a direct construction of a Luzin coherent family of functions with the use of a $2$-construction scheme. No previous direct construction was known.

\begin{theorem}There is a Luzin coherent family of functions supported by an $\omega_1$-tower. 
\begin{proof}Let $\mathcal{F}$ be a $2$-construction scheme. First we define $\omega_1$-tower over the countable set $$N=\bigcup\limits_{k\in\omega} N_k$$
where each $N_k$ is equal to $\{k\}\times k\times r_k\times r_k.$ Given $\alpha\in \omega_1\backslash \omega$ we define $T_\alpha$ as $$\bigcup\limits_{\Xi_\alpha(k)\geq 0} N_k=\{(k,s,i,j)\in \omega^4\,:\,\Xi_\alpha(k)\geq 0\textit{ and }(s,i,j)\in k\times r_k\times r_k\,\}.$$
Note that for each $\alpha<\beta\in\omega_1$ and every $k>\rho(\alpha,\beta)$ we have that $T_\alpha\cap N_k\subseteq T_\beta\cap N_k$ due to the part (c) of Lemma \ref{lemmaxi}. Consequently, $T_\alpha\subseteq ^*T_\beta$. In this way $\mathcal{T}=\langle T_\alpha\rangle_{\alpha\in\omega_1\backslash \omega}$ is a tower.

Now we define a Luzin coherent family of functions supported by $\mathcal{T}$. Let $\alpha\geq \omega$. Note that if $x\in T_\alpha$ then $x=(k,s,i,j)$ where $k>0$, $\lVert \alpha\rVert_k>r_k$ and $i,j<r_k$ (In particular, $(\alpha)_k(i)$ and $(\alpha)_k(j)$ are defined). In this way, we can define $f_\alpha: T_\alpha\longrightarrow \alpha$ as follows:

$$f_\alpha(k,s,i,j)=\begin{cases}(\alpha)_{k}(i)&\textit{ if }\Xi_\alpha(k)=0\\
(\alpha)_{k}(j)&\textit{ if }\Xi_\alpha(k)=1
\end{cases}$$
By definition, the family $\langle f_\alpha\rangle_{\alpha\in\omega_1\backslash \omega}$ satisfies the point (1) of Definition \ref{coherentfamiliesdef}. The two following claims will finish the proof.\\\\
\underline{Claim 1}: $\langle f_\alpha\rangle_{\alpha\in\omega_1\backslash \omega}$ satisfies the point (2)  of Definition \ref{coherentfamiliesdef}.
\begin{claimproof}[Proof of claim] Let $\alpha<\beta\in\omega_1\backslash \omega$ and $(k,s,i,j)\in T_\alpha\cap T_\beta$ be such that $k>\rho(\alpha,\beta)$. By definition of $T_\alpha$ and $T_\beta$ it follows that both $\Xi_\alpha(k)$ and $\Xi_\beta(k)$ are non-negative numbers. In virtue of the part (c) of Lemma \ref{lemmaxi}, it must happen that $\Xi_\alpha(k)=\Xi_\beta(k)$. Furthermore,  $(\alpha)_k\sqsubseteq (\beta)_k$ so $(\alpha)_k(i)=(\beta)_k(i)$ and $(\alpha)_k(j)=(\beta)_k(j)$. Hence, $f_\alpha(k,s,i,j,s)=f_\beta(k,s,i,j)$. Since all but finitely many elements of $T_\alpha\cap T_\beta$ have their first coordinate bigger than $\rho(\alpha,\beta)$, we have shown that $f_\alpha|_{T_\alpha\cap T_\beta}=^*f_\beta|_{T_\alpha\cap T_\beta}$.
\end{claimproof}
\noindent
\underline{Claim 2:} $\langle f_\alpha\rangle_{\alpha\in \omega_1\backslash \omega}$ is Luzin.

\begin{claimproof}[Proof of claim] Let $\xi<\mu\in \omega_1$,  $\beta>\mu$ and $n\in\omega$. We need to prove that the set $$\{\alpha\in \beta\backslash \mu\, :\, |L^\xi_\alpha\cap L^\mu_\beta|\leq n\}$$ is finite. We claim that this set is contained in $(\beta)_l$ where $l=\max(n,\rho^{\{\xi,\mu,\beta\}})$. Indeed, take an $\alpha\in \beta\backslash \mu$ such that $\alpha\notin (\beta)_l$ and let $k=\rho(\alpha,\beta)$ (so obviously $l<k$). According to the part (b) of Lemma \ref{lemmaxi}, $\Xi_\alpha(k)=0$ and $\Xi_\beta(k)=1.$ Furthermore, since $\alpha\notin (\beta)_l$ then $k\geq \rho^{\{\xi,\mu,\beta\}}$. This means that both $\xi$ and $\mu$ belong to $(\alpha)_k$. Hence $$\Xi_\xi(k)\leq \Xi_\mu(k)\leq \Xi_\alpha(k)=0<\Xi_\beta(k).$$By virtue of the part (c) of Lemma \ref{lemmaxi},  $\Xi_\mu(k)=\Xi_\xi(k)=-1$. In other words, $\lVert \xi\rVert_k$ and $\lVert \mu\rVert_k$ are numbers strictly smaller than $r_k$. Moreover, $(\alpha)_k(\lVert\xi\rVert_k)=(\beta)_k(\lVert\xi\rVert_k)=\xi$ and $(\alpha)_k(\lVert\mu\rVert_k)=(\beta)_k(\lVert\mu\rVert_k)=\xi$.  From this it follows that $$\{k\}\times k\times\{\lVert \xi\lVert_k\}\times \{\lVert\mu\lVert_k\}\subseteq f^{-1}_\alpha[\{\xi\}]\cap f^{-1}_\beta[\{\mu\}]=L^\xi_\alpha\cap L^\mu_\beta.$$
As this set has cardinality $k$ (which is bigger than $n$) we are done.
\end{claimproof}

\end{proof}
\end{theorem}
The gap cohomology group of $\omega_1$-towers with a Luzin coherent family of functions can be explicitly calculated using Talayco's ideas. We present here such calculations for the sake of completeness.
\begin{proposition}Let $\mathcal{T}=\langle T_\alpha\rangle_{\alpha\in \omega_1}$ be an $\omega_1$-tower and suppose there is a Luzin coherent family of functions $\langle f_\alpha\rangle_{\alpha\in \omega_1}$ supported by $\mathcal{T}$. Then $|G(\mathcal{T})|=2^{\omega_1}$.
\begin{proof} For each $S\in \mathscr{P}(\omega_1)$, we define a coherent subsystem $g_S:\mathcal{T}\longrightarrow \mathscr{P}(\omega)$ as follows:
$$g_S(T_\alpha)=f^{-1}_\alpha[S].$$
Note that if $S,S'\in \mathscr{P}(\omega_1)$ then $g_{S}\cdot g_{S'}=g_{S\Delta S'}$. Let $\mathcal{S}\subseteq \mathscr{P}(\omega_1)$ be a family of size $2^{\omega_1}$  so that $S\Delta S'$ and $\omega_1\backslash (S\Delta S')$ are infinite for any two distinct $S,S'\in \mathcal{S}$. Given such $S$ and $S'$ we claim that $g_S$ and $g_{S'}$ are not cohomologous. Indeed, let $\xi\in S\Delta S'$ and $\mu\in \omega_1\backslash(S\Delta S')$ with $\xi<\mu$. Then $L^\xi_\alpha\subseteq g_{S\Delta S'}(T_\alpha)$ and $L^\mu_\alpha\subseteq g_{\omega_1\backslash(S\Delta S')}(T_\alpha)=T_\alpha\backslash g_{S\Delta S'}(T_\alpha)$ for any $\alpha>\mu$. By hypothesis we have that $(L^\xi_\alpha,L^\mu_\alpha)_{\alpha\in \omega_1}$ is a Luzin-gap. This implies that $G(\mathcal{T},g_{S\Delta S'})$ is a gap as well. That is, $g_{S\Delta S'}=g_S\cdot g_{S'}$ is a non-trivial coherent subsystem of $\mathcal{T}$. We have proved that for any two distinct $S,S'\in \mathcal{S}$, $g_S$ and $g_{S'}$ are not cohomologous. Since $|S|=2^{\omega_1}$, then $|G(\mathcal{T})|=2^{\omega_1}$.
\end{proof}
\end{proposition}
As a second application of the concept of Luzin representations, we use Theorem \ref{luzinreptheorem} to generalize the results of Talayco, Todor\v{c}evi\'c, Morgan and Farah regarding the existence of \say{big gap cohomology groups}. It is worth noting that, as for now, it is unclear if there is an analogous to the $\aleph_0$ gap theorem for aribrary $*$-lower semi-lattices. Therefore, we can not argue in the same way as Todor\v{c}evi\'c  in order to conclude that the gap cohomology group has always cardinality $2^{\omega_1}$.
\begin{theorem}\label{hausdorffcoherenttheorem}Let $(X,<,\wedge)$ be an $\omega_1$-like lower semi-lattice. Then there is a $*$-lower semi-lattice $\mathcal{T}$ isomorphic to $X$ with $|G(\mathcal{T})|=2^{\omega_1}.$
\begin{proof}Consider $X\times 2$ with the order given by $(x,i)<(y,j)$ if and only if $x<y$ and $i=j$. It is easy to see that $X\times 2$ is also $\omega_1$-like. In virtue of Theorem \ref{luzinreptheorem} we can take $\mathcal{T'}=\langle T'_{(x,i)}\rangle_{(x,i)\in X\times 2}$ and $\mathcal{A}=\langle A_{(x,i)}\rangle_{(x,i)\in X\times 2}$ so that $(\mathcal{T}',\mathcal{A})$ is a Luzin representation of $X\times 2$. For any $x\in X$ let $T_x=T'_{(x,0)}\cup T'_{(x,1)}. $ It is straightforward that the following properties hold for any $x,y \in X$:
\begin{enumerate}[label=$(\alph*)$]
\item $A_{(x,0)}\cup A_{(x,1)}\subseteq T_x$.

\item If $y\not\leq x$ then $A_{(y,0)}\cup A_{(y,1)}\subseteq^*T_y\backslash T_x$. In other words, $(A_{(y,0)}\cup A_{(y,1)})\cap T_x=^*\emptyset$.
\item $T_x\cap T_y=^*T_{x\wedge y}.$

\end{enumerate}
In particular, properties (b) and (c) imply that $\mathcal{T}=\langle T_x\rangle_{x\in X}$ is a $*$-lower semi-lattice isomorphic to $X$. We claim that $|G(\mathcal{T})|=2^{\omega_1}$. Trivially $|G(\mathcal{T})|\leq 2^{\omega_1}$, so we will only prove the other inequality. For this, let $S\in \mathscr{P}(X)$. . We will build recursively (using that $X$ is well-founded) a coherent subsystem $g_S:\mathcal{T}\longrightarrow \mathscr{P}(\omega)$ such that $A_{(x,0)}\subseteq g_S(T_x)$ for any $x\in X$ and:\\
\begin{center}\begin{minipage}{6cm} \begin{center} \textbf{(A)}\end{center} If $x\in S$, then $A_{(x,1)}\cap g_S(T_x)=^*\emptyset.$
\end{minipage}\hspace{2.5cm} \begin{minipage}{6cm}\begin{center} \textbf{(B)}\end{center} If $x\not\in S$, then $A_{(x,1)}\subseteq g_S(T_x)$.
\end{minipage}

\end{center}
Suppose that $y\in X$ and we have defined $g_S(T_x)$ for any $x<y$. By means of Lemma \ref{lemmacoherentsubsystempregap} we know that $\mathcal{G}(\mathcal{T}|_y,g_S|_{\mathcal{T}|_{y}})$ is pregap where $\mathcal{T}|_y=\langle T_x\rangle _{x<y}$. As there are no countable pregaps and $|(-\infty,y)|\leq \omega$, we conclude that there is $C\in \mathscr{P}(\omega)$ separating it.  Since $g_S(T_x)\subseteq T_x\subseteq^*T_y$ for any $x<y$, we may assume without loss of generality that $C\subseteq T_y$. Note that $C\cap T_x=^*g_S(T_x)$ for any $x<y$. We now define $g_S(T_y)$ by cases. If $y\in S$ define $g_S(T_y)$ as $(C\cup A_{(y,0)})\backslash A_{(y,1)}$. Otherwise define $g_S(T_y)$ as $C\cup A_{(y,0)}\cup A_{(y,1)}$. Since both $A_{(y,0)}$ and $A_{(y,1)}$ are almost disjoint with $T_x$ for all $x<y$ due to the point (b), it follows that $g_S(T_y)\cap T_x=^* C\cap T_x=^*g_S(T_x)$ for any such $x$. This finishes the recursion.

Now consider $\mathcal{S}\subseteq [X]^{\omega_1}$  a family of cardinality $2^{\omega_1}$ so that $S\Delta S'$ is uncountable for any two distinct $S,S'\in \mathcal{S}$. The key fact needed to finish the proof is that if $x\in S\Delta S'$ then $A_{(x,1)}\subseteq g_S\cdot g_{S'}(T_x)$ and $A_{(x,0)}\subseteq T_x\backslash (g_S\cdot g_{S'}(T_x))$.  In this way, any separation of the pregap  $\mathcal{G}(\mathcal{T},g_S\cdot g_{S'})$ would also separate the pregap $(A_{(x,1)},A_{(x,0)})_{x\in S\Delta S'}$ which is impossible since $\mathcal{A}$ is a Luzin family and $S\Delta S'$ is uncountable. This shows that $g_S$ and $g_{S'}$ are not cohomologous. Since $|S|=2^{\omega_1}$ then $|G(\mathcal{T})|=2^{\omega_1}$.
\end{proof}
\end{theorem}
\subsection{Destructibility of gaps}
In this subsection we will study $(\omega_1,\omega_1)$-gaps in the context of forcing. The reader can find fairly complete treatments of this subject in \cite{GapsandTowers}, \cite{ScheepersGaps} and \cite{StevoIlias}.

The following fundamental lemma is attributed to Kenneth Kunen. We will prove it for the sake of completeness. For others proofs see \cite{gapsandlimits}, \cite{ScheepersGaps} and \cite{StevoIlias}.
\begin{lemma}\label{gapforcinglemma}Let $(L_\alpha,R_\alpha)_{\alpha\in\omega_1}$ be an $(\omega_1,\omega_1)$-pregap with $L_\alpha\cap R_\alpha=\emptyset$ for any $\alpha$. Then $(L_\alpha,R_\alpha)_{\alpha\in \omega_1}$ is a gap if and only if for any $S\in [\omega_1]^{\omega_1}$ there are $\alpha<\beta\in S$ such that $(L_\alpha\cap R_\beta)\cup (L_\beta\cap R_\alpha)\not=\emptyset$.
\begin{proof}

\begin{claimproof}[Proof of $\Rightarrow$]Assume that $(L_\alpha,R_\alpha)_{\alpha\in\omega_1}$ is a gap and let $S\in [\omega_1]^{\omega_1}$. Suppose towards a contradiction that $(L_\alpha\cap R_\beta)\cup (L_\beta\cap R_\alpha)=\emptyset$ for all $\alpha,\beta\in S$. Let $C=\bigcup\limits_{\beta\in S} L_\beta$. Given $\alpha\in \omega_1$ there is $\beta\in S$ such that $\alpha<\beta$. Note that $L_\alpha\subseteq^*L_\beta\subseteq C$ and $R_\alpha\cap C\subseteq^* R_\beta\cap C=\emptyset$. In this way $C$ separates $(L_\alpha,R_\alpha)_{\alpha\in \omega_1}$ which is a contradiction. Thus, the proof of this implication is over.
\end{claimproof}

 \begin{claimproof}[Proof of $\Leftarrow$] Again, by contradiction. Suppose that $(L_\alpha,R_\alpha)_{\alpha\in \omega_1}$ is not a gap and let $C$ be a set which separates it. By the pigeonhole principle we can find $S\in [\omega_1]^{\omega_1}$ and $a,b\in [\omega]^{<\omega}$ such that $L_\alpha\backslash C=a$ and $R_\alpha\cap C=b$ for any $\alpha\in S$. Note that if $\alpha\in S$ then $L_\alpha\cap b=\emptyset$ because $b\subseteq R_\alpha$ and $L_\alpha\cap R_\alpha=\emptyset$. Therefore, $L_\alpha\cap R_\beta\subseteq ((C\cup a)\backslash b)\cap R_\beta=\emptyset $ whenever $\alpha,\beta\in S$. With this contradiction we finish the proof.
 \end{claimproof}
\end{proof}
\end{lemma}
\begin{rem}For the rest of this section we will assume that if $(L_\alpha,R_\alpha)_{\alpha\in\omega_1}$ is an $(\omega_1,\omega_1)$-pregap then it is normal. That is, $L_\alpha\cap R_\alpha=\emptyset$ any $\alpha.$
\end{rem}
\begin{definition}[Destructible gaps]Let $(\mathcal{L},\mathcal{R})$ be an $(\omega_1,\omega_1)$-pregap. We say that $(\mathcal{L},\mathcal{R})$ is \textit{destructible} if there is a forcing notion $\mathbb{P}$ which preserves $\omega_1$ in such way that $(\mathcal{L},\mathcal{R})$ is not a gap in some generic extension through $\mathbb{P}$. If this does not happen, the gap is said to be \textit{indestructible}.
\end{definition}
\begin{rem}The Hausdorff condition is absolute. Hence any Hausdorff gap is indestructible. Note that under $PID$ any $(\omega_1,\omega_1)$-gap is undestructible. This is due to the Remark \ref{remarkhausdorffgaps1}.
\end{rem}
In Definition \ref{separabilityforcingdonuts} we defined the forcing $\mathbb{P}(\mathcal{D},\mathcal{E})$ where $(\mathcal{D},\mathcal{E})$ is a pregap (not necesarily of type $(\omega_1,\omega_1)$) with $|\mathcal{D}|=|\mathcal{E}|=\omega_1$. According to Proposition \ref{separabilitygapsprop}, this forcing is $ccc$ if and only if the pregap $(\mathcal{D},\mathcal{E})$ can be separated in some $\omega_1$-preserving extension of the universe. Forcings with such properties have already been studied for the case of $(\omega_1,\omega_1)$-gaps. In the following definition we present some well-known reincarnations of them (see \cite{independenceforanalysts}, \cite{ScheepersGaps}, \cite{PartitionProblems} or \cite{yoriokadestructiblegaps}).
\begin{definition}Let $(\mathcal{L},\mathcal{R})$ be an $(\omega_1,\omega_1)$-pregap indexed as $(L_\alpha, R_\alpha)_{\alpha\in X}$ where $X$ is an uncountable set of ordinals. We define the following forcing notions:
\begin{itemize}
 \item $\chi_0(\mathcal{L},\mathcal{R})=\{p\in [X]^{<\omega}\,:\, \big(\bigcup\limits_{\alpha\in p}L_\alpha \big)\cap \big(\bigcup\limits_{\alpha\in p}R_\alpha\big)=\emptyset\}.$
    \item $\chi_1(\mathcal{L},\mathcal{R})=\{p\in [X]^{<\omega}\,:\,\forall \alpha\not=\beta\in p\,\big( (L_\alpha\cap R_\beta)\cup(L_\beta\cap R_\alpha)\not=\emptyset\big)\}.$
\end{itemize}
both ordered by reverse inclusion.
\end{definition}
\begin{lemma}\label{equivalencechi1forcing}Let $(\mathcal{L},\mathcal{R})$ be an $(\omega_1,\omega_1)$-pregap indexed as $(L_\alpha,R_\alpha)_{\alpha\in X}$ where $X$ is an uncountable set of ordinals. Then $(\mathcal{L},\mathcal{R})$ is a gap if and only if $\mathbb{P}=\chi_1(\mathcal{L},\mathcal{R})$ is $ccc$.
\begin{proof}First suppose that $(\mathcal{L},\mathcal{R})$ is a gap and consider $\mathcal{A}\in [\mathbb{P}]^{\omega_1}$. By refining $\mathcal{A}$ we can assume that the following conditions hold:
\begin{itemize}
    \item There is $n\in \omega$ such that $|p|=n$ for any $p\in \mathcal{A}$.
    \item There are $\sigma,\tau:n\times n\longrightarrow [\omega]^{<\omega}$ so that for any $i\leq j\in n$, $L_{p(i)}\backslash L_{p(j)}=\sigma(i,j)=\sigma(j,i)$ and $R_{p(i)}\backslash R_{p(j)}=\tau(i,j)=\tau(j,i)$.
    \item There is $\mu:n\times n\longrightarrow [\omega]^{<\omega}$ so that for any $(i,j)\in n\times n$, $L_{p(i)}\cap R_{q(j)}=\mu(i,j)$.
    \item $\mathcal{A}$ forms a root-tail-tail $\Delta$-system with a root $R$ of cardinality $r<n$. 
    \end{itemize}
    By means of Lemma \ref{gapforcinglemma} we know that there are distinct $p, q\in \mathcal{A}$ so that $$(L_{p(r)}\cap R_{q(r)})\cup (L_{q(r)}\cap R_{p(r)})\not=\emptyset.$$
\noindent
\underline{Claim}: $p\cup q\in  \mathbb{P}$.

\begin{claimproof}[Proof of claim] It is enough to take $i,j<n$ and show that  $(L_{p(i)}\cap R_{q(j)})\cup (L_{q(j)}\cap R_{p(i)})\not=\emptyset.$ If $i<r$ then $p(i)=q(i)$ so  the previous set 
 is equal to $(L_{q(i)}\cap R_{q(j)})\cup (L_{q(j)}\cap R_{q(i)})$ which by hypothesis is non-empty because $q\in \mathbb{P}$. An analogous argument applies when $j<r$ so let us assume that $r\leq i,j$. But if $i\not=j$ then \begin{align*}\emptyset &\not= (L_{p(i)}\cap R_{p(j)})\cup (L_{p(j)}\cap R_{p(i)})=\mu(i,j)\cup \mu(j,i)
 &\subseteq (L_{p(i)}\cap R_{q(j)})\cup (L_{q(j)}\cap R_{p(i)}).
 \end{align*}
 So we might as well assume that $i=j$.
In this case we have that $$(L_{p(r)}\cap R_{q(r)})\cap \tau(r,i)\subseteq(L_{p(r)}\cap R_{q(r)})\cap R_{p(r)}\subseteq R_{p(r)}\cap L_{p(r)}=\emptyset,$$
 $$(L_{p(r)}\cap R_{q(r)})\cap \sigma(r,i)\subseteq R_{q(r)}\cap L_{q(r)}=\emptyset.$$
 Therefore $L_{p(r)}\cap R_{q(r)}\subseteq R_{q(r)}\backslash \tau(r,i)\subseteq R_{q(i)}$ and $L_{p(r)}\cap  R_{q(r)}\subseteq L_{p(r)}\backslash \sigma(r,i)\subseteq L_{p(i)}$. That is, $L_{p(r)}\cap R_{q(r)}\subseteq L_{p(i)}\cap R_{q(j)}$. In the same way we can show that $L_{p(r)}\cap R_{q(r)}\subseteq L_{q(i)}\cap R_{p(i)}$ so $$\emptyset\not=(L_{p(r)}\cap R_{q(r)})\cup (L_{q(r)}\cap R_{p(r)})\subseteq  (L_{p(i)}\cap R_{q(i)})\cup (L_{q(i)}\cap R_{p(i)}).$$
 \end{claimproof}
This finishes the proof of this implication. The other one is direct from Lemma \ref{gapforcinglemma}. 
 \end{proof}
\end{lemma}

\begin{lemma}\label{equivalencechi0forcing} Let $(L_\alpha,R_\alpha)_{\alpha\in\omega_1}$ be an $(\omega_1,\omega_1)$-pregap indexed as $(L_\alpha,R_\alpha)_{\alpha\in X}$ where $X$ is an uncountable set of ordinal. Then $\mathbb{P}=\chi_0(\mathcal{L},\mathcal{R})$ is $ccc$ if and only if for any $S\in [X]^{\omega_1}$ there are $\alpha<\beta$ such that $(L_\alpha\cap R_\beta)\cup (L_\beta\cap R_\alpha)=\emptyset$.
\begin{proof}
\begin{claimproof}[Proof of $\Leftarrow$] Take $S\in [X]^{\omega_1}$. Then $\mathcal{A}=\{\{\alpha\}\,:\,\alpha\in S\}$ is an uncountable subset of $\mathbb{P}$. By hypothesis there are distinct $\alpha,\beta\in S$ so that $\{\alpha\}$ and $\{\beta\}$ are compatible. It is easy to see that $(L_\alpha\cap R_\beta)\cup (L_\beta\cap R_\alpha)=\emptyset$.
\end{claimproof}
\begin{claimproof}[Proof of $\Rightarrow$]Let $\mathcal{A}\in [\mathbb{P}]^{\omega_1}$. We can suppose without loss of generality that $\mathcal{A}$ satisfies the same conditions as in Lemma \ref{equivalencechi1forcing} for some $\sigma,\tau,\mu$ and  $R$. Note that the set $S=\{p(n-1)\,:\,p\in \mathcal{A}\}$ is uncountable so by the hypotheses we can find distinct $p$ and $q$ so that $$(L_{p(n-1)}\cap R_{q(n-1)})\cup (L_{q(n-1)}\cap R_{p(n-1)})=\emptyset.$$
Again, we claim that $p\cup q\in \mathbb{P}$. For this it suffices to prove that if $i,j<n$ then $L_{p(i)}\cap R_{q(j)}=\emptyset.$ Indeed, for any such $i,j$ it happens that $$L_{p(i)}\subseteq L_{p(n-1)}\cup \sigma(i,n-1)\subseteq L_{p(n-1)}\cup L_{q(i)}.$$
Therefore $L_{p(i)}\cap (R_{q(n-1)}\subseteq L_{p(n-1)})\cap R_{q(n-1)}\cup (L_{q(i)}\cap R_{q(n-1)})=\emptyset$. In an analogous way $R_{q(j)}\subseteq R_{q(n-1)}\cup R_{p(j)}$. Hence, we conclude that $$L_{p(i)}\cap R_{q(j)}\subseteq (L_{p(i)}\cap R_{q(n-1)})\cup (L_{p(i)}\cap L_{p(j)})=\emptyset.$$
\end{claimproof}
\end{proof}
\end{lemma}
The following theorem can be found in \cite{independenceforanalysts}, \cite{ScheepersGaps}, \cite{PartitionProblems} and \cite{yoriokadestructiblegaps}.

\begin{theorem}\label{cccdestructibilityequivalence}Let $(\mathcal{A},\mathcal{B})$ be an $(\omega_1,\omega_1)$-pregap:\begin{itemize}
    \item $\chi_1(\mathcal{L},\mathcal{R})$ is $ccc$ if and only if $(\mathcal{L},\mathcal{R})$ is a gap. In this case, there is some condition in $\chi_1(\mathcal{L},\mathcal{R})$ forcing $(\mathcal{L},\mathcal{R})$ to be indestructible.
    \item $\chi_0(\mathcal{L},\mathcal{R})$ is $ccc$ if and only if $(\mathcal{L},\mathcal{R})$ is destructible. In this case, there is some condition in $\chi_0(\mathcal{L},\mathcal{R})$ forcing $(\mathcal{L},\mathcal{R})$ to be separated.
\end{itemize}
\end{theorem}

A particular instance of destrucible gaps are the ones satifying the following condition.

\begin{definition}[Todor\v{c}evi\'c condition]Let $(L_\alpha,R_\alpha)_{\alpha\in\omega_1}$ be an $(\omega_1,\omega_1)$-pregap. We say that $(L_\alpha,R_\alpha)_{\alpha\in \omega_1}$ is \textit{Todor\v{c}evi\'c} if for any $S\in[\omega_1]^{\omega_1}$ there are $\alpha<\beta$ such that $L_\alpha \subseteq L_\beta$ and $R_\alpha\subseteq R_\beta$.
\end{definition}

Fulgencio Lopez and Stevo Todor\v{c}evi\'{c} showed that the existence of Todor\v{c}evi\'{c} gaps follows from $CA_3$ and $CA_2(part)$. Here we give another proof of such result based on the construction given in Theorem \ref{hausdorffgapconstruction}.
\begin{theorem}[Under $CA_2(part)$]\label{todorcevicgapconstruction}Let $\mathcal{F}$ be a $\mathcal{P}$-$2$-capturing $2$-construction scheme with $\mathcal{P}=\{P_0,P_1\}$. Also, consider $(L_\alpha,R_\alpha)_{\alpha\in\omega_1}$ the Hausdorff gap constructed in Theorem \ref{hausdorffgapconstruction}. Let $$C=\bigcup\{\,\{2k,2k+1\}\,:\,k\in P_0\backslash 1\}$$ and define $$L^C_\alpha=L_\alpha\cap C=\{ 2k+\Xi_\alpha(k)\,:\,k\in P_0\backslash 1\textit{ and }\Xi_\alpha(k)\geq 0\},$$
$$R^C_\alpha=R_\alpha\cap C=\{2k+(1-\Xi_\alpha(k))\,:\,k\in P_0\backslash 1,\textit{ and }\Xi_\alpha(k)\geq 0\}$$
for each $\alpha\in \omega_1$. Then $(\mathcal{L}|_C,\mathcal{R}|_C)=(L^C_\alpha,R^C_\alpha)_{\alpha\in \omega_1}$ is a Todor\v{c}evi\'c gap.
\begin{proof}$(L^C_\alpha,R^C_\alpha)_{\alpha\in\omega_1}$ is a pregap by similar reasons as the ones in the proof of Theorem \ref{hausdorffgapconstruction}. That is, $L^C_\alpha\cap R^C_\alpha=\emptyset$ for each $\alpha\in\omega_1$ and if $\alpha<\beta$ then:
\begin{enumerate}[label=$(\alph*)$] 

\item  $L^C_\alpha\backslash L^C_\beta\subseteq \{2k+\Xi_\alpha(k)\,:\,k\leq \rho(\alpha,\beta)\textit{ and }\Xi_\alpha(k)\not=\Xi_\beta(k)\},$
\item $R^C_\alpha\backslash R^C_\beta\subseteq \{ 2k+(1-\Xi_\alpha(k))\,:\,k\leq \rho(\alpha,\beta)\textit{ and }\Xi_\alpha(k)\not=\Xi_\beta(k)\}.$ 
\end{enumerate}
In particular, as a consequence of (a) we have:
\begin{enumerate}
\item[$(c)$] $L^C_\alpha\cap R^C_\beta\subseteq\{2k+\Xi_\alpha(k)\,:\,k\leq \rho(\alpha,\beta)\textit{ and }\Xi_\alpha(k)\not=\Xi_\beta(k)\}.$

\end{enumerate}

In order to prove that $(L^C_\alpha,R^C_\alpha)_{\alpha\in \omega_1}$ is a gap we will use the equivalence provided by Lemma \ref{gapforcinglemma}. Let $S\in [\omega_1]^{\omega_1}$. Since $\mathcal{F}$ is $\mathcal{P}$-2-capturing there are $\alpha<\beta\in S$ so that $\{\alpha,\beta\}$ is captured at some level $l\in P_0$. By Proposition \ref{deltarhoequalityprop} it follows that $\Delta(\alpha,\beta)=l=\rho(\alpha,\beta)$. Hence, $\Xi_\alpha(k)=\Xi_\beta(k)$ for any $k<l$. From this and by the the point $(c)$ above we deduce that $L^C_\alpha \cap R^C_\beta\subseteq \{2l+\Xi_\alpha(l)\}$. Furthermore, since $l=\rho(\alpha,\beta)$ and $\mathcal{F}$ is a construction scheme(or just by point (2) in Proposition \ref{deltarhoequalityprop}) we also have that $\Xi_\alpha(l)=0$ and $\Xi_\beta(l)=1$. Therefore we conclude that $L^C_\alpha\cap R^C_\beta=\{2k\}$, so in particular $(L^C_\alpha\cap R^C_\beta)\cup (L^C_\beta\cap R^C_\alpha)\not=\emptyset.$ As we said before, by Lemma \ref{gapforcinglemma} we are done.

Now we will prove that $(L^C_\alpha,R^C_\alpha)_{\alpha\in\omega_1}$ is Todor\v{c}evi\'c. Let $S\in [\omega_1]^{\omega_1}$. Since $\mathcal{F}$ is $\mathcal{P}$-2-capturing there are $\alpha<\beta\in S$ such that $\{\alpha,\beta\}$ is captured at some level $l\in P_1$. Again, by Proposition \ref{deltarhoequalityprop} and by the points (a) and (b) written at the beginning of the proof we deduce that $L^C_\alpha\backslash L^C_\beta\subseteq \{2l\} $ and $R^C_\alpha\backslash R^C_\beta\subseteq \{2l+1\}$. But $l\in P_1$ so $\{2l,2l+1\}$ has empty intersection with both $L^C_\alpha$ and $R^C_\alpha$ by definition of these two sets. In this way we conclude that $L^C_\alpha\backslash L^C_\beta=\emptyset$ and $R^C_\alpha\backslash R^C_\beta=\emptyset.$ In other words, $L^C_\alpha\subseteq L^C_\beta$ and $R^C_\alpha\subseteq R^C_\beta$. 
\end{proof}
\end{theorem}
It is natural to ask whether destructible can be constructed only by using a $2$-capturing construction scheme. In the next theorem we will show that this is not the case. But first, let us note that the proofs of both Lemma \ref{equivalencechi0forcing} and Lemma \ref{equivalencechi1forcing} actually yield the two following results.

\begin{lemma}\label{lemmachi0mf}Let $(\mathcal{L},\mathcal{R})=(L_\alpha,R_\alpha)_{\alpha\in\omega_1}$ be an $(\omega_1,\omega_1)$-gap. Suppose that $\mathcal{A}\in[\chi_0(\mathcal{L},\mathcal{R})]^{\omega_1}$. Then there are $n\in\omega$ and  $\mathcal{B}\in [\mathcal{A}]^{\omega_1}\cap \mathscr{P}([\omega_1]^n)$ an uncountable root-tail-tail $\Delta$-system such that for any two distinct $p,q\in  \mathcal{B}$, the following conditions are equivalent:
\begin{itemize}
    \item $p$ is compatible with $q$.
    \item $(L_{p(n-1)}\cap R_{q(n-1)})\cup (L_{q(n-1)}\cap R_{p(n-1)})=\emptyset.$ In other words, $\{p(n-1)\}$ is compatible with $\{q(n-1)\}.$
\end{itemize}
\end{lemma}

\begin{lemma}\label{lemmachi1mf}Let $(\mathcal{L},\mathcal{R})=(L_\alpha,R_\alpha)_{\alpha\in\omega_1}$ be an $(\omega_1,\omega_1)$-gap. Suppose that $\mathcal{A}\in[\chi_1(\mathcal{L},\mathcal{R})]^{\omega_1}$. Then there are $n,r\in\omega$ and  $\mathcal{B}\in [\mathcal{A}]^{\omega_1}\cap \mathscr{P}([\omega_1]^n)$  a root-tail-tail $\Delta$-system with root $R$ of cardinality $r$ such that for any two distinct $p,q\in  \mathcal{B}$, the following conditions are equivalent:
\begin{itemize}
    \item $p$ is compatible with $q$.
    \item $(L_{p(r)}\cap R_{q(r)})\cup (L_{q(r)}\cap R_{p(r)})\not=\emptyset.$ In other words, $\{p(r)\}$ is compatible with $\{q(r)\}.$ 
\end{itemize}
    
\end{lemma}
\begin{theorem}\label{nodestructiblemf}Let $\mathcal{F}$ be a $2$-capturing construction scheme. If $\mathfrak{m}_\mathcal{F}>\omega_1$, then there are no destructible gaps.
\begin{proof}Let us assume towards a contradiction that there is a destructible gap $(\mathcal{L},\mathcal{R})=(L_\alpha,R_\alpha)_{\alpha\in\omega_1}.$ The forcing $\chi_0(\mathcal{L},\mathcal{R})$ is never Knaster. Therefore, by Propositions \ref{propequivalenceknasterpreservescheme} and \ref{lemmaequivalencefunctioncapturingpreserving}, there are $X'\in [\chi_0(\mathcal{L},\mathcal{R})]^{\omega_1}$ and an injective function $\zeta':X\longrightarrow \omega_1$ in such way that for any two distinct $p,q\in X'$, if $\{\zeta'(p),\zeta'(q)\}$ is captured, then $p$ and $q$ are incompatible in $\chi_0(\mathcal{L},\mathcal{R})$. By virtue of the Lemma \ref{lemmachi0mf}, we may assume that the elements $X'$ are singletons. In this way, we can define $X=\{\alpha\in \omega_1\,:\,\{\alpha\}\in X'\}$ and $\zeta:X\longrightarrow \omega_1$ given by $\zeta(\alpha)=\zeta'(\{\alpha\})$.
Note that for any two distinct $\alpha,\beta\in X$, if $\{\zeta(\alpha),\zeta(\beta)\}$ is captured, then $$(L_\alpha\cap R_\beta)\cup (L_\alpha\cap R_\beta)\not=\emptyset.$$
\noindent
\underline{Claim}: $\chi_1( L_\alpha,R_\alpha)_{\alpha\in X}$ is $ccc$ and $n$-preserves $\mathcal{F}.$
\begin{claimproof}[Proof of claim]We will prove the claim by appealing to the Lemma  \ref{lemmaequivalencefunctioncapturingpreserving}. Let $\mathcal{A}$ be an uncountable subset of $\chi_1(L_\alpha,R_\alpha)_{\alpha\in X}$ and  $\nu:\mathcal{A}\longrightarrow \omega_1$ be an injective function. According to the Lemma \ref{lemmachi1mf}, we may suppose that there are $n,r\in \omega$ so that $\mathcal{A}\in [X]^n$ and it forms a root-tail-tail $\Delta$-system with root $R$ of cardinality $r$ so that for any two distinct $p,q\in \mathcal{A}$, if $(L_{p(r)}\cap R_{q(r)})\cup (L_{q(r)}\cap R_{p(r)})\not=\emptyset$, then $p$ and $q$ are compatible.  Given $p\in \mathcal{A}$, let us define $$A_p=\ \{\zeta(p(r)),\nu(p)\},$$ $$\alpha_p=\max(A_p).$$ 
Since both $\zeta$ and $\nu$ are injective functions and $p(r)\not=q(r)$ whenever $p\not=q$, we may assume that the $A_p$'s are pairwise disjoint. Even more,we can suppose that there are $k,a\in\omega$ so that the following conditions hold for any two distinct $p,q\in \mathcal{A}$:
\begin{enumerate}[label=$(\alph*)$]
    \item $\rho^{A_p}=k.$
    \item $\lVert\alpha_p\rVert_k=a.$
    \item If $h:(\alpha_p)_k\longrightarrow (\alpha_q)_k$ is the increasing bijection $h(\nu(p))=\nu(q)$ and $h(\zeta(p(r)))=\zeta(q(r))$.
      \item $\rho(\alpha_p,\alpha_q)>k.$
\end{enumerate}
Since $\mathcal{F}$ is assumed to be $2$-capturing, there are distinct $p_0,p_1\in \mathcal{A}$ for which $\{\alpha_{p_0},\alpha_{p_1}\}$ is captured. By virtue of the Lemma \ref{capturedfamiliestosetslemma} and the points $(a)$, $(b)$, $(c)$ and $(d)$ above, are $\{\nu(p_0),\nu(p_1)\}$ and $\{\zeta(p_0(r)),\zeta(p_1(r))\}$ are also captured. In particular, this means that $(L_{p_0(r)}\cap R_{p_1(r)})\cup (L_{p_1(r)}\cap R_{p_0(r)})\not=\emptyset.$ Thus, $p_1$ and $p_1$ are compatible. This finishes the proof. 

\end{claimproof}
By the previous claim $\chi_1(L_\alpha,R_\alpha)_{\alpha\in X}$ is an uncountable $ccc$ forcing which $2$-preserves $\mathcal{F}$. Since $\mathfrak{m}_\mathcal{F}>\omega_1$, it follows that $(L_\alpha,R_\alpha)_{\alpha\in X}$ is an indestructible gap due to the Lemma \ref{cccdestructibilityequivalence} (hence, so is $(\mathcal{L},\mathcal{R})$). This is a contradiction, so the proof is over.
\end{proof}
\end{theorem}

\begin{rem}\label{todorcevicgapremark}Under the hypotheses of the Theorem \ref{todorcevicgapconstruction} , define $L^{\omega\backslash C}_\alpha=L_\alpha\cap (\omega\backslash C)$ and $R^{\omega\backslash C}_\alpha=R_\alpha\cap (\omega\backslash C)$ for each $\alpha\in \omega_1$. A similar arguing yields that $(\mathcal{L}|_{\omega\backslash C},\mathcal{R}|_{\omega\backslash C})=(L^{\omega\backslash C}_\alpha,R^{\omega\backslash C}_\alpha)_{\alpha\in \omega_1}$ is also a Todor\v{c}evi\'c gap. 
\end{rem}

Recall that the product of two forcings, say $\mathbb{P}$ and $\mathbb{Q}$, is $ccc$ if and only if $\mathbb{P}$ is $ccc$ and $\mathbb{P}\Vdash\text{\say{ $\mathbb{Q}\textit{ is ccc }$}}$. Suppose we are in the particular situation where have $(\mathcal{L}_0,\mathcal{R}_0)$ and $(\mathcal{L}_1,\mathcal{R}_1)$ are two gaps such that $\chi_i(\mathcal{L}_0,\mathcal{R}_0)\times \chi_j(\mathcal{L}_1,\mathcal{R}_1)$ is $ccc$ for any $i,j\in 2$. According to the previous theorem, destroying or making indestructible one of these two gaps by any of the previous forcings will not destroy or or make indestructible the remaining one. This motivates the following definition.
\begin{definition}[Independent families of gaps] We say that a family $\langle(\mathcal{L}^c,\mathcal{R}^c)\rangle_{c\in I}$ of $(\omega_1,\omega_1)$-pregaps is \textit{independent} if:
$$\prod\limits_{c\in I}^{FS}\chi_{\phi(c)}(\mathcal{L}^c,\mathcal{R}^c)$$ 
is $ccc$ for any $\phi:I\longrightarrow 2.$ Additionally, we say that a coherent family of functions $\mathfrak{F}$ supported by an $\omega_1$-tower is \textit{independent} if the family $\langle (\mathcal{L}^{c(0)},\mathcal{L}^{c(1)})\rangle_{c\in[\omega_1\backslash \omega]^2}$ is independent.
\end{definition}
\begin{rem}If $\langle (\mathcal{L}^c,\mathcal{R}^c)\rangle_{c\in I}$ is an independent family of $(\omega_1,\omega_1)$-pregaps, then $(\mathcal{L}^c,\mathcal{R}^c)$ is a gap for any $c\in I.$
\end{rem}
\begin{theorem} [Under $CA_2(part)$]There are two $(\omega_1,\omega_1)$-gaps which are Todor\v{c}evi\'c but not independent.
\begin{proof}Let $\mathcal{F}$ be a $2$-capturing $2$-construction scheme and consider the gaps $(\mathcal{L}|_C,\mathcal{R}|_C)$ and $(\mathcal{L}|_{\omega\backslash C},\mathcal{R}|_{\omega\backslash C})$ defined in Theorem \ref{todorcevicgapconstruction} and Remark \ref{todorcevicgapremark} respectively. These two gaps are Todor\v{c}evi\'c, so we only need to prove that they are not independent. This will be done in the following claim.\\

\noindent
\underline{Claim}: $\mathbb{P}=\chi_0(\mathcal{L}|_C,\mathcal{R}|_C)\times \chi_0(\mathcal{L}|_{\omega\backslash C},\mathcal{R}|_{\omega\backslash C})$ is not $ccc$.
\begin{claimproof}[Proof of claim] Suppose towards a contradiction that the forcing $\mathbb{P}$ is $ccc$. By Theorem \ref{cccdestructibilityequivalence}, it follows that there is $G$ a $\mathbb{P}$-generic filter over $V$ so that both of $(\mathcal{L}|_C,\mathcal{R}|_C)$ and $(\mathcal{L}|_{\omega\backslash C},\mathcal{R}|_{\omega\backslash C})$ are separated in $V[G]$ by some $C_0$ and $C_1$ respetively. In order to finish, just note that $C_0\cup C_1$ separates the Hausdorff gap $(L_\alpha,R_\alpha)_{\alpha\in \omega_1}$. This is a contradiction since Hausdorff gaps  undestructible. Thus, the proof is over.  
\end{claimproof}
    
\end{proof}
    
\end{theorem}
In \cite{Shelahabrahamincompactness}, Uri Abraham and Saharon Shelah used large independent families of Suslin trees in order to code certain sets of reals. Independent families of gaps can be used in a similar way.
In \cite{yoriokadestructiblegaps}, Teruyuki Yorioka proved, assuming the $\Diamond$-principle, that there is an independent family of $2^{\omega_1}$ gaps. In here we will construct an independent coherent family of functions from $FCA$. For this purpose, we will need the following proposition (see \cite{StevoIlias}).
\begin{proposition}\label{independentccc}Let $\langle (\mathcal{L}^c,\mathcal{R}^c)\rangle_{c\in I}$ be a finite family of $(\omega_1,\omega_1)$-gaps indexed as $(L^c_\alpha,R^c_\alpha)_{\alpha\in X_c}$ for each $c\in I$. Also let $\phi:I\longrightarrow 2$. Suppose that for any uncountable $\mathcal{B}\subseteq \prod\limits_{i\in I}X_i$ there are distinct $g,h\in \mathcal{B}$ so that for any $c\in I$, if $g(c)\not=h(c)$ then:$$(L^c_{g(c)}\cap R^c_{h(c)})\cup (L^c_{h(c)}\cap R^c_{g(c)})=\emptyset\textit{ if }\phi(c)=0,$$
$$(L^c_{g(c)}\cap R^c_{h(c)})\cup (L^c_{h(c)}\cap R^c_{g(c)})\not=\emptyset\textit{ if }\phi(c)=1.$$
Then $\mathbb{P}=\prod\limits_{c\in I}\chi_{\phi(c)}(\mathcal{L}^c,\mathcal{R}^c)$ is $ccc$.
\begin{proof}
Let $\mathcal{A}$ be an uncountable subset of $\mathbb{P}$. For any $c\in I$ consider $\mathcal{A}^c=\{p(c)\,:\,p\in \mathcal{A}\}$. Note that $\mathcal{A}^c\subseteq \chi_{\phi(c)}(\mathcal{L}^c,\mathcal{R}^c)$. In this way we can suppose that each $\mathcal{A}^c$ satisfies the conditions listed in the proof of Lemma \ref{equivalencechi1forcing} for some $n_c,\sigma_c,\tau_c,\mu_c, R_c$ and $r_c$. For any $p\in \mathcal{A}$ define $g_p\in \prod\limits_{c\in I}X_c$ as:
$$g_p(i)=\begin{cases}p(c)(n_c-1)&\textit{ if }\phi(c)=0\\
p(c)(r_c)&\textit{ if }\phi(i)=1
    
\end{cases}$$
By the hypotheses there are distinct $p,q\in \mathcal{A}$ such that for any $c\in I$, if $g_p(c)\not=g_q(c)$ then
$$(L^c_{g_p(c)}\cap R^c_{g_q(c)})\cup (L^c_{g_q(c)}\cap R^c_{g_p(c)})=\emptyset\textit{ if }\phi(i)=0,$$
$$(L^c_{g_p(c)}\cap R^c_{g_q(c)})\cup (L^c_{g_q(c)}\cap R^c_{g_p(c)})\not=\emptyset\textit{ if }\phi(c)=1.$$

By arguing on each coordinate in the same way as we did in the proof of Lemmas \ref{equivalencechi0forcing}
 an \ref{equivalencechi1forcing} we conclude that the function $p'\in \prod\limits_{c\in I} [X_c]^{<\omega}$ given by
 $p'(c)=p(c)\cup q(c)$ for each $c\in I$, is an element of $\mathbb{P}$ which is obviously smaller than both $p$ and $q$. This finishes the proof.
 \end{proof}
\end{proposition}

\begin{theorem}[FCA]\label{independentcoherentschemefca} There is an independent coherent family of functions supported by an $\omega_1$-tower.
\begin{proof}The construction here is similar to the one in Theorem \ref{hausdorffcoherenttheorem}. Fix a type $\langle m_k,n_{k+1},r_{k+1}\rangle_{k\in\omega}$ such that $n_{k+1}>2^{r_{k+1}^2}$ for all $k\in\omega$, and let $\mathcal{F}$ be a fully capturing construction scheme of that type. For each $k>0$, let $$N_k=\{k\}\times [r_k]^2,$$
$$N=\bigcup\limits_{k>0}N_k.$$
Also enumerate $\mathscr{P}([r_k]^2)$ (possibly with repetitions) as $\langle S^k_i\rangle_{i<n_k}$ in such way that $S^k_0=S^k_1=[r_k]^2$. We start by defining an $\omega_1$-tower over $N$.  
Given $\alpha\geq\omega$ and $k>0$ we define $T^k_\alpha\subseteq N_k$ as follows:
$$T^k_\alpha=\begin{cases}\emptyset &\textit{ if }\Xi_\alpha(k)=-1\\
\{k\}\times S^k_{i}&\textit{ if }\Xi_{\alpha}(k)=i\geq 0
    
\end{cases}$$
Finally, let $T_\alpha=\bigcup_{k\in\omega\backslash 1} T^k_\alpha$. In order to prove that $\mathcal{T}=\langle T_\alpha\rangle_{\alpha\in \omega_1\backslash \omega}$ is in fact an $\omega_1$-tower just note that if $\alpha<\beta\in \omega_1\backslash \omega$ and $k>\rho(\alpha,\beta)$ then  $T^k_\alpha\subseteq T^k_\beta$ by means of the point (c) of Lemma \ref{lemmaxi}.

Now we will construct a coherent family of functions supported by $\mathcal{T}$. For this, let $\beta\geq \omega$. Note that if $x\in T_\beta$ then $x=(k,s)$ where $k>0$, $\Xi_k(\beta)\geq 0$ and $s\in S^k_{\Xi_\beta(k)}\subseteq [r_k]^2$. In this way, we can define $f_\beta:T_\beta\longrightarrow \beta$ as:
$$f_\beta (k, s)=\begin{cases}(\beta)_k(s(0))&\textit{ if }\Xi_\beta(k)=0\\
(\beta)_k(s(1))&\textit{ if }\Xi_\beta(k)>0

\end{cases}$$
It is easy to check that  $\langle f_\alpha\rangle_{\alpha\in\omega_1\backslash \omega}$ is a coherent family of functions supported by $\mathcal{T}$. This is done, again, by appealing to the point (c) of Lemma \ref{lemmaxi}.

The only thing left to show is that $\langle f_\alpha\rangle_{\alpha\in\omega_1\backslash \omega}$ is independent. Recall that a finite support product of forcings is $ccc$ if and only if each finite subproduct is $ccc$. Let $I$ be a non-empty finite subset of $[\omega_1\backslash \omega]^2$  and let  $\phi: I\longrightarrow 2.$ We will finish by proving the following claim.\\\\
\underline{Claim}: $\mathbb{P}=\prod\limits_{c\in I}\chi_{\phi(c)}(\mathcal{L}^{c(0)},\mathcal{L}^{c(1)})$ is $ccc$.

\begin{claimproof}[Proof of claim] The proof of this claim will be performed by appealing to the equivalence provided by Proposition \ref{independentccc}. First note that for any $c\in I$, $(\mathcal{L}^{c(0)},\mathcal{L}^{c(1)})= (L^{c(0)}_\alpha, L^{c(1)}_\alpha)_{\alpha\in X_c}$ where $X_c= \omega_1\backslash (c(1)+1)$. Now, let $\mathcal{B}$ be an uncountable subset of $\prod\limits_{c\in I} X_c$. For any $g\in \mathcal{B}$, define $D_g=\bigcup I\cup im(g)$. We need to prove the following claim.
\begin{enumerate}[label=$(\arabic*)$]
    \item If $g,h\in B$ then $|im(g)|=|im(h)|$. Furthermore, if $\psi:im(g)\longrightarrow im(h)$ is the increasing bijection then $\psi(g(c))=h(c)$ for any $c\in I$.
    \item $\{D_g\,:g\in \mathcal{B}\}$ is a root-tail-tail $\Delta$-system with root $\bigcup I$.
\end{enumerate}
Fix $k>\rho^{\cup I}$. Since $\mathcal{F}$ is a fully capturing construction scheme, there is $l>k$ and there are distinct $g_0,\dots, g_{n_l-1}\in \mathcal{B}$ so that the family $\{ D_{g_0},\dots, D_{g_{n_l-1}}\}$ is a captured at  level $l$.  According to (2) and by Lemma \ref{intersectionrhoisomorphiclemma}, $\Xi_{\cup I}(l)={-1}$. In this way, the set $S=\{ \,\{\lVert c(0)\rVert_l, \lVert c(1)\rVert_l\}\,:\,c\in I\textit{ and }\phi(c)=1\}$ is a subset of $[r_l]^2$. Therefore we can find $i<n_l$ so that $S^l_i=g_i$. The proof will be over once we show the following subclaim.
\\\\
\underline{Subclaim}: $g_0$ and $g_i$ satisfy the hypotheses of Proposition \ref{independentccc}. That is, for any $c\in I$, if $g_0(c)\not=g_i(c)$ then:\footnote{In this case $(\mathcal{L}^c,\mathcal{R}^c)=(L^{c(0)}_\alpha,L^{c(1)}_\alpha)_{\alpha\in X_c}$. Therefore $L^c_{g_0(c)}=L^{c(0)}_{g_0(c)}$, $R^c_{g_i(c)}=L^{c(1)}_{g_i(c)}$,  $L^c_{g_i(c)}=L^{c(0)}_{g_i(c)}$ and  $R^c_{g_0(c)}=L^{c(1)}_{g_0(c)}$.}$$(L^{c(0)}_{g_0(c)}\cap L^{c(1)}_{g_i(c)})\cup (L^{c(0)}_{g_i(c)}\cap L^{c(1)}_{g_0(c)})=\emptyset\textit{ if }\phi(c)=0,$$
$$(L^{c(0)}_{g_0(c)}\cap L^{c(1)}_{g_i(c)})\cup (L^{c(0)}_{g_i(c)}\cap L^{c(1)}_{g_0(c)})\not=\emptyset\textit{ if }\phi(c)=1.$$

\begin{claimproof}[Proof of subclaim] Let $c\in I$. Then  $g(c)\not=h(c)$ due to condition (2). Furthermore,   $\rho(g(c),h(c))=l=\Delta(g(c),h(c))$ for any $c\in I.$ This follows from condition (1) together with the fact that $D_{g_0}$ is strongly isomorphic to $D_{g_i}$. By the definitions of $f_{g_0(c)}$ and $f_{g_i(c)}$ and the points (a) and (c) of Lemma \ref{lemmaxi}, it follows that $$L^{c(0)}_{g_0(c)}\cap L^{c(1)}_{g_i(c)}=\big(f^{-1}_{g_0(c)}[\{c(0)\}]\cap f^{-1}_{g_i(c)}[\{c(1)\}]\big)\cap T^l_{g_0(c)}.$$  
On one hand, since $\Xi_{g_0(c)}(l)=0$ then $$f^{-1}_{g_0(c)}[\{c(0)\}]\cap T^l_{g_0(c)}=\{l\}\times \{ s\in [r_l]^2\,:\, s(0)=\lVert c(0)\rVert_l \,\}.$$
On the other hand, since $\Xi_{g_i(c)}(l)>0$ and $S=S^l_i$ then 
$$f^{-1}_{g_0(c)}[\{c(1)\}]\cap T^l_{g_0(c)}=\{l\}\times \{ s\in S\,:\, s(1)=\lVert c(1)\rVert_l \,\}.$$
Thus,  $L^{c(0)}_{g_0(c)}\cap L^{c(1)}_{g_i(c)}\subseteq \{\,\{\lVert c(0)\rVert_l, \lVert c(1)\rVert_l\}\, \}$. Furthermore, such intersection is non-empty if and only if $\phi(c)=1$ by the definition of $S$.  By arguing in a similar way and using that $\lVert c(0)\rVert_l< \lVert c(1)\rVert_l$ we conclude that $L^{c(0)}_{g_i(c)}\cap L^{c(1)}_{g_0(c)}$ is always empty. This finishes the proof of the subclaim.

\end{claimproof}
\end{claimproof}
 \end{proof}
\end{theorem}

%% file: chapters/Trees_and_lines.tex
\section{Trees and lines}
The purpose of this section is to study uncountable trees and lines (totally ordered sets) from the point of view of construction schemes. Before dealing with any specific notion, it is convenient to recall the definition of a tree. Complete introductions to trees, lines and the duality between them can be found in  \cite{jechsettheory}, \cite{kunensettheory} and \cite{stevotreesandlinearlyorderedsets}.
\begin{definition}[Tree]A partial ordered set $(T,<)$ is called a \textit{tree} if $(-\infty,t)_T$ is well-ordered for any $t\in T$. 
\end{definition}
When dealing with trees, we apply some changes in the notation related to partial orders. In this realm, intervals  $(-\infty,t)_T$ and $(t,\infty)_T$ are denoted as $t_{\downarrow_T}$ and $t_{\uparrow_T}$ respectively. As usual, whenever there is no risk of confusion we will write those sets simply as $t_\downarrow$ and $t_\uparrow$ respectively. 

It is easy to see each tree $T$ is well-founded. Thus, according to the Section \ref{settheoreticnotation}, there is a unique function $rank_T:T\longrightarrow Ord$ so that: $$rank_T(t)=\sup(\,rank_T(s)+1\
,:\,\,s<t\,).$$
Due to uniqueness, it can be shown that $rank_T(t)=ot(t_{\downarrow_T})$ for all $t\in T$. In this way, for each ordinal $\alpha$, the level $\alpha$ of $T$ can be described as $T_\alpha=\{\,t\in T\,:\,ot(t_\downarrow)=\alpha\,\}.$ 
Given an ordinal $\beta$ and $t\in T$, we define $$T|_\beta=\bigcup\limits_{\alpha<\beta} T_\alpha,$$
$$T|_t=\{s\in T\,:\,s\textit{ is comparable with }t\}.$$
Finally, we say that $B\subseteq T$ is a \textit{branch} if it is a maximal chain in $T$.

\subsection{Countryman lines and Aronszajn trees}
Countryman lines are certain kind of linear orders whose existence was proposed by Roger Simmons Countryman in the 1970's. In \cite{decomposingsquares}, Saharon Shelah constructed for the first time Countryman lines without appealing to any extra axioms. Years later, Stevo Todor\v{c}evi\'{c} presented in \cite{partitioningpairs} an easy construction of such a line using walks on ordinals. The reader may find some variations of Todor\v{c}evi\'c original proof in \cite{hudson2007canonical} and \cite{Walksonordinals}. 

\begin{definition}[Countryman line] We say that a totally ordered set $(C,<)$ is a Countryman line if $C$ is uncountable and $C^2$ can be covered by countably many chains\footnote{Here, we consider the order over $C^2$ given by $(x,y)\leq (w,z)$ if and only $x\leq w$ and $y\leq z$.}. 
\end{definition}
In \cite{Afiveelementbasisjustin}, Justin Moore put Countryman lines into the realm of canonical objects inside Mathematics. Specifically, he proved the following theorem.
\begin{theorem}[Under $PFA$]Let $C$ be any Countryman line and $X$ be a fixed subset of $\mathbb{R}$ of cardinality $\omega_1$. Then any uncountable order contains an isomorphic copy of one of the following five orders:
$$X,\,\omega_1,\,\omega_1^*,\,C,\textit{ or }C^*.$$
Where $\omega_1^*$ and $C^*$ denote $\omega_1$ and $C$ but with the reverse order.
\end{theorem}
We will now define a Countryman line using construction schemes. It should be noted that in Chapter 3 of \cite{Walksonordinals}, Todor\v{c}evi\'{c} showed that ordinal metrics can also be used to define Countryman lines. Although working with construction schemes is equivalent to working with their induced ordinal metric, the definition of the line appearing in this thesis differs from the ones in \cite{Walksonordinals}.\\
The following proposition will help by giving us a better picture of the next definitions.
\begin{proposition}\label{propositionintersectiondelta}Let $\alpha,\beta\in\omega_1$ be distinct ordinals and $k<\Delta(\alpha,\beta)$. Then $\Xi_\delta(\Delta(\alpha,\beta))=-1$ for each $\delta\in (\alpha)_k\cap (\beta)_k$. In particular $|(\alpha)_k\cap (\beta)_k|\leq r_{\Delta(\alpha,\beta)}.$
\begin{proof}Let $\delta\in (\alpha)_k\cap (\beta)_k$. Suppose towards a contradiction that $\Xi_\delta(\Delta(\alpha,\beta))\not=-1$. Thus, according to the part (c)
 of Lemma \ref{lemmaxi} and the fact that $\rho(\alpha,\delta),\rho(\beta,\delta)\leq k<\Delta(\alpha,\beta)$, we have that $\Xi_\alpha(\Delta(\alpha,\beta))=\Xi_\delta(\Delta(\alpha,\beta))=\Xi_\beta(\Delta(\alpha,\beta)).$ This is impossible since $\Xi_\alpha(\Delta(\alpha,\beta)\not=\Xi_\beta(\Delta(\alpha,\beta))$ due to the point (d) of the same lemma.
\end{proof}
\end{proposition}

\begin{definition}\label{definitioncalphabeta}Let $\alpha,\beta\in \omega_1$ be distinct ordinals and $k=\Delta(\alpha,\beta)-1$. We define $c(\alpha,\beta)$ as:$$c^\alpha_\beta=\min(\,(\alpha)_k\backslash (\beta)_k\,)$$
\end{definition}
\noindent
\begin{minipage}[b]{0.45\linewidth}
\begin{rem} Note that $c^\alpha_\beta\not=c^\beta_\alpha$ for any two distinct $\alpha,\beta\in\omega_1$. Furthermore, if  $k=\Delta(\alpha,\beta)-1$, then $(c^\alpha_\beta)^-_k=(\alpha)_k\cap (\beta)_k=(c^\beta_\alpha)^-_k.$ In this way, $\lVert c^\beta_\alpha\rVert_k=\lVert c^\alpha_\beta\rVert_k\leq r_{k+1}$ due to the Proposition \ref{propositionintersectiondelta}.
\end{rem}
\end{minipage}\hspace{1.5cm}\begin{minipage}[b]{3.5cm}
\includegraphics[width=7cm, height=3.5cm]{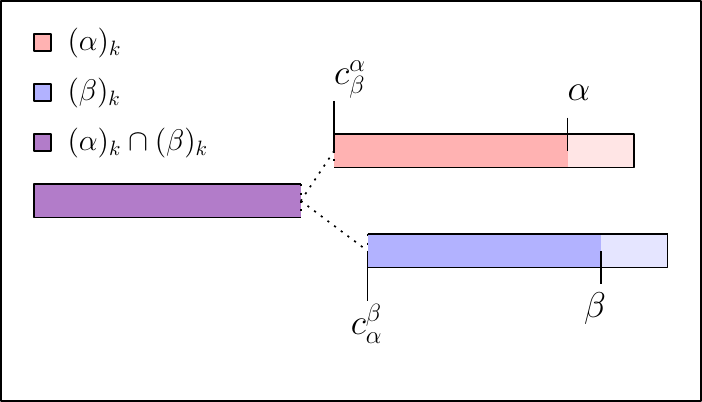}
\end{minipage}\\\\
In the next definition we define the relation that turns $\omega_1$ into a Countryman line.

\begin{definition}\label{definitioncountrymanscheme}Let $\mathcal{F}$ be a construction scheme. Given different $\alpha,\beta\in \omega_1$ and $k=\Delta(\alpha,\beta)-1$, we  recursively decide whether $\alpha<_\mathcal{F} \beta$ when one of the following conditions holds:\begin{enumerate}[label=$(\alph*)$]
    \item $\lVert c^\alpha_\beta\rVert_k=r_{k+1}$ and $\Xi_\alpha(k+1)<\Xi_\beta(k+1).$
   
    \item  $\lVert c^\alpha_\beta\rVert_k<r_{k+1}$ and $c^\alpha_\beta<_\mathcal{F} c^\beta_\alpha.$
\end{enumerate}
Where $k=\Delta(\alpha,\beta)-1.$
\end{definition}
\noindent
\begin{minipage}[b]{7cm}
\centering
   \includegraphics[width=7cm, height=3.5cm]{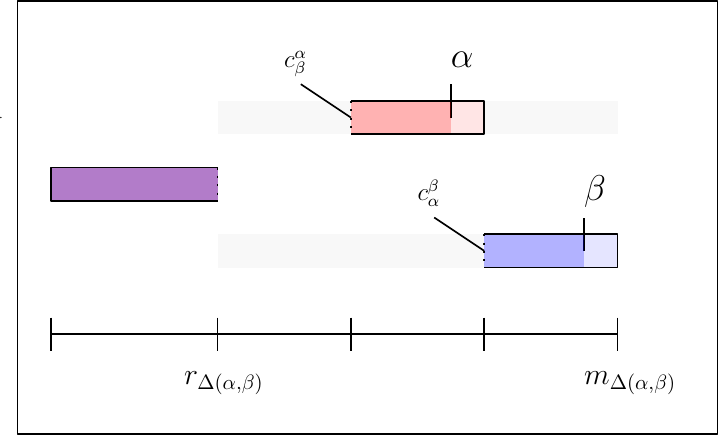}
 \textit{\small In here, $|(\alpha)_k\cap (\beta)_k|=r_{\Delta(\alpha,\beta)}$. Thus, we decide whether $\alpha<_\mathcal{F}\beta $ by comparing $\Xi_\alpha(\Delta(\alpha,\beta))$ and $\Xi_\beta(\Delta(\alpha,\beta)).$}
 
\end{minipage}\hspace{1cm}
\begin{minipage}[b]{7cm}
\centering
   \includegraphics[width=7cm, height=3.5cm]{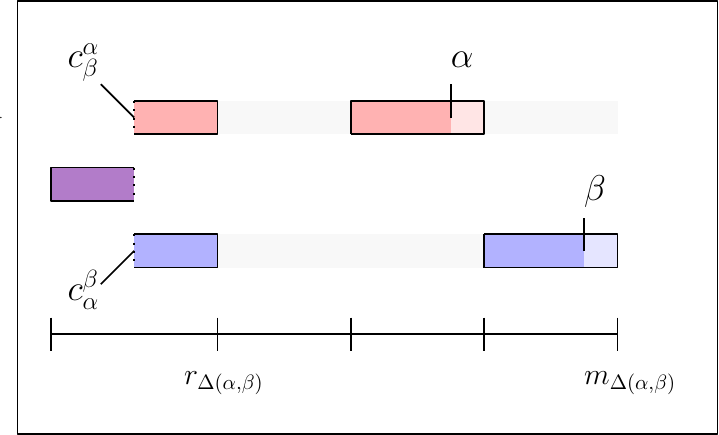}\\
   \textit{\small In here, $|(\alpha)_k\cap (\beta)_k|<r_{\Delta(\alpha,\beta)}$. Thus, we decide whether $\alpha<_\mathcal{F}\beta $ by comparing $c^\alpha_\beta$ and $c^\beta_\alpha.$}
 
\end{minipage}

\begin{lemma}\label{countrymanlemma2} Suppose $\alpha,\beta\in \omega_1$ are distinct and let  $k=\Delta(\alpha,\beta)-1$. Consider $h:(\alpha)_{k}\longrightarrow (\beta)_{k}$ the increasing bijection and take $\gamma\in (\alpha)_k$  such that $h(\gamma)\not=\gamma$. Then the following statements are equivalent:
\begin{enumerate}[label=$(\arabic*)$]
    \item $\alpha<_\mathcal{F}\beta$.
    \item $\gamma<_\mathcal{F}h(\gamma)$.
\end{enumerate}
\begin{proof}
The proof is carried by induction over $\alpha$ and $\beta$. So suppose that we have proved the Lemma for each $\alpha'<\alpha$ and $\beta'<\beta$. We start with some remarks about $\gamma$. First of all, $\Delta(\gamma,h(\gamma))\geq \Delta(\alpha,\beta)$ due to the point (a) of Lemma \ref{lemmainequalitiesdelta}. Furthermore, according to Remark \ref{remarkbijectionidentityballs} it follows that $\gamma\in (\alpha)_k\backslash (\beta)_k$ and $h(\gamma)\in (\beta)_k\backslash (\alpha)_k$. This means that $c^\alpha_\beta\in( \gamma)_k$ and $c^\beta_\alpha\in(h(\gamma))_k$. Lastly,  as $(\alpha)_k\cap (\beta)_k$ is an initial segment of both $(\alpha)_k$ and $(\beta)_k$, we conclude that $(\alpha)_k\cap (\beta)_k=(\delta)_k\cap (h(\delta))_k$. We divide the rest of the proof into two cases.\\

\noindent
\underline{Case 1}: If $\Delta(\alpha,\beta)=\Delta(\gamma,h(\gamma))$.
\begin{claimproof}[Proof of case] In this case we have that $c^\alpha_\beta=c^\gamma_{h(\gamma)}$ and $c^\beta_\alpha=c^{h(\gamma)}_\gamma$. Thus, if $\lVert c^\alpha_\beta\rVert_k< r_{k+1}$ then $\alpha<_\mathcal{F}\beta$ if and only if $c^\alpha_\beta<_\mathcal{F}c^\beta_\alpha$ if and only if $\gamma<_\mathcal{F}h(\gamma)$. On the other hand, if $\lVert c^\alpha_\beta\rVert_k=r_{k+1}$ then $\alpha<_\mathcal{F} \beta$ if and only if $\Xi_\alpha(k+1)<\Xi_\beta(k+1)$, and $\gamma<_\mathcal{F} h(\gamma)$ if and only if $\Xi_\gamma(k+1)<\Xi_{h(\gamma)}(k+1)$. Now, as $k+1=\Delta(\gamma,h(\gamma))$ then both $\Xi_\gamma(k+1)$ and $\Xi_{h(\gamma)}(k+1)$ are non-negative. From this fact and the part (c) of Lemma \ref{lemmaxi}, we conclude that $\Xi_\alpha(k+1)=\Xi_\gamma(k+1)$ and $\Xi_\beta(k+1)=\Xi_{h(\gamma)}(k+1)$. In this way, $\alpha<_\mathcal{F} \beta$ if and only if $\gamma<_\mathcal{F} h(\gamma)$. 
\end{claimproof}

\noindent
\underline{Case 2}:  If $\Delta(\alpha,\beta)<\Delta(\gamma,h(\gamma)).$
\begin{claimproof}[Proof of case] First we argue that $\alpha<_\mathcal{F}\beta$ if and only if $c^\alpha_\beta<_\mathcal{F}c^\beta_\alpha$. For this purpose it is enough to show that, in this case, $\lVert c^\alpha_\beta\lVert_k<r_{k+1}$. Indeed, according to the point (a) of Lemma \ref{lemmaxi}, we have that $\Xi_\gamma(k+1)=\Xi_{h(\gamma)}(k+1)$. On the other hand, by the part (d) of such lemma, we also 
know that $\Xi_\alpha(k+1)\not=\Xi_\beta(k+1)$. Thus,   $\Xi_\gamma(k+1)=-1$ due to the part $(c)$ of Lemma \ref{lemmaxi}. Since $c^\alpha_\beta\in (\gamma)_k$, it follows that $\lVert c^\alpha_\beta\rVert_k<r_{k+1}$. In this way, $\alpha<_\mathcal{F} \beta$ if and only if $c^\alpha_\beta<_\mathcal{F} c^\beta_\alpha$ as we 
wanted. The proof of this case will be over once we show that $\gamma<_\mathcal{F}h(\gamma)$ if and only if $c^\alpha_\beta<_\mathcal{F}c^\beta_\alpha$.  Let $l=\Delta(\gamma,h(\gamma))-1$ and consider $h':(\gamma)_l\longrightarrow (h(\gamma))_l$ be the increasing bijection. As $k<l$ then $c^\alpha_\beta\in (\gamma)_k\subseteq (\gamma)_l$. Furthermore, $$h'(c^\alpha_\beta)=h'(\,(\gamma)_k(\lVert c^\alpha_\beta\rVert_k)\,)=(h(\gamma))_k(\lVert c^\alpha_\beta\rVert_k)=(h(\gamma))_k(\lVert c_\alpha^\beta\rVert_k)=c^\beta_\alpha.$$ 

As $c^\alpha_\beta\not=c^\beta_\alpha$, we may use 
the inductive hypotheses to conclude that $\gamma<_\mathcal{F}h(\gamma)$ if and only if $c^\alpha_\beta<_\mathcal{F}c^\beta_\gamma$. This finishes the proof. 
\end{claimproof}
\end{proof}
\end{lemma}
\begin{corollary}let $\alpha,\beta\in \omega_1$ be distinct ordinals. Then $\alpha<_\mathcal{F}\beta$ if and only if $c^\alpha_\beta<_\mathcal{F} c^\beta_\alpha$. 
\end{corollary}
\begin{corollary}\label{corollarycountrymanbijection}Suppose $\alpha,\beta\in \omega_1$ are distinct and let  $k<\Delta(\alpha,\beta)$. Consider $h:(\alpha)_{k}\longrightarrow (\beta)_{k}$ the increasing bijection and take $\gamma\in (\alpha)_k$  such that $h(\gamma)\not=\gamma$. Then the following statements are equivalent:
\begin{enumerate}[label=$(\arabic*)$]
    \item $\alpha<_\mathcal{F}\beta$.
    \item $\gamma<_\mathcal{F}h(\gamma)$.
\end{enumerate}
\begin{proof}  Let $\psi:(\alpha)_{\Delta(\alpha,\beta)-1}\longrightarrow(\beta)_{\Delta(\alpha,\beta)-1}$ be the increasing bijection. To prove the corollary, just notice that $\psi|_{(\alpha)_z}=\phi.$  
\end{proof}
\end{corollary}
\begin{lemma}\label{lemmadeltainitialsegment}Let $\alpha,\beta,\delta\in \omega_1$ be distinct ordinals and $k\in \omega$. If $(\alpha)_k\cap (\beta)_k\subseteq (\alpha)_k\cap (\delta)_k$ then $(\alpha)_k\cap (\beta)_k=(\beta)_k\cap (\delta)_k.$
\begin{proof}Since $(\alpha)_k\cap (\beta)_k\subsetneq (\alpha)_k\cap (\delta)_k$ then $(\alpha)_k\cap (\beta)_k= (\alpha)_k\cap (\delta)_k\cap (\beta)_k\subseteq (\delta)_k\cap (\beta)_k$. In order to prove that $(\beta)_k\cap (\delta)_k\subseteq (\alpha)_k\cap (\beta)_k$, note that both $(\alpha)_k\cap (\delta)_k$ and $(\beta)_k\cap (\delta)_k$ are initial segments of $(\delta)_k$ by virtue of the point (6) in Proposition \ref{closureprop1}. In this way, either $(\beta)_k\cap (\delta)_k\sqsubseteq(\alpha)_k\cap (\delta)_k$ or $(\alpha)_k\cap (\delta)_k\sqsubseteq (\beta)_k\cap (\delta)_k$. The latter alternative can not happen as it would imply that $(\alpha)_k\cap (\delta)_k\subseteq (\alpha)_k\cap (\beta)_k$. Thus, the first alternative holds.  In this way,  $(\beta)_k\cap(\delta)_k\subseteq (\beta)_k\cap (\alpha)_k\cap (\delta)_k\subseteq (\alpha)_k\cap (\beta)_k.$ This finishes the proof.
\end{proof}
    
\end{lemma}
\begin{proposition}\label{countrymanlinearorderprop}$(\omega_1,<_\mathcal{F})$ is a total order.
\begin{proof}The only non-trivial part of this task is to prove transitivity. This proof will be performed by induction. For this purpose let $\alpha,\beta,\delta\in \omega_1$ be distinct ordinals. Suppose that we have proved that for any  $\alpha'<\alpha$, $\beta'<\beta$ and $\delta'<\delta$, the triplet $\{\alpha',\beta',\delta'\}$ do not form a cycle. That is, neither $\alpha'<_\mathcal{F}\beta'<_\mathcal{F}\delta'<_\mathcal{F}<\alpha'$ nor $\alpha'>_\mathcal{F}\beta'>_\mathcal{F}>\delta'>_\mathcal{F}\alpha'$. We will show that the same holds for $\{\alpha,\beta,\delta\}$. By virtue of Lemma \ref{countrymanlemma3}, we may assume without loss of generality that $\Delta(\alpha,\beta)=\Delta(\alpha,\delta)\leq \Delta(\beta,\delta).$ Let $k=\Delta(\alpha,\beta)-1.$ For each distinct $x,y\in \{\alpha,\beta,\delta\}$ consider $e^x_y=\min(\,(x)_k\backslash (y)_k\,)$ and $h^x_y:(x)_k\longrightarrow (y)_k$ the increasing bijection. Observe that $h^z_x\circ h^y_z=h^y_x$, $h^x_y(e^x_y)=e^y_x$ and $e^x_y\not=e^y_x$. We divide the rest of the proof into the following cases.\\

\noindent
\underline{Case 1}: There are distinct $x,y,z\in \{\alpha,\beta,\delta\}$ for which $(x)_k\cap (y)_k\subsetneq (x)_k\cap (z)_k.$
\begin{claimproof}[Proof of case] By means of the Lemma \ref{lemmadeltainitialsegment},  $(x)_k\cap(y)_k=(z)_k\cap (y)_k$. It then follows that $e^y_z=e^y_x$. Now we will now show that $e^x_y=e^z_y$. First recall that both $(x)_k\cap (y)_k$ and $(x)_k\cap (z)_k$ are initial segments of $(x)_k$. From this fact we can deduce \break that  $(x)_k\cap (y)_k{\mathrlap{\sqsubseteq}{\,_\not}}$ $\; (x)_k\cap (z)_k.$  By minimality, $e^x_y\in (x)_k\cap (z)_k.$ As $e^x_y\notin(y)_k$, it follows that $e^x_y\geq e^z_y$. Moreover, since $(x)_k\cap (z)_k\sqsubseteq (z)_k$ then $e^z_y\in (x)_k\cap (z)_k$. Again, we have that $e^z_y\geq e^x_y$ because $e^z_y\notin (y)_k$. Therefore, $e^z_y=e^x_y$. According to the Corollary \ref{corollarycountrymanbijection}, $y$ and $x$ are ordered (with respect to $<_\mathcal{F}$) in the same way as $e^y_x$ and $h^y_x(e^y_x)=e^x_y$, and $y$ and $z$ are ordered in the same way as $e^y_z=e^y_x$ and $h^y_z(e^y_z)=e^z_y=e^x_y$. We conclude that either $y<_\mathcal{F} x,z$ or $x,z<_\mathcal{F} y.$ This finishes the proof of this case.
\end{claimproof}

\noindent
\underline{Case 2}: $(\alpha)_k\cap (\beta)_k=(\alpha)_k\cap (\delta)_k=(\beta)_k\cap (\delta)_k.$
\begin{claimproof}[Proof of case] We will divide the proof of this case into two subcases, but first, note that $e^\alpha_\beta, e^\beta_\alpha, e^\alpha_\delta$ and $e^\delta_\alpha$ are equal to $c^\alpha_\beta, c^\beta_\alpha, c^\alpha_\delta$ and $c^\delta_\alpha$ respectively. Moreover, $e^x_y=e^x_z$ for all distinct $x,y,z\in\{\alpha,\beta,\delta\}$.\\

\noindent
\underline{Subcase 1}: $|(\alpha)_k\cap (\beta)_k|<r_{k+1}$.
\begin{claimproof}[Proof of subcase] By the remarks made at the start of this case and by Corollary \ref{corollarycountrymanbijection}, we have that $\alpha,\beta$ and $\delta$ are ordered in the same way as $e^\alpha_\beta$, $e^\beta_\alpha$ and $e^\delta_\alpha$ respectively. Thus, by virtue of the inductive hypotheses, it suffices to show that $e^\alpha_\beta<\alpha$, $e^\beta_\alpha<\beta$ and $e^\delta_\alpha<\delta$. Indeed, according to the hypotheses of this subcase, $\lVert e^x_y\rVert_k=|(x)_k\cap(y)_k|<r_{k+1}$ for all $x,y\in\{\alpha,\beta,\delta\}$. On the other hand, since $k+1=\Delta(\alpha,\beta)=\Delta(\alpha,\delta)\leq \Delta(\beta,\delta)$, then  $\Xi_{x}(k+1)\geq 0$ for each $x\in \{\alpha,\beta,\delta\}$. In other words $\lVert x\rVert_k\geq r_{k+1}$. Particularly, $e^\alpha_\beta<\alpha$, $e^\beta_\alpha<\beta$ and $e^\delta_\alpha<\delta$, so we are done. 
\end{claimproof}

\noindent
\underline{Subcase 2}: If $|(\alpha)_k\cap (\beta)_k|=r_{k+1}.$
\begin{claimproof}[Proof of subcase] According to the definition of the order $<_\mathcal{F}$ and since $k+1=\Delta(\alpha,\beta)=\Delta(\alpha,\delta)$, it follows that $\alpha$ and $\beta$ are ordered (with respect to $<_\mathcal{F}$) in the same way as $\Xi_\alpha(k+1)$ and $\Xi_\beta(k+1)$ (with respect to the usual order), and $\alpha$ and $\delta$ are ordered in the same way as $\Xi_\alpha(k+1)$ and $\Xi_\beta(k+1)$. If $\Delta(\beta,\delta)>k+1$ then $\Xi_\beta(k+1)=\Xi_\delta(k+1)$ by the point (a) of Lemma \ref{lemmaxi}. Thus, either $\alpha<\beta,\delta$ or $\beta,\delta<\alpha$. On the other hand, if $\Delta(\beta,\delta)=k+1$ then $e^\beta_\delta=c^\beta_\delta$ and $e^\beta_\delta$.  Therefore, $\beta$ and $\delta$ are ordered in the same way as $\Xi_\beta(k+1)$ and $\Xi_\delta(k+1)$. Evidently, $\Xi_\alpha(k+1)$, $\Xi_\beta(k+1)$ and $\Xi_\delta(k+1)$ do not form a cycle, so we are done.
\end{claimproof}
There are no more subcases by virtue of the Proposition \ref{propositionintersectiondelta}. Thus, the proof is over.
\end{claimproof}
\end{proof}
\end{proposition}
\begin{theorem}\label{countrymantheoremjjj}$(\omega_1,<_\mathcal{F})$ is a Countryman line.
\begin{proof}We have already seen that $(\omega_1,<_\mathcal{F})$ is a total order in Proposition \ref{countrymanlinearorderprop}. The only thing left to show is that $\omega^2$ can be partitioned into countably many chains. For this purpose, it is sufficient to prove the same holds for $D=\{(\alpha,\beta)\in\omega_1^2\;\:;\alpha<\beta\}$. Given $x,y,z\in\omega$, define $$P(x,y,z)=\{(\alpha,\beta)\in D\;:\; \lVert\alpha\rVert_z=x,\;\lVert\beta\rVert_z=y \textit{ and }\rho(\alpha,\beta)=z\}.$$
It is clear that the family of all $P(x,y,z)$ is countable and covers $D$. Thus, the proof will be over once we show each $P(x,y,z)$ is a chain. Indeed, take $(\alpha,\beta),(\delta,\gamma)\in P(x,y,z)$ and suppose  $\alpha<_\mathcal{F}\delta$. Since $\lVert\beta\rVert_z=\lVert\gamma\rVert_z$ then $\Delta(\beta,\gamma)>z$.  In this way, we can consider $h$ 
 the increasing bijection from $(\beta)_{z}$ to $(\gamma)_{z}$. Note that $\alpha\in (\beta)_z$ and $$h(\alpha)=h(\,(\beta)_z(\lVert \alpha\rVert_z)\,)=h(\,(\beta)_z(x))=(\gamma)_z(x)=(\gamma)_z(\lVert \delta\rVert_z)=\delta.$$ So by Lemma \ref{countrymanlemma2},  $\beta<_\mathcal{F}\gamma.$ 
\end{proof}
    
\end{theorem}
Now, we pass to the study of trees which are closely related to Countryman lines.
\begin{definition}[Aronszajn tree] Let $(T,<)$ be a tree. We say that $T$ is an $\omega_1$-tree if $Ht(T)=\omega_1$,  and $T_\alpha$ is countable for each $\alpha<\omega_1.$ Furthermore, if $T$ is an $\omega_1$-tree without uncountable chains, we call it an \textit{Aronszajn tree.}
    \end{definition}
Aronszajn trees were constructed for the first time by Nachman Aronszajn in the 1930's.  An important class of Aronszajn trees are the so called $\textit{special}$. We say that an $\omega_1$-tree, say $(T,<)$, is \textit{special} if it can be written as a countable union of antichains. Since any chain intersect each antichain in at most one point, it follows that each special tree is Aronszajn. It is known that $MA$ implies that any Aronszajn tree is special. On the other hand, the $\Diamond$-principle implies the existence of non-special Aronszajn tree. In the following theorem, we give a simple construction of a special tree from a construction scheme. It is worth pointing out that Aronszajn trees can already be constructed using Countryman lines.
\begin{theorem}There is a special Aronszajn tree.

\begin{proof}
    Let $\mathcal{F}$ be a construction scheme. Given $\beta\in\omega_1$, consider the function $\rho_\beta:\beta+1\longrightarrow \omega$ defined as $\rho_\beta(\alpha)=\rho(\alpha,\beta).$ Let $$T=\{f\in \omega^{<\omega_1}\,:\,\exists \beta\in\omega_1\,(dom(f)=\beta+1\textit{ and }f=^*\rho_\beta)\}.$$
Given $\alpha<\beta\in\omega_1$, it is straightforward that $\rho_\alpha(\xi)=\rho_\beta(\xi)$ for each $\xi\in (\alpha+1)\backslash(\alpha)_{\rho(\alpha,\beta)}$. From this, we conclude that $f|_{\alpha+1}\in T$ for any $f\in T$ and $\alpha\in dom(f).$ Hence, $(T,\subseteq)$ is an $\omega_1$-tree.

Now, we will show that $(T,\subseteq)$ is in fact a special tree. For this, take an arbitrary $f\in \mathcal{F}$ and let $\beta_f$ be such that $dom(f)=\beta_f+1$. Now, consider $k_f\in\omega$ the minimal natural number with the following properties:\begin{enumerate}[label=$(\arabic*)$]
    \item $\rho_{\beta_f}(\xi)=f(\xi)$ for each $\xi\notin (\beta_f)_{k_f}$,
    \item $f(\xi)\leq k_f$ for each $\xi\in (\beta_f)_{k_f}$.
\end{enumerate}
For each $k,s\in\omega$, let $T(k,s)$ be the set of all $f\in T$ for which $k=k_f$ and $\lVert \beta_f\rVert_k=s.$\\

\noindent
\underline{Claim}: $T(k,s)$ is an antichain in $T$.
\begin{claimproof}[Proof of claim]
    
Indeed, take two distinct $f,g\in T(k,s)$. If $\beta_f=\beta_g$ there is nothing to do, so let us assume that $\beta_f<\beta_g$. Since    $\lVert \beta_f\rVert_k=s=\lVert \beta_g\rVert_k$, then $\rho(\beta_f,\beta_g)>k$. By definition of $k_f$ and $k_g$, we get that $f(\beta_f)\leq k<\rho(\beta_f,\beta_g)=g(\beta_f)$. This concludes the claim. 
\end{claimproof}
To finish, just note that $T=\bigcup\limits_{k,s\in\omega}T(k,s)$.
\end{proof}

\end{theorem}

    \subsection{Suslin trees}
    Georg Cantor proved that $\mathbb{R}$ is, up to isomorphism, the unique total order $X$ with the following properties:
    \begin{enumerate}[label=$(\arabic*)$]
        \item $X$ is dense,
        \item $X$ has no end-points,
        \item $X$ is complete. That is, any of its bounded subsets have a supremum,
        \item $X$ is is separable.
    \end{enumerate}
    In 1920, Mikhail Suslin asked whether the same characterization is true when the property (4) is changed by the condition: \begin{enumerate}[label=(4)']
    \item$X$ is $ccc$. That is, any family of disjoint open intervals in $X$ is at most countable.
    \end{enumerate}
Any total order $X$ satisfying the properties (1), (2), (3) and (4)' but not (4) is call a \textit{Suslin line}. Suslin's hypothesis states that there are no Suslin lines. In the 1930's Djuro Kurepa show that the existence of Suslin lines is equivalent to the existence of the trees that we now know as Suslin.

\begin{definition}[Suslin tree]Let $(T,<)$ be an $\omega_1$-tree. We say that $T$ is \textit{Suslin} if $(T,>)$ does not contain uncountable chains or antichains.
\end{definition}

Thomas Jech and Stanley Tennenbaum independently proved the consistency of the existence  Suslin trees in \cite{Jechsuslin} and \cite{Tennenbaumsuslin} respectively. Later, Ronald Jensen showed that the existence of Suslin trees follows from the $\Diamond$-principle. Finally, Robert M. Solovay and Tennenbaum proved that $ZFC$ is consistent with the non existence of Suslin trees (see \cite{Iteratedcohensuslin}).

In this subsection, we will show that the existence of two distinct types of Suslin trees follows from the capturing axioms $FCA(part)$ and $FCA$.
\begin{definition}[Coherent tree]\label{coherenttreedef} Let $(T,<)$ be a Suslin tree. We say that $T$ is \textit{coherent} if there is a family of functions $\langle f_\alpha\rangle _{\alpha\in\omega_1}$ so that the following properties hold for all $\alpha<\beta\in\omega_1$:
\begin{enumerate}[label=$(\arabic*)$]
    \item $f_\alpha:\alpha\longrightarrow \omega,$
    \item $f_\alpha=^*f_\beta|_\alpha.$
\end{enumerate}
Furthermore, 
$T=\{f_\alpha|_\xi\,:\,\xi\leq\alpha<\omega_1\}$ and the order $<$ of $T$ coincides with $\subseteq.$

\end{definition}

Coherent trees have been extensively studied in the past, since they have many interest forcing properties (see \cite{microscopic1}, \cite{Microscopic2},  \cite{variationsouslin}, \cite{katetovproblem}, \cite{chainconditionsmaximalmodels} and \cite{gapstructure}).

\begin{theorem}[Under FCA(part)]\label{coherentsuslinscheme}There is a coherent Suslin tree.
\begin{proof}Let $\tau=\langle m_k,n_{k+1},r_{k+1}\rangle_{k\in\omega}$ be a type so that $n_{k+1}\geq 2^{m_k-r_{k+1}}+1$ for each $k\in\omega$. Furthermore, let $\mathcal{P}=\{P_c,P_a\}$ be a partition of $\omega$ compatible with $\tau$. We will build a coherent Suslin tree by using a $\mathcal{P}$-fully capturing construction scheme $\mathcal{F}$ of type $\tau$.

For each  $k\in\omega$, we first enumerate (possibly with repetitions) the set of all functions from $m_k\backslash r_{k+1}$ into $2$ as $\langle g^k_i\rangle _{0<i<n_{k+1}}$. Now, given $\beta\in\omega_1$, we define $f_\beta:\beta\longrightarrow 2$ as follows:
$$f_\beta(\xi)=\begin{cases}1 &\textit{if } \Xi_\xi(l)= 0,\,\Xi_\beta(l)=1\textit{ and }l\in P_c\\
g^l_i(\lVert \xi\rVert_l)&\textit{if }\Xi_\xi(l)=0\textit{ and }l\in P_a\\
0&\textit{ otherwise }
\end{cases}$$
Here, $l=\rho(\xi,\beta)$.\\
\noindent
Before defining the tree, we will prove that the sequence $\langle f_\alpha\rangle_{\alpha\in\omega_1}$ satisfies the condition (2) of the Definition \ref{coherenttreedef}. This will be done in the following claim.\\

\noindent
\underline{Claim 1}: Let $\xi<\alpha<\beta\in \omega_1$. If $\xi\not\in(\alpha)_{\rho(\alpha,\beta)}$ then $f_\alpha(\xi)=f_\beta(x)$.
\begin{claimproof}[Proof of claim]
Note that $\rho(\xi,\alpha)>\rho(\alpha,\beta)$ we the hypotheses. As $\rho$ is an ordinal metric, we use the previous fact to conclude that $\rho(\xi,\alpha)=\rho(\xi,\beta)$. 
Let us call this number $l$. According to the part (c) of Lemma \ref{lemmaxi}, either $\Xi_\alpha(l)=-1$ or $\Xi_\alpha(l)=\Xi_\beta(l)$. On the other hand,  $0\leq \Xi_\xi(l)<\Xi_\alpha(l)$ due to part (b) of the same lemma. In this way, $\Xi_\alpha(l)=\Xi_\beta(l)$. By definition of $f_\alpha$ and $f_\beta$,  $f_\beta(\xi)=f_\alpha(\xi)$.
\end{claimproof}
Now we define $T$ as expected. That is, $T=\{f_\alpha|_\xi\,:\,\xi\leq \alpha<\omega_1\}$ and consider it ordered by $\subseteq$. In the next two claims, we will prove that $T$ is the coherent Suslin tree we are looking for.\\

\noindent
\underline{Claim 2}: $T$ does not have uncountable chains.
\begin{claimproof}[Proof of claim] Let $S\in[T]^{\omega_1}$. Without loss of generality we can suppose that each element of $S$ is of the form $f_\alpha|_\xi$ for some $\xi<\alpha\in \omega_1$. Consider $$\mathcal{C}=\{\,C\in [\omega_1]^2\,:\,f_{C(1)}|_{C(0)}\in S\}.$$
By refining $\mathcal{C}$, we may assume that its elements are pairwise disjoint. Since $\mathcal{F}$ is $\mathcal{P}$-$3$-capturing, there are distinct $C_0,C_1,C_2\in \mathcal{C}$ so that $\{C_0,C_1,C_2\}$ is captured at some level $l\in P_c$. For convenience, let us denote $C_i(0)$ and $C_i(1)$ simply as $\xi_i$ and $\alpha_i$ for each $i<3$. In order to finish, just note that $\Xi_{\alpha_1}(l)=1$, $\Xi_{\alpha_2}(l)=2$ and $\rho(\xi_0,\alpha_1)=\rho(\xi_0,\alpha_2)=l$. Thus, according to the definition, we have that $f_{\alpha_1}(\xi_0)=f_{\alpha_1}|_{\xi_1}(\xi_0)=1$ and  $f_{\alpha_2}(\xi_0)=f_{\alpha_2}|_{\xi_2}(\xi_0)=0$. Hence, $S$ is not a chain.
\end{claimproof}

\noindent
\underline{Claim 3}: $T$ has no  uncountable antichains. 

\begin{claimproof}[Proof of claim]
Let $A\in [T]^{\omega_1}$. Without loss of generality we can suppose that each element of $A$ is of the form $f_\alpha$ for some $\alpha\in\omega_1$. Consider $$\mathcal{D}=\{\alpha\in\omega_1\,:\,f_\alpha\in A\}.$$ Since $\mathcal{F}$ is $\mathcal{P}$-fully capturing, there is $D\in\text{FIN}(\mathcal{D})$ which is fully-captured at some level $l\in P_a$. For convenience, let us denote $D(i)$ as $\alpha_i$ for each $i<l$. Now, let $g:m_{l-1}\backslash r_l\longrightarrow 2$ be given by:$$g(j)=\begin{cases}f_{\alpha_0}(\,(\alpha_0)_l(j)\,)&\textit{if } j<\lVert\alpha_0\rVert_l\\
0&\textit{otherwise}
\end{cases}$$
 Take $0<i<l$ for which $g=g^l_i$. We claim  that $f|_{\alpha_0}=f_{\alpha_i}|_{\alpha_0}$. For this purpose, take an arbitrary $\xi<\alpha_0$. If $\xi\in (\alpha_0)_l$ and $r_l\leq \lVert \xi\rVert_l$ then $$f_{\alpha_i}(\xi)=g^l_i(\lVert \xi\rVert_l)=g(\lVert \xi\rVert_l)=f_{\alpha_0}((\alpha_0)_l(\lVert\xi\rVert_l))=f_{\alpha_0}(\xi).$$ If $\xi\in (\alpha_0)_l$ and $\lVert \xi\rVert_l<r_l$ then $\xi\in (\alpha_0)_l\cap (\alpha_i)_l$. By means of the Proposition \ref{deltarhoequalityprop}, we know that $\rho(\alpha_0,\alpha_i)=\Delta(\alpha_0,\alpha_i)=l$. In other words, $\alpha_0$ and $\alpha_i$ are strongly $\rho$-isomorphic (see Remark \ref{remarkcapturing2}). From this, it follows that $\rho(\xi,\alpha_0)=\rho(\xi,\alpha_i)$. Thus, $f_{\alpha_i}(\xi)=f_{\alpha_0}(\xi)$ by virtue of the part (a) of Lemma \ref{lemmaxi} (since $\Delta(\alpha_0,\alpha_i)=l$ ). The last case happen when $\xi\not\in(\alpha)_l$. Here,  $f_{\alpha_0}(\xi)=f_{\alpha_i}(\xi)$ as a direct consequence of the Claim 1.
 \end{claimproof}
\end{proof}
\end{theorem}

\begin{definition}[Tree product] Given $k\in\omega$ and $\{(T_i,<_i)\}_{i< k}$ a family of trees we define their \textit{tree product} $$\bigotimes\limits_{i< k} T_i=\{t\in \prod\limits_{i<k}T_i\,:\,\forall i<k\,\big( rank(t(0))=rank(t(i))\big)\}.$$
To this set, we associate a canonical order given by: $$s<t\text{ if and only if }s(i)<_i t(i)\text{ for each }i<k$$
\end{definition}
\begin{rem}It is not hard to see that $\bigotimes\limits_{i<k}T_i$ is always a tree. Furthermore, the tree product of Aronszajn trees is always Aronszajn. Unfortunately, we can not say the same about the tree product of Suslin trees. In fact, the tree product of a Suslin tree with itself is never Suslin.\end{rem}

\begin{definition}[Full Suslin tree] Let $(T,<)$ be a tree. We say that $T$ is \textit{full Suslin} if $\bigotimes\limits_{i<k} T|_{t_i}$
is Suslin for every distinct $t_0,\dots,t_{k-1}\in T$, all of the same rank.
\end{definition}

\begin{theorem}[Under FCA]\label{fullsuslinscheme}There is a full Suslin tree.
\begin{proof}Fix a type $\langle m_k, n_{k+1},r_{k+1}\rangle_{k\in\omega}$ such that $n_{k+1}\geq 2^{m_k}$ for each $k\in\omega$, and let $\mathcal{F}$ be a fully capturing construction scheme of that type. The plan is define an order $<_T$ over $\omega_1$ which turns it into a full Suslin tree. For this purpose, we will recursively define for each $F\in \mathcal{F}$, a tree ordering $\prec_F$. We will ask that the following conditions hold for any two $F,G\in \mathcal{F}$:
\begin{enumerate}
    \item For each $\alpha,\beta\in F$, if $\alpha\prec_F\beta$ then $\alpha<\beta$ in the usual ordering over $\omega_1$.
    \item If $F\subseteq G$, then $\prec_F=\prec_G\cap (F\times F)$.
    \item If $\rho^F=\rho^G$ and $h:F\longrightarrow G$ is the increasing bijection, then $h$ is an isomorphism between $(F,\prec_F)$ and $(G,\prec_G)$.
\end{enumerate}
We proceed to define the orderings by recursion over $k=\rho^F$. \\

\noindent
\underline{Base step}: If $k=0$, there is nothing to do. This is because $|F|=1 $in this case.\\

\noindent
\underline{Recursion step}: Suppose that $0<k\in \omega$ and we have constructed the required orderings over each element of $\mathcal{F}_{k-1}$. Let $F\in \mathcal{F}_k$. For each $i<n_k$ let $h_i:F_0\longrightarrow F_i$ be the increasing bijection.

Given $I\subseteq F_0\backslash R(F)$, let us say that $I$ is \textit{$k$-independent} if for any two distinct $\alpha,\beta\in I$ and each $\xi\in F_0$, if $\xi\prec_{F_0}\alpha,\beta$ then $\xi\in R(F)$. Note that both the empty set and singletons are $k$-independent. Now, let us enumerate (possibly with repetitions) the set of all $k$-independent sets of $F_0\backslash R(F)$ as $\langle I^F_i\rangle _{0<i<n_k}$. Given $\alpha,\beta\in F$, we decide whether $\alpha\prec_F\beta$ if one of the two following cases occur:
\begin{itemize}
    \item Both $\alpha$ and $\beta$ belong to the same $F_i$ and $\alpha\prec_{F_i}\beta$.
    \item $\alpha\in F_0\backslash R(F)$, $\beta\in F_i\backslash R(F)$ for some $i>0$ and there is a (necessarily unique) $\delta\in I^F_i$ so that $\alpha\preceq_{F_0} \delta$ and $h(\delta)$ is comparable with $\beta$.
\end{itemize}
It is straightforward to check that $\prec_F$ is a tree ordering which satisfies the condition (2) of the recursive hypotheses. Furthermore, in the definition we are only adding new relations between elements of $F_0\backslash R(F)$ and $F_i\backslash R(F)$ for $i>0$, it also follows that the condition (1) of the recursive hypotheses hold. It should be also clear that condition (3) holds if we choose the enumerations $\langle I^F_i\rangle_{0<i<n_k}$ always in \say{isomorphic} positions. This finishes the construction.\\

We now define $\prec =\bigcup\limits_{F\in \mathcal{F}}\prec_F$. According to the points (1), (2) and (3) of the recursion hypotheses, it follows that $(\omega_1,<_T)$ is  in fact a tree. Moreover, $\prec\cap (F\times F)=\prec_F$ for any $F\in \mathcal{F}$. We proceed to show that it is full Suslin.\\

\noindent
\underline{Claim 1}: $(\omega_1,<_T)$ does not have uncountable chains.
\begin{claimproof}[Proof of claim] Let $A\in[\omega_1]^{\omega_1}$. As $\mathcal{F}$ is full capturing there is $l\in \omega$ and $D\in [A]^l$ which is captured at level $l$. Consider $F\in \mathcal{F}_l$ so that $D\subseteq F$. $I=\emptyset$ is $l$-independent. Hence,  there is $0<i<n_l$ so that $I=I^F_i$. By definition of $\prec_F$, we have that $D(0)$ is incompatible with $D(1)$. So the claim is over. 
\end{claimproof}

\noindent
\underline{Claim 2}: Let $1\leq k\in\omega$ and $t_0<\dots<t_{k-1}\in \omega_1$ be distinct elements of the same rank.  Then  $\bigotimes\limits_{i<k} \omega_1|_{t_i}$ has no uncountable antichains.
\begin{claimproof}[Proof of claim] Let $\mathcal{A}$ be an uncountable subset of $\bigotimes\limits_{i<k} \omega_1|_{t_i}$. We may assume without loss of generality that $t_{k-1}<s(0)<\dots<s(k-1)$ for any $s\in \mathcal{A}$. Furthermore, we may suppose that if $s$ and $t$ are distinct elements of $\mathcal{A}$,  either $s(k-1)<t(0)$ or $t(k-1)<s(0)$. Given $s\in \mathcal{A}$, let $$D_s=\{t_0,\dots,t_{k-1},s(0),\dots,s(k-1)\}.$$

As $\mathcal{F}$ is fully capturing, there $s_0,\dots, s_{l-1}\in \mathcal{A}$ so that $\mathcal{D}=\{D_{s_0},\dots,D_{s_{n_{l-1}}}\}$ is fully captured at some level  $l\in\omega$. Let $F\in \mathcal{F}_l$ for which $$\bigcup D_{s_i}\subseteq F.$$
Then $D_{s_i}\subseteq F_i$ for any $i<n_l$. Moreover, $F_i\cap R(F)=\{t_0,\dots,t_{k-1}\}$. We will show that there is $i<n_l$ so that $s_0<s_i$ with respect to the tree product ordering. For this aim, we need the following subclaim:\\

\noindent
\underline{Subclaim 1}: $I=\{s_0(0),\dots,s_1(k-1)\}$ is an $l$-independent subset of $F_0\backslash R(F)$.
\begin{proof}[Proof of subclaim] Let $i<j<k$ and $\alpha\in F_0$ be such that $\alpha\prec s_0(i),s_0(j)$. We know that $t_i\prec s_0(i)$ and $t_j\prec s_0(j)$. Furthermore, $t_i$ and $t_j$ are incomparable because they are different elements of the same rank. Hence, it must happen that $\alpha\prec  t_i,t_j$. By the point (1) of the recursive conditions, this implies that $\alpha< t_i,t_j$. As $t_i\in R(F)$, then $\alpha\in R(F)$ too. This completes the proof.
\end{proof}

By virtue of the previous subclaim, there is $i< n_l$ for which $I^F_i=I$. Let us consider $h_i:F_0\longrightarrow F_i$  the increasing bijection. Given $j<k$, we have that $s_0(j)\prec s_i(j)$ because there is $\delta\in I^F_i$, namely $s_0(j)$ so that $s_0(j)\prec \delta$ and $h(\delta)=s_i(j)$ is comparable with $s_i(j)$. In this way, we conclude that $s_0<s_j$ with respect to the tree product ordering. This finishes the proof.
\end{claimproof}
\end{proof}
    
\end{theorem}
We will finish this subsection by showing that $2$-capturing schemes do not suffice to construct Suslin trees. For this task, we will need the following Lemma.

\begin{lemma}\label{lemmacocountablecaptured}Let $\mathcal{F}$ be a $2$-capturing construction scheme. Given $X\in [\omega_1]^{\omega_1}$, the set $\{\alpha\in X\,:\, \{\beta\in X\,:\,\{\alpha,\beta\}\textit{ is captured}\}\textit{ is uncountable}\}$ is co-countable in $X.$
\begin{proof}Let $X$ be as in the hypotheses. Suppose towards a contradiction that there are uncountably many $\alpha\in X$ for which the set $C_\alpha=\{\beta\in X\,:\,\{\alpha,\beta\}\textit{ is captured }\}$. Then we can recursively construct an uncountable $Y\subseteq X$ so that if $\alpha<\beta\in Y$, then $\beta > \sup(C_\alpha\cup\{0\})$. As $\mathcal{F}$ is assumed to be $2$-capturing. There are $\alpha,\beta\in Y$ for which $\{\alpha,\beta\}$ is captured. Since $Y\subseteq X$, then $\beta\in C_\alpha$ which is a contradiction.
    
\end{proof}
    
\end{lemma}

\begin{theorem}\label{nosuslinmf}Let $\mathcal{F}$ be a $2$-capturing construction scheme. If $\mathfrak{m}_\mathcal{F}>\omega_1$, then there are no Suslin trees.
\begin{proof}
Let us assume towards a contradiction that there is a Suslin tree $T.$ It is well known that Suslin trees are never Knaster. Therefore, by Propositions \ref{propequivalenceknasterpreservescheme} and \ref{lemmaequivalencefunctioncapturingpreserving}, there are $X\in [T]^{\omega_1}$ and an injective function $\zeta:X\longrightarrow \omega_1$ in such way that for any two distinct $x,y\in X$, if $\{\zeta(x),\zeta(y)\}$ is captured, then $x$ and $y$ are incompatible. We now define $$\mathbb{P}=\{ p\in[X]^{<\omega}\,:\,p\textit{ is an antichain in T}\}$$
and order it with the reverse inclusion.\\

\noindent
\underline{Claim}: $\mathbb{P}$ is $ccc$ and $2$-preserves $\mathcal{F}$.
\begin{claimproof}[Proof of claim] We will prove the claim by appealing to the Lemma \ref{lemmaequivalencefunctioncapturingpreserving}. Let $\mathcal{A}\in [\mathbb{P}]^{\omega_1}$ and $\nu:\mathcal{A}\longrightarrow \omega_1$ be an injective function. Given $p\in \mathcal{A}$, let us define $$A_p=\zeta[p]\cup \{\nu(p)\},$$ $$\alpha_p=\max(A_p).$$  Without any loss of generality we may assume that  the elements of $\{ A_p\,:\,p\in \mathcal{A}\}$ are pairwise disjoint. Furthermore, we can suppose that there are $n,k,a\in\omega$ so that the following conditions hold for any two distinct $p,q\in \mathcal{A}$:
\begin{enumerate}[label=$(\alph*)$]
    \item $|p|=n$.
    \item $\rho^{A_p}=k.$
    \item $\lVert\alpha_p\rVert_k=a.$
    \item $|A_p|=|A_q|$. Furthermore, if $h:(\alpha_p)_k\longrightarrow (\alpha_q)_k$ is the increasing bijection, then $h[A_p]=A_q$ and $h(\nu(p))=\nu(q).$
    
    \item $\rho(\alpha_p,\alpha_q)>k.$
\end{enumerate}
We now enumerate each $p\in \mathcal{A}$ as $x^p_0,\dots,x^p_{n-1}$ in such way that $\zeta(x^p_i)<\zeta(x^p_j)$ whenever $i<j.$ By virtue of the Lemma \ref{capturedfamiliestosetslemma} and the points $(b)$, $(c)$, $(d)$ and $(e)$ above, if $p,q\in \mathcal{A}$ are such that $\{\alpha_p,\alpha_q\}$ is captured, then so are $\{\nu(p),\nu(q)\}$ and $\{\zeta(x^p_i),\zeta(x^q_i)\}$. The claim follows directly from the next subclaim.\\

\noindent
\underline{Subclaim}: There are distinct $p,q\in \mathcal{A}$ so that $p\cup q$ is an antichain and $\{\alpha_p,\alpha_q\}$ are captured.
\begin{claimproof}[Proof of subclaim] Let us suppose that the subclaim is false. Then, given $p,q\in \mathcal{A}$, if $\{\alpha_p,\alpha_q\}$ is captured, then $p\cup q$ is not an antichain. In this way, there are $i_{p,q},j_{p,q}<n$ for which $x^p_{i_{p,q}}$ and $x^q_{j_{p,q}}$ are different but comparable in $T$. By means of Lemma \ref{lemmacocountablecaptured}, we can construct three sequences $\langle p_s\rangle_{s<n^2+1}\subseteq \mathcal{A}$, $\langle(i_s,j_s)\rangle_{s<n^2+1}\subseteq n\times n$ and $\langle M_s\rangle_{s<n^2+1}\subseteq [\mathcal{A}]^{\omega_1}$ so that the following properties are satisfied for each $s<n^2+1$:
\begin{itemize}
    \item For all $q\in M_s$,  $\{\alpha_{p_s},\alpha_q\}$ is captured and $(i_s,j_s)=(i_{p_s,q},j_{p_s,q}).$
    \item $p_{s+1}\in M_s$.
    \item $M_{s+1}\subseteq M_s$.
\end{itemize}
By the pigenhole principle, there are $s<r<n^2+1$ for which $(i_s,j_s)=(i_r,j_r)=(i,j)$ for some $i,j<n.$ Since $M_r$ is uncountable but ${x^{p_s}_{i}}_\downarrow\cup {x^{p_r}_{i}}_\downarrow$ is countable, there is $q\in M_r$ for which $x^{p_s}_i,x^{p_r}_i<x^q_j$. As $T$ is a tree, then $x^{p_s}_i$ and $x^{p_r}_i$ are comparable. On the other hand, $p_r\in M_s$. Thus, $\{\alpha_{p_s},\alpha_{p_r}\}$ is captured. So by the observations prior to the subclaim, this means that $\{\zeta(x^{p_s}_i),\zeta(x^{p_r}_i)\}$ is captured too. Therefore, $x^{p_s}_i$ and $x^{q_s}_i$ are incompatible. This is a contradiction, so the proof is over.
\end{claimproof}  
\end{claimproof}
 $\mathbb{P}$ is an uncountable $ccc$ forcing which $2$-preserves $\mathcal{F}$.  As $\mathfrak{m}_\mathcal{F}>\omega_1$, such forcing contains an uncountable filter, namely $G$. Note that $\bigcup G$ is an uncountable chain in $T$. This contradicts the fact that $T$ is Suslin. 
\end{proof}
    
\end{theorem}

    \subsection{Suslin lattices}Recall that a pie in a partial order $X$ is a set of pairwise incomparable elements.
\begin{definition}[Suslin lattice]Let $(X,<,\wedge)$ be a lower semi-lattice. We say that $(X,<,\wedge)$ is  \textit{Suslin}  if \begin{itemize}
    \item $(X,<)$ is well founded,
    
    \item $X$ is uncountable,
    \item $X$ does not contain any uncountable chain nor an uncountable pie.
\end{itemize}
\end{definition}

 Suslin lower semi-lattices were first studied by Stephen J. Dilworth, Edward Odell and B\"unyamin Sari in \cite{dilworth2007lattice}, in the context of Banach spaces. In \cite{raghavan2014suslin}, Dilip Raghavan and Teruyuki Yorioka proved that, assuming the $\Diamond$-principle, there is a Suslin lower semi-lattice $\mathbb{S}$ which is a substructure of $(\mathscr{P}(\omega),\subseteq,\cap)$ and such that $\mathbb{S}^n$ does not contain any uncountable pie for each $n\in\omega$. A partial order which satisfies this last property is said to be a \textit{powerful pie}.\\
In the following theorem, we show that $CA_2$ is all that is needed in order to construct a Suslin lower semi-lattice with the above mentioned properties.

\begin{theorem}[Under $CA_2$]\label{suslinlatticescheme} There is $\mathbb{S}\subseteq \mathscr{P}(\omega)$ such that $(\mathbb{S},\subseteq,\cap)$ is a Suslin lower semi-lattice which is powerful pie.
\begin{proof}Let $\mathcal{F}$ be a $2$-capturing morass. For each $k\in\omega$, let $A_k=m_k\times 2^k$. We the $U_k=\{k\}\times m_k\times 2^{k-1}$. Since the last expression has no sense when $k=0$, we let $U_0=\{(0,0,0)\}$. Also, let $\phi_k:A_k\longrightarrow A_{k+1}$ be given as:
$$\phi_k(a,b)=\begin{cases}(a,b)&\textit{if }(a,b)\in r_{k+1}\times 2^k\\
(a+(m_k-r_{k+1}),b)&\textit{in other case}
\end{cases}$$
Notice that $A_k\cup \phi_k[A_k]=m_{k+1}\times 2^k.$ As the final part of the preparation, let $$N_k=\bigcup\limits_{i\leq k}U_k.$$
Our first objective is to construct, for each $k\in\omega$, a family $\langle S^k_x\rangle _{x\in A_k}\subseteq \mathscr{P}(N_k)$ in such way that the following conditions hold for $\mathbb{S}^k=\{\emptyset\}\cup \{S^k_x\,:\,x\in A_k\}$:
\begin{enumerate}[label=(\alph*)]
    \item $(\mathbb{S}^k,\subseteq, \cap)$ is a lower semi-lattice.
    \item $\mathbb{S}^k_0=\{\emptyset\}$ and $\mathbb{S}^k_{i+1}=\langle S^k_x\rangle_{x\in\{i\}\times {2^k}}$ for all $i<m_k$ (where $\mathbb{S}^k_i$ is the set of all elements of $\mathbb{S}^k$ of rank $i$).
    \item For all $x\in A_k$, $S^{k+1}_x=S^k_x$ and $S^{k+1}_{\phi_{k}(x)}\cap N_k=S^k_x.$ In particular, $S^k_x\subseteq S^{k+1}_{\phi_k(x)}.$
    \item The function $\psi_k:\mathbb{S}^k\longrightarrow\mathbb{S}^{k+1}$ given as:$$\psi_k(x)=\begin{cases}\emptyset &\textit{if } x=\emptyset\\
    S^{k+1}_{\phi_k(y)}&\textit{if }x=S^k_y
    \end{cases}$$
    is an lower semi-lattice embedding for each $k\in \omega.$
\end{enumerate}
The construction is carried by recursion over $k$.\\\\
\underline{Base step}: If $k=0$, then $A_k=\{(0,0)\}$. In this case, we let $S^0_{(0,0)}=\{(0,0,0)\}$. Trivially, all the conditions are satisfied.\\\\
\underline{Recursive step}: Suppose that we have defined $\mathbb{S}^k$ for some $k\in\omega$ in such way that  the conditions (a), (b), (c) and (d) are satisfied. In order to define $\mathbb{S}^{k+1}$, we first divide $A_{k+1}$ into three quadrants:
\begin{itemize}
    \item $C_0=A_k$,
    \item $C_1=[m_k,m_{k+1})\times 2^k$,
    \item $C_2=m_{k+1}\times [2^k,2^{k+1}),$
\end{itemize}   
Now, take an arbitrary $x\in A_{k+1}$ and consider the following cases:\begin{enumerate}[label=$(\roman*)$]
    \item If $x\in C_0$, let $S^{k+1}_x=S^k_x$. 
    \item If $x\in C_2$, then $x=(a,b)$ for some $a<m_{k+1}$ and $2^k\leq b<2^{k+1}$. In this case, let $S^{k+1}_x=\{k+1\}\times(a+1)\times\{b-2^k\}$.
    \item If $x\in C_1$, let $z=\phi^{-1}_k(x)$ and consider $D_x=\{b<2^k\,:\,S^k_{(r_{k+1},b)}\subseteq S^k_z\}.$
    Observe that $D_x$ codes the elements of $\mathbb{S}^k_{r_k+1}$ which are below $S^k_z.$ In this case, let $$S^{k+1}_x=S^k_z\cup\big(\bigcup\limits_{b\in D_x}\{k+1\}\times m_k\times \{b\}\big).$$
\end{enumerate}
It follows directly that condition (c) is satisfied. In particular, $\mathbb{S}^k\subseteq \mathbb{S}^{k+1}$. We will prove  that $\mathbb{S}^{k+1}$ is a lower semi-lattice by showing that it is closed under intersections. For this purpose, let $x,y\in A_{k+1}$ and consider the following cases:\\

\noindent
\underline{Case 1}: If $x,y\in C_0$.
\begin{claimproof}[Proof of case]In this case,  $S_x^{k+1},S_y^{k+1}\in \mathbb{S}^{k}$, Thus, $S^{k+1}_x\cap S^{k+1}_y\in \mathbb{S}^{k}\subseteq \mathbb{S}^{k+1}$ due to the recursive hypotheses.
\end{claimproof}
\noindent
\underline{Case 2}: If $x,y\in C_1$.
\begin{claimproof}[Proof of case]Let $x'=\phi^{-1}_k(x)$ and $y'=\phi^{-1}_k(y)$. If $D_x\cap D_y=\emptyset$, thdn $S^{k+1}_x\cap S^{k+1}_y=S^k_{x'}\cap S^k_{y'}\in \mathbb{S}^k$. On the other hand, if the intersection of $D_x\cap D_y\not=\emptyset$, then $S^k_{x'}\cap S^k_{y'}=S^k_w$ for some $w\in[r_{k+1},m_k)\times 2^k$. As $\mathbb{S}^k$ satisfies the conditions (a) and (b), we conclude that  $D_x\cap D_y=D_{\phi_k(w)}$. Thus,   $S^{k+1}_x\cap S^{k+1}_y=S^{k+1}_{\phi_k(w)}.$
\end{claimproof}
\noindent
    \underline{Case 3}: If $x,y\in C_2$.
    
    \begin{claimproof}[Proof of case]We can suppose without loss of generality that the first coordinate of $x$ is smaller or equal than the first coordinate of $y.$ In this case, $S^{k+1}_x\cap S^{k+1}_y=S^{k+1}_x$ if $x$ and $y$ share the second coordinate, or $S^{k+1}_x\cap S^{k+1}_y=\emptyset$ if their second coordinates are distinct.
\end{claimproof}
\noindent
    \underline{Case 4}: If $x\in C_0$ and $y\in C_1$.
    
    \begin{claimproof}[Proof of case]By definition we have that $S^{k+1}_x\cap S^{k+1}_y=S^k_x\cap S^k_{\phi^{-1}(y)}$.
    \end{claimproof}
    \noindent
    \underline{Case 5}: If $x\in C_0$ and $y\in C_2$.
    \begin{claimproof}[Proof of case] In here, $S^{k+1}_x\cap S^{k+1}_y=\emptyset.$
    \end{claimproof}
    \noindent
    \underline{Case 6}: If $x\in C_1$ and $y\in C_2$.
    
    \begin{claimproof}[Proof of case]
    let $a<m_{k+1}$ and $2^k\leq b<2^{k+1}$ such that $y=(a,b)$. If $b-2^k\in D_x$, we have  $S^{k+1}_x\cap S^{k+1}_y=S^{k+1}_{(c,b)}$ where $c=\min(a,r_{k+1})$. On the other hand, if $b-2^k\notin D_x$ then $S^{k+1}_x\cap S^{k+1}_y=\emptyset$.
    \end{claimproof}

In this way, we finish the proof of the satisfaction of condition (a). By carefully looking at the equalities in cases 1, 2 and 3, we also conclude that condition (d) is true for $\mathbb{S}^{k+1}$.\\ 
We now proceed to check condition (b). For this,  first observe that for each $(a,b)\in C_0\cup C_2$ it is trivially true that $rank(S^{k+1}_{(a,b)})=a+1$. To prove the same holds for all $x=(a,b)\in C_1$, we use induction over the first coordinate. The base case is when $a=m_k$. Here, we have that $S^{k+1}_{x}=S^{k+1}_{(r_{k+1},b)}\cup S^{k+1}_{(m_k-1,b+2^k)}$. Furthermore, any element 
of $\mathbb{S}^{k+1}$ which is contained in $S^{k+1}_{x}$, is also contained in either $S^{k+1}_{(r_k,b)}$ or $ S^{k+1}_{(m_k-1,b+2^k)}$. This means $$rank(S^{k+1}_{(m_k,b)})=\max( r_{k+1}+1, m_k+1)=m_k+1.$$ 
For the inductive step, suppose we that have proved  what we want for each $(c,d)\in C_1$ with $c<a$. Let $L$ be the set of all $(a-1,d)$ for which $S^{k+1}_{(a-1,d)}$ is (properly) contained in  $S^{k+1}_x$. By the inductive hypotheses, $rank(S^{k+1}_y)=a$ for each $y \in L$. The key now is to notice that each element of $\mathbb{S}^{k+1}$ which is  properly contained in $S^{k+1}_x$ is also contained in $S^{k+1}_y$ for some $y\in L$.  $rank(S^{k+1}_x)=a+1$ as a direct consequence of this fact. \\\\

Given $k\in\omega$, let us define $f_k:\omega_1\times \omega \longrightarrow \omega^2$  by the formula $f_k(\alpha,b)=(\lVert \alpha\rVert_k,b)$. For each $(\alpha,b)\in \omega_1\times \omega$, let $$S_{(\alpha,b)}=\bigcup\limits_{k> b}S^k_{f_k(\alpha,b)}.$$
Finally, let $\mathbb{S}=\{\emptyset\}\cup\{S_x\,:\,x\in \omega_1\times\omega\}$. We will prove that $\mathbb{S}$ is a Suslin lower semi-lattice. As before, it suffices to show that it is closed under intersections. Indeed, take $x=(\beta,b),y=(\delta,d)\in \omega_1\times \omega$. According to the conditions (c) and (d),  $S_x\cap N_k=S^k_{f_k(x)}$ and $S_y\cap N_k=S^k_{f_k(y)}$ for each $k>\max(b,d)$. Observe that $f_k(x),f_k(y)\in A_k$ where $k=\max(\rho(\beta,\delta),b,d)$. By condition (a), we know  that $S^k_{f_k(x)}\cap S^k_{f_k(y)}\in \mathbb{S}^k$. If this intersection is empty, we can use condition (d) to conclude that $S_x\cap S_y=\emptyset$. On the other hand, if $S^k_{f_k(x)}\cap S^k_{f_k(y)}=S^k_{(a,c)}$ for $(a,b)\in A_k$, then $a<\lVert \beta \rVert_k$. Thus, there is $\alpha\in (\beta)_{k}$ for which $\lVert \alpha\rVert_k=a$. In this way,  $S^k_{f_k(x)}\cap S^k_{f_k(y)}=S^k_{f_k(\alpha,c)}$. From conditions (c) and (d) we conclude that $S_x\cap S_y=S_{(\alpha,c)}$. 
Now, the condition (b) implies $\mathbb{S}$ is well-founded and its rank function satisfies the following properties:\begin{itemize}
    \item $rank(\emptyset)=0,$
    \item $rank(\alpha,b)=\alpha+1$ if $\alpha\in \omega$,
    \item $rank(\alpha,b)=\alpha$ if $\alpha\in \omega_1\backslash \omega.$
\end{itemize}
We do not have to prove that $\mathbb{S}$ has no uncountable chains since it is a substructure of $(\mathscr{P}(N),\,\subseteq,\,\cap)$, where $N=\bigcup\limits_{i\in\omega}N_i$. Thus, the only thing left to do is to prove  $\mathbb{S}$ is powerful pie. For this,  let $n\in\omega$ and $\mathcal{A}$ be an uncountable subset of $(\omega_1\times\omega)^n$. For each $x\in \mathcal{A}$, let $(\alpha^x_i,b^x_i)$ be such that $x(i)=(\alpha^x_i,b^x_i)$ for all $i<n$. Without any loss of generality we may assume that the following conditions hold for any two distinct $x,y\in \mathcal{A}$:
\begin{itemize}
    \item $\alpha^x_0<\dots< \alpha^x_{n-1}$.        
    \item $\{\alpha^x_j\,:\,j<n\}\cap \{\alpha^y_j\,:\,j<n\}=\emptyset$
    \end{itemize}
We also suppose can suppose  there are $b_0,\dots,b_{n-1}\in \omega$ such that $b^x_i=b_i$ for each $i<n$ and $x\in \mathcal{A}$. Since $\mathcal{F}$ is $2$-capturing, there are distinct $x,y\in \mathcal{A}$ for which  $\{\langle\alpha^x_i\rangle_{i<n},\langle\alpha^y_i\rangle_{i<n}\}$ is captured at some level $l>\max(b_i\,:\,i<n)+1$. Note that $$f_{l-1}(x(i)),f_{l-1}(y(i))\in A_{l-1}$$ for each $i<n$. 
    Now, given $i<n$ we have that $\rho(\alpha^x_i,\alpha^y_i)=\Delta(\alpha^x_i,\alpha^y_i)=l$. From this, it follows that $\lVert \alpha_i^x\rVert_{l-1}=\lVert \alpha_i^x\rVert_{l}$ and $\lVert \alpha^y_i\rVert_l=\lVert \alpha^x_i\rVert_{l-1}+(m_{l-1}-r_l)$. Thus, $f_{l-1}(x(i))=f_l(x(i))$ and  $$\phi_{l-1}(f_{l-1}(x(i)))=(\lVert \alpha^x_i\rVert_{l-1}+(m_{l-1}-r_l),b_i)=(\lVert \alpha^y_i\rVert_l,b_i)=f_l(y(i)).$$ \\
By condition (c), we conclude that $S^l_{f_l(x(i))}=S^{l-1}_{f_{l-1}(x(i))}=S^l_{f_l(y(i))}\cap N_k\subseteq S^l_{f_l(y(i))}$.  By a previous argument, this means  $S_{x(i)}\subseteq S_{y(i)}$ for each $i<n$. Consequently, $x$ and $y$ testify that $\{(S_{z(0)},\dots,S_{z(n-1)})\,:\,z\in \mathcal{A}\}$ is not a pie. This finishes the proof.
\end{proof}
\end{theorem}
    
    \subsection{Entangled sets}
One of Georg Cantor's most famous theorems is that any two countable dense total orders without endpoint are isomorphic. Particularly, this result can be applied to any two countable dense subsets of the reals. Thus, it is natural to ask whether Cantor's theorem can be extended to higher infinities (inside $\mathbb{R}$). 
\begin{definition}[$\kappa$-dense sets] Let $\kappa$ be an infinite cardinal. We say that $D\subseteq \mathbb{R}$ is \textit{$\kappa$-dense} if $|(a,b)\cap D|=\kappa$ for all $a<b\in D$.
    
\end{definition}
In \cite{baumgartneraxiom}, James Baumgartner proved that it is consistent that any two $\omega_1$-dense sets of reals are isomorphic. This assertion is now known as the Baumgartner's axiom for $\omega_1$. More generally, if $\kappa$ is an infinite cardinal then we can define:\\

\noindent
{\bf Baumgartner's Axiom for $\kappa$ [$BA(\kappa)$]}: Any two $\kappa$-dense sets of reals are isomorphic.\\

\noindent
One of the ways of proving Cantor's theorem is by defining the forcing of finite approximations of isomorphisms and making us of the Rasiowa Sikorski's Lemma. Hence, it is natural to think that maybe $BA(\kappa)$ has some relation with  $MA$. However, this is not the case. In \cite{MAdoesnotImplyBA}, Uri Abraham and Saharon Shelah showed that $MA$ does not imply $BA(\omega_1)$. For this, they introduced the objects which will study in this subsection. Namely, the entangled sets. Readers interested in learning more about them  may also look at \cite{WhyYcc}, \cite{guzmanentangled}, \cite{AsperoMotaEntangled}, \cite{RemarksonChainConditionsinProducts}, \cite{PartitionProblems} and \cite{EntangledCohen}.    
\begin{definition}[Realization]Let $k\in\omega$ and $t:k\longrightarrow \{<,>\}$. Given a total order $(X,<)$ and $a,b\in [X]^k$, we say that $(a,b)$ \textit{realizes} $t$ if $$a(i)\:t(i)\:b(i)$$
for each $i<k.$ By $T(a,b)$ we denote the unique $t$ which is realized by $(a,b).$
\end{definition}  

\begin{definition}[Entangled set]Let $(X,<)$ be a partial order, $\mathcal{E}\in [X]^{\omega_1}$ and $k\in\omega.$ We say that:\begin{itemize}
    \item $\mathcal{E}$ is \textit{$k$-entangled} if for each  uncountable family $\mathcal{A}\subseteq [\mathcal{E}]^{k}$ of pairwise disjoint sets and  $t:k\longrightarrow\{<,>\}$ there are distinct $a,b\in \mathcal{A}$ for which that $T(a,b)=t.$
    \item$\mathcal{E}$ is \textit{entangled} if it is $k$-entangled for each $k\in\omega$.
\end{itemize}
\end{definition}
\begin{lemma} Let $(X,<)$ be a total order, $k\in\omega$ and $\mathcal{E}\in [X]^{\omega_1}$ injectively enumerated as $\langle r_\alpha\rangle_{\alpha\in \omega_1}$. Then $\mathcal{E}$ is $k$-entangled if and only if for every uncountable family $\mathcal{C}\subseteq [\omega_1]^{k}$ of pairwise disjoint sets and  each $t:k\longrightarrow\{<,>\}$ there are distinct $c,d\in \mathcal{C}$ for which $$r_{c(i)}\,t(i)\,r_{d(i)}$$ 
for each $i<k.$
\begin{proof}As both implications are proved in a completely similar way, we will only show the one from left to right. For this purpose, suppose that $\mathcal{E}$ is $k$-entangled.   Let $\mathcal{C}\subseteq[\omega_1]^k$ be an uncountable family of pairwise disjoint sets and $t:k\longrightarrow \{<,>\}$.  Given $c\in \mathcal{C}$, take $h_c:k\longrightarrow k$ the unique function satisfying that $i<j$ if and only if $r_{c(h(i))}<r_{c(h_c(j))}$ for every $i,j<k$. By refining $\mathcal{C}$, we can suppose without loss of generality that there is $h$ for which $h=h_c$ for all $c\in \mathcal{C}$. The key observation is that $h$ codes the increasing enumeration (with respect to order in $X$) of $\langle r_{c(i)}\rangle_{i<k}$ whenever $c\in \mathcal{C}$. That is, $r_{c(h(0))}<\dots<r_{c(h(k-1))}$. Now, let  $t'=t\circ h$. By the previous observation and since $\mathcal{E}$ is $k$-entangled, there are distinct $c,d\in \mathcal{C}$ for which $$T(\langle\, r_{c(h(i))}\,\rangle_{i<k},\langle\, r_{d(h(i))}\,\rangle_{i<k})=t'.$$ To finish, take an arbitrary $i<k$ and let $i'=h^{-1}(i)$. Then $r_{c(h(i'))}\,t'(i')\, r_{d(h(i'))}$. But $r_{c(h(i'))}=r_{c(i)}$, $t'(i')=t(i)$ and $r_{d(h(i'))}=r_{d(i)}$. In this way, $r_{c(i)}\,t(i)\, r_{d(i)}$. 
\end{proof}
\end{lemma}

\begin{definition}\label{definitionlex} Let $(X,<)$ be a linear order. 
Given functions $f,g\in X^{\omega}$, we say that $f<_{lex} g$ if $f(n)<g(n)$ where $n=\min(\,k\in\omega\,:\,f(k)\not=g(k)\,)$.
\end{definition}
\begin{rem}\label{remarklex} $(X^\omega,<_{lex})$ is a linear order which can be embedded in $\mathbb{R}$ whenever $X$ is countable.
\end{rem}
\begin{theorem}[Under FCA]\label{entangledscheme}There is an entangled set.
\begin{proof}Let $\mathcal{F}$ be a fully capturing construction scheme of type $\langle m_k,n_{k+1},r_{k+1}\rangle_{k\in\omega}$ 
satisfying that $n_{k+1}\geq 2^{m_k}+1$ for each $k\in\omega.$\\ 
Given $k\in\omega$, let us enumerate $\mathscr{P}(m_{k}\backslash r_{k+1})$, possibly with repetitions, as $\langle C^k_i\rangle_{0<i<n_{k+1}}$ Now, for each $\alpha\in \omega_1$ let $f_\alpha:\omega\longrightarrow \mathbb{Z}$ be given as:
$$f_\alpha(k)=\begin{cases}0 &\textit{if }k=0\textit{ or } \Xi_\alpha(k)=-1\\
\Xi_\alpha(k)&\textit{if }k>0,\,\Xi_\alpha(k)\geq 0 \textit{ and }\lVert\alpha\rVert_{k-1}\in C^{k-1}_{\Xi_\alpha(k)}\\
-\Xi_\alpha(k)&\textit{if }k>0,\,\Xi_\alpha(k)\geq 0\textit{ and }\rVert\alpha\lVert_{k-1}\,\notin C^{k-1}_{\Xi_\alpha(k)}\\
\end{cases}$$

Let $\mathcal{E}=\langle f_\alpha\rangle_{\alpha\in\omega_1}$. We claim that $\mathcal{E}$ is an entangled set in $(\mathbb{Z}^\omega,<_{lex})$. Indeed, let $k\in\omega$,  $t:k\longrightarrow \{>,<\}$  and $\mathcal{C}\subseteq [\omega_1]^k$ be an uncountable 
family of pairwise disjoint sets. As $\mathcal{F}$ is fully capturing, there is  
$\{c_1,\dots,c_{n_{l}-1}\}\in[\mathcal{C}]^{n_l}$ which is captured at some level $l>0$. We claim that there is $i<n_l$ so that $T(c_0,c_i)=t$. First of all, note that  $\Delta(c_0(j),c_s(j))$ for each $s<n_l$ and 
$j<k$. This means that $f_{c_0(j)}|_l=f_{c_s(j)}|_l$. Furthermore, since the elements of $\mathcal{C}$ are pairwise disjoint, then $\Xi_l(c_s(j))=\Xi_l(c_s)=s$. This is due to Corollary \ref{corollarydisjointrhoisomorphic}.  From the two previous facts, it is easy to see that, for any $0<s<n_l$ and $j\in\omega$,  the order between $f_{c_0(j)}$ and $f_{c_s(j)}$ is decided by the value of both functions at $l$. Now, let us consider $0<i<n_{l}$ so that $$C^{l-1}_{i}=\{\,\lVert c_0(j)\rVert _{l-1}\,:
\,j<k\textit{ and }t(j)=\,<\,\}.$$ According to the second and third cases of the definitions of the functions, we  have that $0<i=f_{c_i(j)}(l)$ whenever $t(j)=\,<$ and $0>-1=f_{c_i(j)}(l)$ whenever $t(j)=\,>$. This implies $T(c_0,c_{i})=t,$ so we are done.
\end{proof}
\end{theorem}

%% file: chapters/Ramsey.tex
\section{Colorings at the uncountable}
Ramsey theory was implicitly initiated by Frank P. Ramsey in \cite{ramseytheorempaper}. In there, he proved a lemma which is now known as \say{Ramsey's theorem}.
\begin{theorem}For all $2\leq n,k\in \omega$ an each $c:[\omega]^n\longrightarrow k$, there is $A\in [\omega]^{\omega}$ so that $c|_{[A]^n}$ is constant.
\end{theorem}
Informally speaking, Ramsey theory studies relation between order and chaos in the realm of functions. In this context, it is usual to interpret the image of a function $c$ as a set of colors, and the function $c$ itself as a \textit{coloring} of the domain. Classical Ramsey theory is concerned in extending and generalizing the pigeonhole principle into different contexts. This is sometimes interpreted as finding order in chaos. For example, if $c:[\omega]^n\longrightarrow k$ is an arbitrary coloring, the pigeonhole principle would imply that there is an infinite set $X\subseteq [\omega]^n$ so that $c|_X$ is constant. Ramsey's theorem improves this by giving a reduced class of sets which testify the previous conclusion. Namely, the family $\{ [A]^n\,:\,A\in [\omega]^\omega\}.$

Classical Ramsey-type theorems over countable domains tend to have finite counterparts. However, this is not the case for the uncountable. For example, fix an injection $f:\omega_1\longrightarrow \mathbb{R}$ and let  $c:[\omega]^2\longrightarrow 2$ be given by:
$$c(a)=\begin{cases}0&\textit{ if }f(a(0))<f((a(1))\\
1&\textit{ otherwise}
    
\end{cases}$$
Then there is no $A\in [\omega_1]^{\omega_1}$ so that $c|_{[A]^2}$ is constant. Thus, the natural generalization of Ramsey's theorem to $\omega_1$ fails even when $n=k=2$. However, in this case there still exists a useful weakening of theorem. Before stating it, we recall some notation.

\begin{definition}Let $\kappa$ be a (possibly finite) cardinal, $X,Y$ be two arbitrary sets and $c:[X]^\kappa\longrightarrow Y$. Given $A\subseteq X$ and $y\in Y$, we say that $A$ is $y$-monochromatic if $c|_{[A]^\kappa}$ is constant with value $y.$
    
\end{definition}
\begin{definition}Let $\alpha,\beta$ and $\gamma$ be ordinals. The partition relation $\gamma\rightarrow (\alpha,\beta)^2_2$ stands for the following statement:\begin{center}
 For all $c:[\gamma]^2\longrightarrow 2$, there is a $0$-monochromatic subset of $\gamma$ of order type $\alpha$ or there is a $1$-monochromatic subset of $\gamma$ of order type $\beta.$ 
\end{center}
Its negation is written as $\gamma\not\rightarrow (\alpha,\beta)^2_2$. 
\end{definition}
The following  theorem is due to Paul Erd\"os and Richard Rado (see \cite{partitioncalculus}). It is a generalization of a well known theorem of Ben Dushnik, Ernest Miller and Erd\"os (see \cite{partiallyorderedsets}). The reader can find a proof in \cite{jechsettheory} and \cite{kunensettheory}.
\begin{theorem}\label{erdosdusnik} $\omega_1\rightarrow (\omega_1,\omega+1)^2_2$.
\end{theorem}
In \cite{positivepartition}, Stevo Todor\v{c}evi\'c proved that it is consistent to have $\omega_1\rightarrow (\omega_1,\alpha)^2_2$ for each $\alpha<\omega_1$. On the other hand, he showed in \cite{PartitionProblems} that $\mathfrak{b}=\omega_1$ implies that $\omega_1\not\rightarrow(\omega_1,\omega+2)^2_2$. In the following theorem, we will prove that the same is true under $CA_2$. Although a coloring testifying the negation of the previous partition relation can be extracted from some of our other constructions, we decided to include a direct proof due to its simplicity.
\begin{theorem}[Under $CA_2$] \label{notdusnikmillertheorem} $\omega_1\not\rightarrow (\omega_1,\omega+2)^2_2$. 
\begin{proof}Let $\mathcal{F}$ be a $2$-capturing construction scheme of an arbitrary type. Let $c:[\omega_1]^2\longrightarrow 2$ defined as:
$$c(\alpha,\beta)=\begin{cases}1 & \textit{ if }\Delta(\alpha,\beta)=\rho(\alpha,\beta)\\
0 &\textit{ otherwise}  
\end{cases}$$
 Since $\mathcal{F}$ is $2$-capturing, it is easy to see that there are no uncountable $0$-monochromatic subsets of $\omega_1$. Suppose towards a contradiction that there is a $1$-monochromatic set, say $X$, of order type $\omega+2$. Let $\beta$ and $\gamma$ be the last two elements of $X$ and consider $\alpha\in X\backslash (\gamma)_{\rho(\beta,\gamma)}$. Since $\rho$ is an ordinal metric and $\alpha<\beta<\gamma$ then  $\rho(\alpha,\beta)=\rho(\alpha,\gamma)$. In this way, $\Delta(\beta,\gamma)=\rho(\beta,\gamma)< \rho(\alpha,\beta)=\Delta(\alpha,\beta)$. By Lemma \ref{countrymanlemma3}, we conclude that $\Delta(\alpha,\gamma)=\Delta(\beta,\gamma)$. Thus, $\rho(\alpha,\gamma)>\Delta(\alpha,\gamma)$ which means that $c(\alpha,\gamma)=0$. This is a contradiction to $X$ being $1$-monochromatic, so we are done.
\end{proof}
\end{theorem}

In contrast to classical Rasmey theory, polychromatic Ramsey theory is about finding chaos. An example of this can be found in \cite{partitioningpairs}, where Stevo Todor\v{c}evi\'c showed  the complete failure of Ramsey's theorem at $\omega_1$ when we consider an infinite or uncountable amount of colors. Specifically, he proved the following theorem. 
\begin{theorem}\label{stevorainbowtheorem}There is a coloring $c:[\omega_1]^2\longrightarrow \omega_1$ such that $c[[A]^2]=\omega_1$ for each $A\in [\omega_1]^{\omega_1}$. 
\end{theorem}
For this, he first constructed a function $\phi:[\omega_1]^2\longrightarrow \omega_1$ with the property that $\phi[[A]^2]$ contains a closed an unbounded subset (club) of $\omega_1$ for each uncountable $A$. This function was constructed through the use of an ordinal metric. In the following proposition we show how to define such function from a construction scheme. It is worth pointing out that the same proposition holds for arbitrary ordinal metrics. However the proof is slightly different. 
\begin{proposition}Let $\mathcal{F}$ be a construction scheme and $\phi:[\omega_1]^2\longrightarrow \omega_1$ be given as:
$$\phi(\alpha,\beta)=\min(\,(\beta)_{\Delta(\alpha,\beta)}\backslash \alpha ).$$
Then $\phi[[A]^2]$ contains a club for each $A\in [\omega_1]^{\omega_1}$.
\begin{proof}Let $A$ be as in the hypotheses. In order to show that $\phi[[A]^2]$ contains a club, it is enough to show that if $M$ is a countable elementary submodel of some $H(\lambda)$ and $A\in M$, then $\delta_M=M\cap \omega_1\in \phi[[A]^2]$. For each $\alpha\in A\backslash \delta_M$ let $k_\alpha=\rho(\delta_M,\alpha)$. Since $A$ is uncountable there is $B\in [A]^{\omega_1}$ and $k\in\omega$ so that $k_\alpha=k$ for each $\alpha\in B$. By refining $B$ even more, we may assume without loss of generality that $\lVert \alpha\rVert_k=\lVert \beta \rVert_k$ for any two given $\alpha,\beta\in B$. Let us fix $\beta<\gamma \in B$ and consider $l=\Delta(\beta,\gamma)$. Note that $l>k$ because $\lVert \beta\rVert_k=\lVert \gamma\rVert_k$.  By elementarity, there is $\alpha\in M$ for which $\alpha>\max(\,(\beta)_l\cap \delta_M)$ and $\lVert \alpha\rVert_l=\lVert \gamma\rVert_l$. Particularly, this means that $\Delta(\alpha,\beta)=l$. It is straightforward that $\delta_M=\min(\,(\beta)_l\backslash \alpha)$. Thus, the proof is over.
\end{proof}
\end{proposition}
In order to prove Theorem \ref{stevorainbowtheorem}, the trick is to take a function $h:\omega_1\longrightarrow \omega_1$ so that $h^{-1}[\{\xi\}]$ is stationary for any $\xi\in\omega_1$. The function $c=h\circ \phi$ satisfies the conclusion of the theorem.

\subsection{ccc-polychromatic coloring}
In this subsection we build another coloring of the pairs of $\omega_1$ using the capturing axiom $CA_3$. Before discussing the context of this construction, we give some necesary definitions.
\begin{definition}Let $c:[\omega_1]^2\longrightarrow \omega_1$ be a coloring, $A\subseteq \omega_1$ and $\kappa$ be a (possibly finite) cardinal. We say that:\begin{itemize}
    \item $c$ is $\kappa$-bounded if $|\,c^{-1}[\{\xi\}]\,|<\kappa$ for each $\xi\in \omega_1.$
    \item $A$ is injective if $c|_{[A]^2}$ is injective.
\end{itemize}
\end{definition}
The problem whether every $2$-bounded coloring $c:[\omega_1]^2\longrightarrow \omega_1$ has an uncountable injective sets was first asked by F. Galvin in the early 1980's, who proved that $CH$ implies a negative answer to that question. In \cite{positivepartition}, Stevo Todor\v{c}evic showed that it is consistent, and in particular that $PFA$ implies that every $\omega$-bounded $c:[\omega_1]^2\longrightarrow \omega_1$ has an uncountable injective set.\\
In \cite{abrahampolychromatic}, U. Abraham, J. Cummings and C. Smyth proved that it is consistent that there is a $2$-bounded coloring $c:[\omega_1]^2\longrightarrow \omega_1$ without uncountable injective sets in any $ccc$ forcing extension. After hearing this theorem, S. Friedman asked for a concrete example of a $2$-bounded coloring without an uncountable injective set, but which adquires one in a $ccc$ forcing extension ($ccc$-destructible). Such example was produced in \cite{abrahampolychromatic} assuming $CH$ and the existence of a Suslin tree. Here, we construct one using $3$-capturing construction schemes.
\begin{theorem}[Under $CA_3$]\label{policromaticscheme}There is a coloring $c:[\omega_1]^2\longrightarrow \omega_1$ with the following properties:
\begin{enumerate}[label=$(\arabic*)$]
    \item $c$ is $2$-bounded,
    \item $c$ has no uncountable injective sets,
    \item $c$ is $ccc$-destructible.
\end{enumerate}
\begin{proof}
Let $\mathcal{F}$  be a $3$-capturing construction scheme of an arbitrary type. Let $\psi:\omega_1\times\omega\times\omega\longrightarrow \omega_1$ be a bijection. We define 
$c:[\omega_1]^2\longrightarrow \omega$ as follows:
$$c(\alpha,\beta)=\begin{cases}\psi(\beta,\rho(\alpha,\beta),\lVert \alpha\rVert_{\rho(\alpha,\beta)})&\textit{if }\alpha<\beta\textit{ and }\Xi_\beta(\rho(\alpha,\beta))>2\\
\psi(\beta,\rho(\alpha,\beta),\lVert\alpha\rVert_{\rho(\alpha,\beta)-1})&\textit{if }\alpha<\beta \textit{ and } \Xi_\beta(\rho(\alpha,\beta))\leq 2
\end{cases}$$
In the following three claims we prove that $c$ satisfies the conclusions of the theorem.\\

\noindent
\underline{Claim 1}: $c$ is $2$-bounded.
\begin{claimproof}Let $\xi\in \omega_1$ and suppose that $\{\alpha_0,\beta_0\},\{\alpha_1,\beta_1\}$ and $\{\alpha_2,\beta_2\}$ are elements of $c^{-1}[\{\xi\}]$. We will show that two of these pairs are equal. For this purpose, take $\beta\in \omega_1$ and $k,a\in\omega$ for which $\phi^{-1}(\xi)=(\beta,k,a)$. Given $i<3$ we have that $\psi^{-1}\circ c(\alpha_i,\beta_i)=(\beta,k,a)$. This means that $\beta_i=\beta$ and $\rho(\alpha_i,\beta_i)=k$. Since $\rho$ is an ordinal metric,  $$\rho(\alpha_i,\alpha_j)\leq\max(\rho(\alpha_i,\beta),\rho(\alpha_j,\beta))=k$$ for any $i,j<3$. In order to finish the proof of the claim we have to consider two cases. The first one is when $\Xi_\beta(k)> 2$. Here, $\lVert\alpha_i\rVert_k=a$ for each $i$. Particularly,  $\Delta(\alpha_0,\alpha_1)>k$. But $\rho(\alpha_0,\alpha_1)\leq k$. The only way in which this is possible is if $\alpha_0=\alpha_1$. So this case is over. The remaining case is when $\Xi_\beta(k)<2$. According to the part $(b)$ of Lemma \ref{lemmaxi},  $0\leq \Xi_{\alpha_i}(k)<\Xi_\beta(k)\leq 2$ for each $i$. Therefore, there are $i<j<3$ for which $\Xi_{\alpha_i}(k)=\Xi_{\alpha_j}(k)$. Hence, by the point $(d)$ of Lemma \ref{lemmaxi}, $\Delta(\alpha,\beta)\not=k$. As $\lVert \alpha_i\rVert_{k-1}=a=\lVert \alpha_j\rVert_{k-1}$, we conclude that $\Delta(\alpha_i,\alpha_j)>k\geq \rho(\alpha_i,\alpha_j)$. In this way, $\alpha_i=\alpha_j$.
\end{claimproof}

\noindent
\underline{Claim 2}: $c$ has no uncountable injective set.
\begin{claimproof}

Now, we prove $c$ has no uncountable injective set. For this, let $S\in[\omega_1]^{\omega_1}$. Since $\mathcal{F}$ is $3$-capturing, there is  $\{\alpha_0,\alpha_1,\alpha_2\}\in [S]^3$ which is captured at some level $l\in\omega$. In particular, $\Xi_{\alpha_2}(l)=2$ and for each $i<2$ the following properties hold:
\begin{itemize}
    \item $\rho(\alpha_i,\alpha_2)=l,$
    \item $ \lVert\alpha_i\rVert_{l-1}=\lVert\alpha_2\rVert_{l-1}.$
\end{itemize}
In other words, $c(\alpha_0,\alpha_2)=c(\alpha_1,\alpha_2)$.
\end{claimproof}
\noindent
\underline{Claim 3}: $c$ is $ccc$-destructible.
\begin{claimproof}[Proof of claim]Let $\mathbb{P}=\{p\in[\omega_1]^{<\omega}\,:\, p\textit{ is injective}\}$ ordered by reverse inclusion. We claim that $\mathbb{P}$ is $ccc$. Since $\mathcal{F}$ is $3$-capturing ( thus, $2$-capturing ), if $\mathcal{A}\in[\mathbb{P}]^{\omega_1}$ there are $p,q\in \mathcal{A}$, $l\in\omega$ and $F\in\mathcal{F}_l$ capturing $p$ and $q$. By definition of $c$, it is easy to see that $p\cup q$ is injective. Hence, $\mathcal{A}$ is not an antichain, and since $\mathcal{A}$ was arbitrary,  $\mathbb{P}$ is $ccc.$ Finally, since $\mathbb{P}$ is $ccc$ and uncountable, there is $p\in \mathbb{P}$ which forces the generic filter to be uncountable. From this it follows that if $G$ is a $\mathbb{P}$-generic filter over $V$ containing $p$, then $\bigcup G$ is an uncountable injective set. Thus, the proof is over.
\end{claimproof}
\end{proof}
\end{theorem}

%% file: chapters/Oscillation.tex
\section{Oscillation theory of $2$-capturing schemes}
In \cite{PartitionProblems}, Stevo Todor\v{c}evic developed powerful Ramsey-type results by analyzing the behaviour of unbounded families of functions of $\omega^\omega$. From these results, he deduced very interesting theorems in a great variety of areas. His analysis was based in the following concept.
\begin{definition}[Oscillation of functions]Let $f,g\in\omega^\omega$ and $k\in\omega$. We define the \textit{$k$-oscillation set of $f$ to $g$} as $$\overline{osc}_k(f,g)=\{s\in\omega\backslash k\,:\,f(s)\leq g(s)\textit{ and }f(s+1)>g(s+1)\,\},$$
Additionally, we define the \textit{oscillation number of $f$ to $g$} as $osc_k(f,g)=|\overline{osc}_k(f,g)|. $
\end{definition}
\begin{rem}Suppose that $f,g\in \omega^\omega$ are so that $f\not=^*g$. If $f$ and $g$ are not comparable with respect to $<^*$ then $osc_k(f,g)=\omega=osc_k(g,f)$. However, if $f$ and $g$ are comparable, it is not necessarily true in general that $osc(f,g)=osc(g,f)$. 
\end{rem}

In this section we will prove that $2$-capturing construction schemes can be used to define bounded families of functions whose properties reassemble the ones from unbounded families. In particular, we will be able to show that $CA_2$ implies a lot of things that are also implied by the hypothesis $\mathfrak{b}=\omega_1$. This is an interesting phenomenom as $CA_2$ is independent from the previous assumption. The reader interested in knowing about other oscillation theories may look for \cite{oscillationsintegers}.

For the rest of this section we fix $\mathcal{F}$ a $2$-capturing construction scheme of some type $\langle m_k,n_{k+1},r_{k+1}\rangle_{k\in\omega}$. We will analyse the behaviour of the oscillation number associated to the $k$-cardinality functions associated to $\mathcal{F}$. Given $\alpha\in\omega_1$ let $f_\alpha:\omega\longrightarrow \omega$ be given as:
$$f_\alpha(l)=\lVert\alpha\rVert_l.$$
We define $\mathcal{B}_\mathcal{F}$ as $\langle f_\alpha\rangle_{\alpha\in\omega_1}.$  The following lemma is a direct consequence of the definitions of $\rho$, $\Delta$ and Lemma \ref{lemmaxi}.
\begin{lemma}\label{lemmaf}Let $\alpha<\beta\in \omega_1$. Then:\begin{enumerate}[label=$(\arabic*)$]
\item $f_\alpha(i)=f_\beta(i)$ if $i<\Delta(\alpha,\beta)$,
\item$f_\alpha(j)<f_\beta(j)$ whenever  $j\geq \rho(\alpha,\beta),$
\item $f_\alpha<f_\beta$ provided that $\Delta(\alpha,\beta)=\rho(\alpha,\beta)$
\item In particular, $f_\alpha<^*f_\beta$.
\end{enumerate}
Furthermore, $\mathcal{B}_\mathcal{F}$ is bounded by the function in $\omega^\omega$ which sends each $i$ to $m_i.$
\end{lemma}
It is interesting that even though $\mathcal{B}_\mathcal{F}$ is bounded, its oscillation theory mirrors the oscillation theory of \cite{PartitionProblems} for unbounded families. Since $\mathcal{F}$ is 2-capturing, given any $\mathcal{A}\in [\mathcal{B}_\mathcal{F}]^{<\omega}$, there are $\alpha<\beta\in \mathcal{A}$ with $\Delta(\alpha,\beta)=\rho(\alpha,\beta).$ Thus, we have the following corollary.
\begin{corollary}\label{corollaryBuncountableantichains}$(\mathcal{B}_\mathcal{F}, \leq )$ has no uncountable pies\footnote{Recall that by pie, we mean a set of pairwise incomparable elements}.
\end{corollary}

Given $\alpha,\beta\in \omega_1$ and $k\in\omega $ we will write $osc_k(\alpha,\beta)$ and $\overline{osc}_k(\alpha,\beta)$ instead of $osc_k(f_\alpha,f_\beta)$ and $\overline{osc}_k(f_\alpha,f_\beta)$ respectively. These two objects will be written as $osc(\alpha,\beta)$ and $\overline{osc}(\alpha,\beta)$ whenever $k=0.$
\begin{definition}Let $a,b\in \text{FIN}(\omega_1)$ and $k\in \omega$. We define the \text{oscillation from $a$ to $b$} as
$$osc_k[a,b]=\{osc_k(\alpha,\beta)\,|\,\alpha\in a\textit{ and }\beta\in b\}.$$
\end{definition}
\begin{proposition}\label{fulloscilation}Let $n,k\in\omega$ and $\mathcal{A}\in[\omega_1]^{n}$ be an uncountable family of pairwise disjoint sets such that  $\rho^a=k$ for each $a\in \mathcal{A}$. Given $l\in\omega$, there are $a<b\in \mathcal{A}$ such that $osc_k[a,b]\subseteq [l,2l].$
\begin{proof}The proof is by induction over $l.$\\\\
\underline{Base step}: Suppose that $l=0$. Since $\mathcal{F}$ is $2$-capturing, there are $a<b\in \mathcal{A}$ so that the pair $\{a,b\}$ is captured at some level $s>k$. In particular, for each $\rho(a(i),b(i))=s=\Delta(a(i),b(i))$  and $\rho(a(i),b(j))=s$ for all distinct $i,j<k$. According to by Lemma \ref{lemmaf} and the previous observation, the following properties hold for all $i,j<k$ :
\begin{enumerate}[label=$(\arabic*)$]
    \item $f_{a(i)}|_{[k,s)}=f_{b(i)}|_{[k,s)},$
    \item $f_{a(i)}|_{\omega\backslash s}<f_{b(j)}|_{\omega\backslash s},$
    \item $f_{a(i)}|_{[k,\rho^F)}<f_{a(j)}|_{[k,\rho^F)}$ provided that $i<j$.
\end{enumerate}
From the previous facts, we conclude that $osc[a,b]=\{0\}=[l,2l].$\\\\
\underline{Inductive step}: Suppose that we have proved the proposition for some $l\in\omega$ and let $\mathcal{A}$ be as in the hypotheses. Using the inductive hypotheses, we can recursively construct an uncountable $\mathcal{C}\subseteq [\mathcal{A}]^{2}$ of pairwise disjoint sets so that for each $\{a,b\}\in \mathcal{C}$ the following condititions hold:
\begin{itemize}
    \item Either $a<b$ or $b<a$.\footnote{Recall that $a<b$ means that $\max(a)<\min(b)$.}
    \item if $a<b$ then $osc_k[a,b]\subseteq [l,2l].$ 
    \end{itemize}
    From this point on, whenever say that $\{a,b\}\in \mathcal{C}$ we will assume that $a<b$.
    Since $\mathcal{F}$ is $2$-capturing, we can find an uncountable family $\mathcal{D}\subseteq [\mathcal{C}]^2$ and $r>k+1$ with the following properties:
\begin{itemize}
\item Whenever $\{\{a,b\},\{c,d\}\}\in \mathcal{D}$, the pair $\{a\cup b, c\cup d\}$ is captured at level $r$.  In particular, this implies that $r>\rho^{a\cup b}$ and $a\cup b<c\cup d.$
\item For each $x,y\in D$, $\bigcup x\cap \bigcup y=\emptyset.$
\end{itemize}
Using once again that $\mathcal{F}$ is $2$-capturing, we can get $s>r$ and two elements of $ \mathcal{D}$, say $\{\{a_0,b_0\},\{c_0,d_0\}\},\{\{a_1,b_1\},\{c_1,d_1\}\}$, for which the pair $$\{(a_0\cup b_0)\cup (c_0\cup d_0),\, (a_1\cup b_1)\cup (c_1\cup d_1)\}$$
is captured at level $s$. We will finish by proving the following claim.\\

\noindent
\underline{Claim}: $osc_k[c_0,b_1]\subseteq [l+1,2(l+1)]$
\begin{claimproof}[Proof of claim]For this, take $i,j<n$.  The following properties follow from Lemma \ref{lemmaf}:
\begin{enumerate}[label=$(\arabic*)$]
\item $f_{b_1(j)}|_{[k,s)}=f_{b_0(j)}|_{[k,s)}$. This is because, in particular, $\{b_0,b_1\}$ is captured at level $s$.
\item $f_{c_0(i)}|_{[k,r)}=f_{a_0(i)}|_{[k,r)}.$ This is due to the fact that $\{\{a,b\},\{c,d\}\}\in \mathcal{D}.$ That is, the pair $\{a_0\cup b_0,c_0\cup d_0\}$ is captured at level $r$.
\item $f_{b_1(j)}|_{\omega\backslash s}>f_{c_0(i)}|_{\omega\backslash s}$. This is true since $\rho(b_1(j),c_0(i))=s$ and $c_0<b_1$.
\item $f_{b_0(j)}|_{[r,\rho^F)}<f_{c_0(i)}|_{[r,\rho^F)}$. Similarly to the previous point. This inequality holds because $\rho(b_0(j),c_0(i))=r$ and and $b_0<c_0.$
\end{enumerate} 
We will use these properties to calculate the oscillation. First observe that  we can use the part $(2)$ of Lemma \ref{lemmaf} to conclude $\overline{osc}_k(a_0(i),b_0(j))\subseteq [k,\rho^{a_0\cup b_0})$ and $$\overline{osc}_k(c_0(i),b_1(j))\subseteq [k,\rho^F).$$ According the properties (1) and 
(2) written above and since $r>\rho^{a_0\cup b_0}$,  $$\overline{osc}_k(c_0(i),b_1(j))\cap [k,r-1)=\overline{osc}_k(a_0(i),b_0(j)).$$ Due to properties (1), (3) and (4) we also have that $\rho^F-1\in \overline{osc}_k(c_0(i),b_1(j))$. In fact, properties (1) and (4) also imply that $\rho^F-1$ is the only element in the interval $[r,\rho^F)$ which belong to $\overline{osc}_k(c_0(i),b_1(j)).$ By joining all the previous observations, we get that: $$\overline{osc}_k(a_0(i),b_0(j))\cup\{\rho^F-1\}\subseteq \overline{osc}_k(c_0(i),b_1(j))\subseteq \overline{osc}_k(a_0(i),b_0(j))\cup\{\rho^F-1\}\cup\{r-1\}.$$
This means that $l+1\leq osc_k(c_0(i),b_1(j))\leq 2l+2$. Thus,  the proof is complete.
\end{claimproof}
\end{proof}
\end{proposition}
By a careful analysis of the argument of the preceding theorem, one can show that whenever $\mathcal{A}\in [\omega_1]^{\omega}$ and $l\in\omega$, then there are $\alpha<\beta\in \mathcal{A}$ for which $osc(\alpha,\beta)=l$. Unfortunately, this property does not hold for arbitrary uncountable families of finite sets. Nevertheless, the previous result is enough to redefine a \say{corrected} oscillation.\\\\
The following lemma is easy.
\begin{lemma}\label{lemmapartition}There is a partition $\langle P_n\rangle_{n\in\omega}$ of $\omega$ such that for every $k,n\in\omega$ there is $l\in\omega$ such that $[l,2l+k]\subseteq P_n.$
\end{lemma}
\begin{theorem}[Under $CA_2$]\label{coloringca2}There is a coloring $o:[\omega_1]^2\longrightarrow \omega $ such that for every uncountable family $\mathcal{A}\subseteq [\omega_1]^{<\omega}$ of pairwise disjoint sets and each $n\in\omega$, there are $a<b\in \mathcal{A}$ for which $\{o(\alpha,\beta)\,:
\,\alpha\in a\textit{ and }\beta\in b\}=\{n\}$.
\begin{proof}
Let $\langle P_n\rangle_{n\in\omega}$ be a partition of $\omega$ as in Lemma \ref{lemmapartition}. Let $o:[\omega_1]^2\longrightarrow \omega$ be defined as:$$o(\alpha,\beta)=n\textit{ if and only if }osc(\alpha,\beta)\in P_n.$$
We claim  that $o$ satisfies the conclusion of the theorem. Indeed, let $\mathcal{A}$ be an uncountable family of pairwise disjoint finite subsets of $\omega_1$ and $n\in\omega$. By refining $\mathcal{A}$ we may suppose that there is $k\in\omega$ such that $\rho^a=k$ for every $a\in \mathcal{A}$. Let $l\in\omega$ be such that $[l,2l+k]\subseteq P_n$. Due to Proposition \ref{fulloscilation}, there are $a<b\in \mathcal{A}$ such that $osc_k(a,b)\subseteq[l,2l].$ Given $\alpha\in a$ and $\beta\in b$, it is easy to see that $\overline{osc}(\alpha,\beta)$ has at most $k$ more elements than $\overline{osc}_k(\alpha,\beta)$. In this way, $osc(\alpha,\beta)\in [l,2l+k]\subseteq P_n$. In other words, $o(\alpha,\beta)=n$. This finishes the proof.
\end{proof}
\end{theorem}
The existence of a coloring with the properties stated above, already implies the existence of a much more powerful coloring. As we shall mention later, such a coloring can be used to build topological spaces with important properties. 
\begin{corollary}[$CA_2$]\label{coloringo*} There is a coloring $o^*:[\omega_1]^2\longrightarrow \omega$ such that for all $n\in\omega$, $h:n\times n\longrightarrow \omega$ and any uncountable family $\mathcal{A}\subseteq[\omega_1]^{n}$ of pairwise disjoint sets, there are $a<b\in \mathcal{A}$ for which $$o^*(a(i),b(j))=h(i,j)$$
for all $i,j<n.$
\begin{proof}Let $\langle h_n\rangle_{n\in\omega}$ be an enumeration of all $h:X\longrightarrow \omega$ for which $X\subseteq \omega^{<\omega}$ is finite and its elements are pairwise incomparable. Let us call $X_n$ the domain of $h_n$. Note that for each $f\in\omega^\omega$ and every $n\in\omega$ there is at most one $\sigma \in X_n$ which is extended by $f.$ Take a coloring $o$ satisfying the conclusion of Theorem \ref{coloringca2}. We define $o^*:[\omega_1]^2\longrightarrow \omega$ as follows: Given distinct $\alpha,\beta\in \omega_1$, if there are $\sigma_\alpha,\sigma_\beta\in X_{o(\alpha,\beta)}$ for which $\sigma_\alpha\subseteq f_\alpha$ and $\sigma_\beta\subseteq f_\beta$, put $$o^*(\alpha,\beta)=h_{o(\alpha,\beta)}(\sigma_\alpha,\sigma_\beta).$$
In any other case, let $o^*(\alpha,\beta)=17$.
In order to prove that $o^*$ satisfies the conclusion of the corollary, let $n\in\omega$, $h:n\times n\longrightarrow \omega$ and  $\mathcal{A}\subseteq[\omega_1]^n$ be an uncountable family of pairwise disjoint sets. By refining $\mathcal{A}$ we may suppose there is $k\in\omega$ with the following properties:\begin{enumerate}[label=$(\arabic*)$]
    \item $\forall a\in \mathcal{A}\,\forall i\not=j<n\, (f_{a(i)}|_k\not=f_{a(j)}|_k)$,
    \item $\forall a,b\in \mathcal{A}\,\forall i<n\,(f_{a(i)}|_k=f_{b(i)}|_k).$

\end{enumerate}
Fix $a_0\in \mathcal{A}$. Let $X=\{f_{a_0(i)}|_k\,:\,i<n\}$ and define $h:X\times X\longrightarrow \omega$ as:
$$h(f_{a_0(i)}|_k,f_{a_0(j)}|_k)=h(i,j).$$
We know that there is $m\in\omega$ for which $X=X_m$ and $h=h_m$. For such $m$, there are $a<b\in \mathcal{A}$ for which $o(a(i),b(j))=m$ for all $i,j<n.$ For all such $i$ and $j$,  $f_{a_0(i)}|_k\subseteq f_{a(i)}$ and $f_{a_0(j)}|_k\subseteq f_{a(j)}$. In this way, $o^*(a(i),b(j))=h(f_{a_0(i)}|_k,f_{a_0(j)}|_k)=h(i,j)$. So we are done.

\end{proof}
\end{corollary}
As an application, we get the following:
\begin{corollary}[$CA_2$]\label{cccnotproductive}ccc is not productive.
\begin{proof}Let $o$ be a coloring of Theorem \ref{coloringca2}. For each $n\in\omega$, let $$\mathbb{P}_n=\{p\in [\omega_1]^{<\omega_1}\,:\,\forall\alpha,\beta\in p\,(\textit{if }\alpha\not=\beta\textit{ then }o(\alpha,\beta)=n\,)\}.$$ In particular, $\mathbb{P}_0$ and $\mathbb{P}_1$ are $ccc$ but $\mathbb{P}_0\times\mathbb{P}_1$ is not. The set $\{(\alpha,\alpha)\,|\,\alpha\in\omega_1\}$ testifies this last fact.
\end{proof}
    
\end{corollary}

\subsection{A sixth Tukey type}

The purpose of this subsection is to show that the family $\mathbb{F}_\mathcal{F}$ is a so called sixth Tukey type. Before we begin,  recall that a partial order $(D,\leq)$ is \textit{(upwards) directed} if for every $x,y\in X$ there is $z\in D$ bigger than $x$ and $y$.
\begin{proposition} $(\mathcal{B}_\mathcal{F},\leq)$ is directed.
\begin{proof}Let $\alpha<\beta\in \omega_1$ and let $F\in \mathcal{F}_{\rho(\alpha,\beta)}$ be such that $\{\alpha,\beta\}\subseteq F$. Observe that if $\delta=\max F$ then  $f_\delta(i)=m_i$ for each $i\leq \rho(\alpha,\beta)$. Moreover, $\rho(\delta,\alpha),\rho(\delta,\beta)\leq \rho(\alpha,\beta)$ as testified by $F$.  Hence, for each $i>\rho(\alpha,\beta)$, $f_\alpha(i),f_\beta(i)\leq f_\delta(i)$ by means of the point (2) of the Lemma \ref{lemmaf}. From the previous two observations we conclude that $f_\alpha,f_\beta<f_\delta.$
\end{proof}
\end{proposition}
The following concept was introduced by John W. Tukey in \cite{Tukey}. Among other things, he proved that the sets $1$, $\omega$, $\omega_1$, $\omega\times \omega_1$ and $[\omega_1]^{<\omega}$ are non Tukey equivalent when equipped with their natural orderings.
\begin{definition}[The Tukey ordering]Let $(D,\leq_D)$ and $(E,\leq_E)$ be directed partial orders. We say that $E$ is \textit{Tukey below} $D$, and write it as $E\leq_T D$ if there is $\phi:D\longrightarrow E$ such that $\phi[X]$ is cofinal in $E$ for each cofinal $X\subseteq D.$ Furthermore, we say that 
$E$ is \textit{Tukey equivalent} to $D$, and write it  as $E\equiv_T D$, if $E\leq_T D$ and $D\leq_T E$. 
\end{definition}
\begin{rem}In view of Tukey's result We will refer to the orders $1$, $\omega$, $\omega_1$, $\omega_1\times \omega_1$ and $[\omega_1]^{<\omega}$ as the \textit{five canonical Tukey types}.
\end{rem}
  In \cite{CategorycofinaltypesII}, John R. Isbell showed that under $CH$, there is at least one directed partial order of cardinality $\omega_1$ which is non Tukey equivalent to none of the five canonical Tukey types. He later improved his result in \cite{sevencofinaltypes}.  In \cite{PartitionProblems}, Stevo Todor\v{c}evi\'c proved the existence of such a directed partial order under the hypothesis $\mathfrak{b}=\omega_1$. In \cite{directedsetscofinaltypes}, he proved that consistently every directed partial order of cardinality $\omega_1$ is Tukey equivalent to one of the five canonical Tukey types.
  
\begin{theorem}[Under $PFA$]Let $(D.\leq)$ be a directed set of cardinality $\omega_1$. Then $D$ is Tukey equivalent to $1$, $\omega$, $\omega_1$, $\omega\times \omega_1$ or $[\omega_1]^{<\omega}$.
    
\end{theorem}
 \begin{rem}From now on, we will call such an order, a \textit{sixth Tukey type}.
 \end{rem}
 The reader interested in learning more about the Tukey ordering and related topics is invited to search for \cite{surveytukeytheory}, \cite{tukeytypesofultrafilters}, \cite{tukeyorderamong}, \cite{Combinatoricsoffiltersandideals}, \cite{analyticidealscofinal}, \cite{idealcompacttukey}, \cite{CombinatorialDichotomiesandCardinalInvariants},  \cite{cofinaltopological} and \cite{avoidingfamilies}.\\
 
 The following Proposition can be found in \cite{Tukey}.
\begin{proposition} Let $(D,\leq_D)$ and $(E,\leq_E)$ be directed partial orders.\begin{itemize}
    \item $E\leq_T D$ if and only if there is $\phi:E\longrightarrow D$ such that $\phi[X]$ is unbounded in $D$ for each unbounded $X\subseteq E.$
    \item $E\equiv_T D$ if and only if there is a partially ordered set $C$ in which both $D$ and $E$ can be embedded as cofinal subsets.
\end{itemize}
\end{proposition}
\begin{definition}Let $(D,\leq)$ be a directed partial order. We say $S\subseteq D$ is $\omega$-bounded if every countable subset of $S$ is bounded in $D.$
\end{definition}
The following proposition appeard in \cite{directedsetscofinaltypes}.
\begin{proposition}\label{propotukeyequivalence}Let $(D,\leq)$ be a directed set with $|D|=\omega_1.$ Then:
\begin{enumerate}[label=$(\arabic*)$]
    \item $D\leq_T 1$ if and only if $D$ has a greatest element.
    \item $D\leq_T \omega$ if and only if $cof(D)\leq \omega.$
    \item $D\leq \omega_1$ if and only if $D$ is $\omega$-bounded.
    \item $D\leq_T \omega\times\omega_1$ if and only if $D$ can be covered by countably many $\omega$-bounded sets.
    \item $[\omega_1]^{<\omega}\leq_T D$ if and only if there is $A\in[D]^{\omega_1}$ for which every  $X\in[A]^{\omega}$  is unbounded in $D.$
\end{enumerate}
\end{proposition}
As an easy consequence of the previous result, we have the following. 
\begin{corollary}\label{corollarytukey}Let $(D,\leq)$ be a directed set with $|D|=\omega_1.$ Then either $D\equiv_T 1$, $D\equiv_T \omega$, $D\equiv_T \omega_1$ or $\omega\times \omega_1\leq_T D\leq_T [\omega_1]^{\omega}$.
\end{corollary}

Now, we are ready to prove that $\mathbb{B}_\mathcal{F}$ is a sixth Tukey type. It is worth mentioning that this proof is completely similar to the one under $\mathbb{b}=\omega_1$. 
\begin{theorem}\label{sixthtukeytypethm}$(\mathcal{B}_\mathcal{F},\leq)$ is a sixth Tukey type.
\begin{proof}By Corollary \ref{corollarytukey}, it is enough to show $\mathcal{B}\not\leq_T\omega\times\omega_1$ and $[\omega_1]^{<\omega}\not\leq_T \mathcal{B}_\mathcal{F}$. This will be a consequence of the next two claims, and due to  the points (4) and (5) of Proposition \ref{propotukeyequivalence}.\\

\noindent
\underline{Claim}: $\mathcal{B}_\mathcal{F}$ does not contain any uncountable $\omega$-bounded set. 

\begin{claimproof}[Proof of claim]For this, we argue by contradiction.  Assume there is $A\in[\omega_1]^{\omega_1}$ for which $\langle f_\beta\rangle_{\beta\in A}$ is $\omega$-bounded. Recursively,  we can build a sequence  $\langle \alpha_\xi,\beta_\xi \rangle_{\xi<\omega_1}$ of pairs of countable ordinals satisfying the following properties:
\begin{enumerate}[label=$(\arabic*)$]
    \item $\langle \alpha_\xi\rangle _{\xi\in\omega_1}$ and $\langle \beta_\xi\rangle_{\xi\in\omega_1}$ are increasing,
    \item  $\beta_\xi<\alpha_\xi$ for any $\xi\in \omega_1$, 
    \item $\langle \alpha_\xi\rangle_{\xi\in\omega_1}\subseteq A,$
    \item for each $\xi\in\omega_1$, $f_{\beta_\xi}$ is an upper bound of the set $\{ f_{\alpha_\nu} \,:\,\nu<\xi\}$.
\end{enumerate}
Since $\mathcal{F}$ is $2$-capturing, $\delta<\gamma$ so that  $\{\{\beta_\delta,\alpha_\delta\},\{\beta_\gamma,\alpha_\gamma\}\}$ is captured at some level $l\in\omega$. It follows that $f_{\beta_\gamma}(l-1)=f_{\beta_\delta}(l-1)<f_{\alpha_\delta}(l-1).$ But this is a contradiction since $f_{\beta_\gamma}$ was supposed  to bound $f_{\alpha_\delta}$.
\end{claimproof}
\noindent
\underline{Claim}: Let $A\in [\omega_1]^{\omega_1}$. Then $\langle f_\alpha\rangle_{\alpha\in A}$ contains an infinite bounded set.

\begin{claimproof}[Proof of claim]
    
Let us define the coloring $d:[A]^2\longrightarrow 2$ as:$$d(\alpha,\beta)=\begin{cases}0 & \textit{ if } f_\alpha\not\leq f_\beta\textit{ and } f_\beta\not\leq f_\alpha\\ 1 & \textit{ otherwise }
\end{cases}$$

By Theorem \ref{erdosdusnik}, there are two possibilities:\\
\begin{center}\begin{minipage}{6cm} \begin{center} \textbf{(A)}\end{center} $A$ contains a $0$-monochromatic uncountable set.\end{minipage}\hspace{2.5cm} \begin{minipage}{6cm}\begin{center} \textbf{(B)}\end{center} $A$ contains a $1$-monochromatic set of order type $\omega+1.$
\end{minipage}
\end{center}
\bigskip
Every $0$-monochromatic set is an antichain in $\mathcal{B}_\mathcal{F}$, so by Corollary \ref{corollaryBuncountableantichains} there can not be uncountable  $0$-monochromatic sets. Hence, there is  $1$-monochromatic subset of $A$, say $X$, of order type $\omega+1$. Observe that $f_\beta$ bounds $\langle f_\alpha\rangle_{\alpha\in X}$ where $\beta=\max X.$
\end{claimproof}
\end{proof}
\end{theorem}

\subsection{A Suslin tower}
\begin{definition}Let $\mathcal{T}$ be an $\kappa$-tower. We say $\mathcal{T}$ is \textit{Suslin} if for every  uncountable $\mathcal{A}\subseteq \mathcal{T}$ there are distinct $A,B\in \mathcal{A}$ with $A\subseteq B$. 
\end{definition}
Suslin towers were studied in \cite{GapsandTowers}. There, Piotr Borodulin-Nadzieja and David Chodounsk\'y proved, in particular, that Suslin $\omega_1$-towers 
exist under $\mathfrak{b}=\omega_1.$

Note that whenever $\mathcal{B}$ is an increasing family of functions in $\omega^\omega$ with respect to $<^*$ of order type $\omega_1$, then the family $\langle T_f\rangle_{f\in \mathcal{B}}$, where $$T_f=\{(n,m)\,:\,m\leq f(n)\}$$
is a Suslin tower. Furthermore,  if $f,g\in \mathcal{B}$ are such that $f\leq g$, then $T_f\subseteq T_g$. Consequently, if $(\mathcal{B},\leq)$ has no uncountable pies then $\langle T_f\rangle_{f\in \mathcal{B}}$ is Suslin. Thus, we have the following corollary.
\begin{corollary}[$CA_2$]\label{suslintowercoro}There is a Suslin $\omega_1$-tower.   
\end{corollary}

\subsection{S and L spaces}
In this subsection we will study two of the most famous problems in general topology. 
\begin{definition}[$S$-space and $L$-space] Let $(X,\tau)$ be a topological space. We say that:
\begin{itemize}
    \item $X$ is an \textit{$S$-space} if it is $T_3$, hereditarily separable and not hereditarily Lindel\"of. Additionally, we say that $X$ is a \textit{strong $S$-space} if $X^n$ is an $S$-space for each $1\leq n\in\omega.$
    \item $X$ is an \textit{$L$-space} if it is $T_3$, hereditarily Lindel\"of and not hereditarily separable. Moreover, we say that $X$ is a \textit{strong $L$-space} if $X^n$ is an $L$-space for all $1\leq n\in \omega.$
    \end{itemize}
\end{definition}
The existence of an $S$-space and an $L$-space used to be one of the main open problems in set-theoretic topology. Such spaces exist under a large variety of axioms, like $CH$ and some parametrized diamonds of \cite{ParametrizedDiamonds}. The question regarding the existence of an $S$-space  settled when Stevo Todor\v{c}evic proved the following theorem.

\begin{theorem}[Under $PFA$] There are no $S$-spaces.
\end{theorem}For some time, people thought that the existence of $S$-spaces was equivalent to the existence of $L$-spaces. That is, there is an $S$-space if and only if there is an $L$-space. However, this is not the case. By using the technic of walks on ordinals, and some number theory, Justin Moore proved that the existence of an $L$-space already follows from the usual axioms of Set Theory. This was done in \cite{solutionlspace}.
\begin{theorem}There is an $L$-space.
    
\end{theorem}

To learn more about $S$-spaces (and $L$-spaces) the reader may consult \cite{martinaxiomfirstcountable}, \cite{reformulationsandl}, \cite{basicsandl}and \cite{PartitionProblems}. 

\begin{definition}[Right and Left separated spaces] Let $(X,\tau)$ be a topological space with $|X|=\omega_1$ and $\langle x_\alpha\rangle_{\alpha\in\omega_1}$ be an enumeration of $X.$ We say that:
\begin{itemize}
    \item $X$ is \textit{left-open (right-separated)} if $\{x_\xi\,:\,\xi\leq \alpha\}$ is open for every $\alpha\in \omega_1$.
\item $X$ is \textit{right-open (left-separated)} if $\{x_\xi\,:\,\xi\geq \alpha\}$ is open for every $\alpha\in \omega_1$.

\end{itemize}
\end{definition}

The following is lemma is well-known. We prove it here for convenience of the reader.

\begin{lemma}Suppose that $(X,\tau)$ is a topological space with the following properties:
\begin{enumerate}
    \item $X$ is $T_3$,
    \item $X$ is locally countable,
    \item $X$ does not have uncountable discrete sets.
\end{enumerate}
Then $X$ is an $S$-space
\begin{proof}For each $x\in X$ let $U_x$ be a countable open neighborhood of $x$. Since $X$ is uncountable it follows that the cover $\{U_x\,:\,x\in X\}$ does not have a countable subcover. Hence, $X$ is not Lindel\"of. The only thing left to prove is that $X$ is hereditarily separable. For this aim, let $Y\subseteq X$. Suppose towards a contradiction that $Y$ is not separable. It is easy to see that $Y\backslash \overline{D}$ is uncountable for each $D\in [Y]^{\leq \omega}$. By using this fact, we can recursively construct a sequence $\langle y_\alpha\rangle_{\alpha\in \omega_1}$ so that $$y_\beta\not\in\overline{\langle y_\alpha\rangle_{\alpha<\beta}}\cup\big(\bigcup\limits_{\alpha<\beta} U_{y_\alpha}\big) $$
for each $\beta\in \omega_1$. Note that for any such $\beta$, the set $V_\beta=U_{y_\beta}\backslash \overline{\{y_\alpha\,:\,\alpha<\beta\}}$ is an open neighborhood of $y_\beta$ so that $V_\beta\cap \{y_\alpha\}_{\alpha\in \omega}=\{y_\beta\}.$ In this way, $\langle y_\alpha\rangle_{\alpha\in \omega_1}$ is an uncountable discrete set. This contradiction ends the proof.
\end{proof}   
\end{lemma}

\begin{corollary}\label{lemmaspacequivalence}Let $(X,\tau)$ be a $T_3$ topological space and $\langle x_\alpha\rangle _{\alpha\in\omega_1}$ be an enumeration of $X.$ If $X$ is left-open and does not have have uncountable descrete sets then it is is an $S$-space. Furthermore, if $X^n$ does not have uncountable discrete sets for any $1\leq n\in \omega$, then $X$ is a strong $S$-space.
\end{corollary}
\begin{theorem}\label{sspace1}Suppose that there is a coloring $o^*:[\omega_1]^2\longrightarrow \omega$ satisfying the conclusions of corollary \ref{coloringo*}. Then there is a strong $S$-space.
\begin{proof}
Let $o^*$ be a coloring as in the hypotheses. For each $\alpha\in\omega_1$ define $x_\alpha:\omega_1\longrightarrow 2$ as follows:
$$x_\alpha(\beta)=\begin{cases}\min(o^*(\alpha,\beta),1) &\textit{ if }\alpha<\beta\\
0 &\textit{ if }\alpha>\beta\\
1 &\textit{ if }\alpha=\beta
\end{cases}$$
Consider $X=\langle x_\alpha\rangle_{\alpha\in\omega_1}$ endowed with the product topology inherited by $2^{\omega_1}$. It is clear that $X$ is $T_3$. Furthermore, for each $\alpha\in \omega_1$, the set $\{x_\beta\,:\,\beta\in \omega_1\textit{ and }x_\beta(\alpha)=1\}$ is an open set contained in $\{x_\xi\,:\,\xi\leq \alpha\}$ and having $x_\alpha$ as an element. Hence, $X$ is also left-open. According to the Corollary \ref{lemmaspacequivalence}, the only  thing left to do is to show that none the finite powers of $X$ contain an uncountable discrete set. As we will construct another strong $S$-space later in this subsection, we will only prove the previous fact for $X$.\\

\noindent
\underline{Claim}: $X$ has no uncountable discrete set.
\begin{claimproof}[Proof of claim] Suppose towards a contradiction that there is $A\in [\omega_1]^{\omega_1}$ so that $\langle x_\alpha\rangle_{\alpha\in \omega_1}$ is discrete. For each $\beta\in A$ let $h_\beta\in \omega_1^{<\omega}$ so that $\beta\in dom(h_\beta)$, $x_\beta\in [h_\beta]=\{x_\alpha\,:\,\alpha\in \omega_1\textit{ and }h\subseteq x_\beta\}$ and $[h]_\beta\cap\langle x_\alpha\rangle_{\alpha\in \omega_1}= \{x_\beta\}$. Put $D_\beta=dom(x_\beta)$. By applying usual refining arguments we can assume that for each $\alpha,\beta\in A$ the following conditions hold:
\begin{itemize}
    \item $D_\alpha<D_\beta$ whenever $\alpha<\beta$,
    \item  $|D_\alpha|=|D_\beta|$. Furthermore, if $\phi:D_\beta\longrightarrow D_\alpha$ is the increasing bijection, then $\phi(\beta)=\alpha$ and $h_\alpha(\phi(\xi))=h_\beta(\xi)$ for each $\xi\in D_\beta.$
\end{itemize}
Let $m$ be the common cardinality of each $D_\beta$. By the second point above, there is $k\in m$ and $h':m\longrightarrow 2$ so that $D_\alpha(k)=\alpha$ and $h'(i)=h_\alpha(D_\alpha(i))$ for each $\alpha\in A$. We now  define $h:m\times m\longrightarrow 2$ as follows:
$$h(i,j)=\begin{cases}h'(j)&\textit{ if }i=k\\
0 &\textit{ otherwise}
    
\end{cases}$$
Since $o^*$ satisfies the conclusions of corollary \ref{coloringo*}, there are $\alpha<\beta\in A$ so that $$o^*(D_\alpha(i),D_\beta(j))=h(i,j)$$
for each $i,j\in m$. Particularly, $$h_\beta(D_\beta(j))=h'(j)=h(k,j)=o^*(D_\alpha(k),D_\beta(j))=o^*(\alpha,D_\beta(j))=x_\alpha(D_\beta(j))$$
for any given $j\in m$. Therefore, $x_\alpha\in [h_\beta]$. This is a contradiction, so the claim  is over.
    
\end{claimproof}
\end{proof}

\end{theorem}
In the same way, we can define for each $\alpha\in\omega_1$, $y_\alpha$ as follows:
$$y_\alpha(\beta)=\begin{cases} \min(o^*(\alpha,\beta),1) &\textit{ if }\alpha>\beta\\
0 &\textit{ if } \alpha<\beta\\
1 &\textit{ if }\alpha=\beta
\end{cases}$$
By a similar argument, one can show that $Y=\langle y_\alpha\rangle_{\alpha\in \omega_1}$ is a strong $L$-space.

Now, we will build another $S$-space using the family $\mathcal{B}_\mathcal{F}$. 
\begin{definition}[The topology $\tau_S$] Given $\alpha\in\omega_1$, let $C(\alpha)=\{f_\xi\,:\,f_\xi\leq f_\alpha\}$. We define $\tau_S$ to be the topology over $\mathcal{B}_\mathcal{F}$ obtained by refining the canonical Baire topology of $\omega^\omega$ restricted to $\mathcal{B}_\mathcal{F}$ by declaring the sets $C(\alpha)$ open.
\end{definition}
\begin{rem}It is straight forward that each $\alpha\in\omega_1$ has as a local base the following family:
$$\{C(\alpha)\cap [s]\,:\,s\in\omega^{<\omega}\textit{ and }f_\alpha\in[s]\}.$$
\end{rem}
Here, $[s]=\{f\in\omega^\omega\,:\,s\subseteq f\}$
The following is based on the third author's proof that $\mathfrak{b}=\omega_1$ implies that there is an $S$-space.
\begin{proposition}\label{sspace2} $(\mathcal{B}_\mathcal{F}, \tau_S)$ is an $S$-space.
\begin{proof}
This proposition will be proved by appealing to the Lemma \ref{lemmaspacequivalence}. For each $\alpha\in\omega_1$ we have that $C(\alpha)$ is closed in the Baire topology. Therefore, such set is clopen in $\tau_S.$ From this it easily follows that  $\mathcal{B}_\mathcal{F}$ is $0$-dimensional. Consequently, $\mathcal{B}_\mathcal{F}$  is also regular . Moreover, $\mathcal{B}_\mathcal{F}$ is left-open by definition. Thus, the only thing left to show is that $\mathcal{B}_\mathcal{F}$ does not contain any uncountable discrete set.

Let $S\in[\omega_1]^{\omega_1}$ and assume towards a contradiction that $\langle f_\alpha\rangle_{\alpha\in S}$ is discrete. In this way, for each $\alpha\in S$ we can find $s_\alpha\in \omega^{<\omega}$ so that, for $U_\alpha=C(\alpha)\cap[s_\alpha]$, we have that $U_\alpha\cap Y=\{\alpha\}$. Let $W\in [S]^{\omega_1}$ and $s\in\omega^{<\omega}$ for which $s_\alpha=s$ for all $\alpha\in W.$ According to the Theorem \ref{fulloscilation}, there are $\alpha<\beta\in W$ for which $osc(\alpha,\beta)=0$. In other words, $f_\alpha<f_\beta$. By definition of $C(\beta)$, $f_\alpha\in C(\beta)$. Thus, $f_\alpha\in U_\beta$ which is a contradiction.
\end{proof}
\end{proposition}
\begin{corollary}[$CA_2$]There is a  first countable $S$-space. 
\end{corollary}
Now, we present the construction of a distinct  $S$-space. For this construction, we adapt the ideas from Chapter 2 of \cite{PartitionProblems}.
\begin{definition} Given $\beta\in \omega_1$, we define $H(\beta)=\{\alpha<\beta\,:\,\rho(\alpha,\beta)=\Delta(\alpha,\beta)\,\}.$ 
\end{definition}
\begin{rem}Given $\beta\in \omega_1$, the set of all $\alpha<\beta$ so that $\{\alpha,\beta\}$ is captured, is contained in $H(\beta)$. Furthermore, if $\mathcal{F}$ is a morass then these two sets are equal. 
\end{rem}

As a consequence of the previous remark, and since $\mathcal{F}$ is $2$-capturing, we also have the following.
\begin{lemma}\label{Hbetalemma}Let $\mathcal{S}\subseteq \text{FIN}(\omega_1)$ be an uncountable family of pairwise disjoint sets. Then
there are $a,b\in \mathcal{S}$ of the same cardinality $n$ such that $a<b$ and 
$a(i)\in H(b(i))$ for all $i<n.$
\end{lemma}
\begin{definition}\label{Cbeta}For each $\beta\in\omega_1$, we recursively define $C(\beta)\subseteq \beta+1$ as the set containing $\beta$ and all $\alpha<\beta$ for which there is $\gamma\in H(\beta)$ such that:\begin{enumerate}[label=$(\alph*)$]
\item $\alpha\in C(\gamma)$, 
\item for all $\gamma\not=\xi\in H(\beta)\cup \{\beta\}$, $\Delta(\alpha,\gamma)>\Delta(\alpha,\xi).$
\end{enumerate}
Finally, we define $C_k(\beta)=\{\alpha\in C(\beta)\,:\,\Delta(\alpha,\beta)\geq k\}$ for each $k\in\omega.$
\end{definition}
\begin{rem}\label{deacreasingCl}Note that $C_l(\beta)\subseteq C_k(\beta)$ for all $k<l\in \omega.$
\end{rem}
\begin{rem}If $\gamma\in H(\beta)$, evidently $\gamma\in C(\gamma)$. Furthermore, if $\gamma\not=\xi\in H(\beta)\cup \{\beta\}$ then $\Delta(\gamma,\gamma)=\omega>\Delta(\gamma,\xi)$. From this we conclude that $H(\beta)\subseteq C(\gamma)$.
    
\end{rem}
\begin{lemma}\label{lemmaCbetaHbeta}Let $\beta\in \omega_1$ and $\gamma\in H(\beta)$.  Then $C_l(\gamma)\subseteq C(\beta)$ for $l=\Delta(\gamma,\beta)+1.$
\begin{proof}

Let $\alpha\in C_l(\gamma)$. According to the Definition \ref{Cbeta}, it suffices to show that $\Delta(\alpha,\xi)<\Delta(\alpha,\gamma)$ for each $\xi\in H(\beta)\cup \{\beta\}$ distinct from $\gamma$. If $\xi=\beta$ then  $\Delta(\gamma,\beta)<l\leq \Delta(\alpha,\gamma)$. Therefore, $\Delta(\alpha,\beta)=\Delta(\gamma,\beta)$ due to the Lemma \ref{countrymanlemma3}. In particular, $\Delta(\alpha,\beta)<\Delta(\alpha,\gamma)$. Now, if $\xi\in H(\beta)$ let us suppose towards a contradiction that $\Delta(\alpha,\gamma)\leq \Delta(\alpha,\xi)$. From the definition of $l$ and the equality obtained in the case where $\xi=\beta$, we get that $\Delta(\alpha,\beta)=\Delta(\gamma,\beta)<l\leq \Delta(\alpha,\gamma)\leq \Delta(\alpha,\xi)$. Therefore, $$\Delta(\beta,\xi)=\Delta(\alpha,\beta)=\Delta(\gamma,\beta)$$ by means of Lemma \ref{countrymanlemma3}. The definition of $H(\beta)$ yields that $$\Delta(\xi,\gamma)\leq\rho(\xi,\gamma)\leq \max(\rho(\xi,\beta),\rho(\gamma,\beta)\,)=\max(\Delta(\xi,\beta),\Delta(\gamma,\beta)\,)<\Delta(\alpha,\gamma).$$
 In this way, $\Delta(\xi,\gamma)=\Delta(\alpha,\xi)$. This implies that $\Delta(\alpha,\xi)<\Delta(\alpha,\gamma)$ which is a contradiction.
\end{proof}   
\end{lemma}

\begin{lemma}\label{Calphabaselemma}Let $\beta\in\omega_1$. For each $k\in\omega$ and $\alpha\in C_k(\beta)$ there is $l\in\omega$ such that $C_l(\alpha)\subseteq C(\beta).$
\begin{proof}The proof is carried by induction over $\beta$. So suppose that we have proved the lemma for each ordinal smaller than some $\beta\in \omega_1$. Let $\alpha\in C_k(\beta)$. According to the Definition \ref{Cbeta}, there is $\gamma\in H(\beta)$ so that $\alpha\in C(\gamma)$ and $\Delta(\alpha,\gamma)>\Delta(\alpha,\xi)$ for each $\xi\in H(\beta)\cup\{\beta\}$. Particularly, $\Delta(\alpha,\gamma)>\Delta(\alpha,\beta)$. Lemma \ref{countrymanlemma3} implies that $\Delta(\alpha,\beta)=\Delta(\gamma,\beta)$. Now, by the Lemma \ref{lemmaCbetaHbeta}, $C_{l_0}(\gamma)\subseteq C(\beta)$ where $l_0=\Delta(\gamma,\beta)+1$. Furthermore, $C_{l_1}(\alpha)\subseteq C_{\Delta(\alpha,\gamma)}(\gamma)$ for some ${l_1}\in \omega$ due to the inductive hypotheses. Note that $C_{\Delta(\alpha,\gamma)}(\gamma)\subseteq C_{l_0}(\gamma)$. In this way, $C_{l_1}(\alpha)\subseteq C(\beta)$.

In order to finish, let $l=\max(\Delta(\alpha,\beta),l_1)+1$. We claim that $C_l(\alpha)\subseteq C_k(\beta)$. For this purpose, take $\xi\in C_l(\alpha)$. Then $\xi\in C(\beta)$ because $l>l_1$. The only thing missing is to show that $\Delta(\xi,\beta)>k$. Indeed, $\Delta(\xi,\alpha)\geq l>\Delta(\alpha,\beta)$. Thus, by the Lemma \ref{countrymanlemma3}, $\Delta(\xi,\beta)=\Delta(\alpha,\beta)>k$.    
\end{proof}

\end{lemma}

\begin{rem}By the previous corollary it is easy to see that the set  $\{C_k(\beta)\,:\,k\in\omega\textit{ and }\beta\in \omega_1\}$ forms a base for a topology in $\omega_1$. It turns out that this defines a first countable locally compact 
strong $S$-space. It is convenient to transfer such topology to the family $\mathcal{B}_\mathcal{F}$.     
\end{rem}
\begin{definition}[The  topology $\tau_C$] 
Let $\beta\in\omega_1$ and $k\in\omega$. We define $\hat{C}_k(\beta)=\{f_\alpha\,|\,\alpha\in C_k(\beta)\}$ and $\hat{C}(\beta)=\{f_\alpha\,|\,\alpha\in C(\beta)\}.$ Note that $\{\hat{C}_k(\beta)\,:\,k\in\omega\textit{ and }\beta\in\omega_1\}$ forms a base for a topology over $\mathcal{B}_\mathcal{F}.$ We will call this topology  $\tau_C.$
\end{definition}
The following lemma follows directly from the fact that $\hat{C}_k(f_\beta)\subseteq [f_\beta|_k]$ for each $\beta\in\omega_1$ and $k\in\omega.$
\begin{lemma}\label{bairerefinement}Let $s\in\omega^{<\omega}$. Then $[s]\cap \mathcal{B}_\mathcal{F}$ is open in $(\mathcal{B}_\mathcal{F},\tau_C).$ In particular, $(\mathcal{B}_\mathcal{F},\tau_C)$ is Hausdorff.
\end{lemma}
\begin{lemma}\label{locallycompactlemma}Let $\beta\in\omega_1$. $\hat{C}_k(\beta)$ is compact for each $k\in\omega$.
\begin{proof}The proof is carried by induction over $\beta$. Let $\langle \alpha_n\rangle_{n\in\omega}\subseteq C_k(\beta)$. We need to show that $S=\{f_{\alpha_n}\,:\,n\in\omega\}$ has an accumulation point in $\hat{C}_k(\beta)$. If the set of all $\Delta(\alpha_n,\beta)$ is unbounded in $\omega$, clearly $\beta$ is an accumulation point of $S$. So let us assume that such collection is bounded. Then there is $A\in [\omega]^{\omega}$ and $k\leq z\in\omega$ so that $\Delta(\alpha_n,\beta)=z$ for any $n\in A$. Given any such  $n$, let $\gamma_n\in H(\beta)$ be so that $\alpha_n\in C(\gamma_n)$ and $\Delta(\alpha_n,\gamma_n)>\Delta(\alpha_n,\xi)$ for each $\xi\in H(\beta)\cup\{\beta\}$ distinct from $\gamma_n$. Particularly, we have that $\Delta(\alpha_n,\gamma_n)>\Delta(\alpha_n,\beta)$.  Thus, $\Delta(\gamma_n,\beta)=\Delta(\alpha_n,\beta)=z$ due to the Lemma \ref{countrymanlemma3}. As $\gamma_n\in H(\beta)$ then $\rho(\gamma_n,\beta)=z.$ As $(\beta)_z$ is finite, there is $B\in [A]^\omega$ and $\gamma\in H(\beta)$ so that $\gamma_n=\gamma$ for all $n\in B$. In this way, $\alpha_n\in C(\gamma)$ and $\Delta(\alpha_n,\gamma)>\Delta(\alpha_n,\beta)=\Delta(\gamma,\beta)$ for each $n\in B$. In other words, $$\langle \alpha_n\rangle_{n\in B}\subseteq C_l(\gamma)$$
where $l=\Delta(\gamma,\beta)+1$. According to the inductive hypotheses, there is $\delta\in C_l(\gamma)$ for which $f_\delta$ is an accumulation point of $\langle \alpha_n\rangle_{n\in \beta}$. By means of Lemma \ref{lemmaCbetaHbeta}, $\delta\in C(\beta)$. To finish, just note that $\Delta(\gamma,\beta)<l\leq\Delta(\gamma,\delta)$. Therefore, $\Delta(\beta,\delta)=\Delta(\gamma,\beta)=z\geq k$. Consequently, $\delta\in C_k(\beta)$.
\end{proof}
 \end{lemma}

 \begin{proposition}$(\mathcal{B}_\mathcal{F},\tau_C)$\label{sspace3} is a locally compact strong $S$-space.
 \begin{proof}By definition, $\mathcal{B}_\mathcal{F}$ is left-open and it is locally compact due to Lemma \ref{locallycompactlemma}. By Lemma \ref{bairerefinement}, $\mathcal{B}_\mathcal{F}$ is also Hausdorff.  The last two properties imply that the space is $T_3$.  Fix $n\in \omega$. It remains to prove that $\mathcal{B}_{\mathcal{F}}^{n+1}$ has no uncountable discrete subspaces. For this, let $S\subseteq \omega_1^n$ be uncountable, and assume towards a contradiction that $\langle(f_{x(0)},\dots,f_{x(n)})\rangle_{x\in S}$ is discrete. Without loss of generality we can suppose $x(i)<x(j)$ whenever $i<j\leq n$ and $\langle x(i)\rangle_{i\leq n}\cap \langle y(i)\rangle_{i\leq n}=\emptyset$ for all $x,y\in S$ with $x\not=y$. Furthermore, be a refining argument we can also suppose there is $k\in\omega$ such that for all $x,y\in S$, the following happens: \\\begin{enumerate}
     \item $\big(\prod\limits_{i\leq n}C_k(x(i))\big)\cap S=\{x\}$,
     \item $f_{x(i)}|_k=f_{y(i)}|_k$ for every $i\leq n$.\\
 \end{enumerate}
     Due to Lemma \ref{Hbetalemma}, we know there are distinct $x,y\in S$ such that $x(i)\in H(y(i))$ for every $i\leq n$. For any such $i$, we know  $\Delta(x(i),\xi)=\omega$ if and only if $\xi=x(i)$. Since $x(i)$ clearly belongs to $C(x(i))$, it follows from the definition that $x(i)\in C(y(i))$. But $f_{x(i)}|_k=f_y(i)|_k$, so $x(i)$ in fact is an element of $ C_k(x(i))$. In this way, we conclude that $x\in \big(\prod\limits_{i\leq n}C_k(y(i))\big)\cap S$, which is a contradiction. 
     \end{proof}
  \end{proposition}
  \begin{proposition}The Alexandroff compactification of $(\mathcal{B}_\mathcal{F},\tau_C)$ is a compact strong $S$-space.
   \end{proposition}
   \begin{corollary}[$CA_2$]There is a scattered compact strong $S$-space $K$ whose function space $C(K)$
is hereditarily weakly Lindel\"of and whose space $P(K)$ of all probability
measures is also a strong $S$-space.
\end{corollary}

\subsection{Failure of Baumgartner's Axiom}

As we saw on Theorem \ref{entangledscheme}, FCA implies the existence of entangled sets, which means FCA also implies the failure of $BA(\omega_1)$. Although we do not know if $CA_2$ implies the existence of an entangled set, we will prove in this subsection that it does imply the negation of $BA(\omega_1).$ Our proof is based on Stevo Todor\v{c}evi\'c proof that $\mathfrak{b}=\omega_1$ implies the failure of $BA(\omega)$ (see page 308 of \cite{homeomorphismsofmanifolds}).  Remember that at the beginning of this section we fixed $2$-capturing construction scheme, namely $\mathcal{F}$.     \begin{lemma}\label{countableZ}For any $\alpha\in \omega_1$ there are infinitely many $k\in\omega$ for which $0\leq \Xi_\alpha(k)<n_k-1.$
   \begin{proof}Let $\alpha\in \omega_1$. Assume towards a contradiction that there is $k\in\omega$ so that for any $l>k$ we have that $\Xi_\alpha(l)=n_l-1$. Since $\rho$ is an unbounded metric, there are $\beta<\gamma\in \omega_1$ so that $\alpha<\beta$ and $l=\rho(\beta,\gamma)>k$. According to the Lemma \ref{lemmaxi} we have that $\Xi_\beta(l)\leq \Xi_\gamma(l)\leq n_l-1$. Let $F\in \mathcal{F}_l$ be such that $\{\beta,\gamma\}\subseteq F$ and $a=\lVert \alpha\rVert_l$. Observe that since $\alpha<\beta$ and $\Xi_\beta(l)<n_l-1$ then $F_{n_l-1}\backslash R(F)\subseteq \omega_1\backslash (\alpha+1)$. In this way, $\alpha<F(a)$. This is because $\Xi_\alpha(l)=n_l-1$. To finish, just note that $\lVert \alpha\rVert_l=\lVert F(a)\rVert_l$. This means that $$\rho(\alpha,F(a))\geq \Delta(\alpha,F(a))\geq l>k.$$
   By the point (b) of Lemma \ref{lemmaxi} we get that $0\leq \Xi_\alpha(l')<\Xi_{F(a)}(l')$ where $l'=\rho(\alpha,F(a))$. This is a contradiction, so the proof is over. 
   \end{proof}
   \end{lemma}
   Now, let us fix $M$ a countable elementary submodel of some large enough $H(\theta)$ with $\mathcal{F},\mathcal{B}_\mathcal{F}\in M.$ Define $\mathcal{A}=\mathcal{B}_\mathcal{F}\backslash M.$ The following lemma is a direct consequence of elementarity.
   \begin{lemma}\label{lemmaelementarysubmodel}Let $s\in\omega^{<\omega}.$ If $\mathcal{A}\cap[s]\not=\emptyset,$ then $\mathcal{A}\cap[s]$ is uncountable.
   \end{lemma}
   Remember that, by Definition \ref{definitionlex} and Remark \ref{remarklex}, we can think of $(\mathcal{A},<_{lex})$ as a subset of $\mathbb{R}.$
   \begin{lemma} $(\mathcal{A},<_{lex})$ is $\omega_1$-dense.
   \begin{proof}Let $f_\alpha,f_\beta\in \mathcal{A}$ with $f_\alpha<_{lex} f_\beta$, and let $l=\Delta(\alpha,\beta)$. According to the Lemma \ref{countableZ}, there is $k>l$ so that $0\leq \Xi_k(\alpha)<n_k-1$. Consider $F\in \mathcal{F}_k$ such that $\alpha\in F$. Then $\alpha\in F\backslash F_{n_k-1}$. Thus, the unique $\gamma\in F_{n_k-1}$ for which $\lVert\gamma\rVert_{k-1}=\lVert \alpha\rVert_{k-1}$ is greater than $\alpha$. Furthermore, $\lVert \gamma\rVert_k>\rVert\alpha\rVert_k$ and $f_\gamma\in \mathcal{A}$ because $f_\alpha\in \mathcal{A}$. In this way, $S=\mathcal{A}\cap [f_\gamma|_k]$ is uncountable due to the Lemma \ref{lemmaelementarysubmodel}. In order to finish, we will show that $S$ is contained in the open interval given by $f_\alpha$ and $f_\beta$. Indeed, let $f_\xi\in S$. By definition, $f_\xi|_k=f_\alpha|_k$. This implies that $f_\xi<_{lex}f_\beta$ because $k>l$. Finally, since  $f_\alpha(k)=\lVert \alpha\rVert_k<\lVert \gamma\rVert_k=f_\gamma(k)=f_\xi(k)$ then $f_\alpha<_{lex}f_\xi$. 
   \end{proof}
   \end{lemma}
  \begin{definition}
      Let $\alpha\in\omega_1.$ We define $h_\alpha:\omega\longrightarrow \omega$ as:$$h_\alpha(i)=m_i-f_{\alpha}(i.)$$
Additionally, we let $-\mathcal{A}=\{h_\alpha(i)\,|\,\alpha\in \omega_1\backslash M\}$. 
  \end{definition}
It is easy to see that $f_\alpha<_{lex} f_\beta$ if and only if $h_\beta<_{lex}h_\alpha.$ Hence, we have the following corollary.
  
  \begin{corollary}$(-\mathcal{A},<_{lex})$ is $\omega_1$-dense.
  
 \end{corollary}
 \begin{proposition}\label{failurebaomega1scheme} There is no increasing function from $\mathcal{A}$ to $-\mathcal{A}.$
 \begin{proof}Let us assume towards a contradiction that there is an increasing function $\Psi:\mathcal{A}\longrightarrow-\mathcal{A}$ increasing. Let $\psi:\omega_1\backslash M\longrightarrow \omega_1\backslash M$ be the unique function so that $\Psi(f_\alpha)=h_{\psi(\alpha)}$.\\

 \noindent
 \underline{Claim}: $\psi$ has at most on fixed point.
 \begin{claimproof}
 
Suppose that this is not true and let $\alpha,\beta\in \omega_1$ be distinct  fixed points  of $\psi$ such that $f_\alpha<_{lex}f_\beta$. As $\Psi$ is increasing, then  \begin{align*}h_\alpha= h_{\psi(\alpha)}=\Psi(f_\alpha)<_{lex} \Psi(f_\beta)=h_{\psi(\beta)}=h_\beta.
 \end{align*}
 But this is a contradiction since, in fact, $h_\beta<_{lex}h_\alpha.$
 \end{claimproof}
 Now, let $X=\{\alpha\in\omega_1\,:\,\alpha<\psi(\alpha)\}$ and $Y=\{\alpha\in\omega_1\,:\,\alpha>\psi(\alpha)\}$. By the previous claim, one of this sets is  uncountable.
 Suppose without loss of generality that $X$ is uncountable. For each $\alpha\in X$, let $b_\alpha=\{\alpha,\phi(\alpha)\}$. Since $\mathcal{F}$ is $2$-capturing, we can find distinct $\alpha,\beta\in X$ for which $\{b_\alpha,b_\beta\}$ is captured. Observe that $f_\alpha<_{lex} f_\beta$ and $f_{\psi(\alpha)}<_{lex}f_{\psi(\beta)}$, or equivalently, $\Psi(f_\alpha)=h_{\psi(\alpha)}>_{lex}h_{\psi(\beta)}=\Psi(f_\beta)$. Note that this is a contradiction. Thus, the proof is over. 
 \end{proof}
 \end{proposition}
 \begin{corollary}[$CA_2$] There are two $\omega_1$-dense sets of reals which are not isomorphic.
 \end{corollary}
 

%% file: chapters/Onsigmamonotone.tex
\section{On $\sigma$-monotone spaces}
Monotone and $\sigma$-monotone spaces are a particular kind of metric spaces which were defined in \cite{monotonemetricspaces} by Ale\v{s} Nekvinda and Ond\v{r}ej Zindulka. In \cite{universalzindulka}, Zindulka used such spaces to prove the existence of universal measure zero sets of large Hausdorff dimension.
\begin{definition}Let $(X,d)$ be a metric space and $c>0$. We say that $X$ is \textit{$c$-monotone} if there is a linear order $<$ on $X$ such that for all $x<y<z\in X$, the following inequality holds: $$d(x,y)\leq c\cdot d(x,z).$$
We say that $X$ is \textit{monotone} if it is $c$-monotone for some $c>0$. Finally, we say that $X$ is \textit{$\sigma$-monotone} if it is a countable union of monotone subspaces.
    
\end{definition}
In \cite{hrusakzindulka}, Zindulka and Michael Hru\v{s}\'ak investigated the (possibly trivial) ideal $Mon(X)$ of all $\sigma$-monotone subspaces of a given metric space $(X,d)$. Particularly, they showed that every separable metric space of size less than the cardinal invariant $\mathfrak{m}_{\sigma-linked}$ is $\sigma$-monotone. In other words, for all such spaces $X$, the ideal $Mon(X)$ is trivial.  In Question 6.7 of that same paper, they asked whether the previous result remains true for non-separable spaces. That is: \begin{center}Is there a metric space of cardinality $\omega_1$ that is not $\sigma$-monotone?
\end{center}
In this section we will show that the capturing axiom $CA$ implies the existence of a metric space of cardinality $\omega_1$ which has no uncountable $\sigma$-monotone subspaces. Therefore, it is consistent with arbitrarily large values of the continuum that the previous question has an affirmative answer. 
\begin{rem}\label{remarkdiameter}Let $r, c>0$ and $(Y,d)$ be a $c$-monotone space. Then $(Y,r\cdot d)$ is also $c$-monotone.
    
\end{rem}

\begin{rem}\label{remarksubspacemonotone}Let $c>0$ and $(Y,d)$ be a $c$-monotone space. If $Z\subseteq Y$ then $Z$ is also $c$-monotone.
\end{rem}

\begin{rem}\label{remarkfinalmonotone}Let $c>0$ and $(Y,d)$ be a $c$-monotone space. Then $Y$ is $c'$-monotone for each $c'>c$.    
\end{rem}

By the previous remark, it is easy to see that a space $Y$ is not monotone if and only if it is not $\frac{1}{n}$ for each $0<n\in \omega$. The following result is implicit in Lemma 4.4 and Proposition 4.5 of \cite{monotonemetricspaces}.
\begin{lemma}\label{lemmamonotone}Let $n\in\omega$. There is a finite metric space $(Z_n,d)$ so that $Z_n$ is not $\frac{1}{n}$-monotone.
\end{lemma}
\begin{rem}\label{remafterlemmamonotone}Note that if $m$ is a natural number for which there is a space $Z_n$ of size $m$ which is not $\frac{1}{n}$-monotone, then the same is true for each $m'\geq m$.
\end{rem}
\begin{theorem}[Under $CA$]\label{nonmonotonescheme}There is a metric $d$ over $\omega_1$ so that $(\omega_1,d)$ has no uncountable monotone subspaces.
\begin{proof}
Let $\tau=\langle m_k,n_{k+1},r_{k+1}\rangle_{k\in\omega} $ be a type so that for any $k>0$, there is a metric $d_k$ over $n_k$ for which the following conditions holds:
\begin{center}For each $0<k'\leq k$, $n_{k'}$ is not a  $\frac{1}{k'}$-monotone subspace of $(n_k,d_k)$.\end{center}  Such type exists due to the Lemma \ref{lemmamonotone} and the Remarks \ref{remarksubspacemonotone} and \ref{remafterlemmamonotone}. Furthermore, by Remark \ref{remarkdiameter} we may assume that $diam(n_k)=1$ for each $k$. Now, let $\mathcal{F}$ be a capturing construction scheme of type $\tau$. We will proceed to define for each $F\in \mathcal{F}$ a metric $d_F$ over $F$ in such way that the following properties hold for any two $F,G\in \mathcal{F}$:
\begin{enumerate}[label=$(\alph*)$]
\item If $\rho^F=\rho^G$ and $h:F\longrightarrow G$ is the increasing bijection, then $h$ is also an isometry between $(F,d_F)$ and $(G,d_F)$.
\item If $F\subseteq G$, then $(F,d_F)$ is subspace of $(G,d_F)$.
\item If $k=\rho^F>0$, there is $0<s\in \mathbb{R}$ so that for any two $j,l<n_k$ and each $r_k\leq i<m_{k-1}$, $d_F(F_j(i),F_l(i))=s\cdot d_k(j,l)$. 
\end{enumerate}
The proof is carried by recursion over $k=\rho^F$.\\

\noindent
\underline{Base step}: If $k=1$, then $n_k=m_k=|F|$. In this case, we let $d_F:F^2\longrightarrow \mathbb{R}$ be defined as: $$d_F(F(i),F(j))=d_1(i,j)$$
for all $i,j<F$.\\

\noindent
\underline{Recursive step}: Let $k\leq 1$ and suppose that we have defined, for all $1\leq k'\leq k$,  metrics $d_G$ over each $G\in \mathcal{F}_{k'}$ in such a way that the properties $(a)$, $(b)$ and $(c)$ hold. Given $G\in \mathcal{F}_k$, let $$s^G_k=\min(d_G(\alpha,\beta)\,:\,\alpha\not=\beta\in G).$$
By means of the property $(a)$, the number $s^G_k$ does not depend on depend on $G$, so let just call $s_k$. Now, let $F\in\mathcal{F}_{k+1}$. Note that each element of $F$ is of the form $F_j(i)$ for some $j<n_{k+1}$ and $i<m_k$. Having this in mind, we define $d_F:F^2\longrightarrow F$ as:
$$d_F(F_j(i),F_l(t))=\begin{cases}d_{F_j}(F_j(i),F_j(t))&\textit{ if }i\not=t\\
\frac{s_k}{2}\cdot d_{k+1}(j,l)&\textit{ if }i=t \textit{ and }F_j(i)\not=F_l(t)\\
0&\textit{ if }F_j(i)=F_l(t)
    
\end{cases}$$
It is not hard to see that $d_F$ is well defined and is, in fact, a metric. Furthermore, from the definition, it should be clear that for any $F,G\in \mathcal{F}$ with $\rho^F,\rho^G\leq k+1$, the properties $(a)$, $(b)$ and $(c)$ are still true. This finishes the recursion.\\

Now, we define $d=\bigcup\limits_{F\in \mathcal{F}} d_F$. By the property $(b)$, it follows that $d$ is a (well defined) metric over $\omega_1$. We will now show that $(\omega_1,d)$ has no uncountable monotone subspaces. For this, let $S\in [\omega_1]^{\omega_1}$ According to the Remark \ref{remarkfinalmonotone}, it is enough to show that for each $0<k\in \omega$, $S$ is not $\frac{1}{k}$-monotone. Indeed, since $\mathcal{F}$ is capturing (and thus, it is $n_k$-capturing) there is $D\in [S]^{n_k}$ which is captured at some level $l>k$. Let $F\in \mathcal{F}_l$ be such that $D\subseteq F$. Then $D=\{ F(i)\,:\,i<n_k\}$ for some $r_l\leq i<m_{l-1}$. In this way, $D$ is isomorphic to $n_k$ seen as a subspace of $(n_l,\frac{s_l}{2}\cdot d_l)$.  This is due to the property $(c)$. But $n_k$ is not an $\frac{1}{k}$-monotone
 subspace of $n_l$. We conclude, using the Remark \ref{remarkdiameter}, that $D$ is not $\frac{1}{k}$-monotone. As $D\subseteq S$, the same holds for $S$. Thus, the proof is over.
 
 \end{proof}    
\end{theorem}

%% file: chapters/Fragments_of_Martins_axiom_and_ultrafilters.tex
\chapter{Fragments of Martin's axiom}\label{martinschapter}

In this chapter we will study the relation between construction schemes and forcing notions which preserve some of their \say{capturing properties}. Given an $n$-capturing construction scheme $\mathcal{F}$, we show that there is a filter $\mathcal{U}_n(\mathcal{F})$ over $\omega$ so that, for any partition $\mathcal{P}$, $\mathcal{F}$ is $\mathcal{P}$-$n$-capturing if and only if $\mathcal{P}\subseteq \mathcal{U}_n(\mathcal{F})^+$. Later, we will establish the consistency of the statement $\mathfrak{m}^n_\mathcal{F}>\omega_1$ and prove that under this assumption, $\mathcal{U}_n(\mathcal{F})$ is an ultrafilter. In particular, this implies that there may by construction schemes which are $n$-capturing but not $\mathcal{P}$-$n$-capturing for any non-trivial partition. Finally, we will show that $\mathcal{U}_n(\mathcal{F})$ is in fact a Ramsey ultrafilter and give an example of an ultrafilter over $\omega_1$ which can be explicitly defined from $\mathcal{F}$, and which can be projected to another Ramsey ultrafilter over $\omega$.

Before stating and proving the results that are due to the author, it will be convenient to recall some of the previous work done by  Damjan Kalajdzievski and Fulgencio Lopez regarding the preservation of capturing properties.

\begin{definition}[$n$-preserving property]Let $\mathbb{P}$ be a forcing notion and which preserves $\omega_1$ and $\mathcal{F}$ be a construction scheme. Given $2\leq n\in \omega$, we say that $\mathbb{P}$ \textit{$n$-preserves $\mathcal{F}$} if  $$\mathbb{P}\Vdash\text{\say{ $\mathcal{F}$ is $n$-capturing }}.$$
\end{definition}
\begin{rem}If $\mathbb{P}$ and $\mathcal{F}$ are as in the previous definition, then $\mathcal{F}$ is $n$-capturing.
\end{rem}
 Usually, proving that a $ccc$ forcing $\mathbb{P}$ $n$-preserves a scheme $\mathcal{F}$, is similar to proving that $\mathbb{P}$ is $ccc$. This can be exemplified in the proof of Theorem \ref{mfstronglydonuttheorem}. The following lemma gives us a useful way of handling this situation. 

\begin{lemma}\label{lemmaequivalencefunctioncapturingpreserving}Let $\mathcal{F}$ be an $n$-capturing construction scheme and $\mathbb{P}$ be a forcing. The two following statements are equivalent:
\begin{enumerate}[label=$(\arabic*)$]
    \item $\mathbb{P}$ is $ccc$ and $n$-preserves $\mathcal{F}.$
    \item For any $\mathcal{A}\in[\mathbb{P}]^{\omega_1}$ and each injective function $\nu:\mathcal{A}\longrightarrow \omega_1$, there is $\{p_0,\dots, p_{n-1}\}\in [\mathcal{A}]^n$ for which:
    \begin{itemize}
        \item $\{p_0,\dots, p_{n-1}\}$ is centered. That is, there is $p\in \mathbb{P}$ so that $p\leq p_i$ for any $i<n$.
        \item $\{\nu(p_0),\dots,\nu( p_{n-1})\}$ is captured.
    \end{itemize}
\end{enumerate}
\begin{proof}\begin{claimproof}[Proof of $\Rightarrow$]Let $\mathcal{A}$ be an uncountable subset of $\mathbb{P}$ and consider $\nu:\mathcal{A}\longrightarrow \omega_1$ an injective function. As $\mathbb{P}$ is $ccc$, there is $G$ a $\mathbb{P}$-generic filter over $V$ for which $G\cap \mathcal{A}$ is uncountable. In $V[G]$, let $X=\nu[\mathcal{A}\cap G]$. Note that $X$ is uncountable since $\nu$ is injective. As $\mathbb{P}$ $n$-preserves $\mathcal{F}$, there is $\{\alpha_0,\dots,\alpha_{n-1}\}\in [X]^n$ which is captured. Given $i<n$, let $p_i\in \mathcal{A}\cap G$ be such that $\nu(p_i)=\alpha_i$. Then $\{p_0,\dots,p_{n-1}\}$ is centered because it is included in $\mathcal{A}$, and $\{\nu(p_0),\dots,\nu(p_{n-1})\}$ is captured. 
\end{claimproof}
\begin{claimproof}[Proof of $\Leftarrow$] In order to show that $\mathbb{P}$ is $ccc$ let $\mathcal{A}\in [\mathbb{P}]^{\omega_1}$. Let us consider $\nu:\mathcal{A}\longrightarrow \omega_1$ an aribtrary inyective function. According to the hypotheses, there is $\{p_0,\dots p_{n-1}\}\in [\mathcal{A}]^n$ which is centered. Particularly, $p_0$ and $p_1$ are two distinct compatible elements of $\mathcal{A}$. Thus, $\mathcal{A}$ is not an antichain.\\
Now we will show that $\mathbb{P}$ $n$-preserves $\mathcal{F}$. For this purpose, let $G$ be a $\mathbb{P}$-generic filter over $V$ and  $X\in V[G]\cap [\omega_1]^{\omega_1}$. Let $\dot{X}$ be a name for $X$ which is forced to be an uncountable subset of $\omega_1$ by $\mathbb{1}_\mathbb{P}$. To finish, it is enough to show that the set of all $p\in \mathbb{P}$ for which there is $D\in [\omega_1]^{n}$ so that $D$ is captured and $p\Vdash\text{\say{ $D\subseteq \dot{X}$}}$, is dense in $\mathbb{P}$. For this purpose, let $q\in \mathbb{P}$. If there is $p\leq q$ for which the set $Y_p=\{\alpha\in \omega_1\,:\,p\Vdash\text{\say{$\alpha\in \dot{X}$} } \}$
is uncountable, we are done.  Therefore, we may assume that $Y_p$ is at most countable for each $p\leq q$. From this fact, it is easy to see that there is $\mathcal{A}\in[\mathbb{P}]^{\omega_1}$ and an injective function $\nu:\mathcal{A}\longrightarrow \omega_1$ so that $p\leq q$ and $p\Vdash \text{\say{ $\nu(p)\in \dot{X}$}}$ for any $p\in \mathcal{A}$. According to the hypotheses, there is $\{p_0,\dots,p_{n-1}\}\in [\mathcal{A}]$ which is centered and for which $D=\nu[\{p_0,\dots,p_{n-1}\}]$ is captured. Let $p$ be such that $p\leq p_i$ for each $i<n$. Then $p\leq q$ and $p\Vdash\text{\say{$D\subseteq \dot{X}$}}$. This finishes the proof.
\end{claimproof}
    
\end{proof}
\end{lemma}

The following lemma was proved in \cite{forcingandconstructionschemes}.
\begin{lemma}\label{nknastercapturinglemma} Let $n\in \omega$, $\mathcal{P}$ a partition of $\omega$ and suppose that $\mathcal{F}$ is a $\mathcal{P}$-$n$-capturing construction scheme. If $\mathbb{P}$ is an $n$-Knaster forcing then $\mathbb{P}$ $\mathcal{P}$-$n$-preserves $\mathcal{F}$.
\begin{proof}Let $\dot{X}$ be a name for an uncountable subset of $\omega_1$. We will prove the lemma by appealing to the equivalence of $n$-capturing provided by Lemma \ref{equivalencecapturing}. Consider an arbitrary $q\in \mathbb{P}$. Let $\langle p_\xi,\alpha_\xi\rangle_{\xi\in \omega_1}\subseteq \mathbb{P}\times \omega_1$ be a sequence so that the following properties hold:
\begin{itemize}
\item For all $\xi\in \omega_1,$ $p_\xi\leq q$ and $p_\xi\Vdash \text{\say{ $\alpha_\xi\in \dot{X}$ }}$. 
\item For any two distinct $\xi,\mu\in \omega_1$, $\alpha_\xi\not=\alpha_\mu$.
\end{itemize}
Since $\mathbb{P}$ is $n$-Knaster, there is $S\in [\omega_1]^{\omega_1}$ such that $\{ p_\xi\,:\,\xi\in S\,\}$ is $n$-linked. Now, as $\mathcal{F}$ is $\mathcal{P}$-$n$-capturing and $S$ is uncountable, for any $P\in \mathcal{P}$ there are distinct $\xi_0,\dots,\xi_n\in S$ for which $D=\{ \alpha_{\xi_0},\dots,\alpha_{\xi_{n-1}}\}$ is captured at some level $l\in P$. Let $p\in \mathbb{P}$ be such that $p\leq p_{\xi_i}$ for each $i<n$. Then $p\leq q$ and $p\Vdash\text{\say{ $D\in [\dot{X}]^n$ and is captured at $l$ }}$. This finishes the proof.
\end{proof}
\end{lemma}
In \cite{forcingandconstructionschemes}, Damjan Kalajdzievski and Fulgencio Lopez used the previous lemma to show that $MA(K_n)$ is consistent with $CA_n$. 
    \begin{theorem}Let $2\leq n\in \omega$ and $V$ be a model of $CH$. For any cardinal $\kappa>\omega_1$ of uncountable cofinality, there is an $n$-Knaster forcing $\mathbb{K}$ so that $$\mathbb{K}\Vdash\text{\say{ $\mathfrak{c}=\kappa+MA(K_n)+CA_n$ }}.$$
\begin{proof}Let $(\langle \mathbb{K}_\xi\rangle_{\xi\leq \kappa},\langle \dot{\mathbb{Q}}_\xi\rangle_{\xi<\kappa})$ be a finite support iteration of $n$-Knaster forcings of length $\kappa$. If an appropriate bookkeeping is used, we can summon both Lemma \ref{nknastercapturinglemma} and Theorem \ref{Cohenforcingfcatheorem} to show that $\mathbb{K}=\mathbb{K}_\kappa$ satisfies the desired conclusions.
\end{proof}
\end{theorem}
In that same paper, they showed that $\mathfrak{m}_{K_n}>\omega_1$ implies that there are no $n+1$-capturing construction schemes. For that purpose, they used  the following property.\\\\
{\bf The Property $(\star)_n$:} For any uncountable $\Gamma\subseteq \omega^\omega$ there is $\Gamma_0\in[\Gamma]^{\omega_1}$ such that there are no $g_0,\dots,g_n\in \Gamma_0$ and $k\in \omega$ with $g_0|_k=\dots=g_n|_k$ and $g_i(k)\not=g_j(k)$ for any $i<j\leq n$.\\\\
The property $(\star)_n$ was considered by Stevo Todor\v{c}evi\'c in \cite{remarkscellularityproducts}. The following theorem follows from Lemma 6 of such paper.
\begin{theorem}Let $2\leq n\in \omega$. Then $\mathfrak{m}_{K_n}>\omega_1$ implies $(\star)_n$.
\end{theorem}
In contrast, $n+1$-capturing construction schemes imply the failure of $(\star)_n$ (see Theorem 2.4 in \cite{forcingandconstructionschemes}).
\begin{theorem}\label{starfailscapturing}Let $1\leq n\in \omega$. If there is an $n+1$-capturing construction scheme,  then $(\star)_n$ fails.
\begin{proof}Let $\mathcal{F}$ be an $n+1$-capturing construction scheme. For any $\alpha\in \omega_1$ let $f_\alpha:\omega\longrightarrow \omega$ be defined as:
$$f_\alpha(k)=\lVert \alpha\rVert_k.$$
Now, let $\Gamma=\{\,f_\alpha\,:\,\alpha\in \omega_1\}$. Since  $\mathcal{F}$ is $n+1$-capturing, for any $S\in [\omega_1]^{\omega_1}$ there is $D\in [S]^{n+1}$ which is captured at some level $l.$ Then $g_{D(0)}|_l=\dots=g_{D(n)}|_l$ and $g_{D(i)}(l)<g_{D(j)}(l)$ for any $i<j\leq n$. In this way, $\Gamma$ testifies the failure of $(*)_n$.
\end{proof}
\end{theorem}
\section{The principle $\mathfrak{m}^n_\mathcal{F}$}
The main goal of this small section is to show that for a given $n$-capturing construction scheme $\mathcal{F}$, it is consistent that $\mathfrak{m}^n_\mathcal{F}>\omega_1$. Let us first recall the exact definition of $\mathfrak{m}^n_\mathcal{F}$.

\begin{definition}Let $\mathcal{F}$ be a construction scheme and $n\in\omega$. We define $\mathfrak{m}^n_\mathcal{F}$ as follows:
$$\mathfrak{m}^n_\mathcal{F}=\begin{cases}
\omega &\textit{if }\mathcal{F}\textit{ is not }n\textit{-capturing}\\
\min(\mathfrak{m}(\mathbb{P})\,:\,\mathbb{P}\textit{ is }ccc\textit{ and }\mathbb{P}\textit{ $n$-preserves $\mathcal{F}$})&\textit{if }\mathcal{F}\textit{ is }n\textit{-capturing}
\end{cases}$$
$\mathfrak{m}^2_\mathcal{F}$ is denoted simply as $\mathfrak{m}_\mathcal{F}.$ 
\end{definition}

\begin{rem} Suppose that $\mathcal{F}$ is an $n$-capturing construction scheme. The following diagram  represents the basic relations between the cardinals $\mathfrak{m}^i_\mathcal{F}$ and $\mathfrak{m}_{K_i}$ for each $i<n$. Note that the non-trivial relations hold due to Lemma \ref{nknastercapturinglemma}.
\begin{center}
\begin{tikzcd}
 &\mathfrak{m}_{K_2} \arrow[d, rightarrow] &\mathfrak{m}_{K_3} \arrow[l, rightarrow] \arrow[d, rightarrow]&\dots \arrow[l, rightarrow] \arrow[d, rightarrow] & \mathfrak{m}_{K_n}\arrow[l, rightarrow]\arrow[d, rightarrow] & \mathfrak{c} \arrow[l, rightarrow]\\
 \omega_1 & \mathfrak{m}^2_\mathcal{F} \arrow[l,rightarrow] & \mathfrak{m}^3_\mathcal{F} \arrow[l, rightarrow] &\dots\arrow[l, rightarrow] & \mathfrak{m}^n_\mathcal{F} \arrow[l, rightarrow] &
\end{tikzcd}
\end{center} 
\end{rem}
As the following proposition suggests,  having $n$-Knaster property may be optimal in terms $n$-preserving a construction scheme.
\begin{proposition}\label{propequivalenceknasterpreservescheme}Let $2\leq n\in \omega$ and $\mathcal{F}$ be an $n$-capturing construction scheme for which $\mathfrak{m}_\mathcal{F}^n>\omega_1$. Given a $ccc$ forcing $\mathbb{P}$ and $m\geq n$, the following statements are equivalent:
\begin{enumerate}[label=$(\alph*)$]
\item $\mathbb{P}$ is $m$-Knaster.
    \item $\mathbb{P}$ is $n$-Knaster.
    \item $\mathbb{P}$ $n$-preserves $\mathcal{F}$.
    \item $\mathbb{P}$ has precaliber $\omega_1$.
\end{enumerate}
    \begin{proof} $(a)$ implies $(b)$ and $(d)$ implies $(a)$ are obvious, and $(b)$ implies $(c)$ is just Lemma \ref{nknastercapturinglemma} (which does not require $\mathfrak{m}^n_\mathcal{F}>\omega_1$). In order to prove that $(c)$ implies $(d)$, let $\mathcal{A}\in [\mathbb{P}]^{\omega_1}$. Since $\mathfrak{m}_\mathcal{F}^n>\omega_1$, so is $\mathfrak{m}(\mathbb{P})$. As $\mathbb{P}$ is $ccc$, this means that there is a filter $G$ over $\mathbb{P}$ so that $\mathcal{A}\cap G$ is uncountable. 
    \end{proof}
    \end{proposition}

\begin{corollary}Let $2\leq n\in \omega$ and $\mathcal{F}$ be an $n$-capturing construction scheme. If $\mathfrak{m}^n_\mathcal{F}>\omega_1$, then $\mathfrak{m}^n_\mathcal{F}=\mathfrak{m}_{K_m}$ for each $m\geq n$.
    
\end{corollary}
A a corollary of Theorem \ref{starfailscapturing}, we also have the following.
\begin{corollary}Let $2\leq n\in \omega$ and $\mathcal{F}$ be an $n$-capturing construction scheme. If $\mathfrak{m}^n_\mathcal{F}>\omega_1$ then $\mathfrak{m}_{K_m}=\omega_1$ for each $m<n$ and $\mathfrak{m}^m_\mathcal{F}=\omega$ for each $m>n.$
    
\end{corollary}
In Corollary \ref{corollarymanotca2} we will explicitly show that the inequality $\mathfrak{m}>\omega_1$ is inconsistent with the existence of a $2$-capturing construction scheme. As a consequence of such result, we have that the cardinals $\mathfrak{m}^n_\mathcal{F}$ and $\mathfrak{m}$ are incomparable in a strong sense.
\begin{corollary}Let $2\leq n\in \omega$ and $\mathcal{F}$ be a construction scheme. Either $\mathfrak{m}^n_\mathcal{F}\leq \omega_1$ or $\mathfrak{m}=\omega_1$.
\end{corollary}
The following Lemma is all that we need to prove the main result of this section.

\begin{lemma}\label{preservincapturinglemma2}Let $\mathcal{F}$ be an $n$-capturing construction scheme and $(\langle \mathbb{P}_\xi\rangle_{\xi\leq \gamma},\langle \dot{\mathbb{Q}}_\xi\rangle_{\xi<\gamma})$ be a finite support iteration of $ccc$ forcings so that $$\mathbb{P}_\xi\Vdash\text{\say{ $\dot{\mathbb{Q}}_\xi$ $n$-preserves $\mathcal{F}$ }}.$$
Then $\mathbb{P}_\gamma$ also $n$-preserves $\mathcal{F}$.
\begin{proof}The proof is carried by induction over $\gamma$ by appealing to the equivalence of $n$-capturing provided by Lemma \ref{equivalencecapturing}. Both the base and the successor steps of the induction are trivial to show. Hence, we will only do the limit case. For this, let us assume that $\gamma$ is limit and we have already showed that $\mathbb{P}_\alpha$ $n$-preserves $\mathcal{F}$ for each $\alpha<\gamma$. Let $\dot{X}$ be a $\mathbb{P}_\gamma$-name for an uncountable subset of $\omega_1$ and consider an arbitrary $q\in \mathbb{P}_\gamma$. We can take a sequence $\langle p_\xi,\alpha_\xi\rangle_{\xi\in\omega_1}\subseteq \mathbb{P}_\gamma\times \omega_1$ such that for any two distinct $\xi,\mu\in \omega_1$, the following properties hold:
\begin{itemize}
    \item $p_\xi\leq q$ and $p_\xi\Vdash\text{\say{$\alpha_\xi\in \dot{X}$ }}.$
    \item $\alpha_\xi\not=\alpha_\mu$.
\end{itemize}
By refining the sequence if necessary, we may assume that $\{ dom(p_\xi)\,:\,\xi\in\omega_1\}$ forms a $\Delta$-system with root $R\subseteq \gamma$. Now, let $\alpha<\gamma$ be such that $R\subseteq \alpha$. Observe that $q\in \mathbb{P}_\alpha$. As $\mathbb{P}_\alpha$ is $ccc$, there is $G$ a $\mathbb{P}_\alpha$-generic filter over $V$ so that $S=\{\xi\in \omega_1\,:\,p_\xi|_\alpha\in G \}$ is uncountable. In particular, this implies that for any $\xi,\mu\in S$, the conditions $p_\xi|_\alpha$ and $p_\mu|_\alpha$ are compatible. Since $R=dom(p_\xi)\cap dom(p_\mu)\subseteq \alpha$, it follows that $p_\xi$ and $p_\mu$ are compatible. Now, according to the inductive hypotheses, $\mathcal{F}$ is $n$-capturing inside $V[G]$. Therefore, there is $D\in [S]^{n-1}$ so that the family $\{\alpha_\xi\,:\,\xi\in D\}$ is captured. To finish, let $q\in \mathbb{P}_\gamma$ be such that $p\leq p_\xi$ for any $\xi\in D$. It is straightforward that $p\leq q$ and $$p\Vdash\text{\say{ $\{\alpha_\xi\,:\,\xi\in D\}\in [\dot{X}]^n$ and it is captured }}.$$
Thus, the proof is over.
\end{proof}
\end{lemma}

\begin{theorem}\label{consistencymnfcontinuum}Let $\kappa>\omega_1$ be a regular cardinal such that $2^{<\kappa}=\kappa$. Given $\mathcal{F}$ an $n$-capturing construction scheme, there is a $ccc$-forcing $\mathbb{P}$ for which $$\mathbb{P}\Vdash\text{\say{ $\mathfrak{m}_{\mathcal{F}}^n=\mathfrak{c}=\kappa$ }}.$$
\begin{proof}We can construct a finite support iteration $(\langle \mathbb{P}_\xi\rangle_{\xi\leq \kappa},\langle \dot{\mathbb{Q}}_\xi\rangle_{\xi<\kappa})$ of $ccc$ forcings so that following properties hold for each $\xi<\kappa$:
\begin{itemize}
    \item $\mathbb{P}_\xi\Vdash \text{\say{$|\dot{\mathbb{Q}}_\xi|<\kappa$}},$
    \item $\mathbb{P}_\xi\Vdash \text{\say{ $\dot{\mathbb{Q}}_\xi$ $n$-preserves $\mathcal{F}$ }}$.
\end{itemize}
If the iteration is constructed with an appropriate bookkeeping, we can arrange that whenever $G$ is $\mathbb{P}_\kappa$-generic over $V$ and $\mathbb{Q}$ is a $ccc$ forcing inside $V[G]$ of size less than $\kappa$ which $n$-preserves $\mathcal{F}$, then there are cofinally many $\xi<\kappa$ for which $\mathbb{Q}$ is order isomorphic to $\dot{\mathbb{Q}}^G_\kappa.$ By standard arguments, all of these properties imply that $\mathbb{P}=\mathbb{P}_\kappa$ satisfies the desired conclusion.
\end{proof}
\end{theorem}

\section{The $n$-projection filter}

Our next goal is to show construction schemes which are $n$-capturing may not be $\mathcal{P}$-$n$-capturing for any non-trivial partition $\mathcal{P}$ of $\omega$. For this, we introduce a natural filter over $\omega$ which is definable from a construction scheme.
\begin{definition}Let $2\leq n \in \omega$ and $\mathcal{F}$ be a construction scheme. Given $S\subseteq \omega_1$, we define the \textit{$n$-projection} of $S$ as $$\pi_n(S)=\{ \rho^D\,:\,D\in [S]^n\textit{ and }D\textit{ is captured }\}$$
    
\end{definition}
\begin{rem}Given $S\subseteq \omega_1 $, the $n$-projection of $S$ can be recovered from any cofinal subset of $[S]^{<\omega}$. That is, if $G$ is a cofinal subset of $[S]^{<\omega}$, then  $$\pi_n(S)=\bigcup_{D\in G}\pi_n(D).$$
    
\end{rem}

\begin{definition}[The $n$-projection filter]Let  $n\in \omega$ and $\mathcal{F}$ be a construction scheme. We define $\mathcal{U}_n(\mathcal{F})$ as the set of all $A\subseteq \omega$ for which there is $S\in [\omega_1]^{\omega_1}$ such that $$\pi_n(S)\subseteq A.$$ 
\end{definition}
\begin{rem}Note that if $\mathcal{F}$ is not $n$-capturing then $\emptyset\in \mathcal{U}_n(\mathcal{F})$. Therefore, $\mathcal{U}_n(\mathcal{F})$ is not actually a filter in this case. 
\end{rem}

\begin{lemma}\label{lemmafilteruf}Let $2\leq n \in \omega$ and $\mathcal{F}$ be a construction scheme. If $\mathcal{F}$ is $n$-capturing, then $\mathcal{U}_n(\mathcal{F})$ is a non-principal filter over $\omega$.
\begin{proof}Since $\mathcal{F}$ is $n$-capturing, it follows that each member of $\mathcal{U}_n(\mathcal{F})$ is non-empty. Thus, in order to prove that $\mathcal{U}_n(\mathcal{F})$ is a filter, it is enough to show that the family $\{\,\pi_n(S)\,:\,S\in[\omega_1]^{\omega_1}\,\}$ is downwards directed with respect to $\subseteq$. Let $S,S'\in [\omega_1]^{\omega_1}$. We will prove that there is $A\in [\omega_1]^{\omega_1}$ such that $\pi_n(A)\subseteq \pi_n(S)\cap \pi_n(S')$. For this, first note that we can recursively construct a sequence $\langle \alpha_\xi,\beta_\xi\rangle\subseteq S\times S'$ so that  $\alpha_\xi<\beta_\xi<\alpha_{\xi+1}$
for each $\xi\in \omega_1$. By refining the sequence if necessary, we may assume without loss of generality that there is $k\in\omega$ so that for any two distinct $\xi,\mu\in \omega_1$, the following conditions hold:
\begin{enumerate}[label=$(\arabic*)$]
    \item $\rho(\alpha_\xi,\beta_\xi)<k,$
    \item $\lVert \alpha_\xi\rVert_k=\lVert \alpha_\mu\rVert_k$ and $\lVert \beta_\xi\rVert_k=\lVert \beta_\mu\rVert_k$. In particular, $\rho(\beta_\mu,\beta_\xi)>k$ by means of Lemma \ref{lemmadelta1}.
\end{enumerate}
Let $A=\{\beta_\xi\,:\,\xi\in \omega_1\}$. Then $A\in[ S']^{\omega_1}$, so $\pi_n(A)\subseteq \pi_n(S')$. We claim that $\pi_n(A)\subseteq \pi_n(S)$. For this purpose, let $l\in \pi_n(A)$ and consider $D\in[\omega_1]^n$ so that $\{ \beta_\xi\,:\xi\in D\}$ is captured at level $l$. Note that  $l>k$ due to the condition (2).\\
 
\noindent
\underline{Claim}: If $\xi,\mu\in D$ are distinct then $\Delta(\alpha_\xi,\alpha_\mu)=l=\rho(\alpha_\xi,\alpha_\mu)$.
\begin{claimproof}[Proof of claim] Let $h:(\beta_\xi)_k\longrightarrow (\beta_\mu)_k$ be the increasing bijection. By the point (1),  we have that $\alpha_\xi\in (\beta_\xi)_k$ and $\alpha_\mu\in (\beta_\mu)_k$. Furthermore, $h(\alpha_\xi)=\alpha_\mu$ by the point (2).  Since $\alpha_\xi\not=\alpha_\mu$, we may use Lemma \ref{lemmahdeltarhoinequalities} to conclude that $$l=\rho(\beta_\xi,\beta_\mu)\geq \rho(\alpha_\xi,\alpha_\mu)\geq \Delta(\alpha_\xi,\alpha_\mu)\geq \Delta(\beta_\xi,\beta_\mu)=l.$$
So we are done.
\end{claimproof}
Note that $\Xi_{\alpha_\xi}(l)\geq 0$ for each $\xi\in D$. This is due to the point (b) of Lemma \ref{lemmaxi}. By the point (c)  of such lemma, $\Xi_{\alpha_\xi}(l)=\Xi_{\beta_\xi}(l)$. Thus, $\{\alpha_\xi\,:\xi\in D\}$ is captured at level $l$ by means of Proposition \ref{deltarhoequalityprop}. That is, $l\in \pi_n(S)$. As $l\in \pi_n(A)$ was arbitrary, we get that  $\pi_n(A)\subseteq \pi_n(S)$.\\

 Now, we will show that $\mathcal{U}_n(\mathcal{F})$ is non-principal. Let $S\in [\omega_1]^{\omega_1}$ and  $k\in \omega$. Then there is $S'\in [S]^{\omega_1}$ such that $\lVert \alpha\rVert_k=\lVert \beta\rVert_k$ for all $\alpha,\beta\in S$. This implies that $\rho^D>k$ for any $D\in [S']^n$. Therefore, $\pi_n[S']\subseteq \pi_n[S]\backslash k$. This finishes the proof.
    
\end{proof}
    
\end{lemma}
\begin{rem}\label{remarkpartitionpositive}Suppose that $\mathcal{F}$ is an $n$-capturing construction scheme. Note that $A\in \mathcal{U}_n(\mathcal{F})^+$ if and only if for any $S\in [\omega_1]^{\omega_1}$ there are infinitely many $l\in A$ for which there is $D\in [S]^n$ which is captured at level $l$. In this way,  $\mathcal{F}$ is $\mathcal{P}$-$n$-capturing for a partition $\mathcal{P}$ of $\omega$ if and only if $\mathcal{P}\subseteq \mathcal{U}_n(\mathcal{F})^+$.
\end{rem}
As a direct corollary of the previous remark we have that:
\begin{corollary}Let $2\leq n\in \omega$ and $\mathcal{F}$ be a $n$-capturing construction scheme. There is a non-trivial partition $\mathcal{P}$ of $\omega$ for which $\mathcal{F}$ is $\mathcal{P}$-$n$-capturing if and only if $\mathcal{U}_n(\mathcal{F})$ is not an ultrafilter.
\end{corollary}
\section{$n$-capturing with partitions}
In this subsection we will prove Theorem \ref{theoremofmF}. Furthermore, we will show that it is consistent the existence of a construction scheme which is $n$-capturing and $\mathcal{P}$-$(n-1)$-capturing for some non-trivial partition of $\omega$, but it is not $\mathcal{P}'$-$n$-capturing for any non-trivial partition $\mathcal{P}'$.  This will be done by  starting with a $\mathcal{P}$-$n$-capturing construction scheme, and building a suitable finite support iteration of $ccc$ forcings which force $\mathcal{U}_{n}(\mathcal{F})$ to be an ultrafilter. In the following definition, we describe the forcings that will be used for this task.

\begin{definition}Let $\mathcal{F}$ be a construction scheme. Given $2\leq n\in \omega$ and $A\subseteq \omega$, we define the forcing $\mathbb{D}_n(\mathcal{F},A)$ as the family of all $p\in \text[\omega_1]^{<\omega}$ with the following property:
\begin{center}
There is no $D\in[p]^n$ such that $D$ is captured at some level $l\in A$.
\end{center}
We order $\mathbb{D}_n(\mathcal{F},A)$ with respect to $\subseteq$.
\end{definition}
The previous forcing was considered in \cite{forcingandconstructionschemes} for the particular case where $A=\omega$.

\begin{lemma}\label{Dnccc}Let  $2\leq n\in \omega$ and $\mathcal{F}$ be a $2$-capturing construction scheme. Then $\mathbb{D}_n(\mathcal{F},A)$ is $ccc$ for any $A\subseteq \omega$.
\begin{proof}Let $\mathcal{A}$ be an uncountable subset of $\mathbb{D}_n(\mathcal{F},A)$. Consider  $\mathcal{A}'\subseteq \mathcal{A}$ an uncountable root-tail-tail $\Delta$-system with root $R$ so that any two elements of $\mathcal{A}'$ have the same cardinality. Note that  we can enumerate $\mathcal{A}'$ as $\langle p_\alpha\rangle_{\alpha\in \omega_1}$ in such way that $\max(p_\alpha\backslash R)<\min(p_\beta\backslash R)$ whenever $\alpha<\beta.$ Since $\mathcal{F}$ is $2$-capturing, we can find $\delta<\gamma\in \text{Lim}$ so that the set $\{p_\delta\cup p_{\delta+1},p_\gamma\cup p_{\gamma+1}\}$ is captured at some level $l$. The proof will end by showing the following claim.\\

\noindent
\underline{Claim:} $p=p_\delta\cup p_{\gamma+1}\in \mathbb{D}_n(\mathcal{F},A)$.
\begin{claimproof}[Proof of claim] Let $D\in [p]^n$ which is captured and $h:(p_\delta\cup p_{\delta+1})_{l-1}\longrightarrow (p_\gamma\cup p_{\gamma+1})_{l-1}$ be the increasing bijection. It is easy to see that $h[p_\delta]=p_\gamma$ and $h[p_{\delta+1}]=p_{\gamma+1}$. Therefore, $$\lVert \alpha\rVert_{l-1}=\lVert h(\alpha)\rVert_{l-1}<\lVert \beta\rVert_{l-1} $$
for each $\alpha\in p_\delta\backslash R$ and $\beta\in p_{\gamma+1}\backslash R.$ Thus, $\rho^D<l$ by means of the point (1) of Proposition \ref{deltarhoequalityprop}. In this way,  either $D\subseteq p_\delta$ or $D\subseteq p_{\gamma+1}$. As both $p_\delta$ and $p_{\gamma+1}$ belong to $\mathbb{D}_n(\mathcal{F},A)$, it follows that $\rho^D\notin A$. This finishes the proof.
 \end{claimproof}
\end{proof} 
\end{lemma}
The following corollary will help us turn $\mathcal{U}_n(\mathcal{F})$ into an ultrafilter.
\begin{corollary}\label{notpositiveforcinglemma} Let $2\leq n\in \omega$, $\mathcal{F}$ be a $2$-capturing construction scheme and $A\subseteq \omega$. Then there is a condition $p\in \mathbb{D}_n(\mathcal{F},A)$ so that $$p\Vdash \text{\say{ $A\not\in \mathcal{U}_n(\mathcal{F})^+$ }}.$$
\begin{proof}Since $\mathcal{F}$ is $2$-capturing, $\mathbb{D}_n(\mathcal{F},A)$ is a $ccc$-forcing. Furthermore, it is uncountable. Thus, there is $p\in \mathbb{D}_n(\mathcal{F},A)$ which forces the generic filter to be uncountable. It follows that if $G$ is a generic filter over $V$ so that $p\in G$, then $S_G=\bigcup G$ is an uncountable subset of $\omega_1$ with $\pi_n(S_G)\cap A=\emptyset.$ Thus, in $V[G]$, $A\notin \mathcal{U}_{n+1}(\mathcal{F})$.
\end{proof}
\end{corollary}
Many of the results in the previous chapter already imply that $\mathfrak{m}>\omega_1$ is inconsistent with $CA_2$. We give a direct proof here for convenience of the reader.
\begin{corollary}[Under $\mathfrak{m}>\omega_1$]\label{corollarymanotca2}There are no $2$-capturing construction schemes for any type.
\begin{proof}Suppose towards a contradiction that there is a $2$-capturing construction scheme, namely $\mathcal{F}$. Then $\mathbb{D}_n(\mathcal{F},\omega)$ is an uncountable $ccc$-forcing. As $\mathfrak{m}>\omega_1$ such forcing contains an uncountable filter, namely $G$. Then $S_G=\bigcup G$ is an uncountable subset of $\omega_1$ so that $\pi_2(S_G)$ is empty. This contradicts the fact that $\mathcal{F}$ is $2$-capturing, so we are done.
    
\end{proof}
\end{corollary}

\begin{lemma}\label{Dncapturingpreservinglemma}Let $2\leq n\in \omega$ and $\mathcal{F}$ be an $n$-capturing construction scheme. If $A\in \mathcal{U}_n(\mathcal{F})^+$ then $\mathbb{D}_n(\mathcal{F},\omega\backslash A)$ $n$-preserves $\mathcal{F}.$
\begin{proof}We will prove this lemma by appealing to the equivalence provided by Lemma \ref{lemmaequivalencefunctioncapturingpreserving}. Let $\mathcal{A}$ be an uncountable subset of $\mathbb{D}_n(\mathcal{F},\omega\backslash A)$ and $\nu:\mathcal{A}\longrightarrow \omega_1$ be an injective function.
Given $p\in \mathcal{A}$, let $D_p=p\cup \{\nu(p)\}$ and $\alpha_p=\max(D_p)$. By refining $\mathcal{A}$ if necessary, we may assume that are $j\in \omega$, $a<m_j$, $C\subseteq a+1$ and $b\in C$ so that the following conditions hold for any $p\in \mathcal{A}$:
\begin{enumerate}[label=$(\alph*)$]
    \item $\rho^{D_p}\leq j,$
    \item $\lVert \alpha_p\rVert_j=a$,
    \item $(\alpha_p)_j[C]=p$ and $(\alpha_p)_j(b)=\nu(p)$.
\end{enumerate}
As $A\in \mathcal{U}_n(\mathcal{F})^+$, there are $l\in A$ and $\{p_0,\dots,p_{n-1}\}\in [\mathcal{A}]^n$ for which $\{\,\alpha_{p_i}\,:i<n\,\}$ is captured at level $l$. Note that $l>j$ due to the point (b) above. According to Lemma \ref{capturedfamiliestosetslemma}, the family $\{D_{p_i}\,:\,i<n\}$ is  captured at level $l$. From this fact and by the point (c), the same holds for both $\{p_0,\dots,p_{n-1}\}$ and $\{\nu(p_0),\dots,\nu(p_{n-1})\}$. The proof follows from the following claim. \\

\noindent
\underline{Claim}: $p=\bigcup\limits_{i<n}p_i\in \mathbb{D}_n(\mathcal{F},\omega\backslash A).$
\begin{claimproof}[Proof of claim]Let  $D\in [p]^n$ which is captured and consider $F\in \mathcal{F}_l$ such that $p\subseteq F$. As $D\subseteq p$, then $\rho^D\leq \rho^p=l$. If $\rho^D=l$, we are done because $l\not\in \omega\backslash A.$ On the other hand, if  $\rho^D<l$ then $D\subseteq F_{\Xi_D(l)}$.  Since  the family $\{p_0,\dots,p_n\}$ is captured at level $l$, there is $i<n$ so that   $F_{\Xi_D(l)}\cap p=p_i$. In this way, $D\subseteq p_i$. Thus, as $p_i\in \mathbb{D}_n(\mathcal{F},\omega\backslash A)$ then $\rho^D\not\in \omega\backslash A$.
\end{claimproof}
\end{proof}
\end{lemma}

\begin{lemma}\label{preservingpositivelemma1}Let $3\leq n\in \omega$, $\mathcal{F}$ be an $n$-capturing construction scheme and $A \subseteq \omega$. If $B\in \mathcal{U}_{n-1}(\mathcal{F})^+$, then  $$\mathbb{D}_n(\mathcal{F},A)\Vdash \text{\say{ $B\in \mathcal{U}_{n-1}(\mathcal{F})^+$ }}.$$

\begin{proof}Let $\dot{X}$ be a name for an  uncountable subset of $\omega_1$ and $q\in \mathbb{D}_{n}(\mathcal{F}, A)$. Using that $B\in \mathcal{U}_{n-1}(\mathcal{F})^+$ and by arguing in a similar way as in the previous lemmas, we can find a sequence $\langle p_i,\alpha_i\rangle_{i<{n-1}}\subseteq \mathbb{D}_n(\mathcal{F},A)\times \omega_1$ so that:
\begin{enumerate}[label=$(\arabic*)$]
    \item both $\{p_0,\dots,p_{n-2}\}$ and $\{\alpha_0,\dots,\alpha_{n-2}\}$ are captured at some level $l\in B$,
    \item $\alpha_i\not=\alpha_j$ whenever $i\not=j$,
    \item $p_i\leq q$ and $p_i\Vdash\text{\say{ $\alpha_i\in \dot{X}$ }}$ for each $i<n-1.$ 
    \end{enumerate}
    It is easy to see that $p=\bigcup\limits_{i<n-1}p_i$ is a condition in the forcing that we are considering. Furthermore, $p\leq q$ and $p\Vdash \text{\say{ $\{\alpha_0,\dots,\alpha_{n-1}\}\in [\dot{X}]^{n-1}$ }}.$ This means that $p\Vdash \text{\say{ $\pi_{n-1}(\dot{X})\cap B\not=\emptyset$ }}.$ Thus, the proof is over.
\end{proof}
\end{lemma}

The following lemma is proved in the exact same way as the previous one. For that reason, we leave the proof to the reader.
\begin{lemma}\label{preservingpositiveslemma2}Let $\mathcal{F}$ be an $n$-capturing construction scheme, $B\in \mathcal{U}_n(\mathcal{F})^+$ and $(\langle \mathbb{P}_\xi\rangle_{\xi\leq \gamma},\langle \dot{\mathbb{Q}}_\xi\rangle_{\xi<\gamma})$ be a finite support iteration of $ccc$ forcings so that $$\mathbb{P}_\xi\Vdash\text{\say{ $\dot{\mathbb{Q}}_\xi\Vdash\text{\say{$B\in \mathcal{U}_{n}(\mathcal{F})^+$}}$}}.$$  Then $\mathbb{P}_\gamma\Vdash\text{\say{$B\in \mathcal{U}_{n}(\mathcal{F})^+$}}$.
\end{lemma}

By combining all the results we have so far, we get the following theorem.
\begin{theorem}Let $\mathcal{F}$ be $\mathcal{P}$-$n$-capturing construction scheme. There is a $ccc$ forcing $\mathbb{P}$ satisfies the following properties: $$\mathbb{P}\Vdash\text{\say{$\mathcal{F}$ is $n$-capturing and $\mathcal{P}$-$(n-1)$-capturing}},$$
$$\mathbb{P}\Vdash\text{\say{$\mathcal{U}_n(\mathcal{F})$ is an ultrafilter }}.$$
In particular, $\mathbb{P}$ forces that $\mathcal{F}$ is not $\mathcal{P}'$-$n$-capturing for any non-trivial partition $\mathcal{P}'$ of $\omega$.
\begin{proof}Let $\kappa=\mathfrak{c}$. We can construct a finite support iteration  $(\langle \mathbb{P}_\xi\rangle_{\xi\leq \kappa},\langle \dot{\mathbb{Q}}_\xi\rangle_{\xi<\kappa})$  of forcings so that given $\xi<\kappa$, the following condition holds:
$$\mathbb{P}_\xi\Vdash\text{\say{$\dot{\mathbb{Q}_\xi}=\mathbb{D}_n(\mathcal{F},\dot{A})$ for some $\dot{A}\in[\omega]^{\omega}$ with $\dot{A},\omega\backslash \dot{A}\in \mathcal{U}_n(\mathcal{F})^+$}}.$$ 
By means of Lemma \ref{Dnccc}, it follows that $\mathbb{P}_\kappa$ is $ccc$. According to Lemma \ref{Dncapturingpreservinglemma} and \ref{preservincapturinglemma2}, we have that $\mathbb{P}$ $n$-preserves $\mathcal{F}$. Furthermore, $\mathbb{P}\Vdash\text{\say{ $\mathcal{P}\subseteq\mathcal{U}_{n-1}(\mathcal{F})^+$}}$ due to Lemma \ref{preservingpositivelemma1} and \ref{preservingpositiveslemma2}. If we construct the iteration with an appropriate bookkeeping and by making use of Corollary \ref{notpositiveforcinglemma}, we can arrange $\mathbb{P}_\kappa$ to force   $\mathcal{U}_n(\mathcal{F})$ to be an ultrafilter.
This finishes the proof.
\end{proof}
\end{theorem}
\begin{corollary}\label{partitionncapturingcorolary}Let $1\leq n\in \omega$. It is consistent that there is a construction scheme $\mathcal{F}$ such that:\begin{itemize}
    \item There is a non-trivial partition of $\omega$, namely $\mathcal{P}$, such that $\mathcal{F}$ is $\mathcal{P}$-$n$-capturing.
    \item $\mathcal{F}$ is $n+1$-capturing but it is not $\mathcal{P}'$-$n+1$-capturing for any non-trivial partition $\mathcal{P}'$.
\end{itemize}
\end{corollary}

    \section{Ramsey ultrafilters from construction schemes}
Let $\mathcal{U}$ be an ultrafilter over $\omega$. Recall that $\mathcal{U}$ is \textit{Ramsey} if for each infinite partition  $\mathcal{P}\subseteq \mathcal{U}^*$ of $\omega$, there is $A\in \mathcal{U}$ so that $|A\cap P|\leq 1$
    for all $P\in \mathcal{P}.$ Ramsey ultrafilters, also called selective, play a highly important roll in modern set theory. In \cite{mapsontrees} Stevo Todor\v{c}evi\'c presented a clever ways of defining filters over $\omega_1$ by using walks on ordinals. In \cite{walksultrafilters}, he used such filters in order to show that the existence of Ramsey ultrafilters follows from $\mathfrak{m}>\omega_1$. This result is quite interesting since the existence of ultrafilters with particular properties usually follows from equalities of the form \say{$\textit{cardinal invariant}=\mathfrak{c}$} instead of inequalities of the form \say{$\textit{cardinal invariant}>\omega_1$}. In this section we will show that given a $2$-capturing construction scheme $\mathcal{F}$, the filter $\mathcal{U}_n(\mathcal{F})$ is a Ramsey ultrafilter under $\mathfrak{m}^n_\mathcal{F}>\omega_1$. Furthermore, we will define the filter $\mathcal{V}(\mathcal{F})$ and show that, under $\mathfrak{m}_\mathcal{F}>\omega_1$, this filter is an ultrafilter over $\omega_1$ which projects into a Ramsey ultrafilter. Beyond the hypothesis imposed over the cardinal invariant $\mathfrak{m}^n_\mathcal{F}$, the most important feature of the following constructions is that both ultrafilters are explicitly definable using a combinatorial structure over $\omega_1$.

    \begin{theorem}\label{UnFramseytheorem}Let $2\leq n\in \omega$ and $\mathcal{F}$ be an $n$-capturing construction scheme for which $\mathfrak{m}^n_\mathcal{F}>\omega_1$. Then $\mathcal{U}_n(\mathcal{F})$ is a Ramsey ultrafilter.
    \begin{proof} Suppose towards a contradiction that $\mathcal{U}_n(\mathcal{F})$ is not an ultrafilter and let $A\subseteq \omega$ be such that $A, \omega\backslash A\in \mathcal{U}_n(\mathcal{F})^+$. According to the Lemmas \ref{Dnccc} and \ref{Dncapturingpreservinglemma}, $\mathbb{D}_n(\mathcal{F},A)$ is a $ccc$ forcing so that $\mathbb{D}_n(\mathcal{F},A)\Vdash\text{\say{ $\mathcal{F}$ is $n$-capturing }}$. Hence, its Martin's number is greater than $\omega_1$ due to the hypotheses. In this way, there is an uncountable filter $G$ over $\mathbb{D}_n(\mathcal{F},A)$. Let $S_G=\bigcup G$. Then $S_G$ is an uncountable subset of $\omega_1$ satisfying that $\pi_n(S_G)\cap A=\emptyset.$ This contradicts the fact that $A\in \mathcal{U}_n(\mathcal{F})^+$.

    Now, we will show that $\mathcal{U}_n(\mathcal{F})$ is Ramsey. For this purpose, let $\langle P_i\rangle_{i\in \omega}$ be a partition of $\omega$ so that $P_i\notin \mathcal{U}_n(\mathcal{F})$ for every $i\in \omega$. It suffices to prove that there is $S\in [\omega_1]^{\omega_1}$ with $|\pi_n(S)\cap P_i|\leq 1$ for any $i\in \omega$. Let $\mathbb{P}$ be the forcing consisting of the empty set and all $p\in \text{FIN}(\omega_1)$ so that: $$|\pi_n(p)\cap A_i|\leq 1$$
    for each $i\in \omega$. It is easy to see that if $G$ is an uncountable filter over $\mathbb{P}$ then  $S_G=\bigcup G$ satisfies the required property.  As we are assuming that $\mathfrak{m}_\mathcal{F}^n>\omega_1$, the existence of such $G$ will follow if $\mathbb{P}$ is $ccc$ and $\mathbb{P}\Vdash\text{\say{ $\mathcal{F}$ is $n$-capturing }}$. The proof this two facts is similar to ones in Lemmas \ref{Dnccc} and \ref{Dncapturingpreservinglemma}. For that reason, we leave the details to the reader.

    \end{proof} 
    \end{theorem}
    \begin{definition}[A square bracket operation]Let $\mathcal{F}$ be a construction scheme. We define $\llbracket\cdot,\cdot\rrbracket:[\omega_1]^2\longrightarrow \omega$ as:
    $$\llbracket \alpha,\beta\rrbracket=\min(\,(\beta)_{\rho(\alpha,\beta)-1}\backslash \alpha\,)$$
    for each $\alpha<\beta$. Given $S\subseteq \omega_1$, we also define $$\llparenthesis S\rrparenthesis=\{\,\llbracket \alpha,\beta\rrbracket \,:\,\alpha,\beta\in S\textit{ and }\{\alpha,\beta\} \textit{ is captured }\}$$

    \end{definition}

 \begin{multicols}{2}
\begin{rem}\label{remarksqbracket1} Suppose that $\alpha<\beta\in \omega_1$ and let $l=\rho(\alpha,\beta)$. If $F\in \mathcal{F}_{l}$ is such that $\{\alpha,\beta\}\subseteq F$, then $$\llbracket \alpha,\beta\rrbracket= \min(\,F_{\Xi_\beta(l)}\backslash R(F)\,).$$
\end{rem}
\columnbreak
\begin{center}
\begin{minipage}[b]{0.8\linewidth}
\centering
\includegraphics[width=6.6 cm, height=2.4cm]{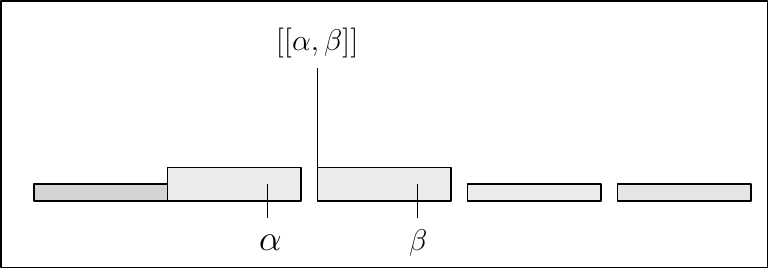}
\end{minipage}
\end{center}

\end{multicols}
\begin{rem}\label{remarkbracketcaptured}As a consequence of the previous remark we have that if $D_0,D_1\in \text{FIN}(\omega_1)$ are disjoint sets for which $\{D_0,D_1\}$ is captured, then $$\llbracket D_0(i),D_1(i)\rrbracket=\llbracket D_0(j),D_1(j)\rrbracket$$
for all $i,j<|D_0|.$
\end{rem}
\begin{definition}[The square-bracket filter]Let $\mathcal{F}$ be a construction scheme. We define $\mathcal{V}(\mathcal{F})$ as the set of all $A\subseteq \omega_1$ for which there is $S\in [\omega_1]^{\omega_1}$  such that $$\llparenthesis S\rrparenthesis \subseteq A.$$
\end{definition}
Recall that an filter over $\omega_1$ is uniform whenever all of its elements are uncountable.
\begin{lemma}Let $\mathcal{F}$ be a $2$-capturing construction scheme. Then $\mathcal{V}(\mathcal{F})$ is a uniform filter.
\begin{proof}Since $\mathcal{F}$ is $2$-capturing, each element of $\mathcal{V}(\mathcal{F})$ is non-empty. In order to show that $\mathcal{V}(\mathcal{F})$ is a filter, it is enough to show that the family $\{\llparenthesis S\rrparenthesis\,:\,S\in[\omega_1]^{\omega_1}\}$ is downwards directed. Let $S,S'\in [\omega_1]^{\omega_1}$. We will prove that there is $A\in [\omega_1]^{\omega_1}$ so that $\llparenthesis A\rrparenthesis\subseteq \llparenthesis S\rrparenthesis \cap \llparenthesis S'\rrparenthesis.$ Just as in the Lemma \ref{lemmafilteruf}, we can consider a sequence $\langle \alpha_\xi,\beta_\xi\rangle_{\xi\in\omega_1}\subseteq S\times S'$ so that $\alpha_\xi<\beta_\xi<\alpha_{\xi+1}$ for each $\xi\in\omega_1$. We may assume without loss of generality that there is $k\in\omega$ so that the following properties hold for any two distinct $\xi,\mu\in\omega_1:$
\begin{enumerate}[label=$(\arabic*)$]
    \item $\rho(\alpha_\xi,\beta_\xi)<k$,
    \item $\lVert \alpha_\xi\rVert_k=\lVert \alpha_\mu \rVert_k$ and $\lVert \beta_\xi\rVert_k=\lVert \beta_\mu\rVert_k$. In particular, this means that $\rho(\beta_\xi,\beta_\mu)>k.$
\end{enumerate}
 We claim that $A=\{\beta_\xi\,:\,\xi\in\omega_1\}$ works. Trivially $\llparenthesis A\rrparenthesis \subseteq \llparenthesis S'\rrparenthesis$  because $A\subseteq S'$. In order to prove that $\llparenthesis A\rrparenthesis \subseteq \llparenthesis S\rrparenthesis$, let $\delta\in \llparenthesis A\rrparenthesis.$ Then there are $\xi<\mu\in\omega_1$ for which $\{\beta_\xi,\beta_\mu\}$ is captured at some level $l$ and $\llbracket \beta_\xi,\beta_\mu\rrbracket=\delta$. The properties (1) and (2) stated above imply that the family $\{\{\alpha_\xi,\beta_\xi\},\{\alpha_\mu,\beta_\mu\}\}$ is also captured at level $l$ (see Lemma \ref{capturedfamiliestosetslemma}). In particular, $\{\alpha_\xi,\alpha_\mu\}$ is captured. Furthermore, $\delta=\llbracket \beta_\xi,\beta_\mu\rrbracket=\llbracket \alpha_\xi,\alpha_\mu\rrbracket$ by means of the Remark \ref{remarkbracketcaptured}. As $\alpha_\xi,\alpha_\mu \in S$, we have shown that $\delta\in \llparenthesis S\rrparenthesis.$ Thus, the proof is over.
\end{proof}    
\end{lemma}

\begin{theorem}Let $\mathcal{F}$ be a $2$ capturing construction scheme for which $\mathfrak{m}_\mathcal{F}>\omega_1$. Then $\mathcal{V}(\mathcal{F})$ is an ultrafilter.  
\begin{proof}Let $A\in \mathcal{V}(\mathcal{F})^+$. We will show that $A\in \mathcal{V}(\mathcal{F})$. For this purpose, let us consider the forcing $$\mathbb{P}=\{\,p\in[\omega_1]^{<\omega}\,:\,\llparenthesis p\rrparenthesis\subseteq A\,\}$$
ordered by reverse inclusion. It is straightforward that if $G\subseteq \mathbb{P}$ is an uncountable filter, then $S_G=\bigcup G$ is an uncountable subset of $\omega_1$ for which $\llparenthesis S_G\rrparenthesis \subseteq A$. In order to guarantee the existence of such filter, it is enough to prove that $\mathbb{P}$ is $ccc$ and $2$-preserves $\mathcal{F}$. This is because we are assuming that $\mathfrak{m}_\mathcal{F}>\omega_1$.\\

\noindent
\underline{Claim}: $\mathbb{P}$ is $ccc$ and $2$-preserves $\mathcal{F}$.
\begin{claimproof}[Proof of claim]Let $\mathcal{A}\in [\mathbb{P}]^{\omega_1}$ and $\nu:\mathcal{A}\longrightarrow \omega_1$ bi an injective function. For each $p\in \mathcal{A}$ put $D_p=p\cup\{\nu(p)\}$$\alpha_p=\max(D_p)$. By refining $\mathcal{A}$ if necessary, we may assume that there are $j\in\omega$, $a<m_j$, $C\subseteq a+1$ and $b\in C$ so that the following conditions hold for any $p\in \mathcal{A}$:
\begin{enumerate}[label=$(\alph*)$]
    \item $\rho^{D_p}\leq j$,
    \item $\lVert \alpha_p\rVert_j=a, $
    \item $(\alpha_p)_j[C]=p$ and $(\alpha_p)_j(b)=\nu(p)$.
\end{enumerate}
Since $A\in \mathcal{V}(\mathcal{F})^+$, there are $p_0,p_1\in \mathcal{A}$ for which $\{\alpha_{p_0},\alpha_{p_1}\}$ is captured at some level $l$ and $\llbracket \alpha_{p_0},\alpha_{p_1}\rrbracket\in A$. According to Lemma \ref{capturedfamiliestosetslemma}, $\{D_{p_0},D_{p_1}\}$ is also captured at level $l.$ From this fact and by point (c), the same holds for $\{p_0,p_1\}$ and $\{\nu(p_0),\nu(p_1)\}$. In order to finish the claim, we will show that $q=p_0\cup p_1\in \mathbb{P}$. Indeed, let $\alpha<\gamma\in q$ be such that $\{\alpha,\gamma\}$ is captured.  If $\rho(\alpha,\gamma)<l$, then there is $i\in 2$ for which $\alpha,\gamma\in p_i$. Thus, $\llbracket \alpha,\gamma\rrbracket \in A$. On the other hand, if $\rho(\alpha,\gamma)=l$, let $F\in \mathcal{F}_l$ be such that $q\subseteq F$. By Remark \ref{remarksqbracket1}, we have that $$\llbracket \alpha,\gamma\rrbracket= F_1\backslash R(F)=\llbracket \alpha_{p_0},\alpha_{p_1}\rrbracket\in A.$$

\end{claimproof}

\end{proof}    
\end{theorem}
Suppose that $\mathcal{F}$ is as in the previous Theorem. Since $\omega_1$ is not measurable, there is a surjective function $\pi:\omega_1\longrightarrow \omega$ so that $\pi^{-1}[\{n\}]\not\in \mathcal{V}(\mathcal{F})$
for each $n\in\omega$. From this, it follows that the family $$\pi \mathcal{V}(\mathcal{F})=\{A\in \omega\,:\,\pi^{-1}[A]\in \mathcal{V}(\mathcal{F})\}$$
is a non-principal ultrafilter over $\omega$.
We end this chapter by proving that it is a Ramsey ultrafilter as well.
\begin{theorem}\label{pivframseytheorem}Let $\mathcal{F}$ be a $2$-capturing construction scheme for which $\mathfrak{m}_\mathcal{F}>\omega_1$. then $\pi \mathcal{V}(\mathcal{F})$ is a Ramsey ultrafilter.
\begin{proof}Let $\langle P_i\rangle_{i\in\omega}$ be a partition of $\omega$ so that $P_i\not\in \pi\mathcal{V}(\mathcal{F}) $ for each $\in \omega$. We define the forcing $\mathbb{P}$ as the set of all $p\in [\omega_1]^{\omega}$ so that $$|\pi[ \llparenthesis p\rrparenthesis]\cap P_i|\leq 1$$
for each $i\in\omega$. We order $\mathbb{P}$ with the reverse inclusion.
Note that if $G\subseteq\mathbb{P}$ is an uncountable filter, then $S_G$ is an uncountable subset of $\omega_1$ such that $|\pi[\llparenthesis S_G\rrparenthesis]\cap P_i|\leq 1$ for each $i$. The existence of such filter follows from the following claim.\\\\
\noindent\underline{Claim}: $\mathbb{P}$ is $ccc$ and $2$-preserves $\mathcal{F}$.
\begin{claimproof}[Proof of claim]
Let $\mathcal{A}$ be an uncountable subset of $\mathbb{P}$ and $\nu:\mathcal{A}\longrightarrow \omega_1$ be an injective function. Given $p\in \mathcal{A}$, let $D_p=p\cup\{\nu(p)\}$ and $\alpha_p=\max(p)$.  By refining $\mathcal{A}$ if necessary, we may assume that there are $j\in\omega$, $a<m_j$, $C\subseteq a+1$ and $b\in C$ so that the following conditions hold for any $p\in \mathcal{A}$:
\begin{enumerate}[label=$(\alph*)$]
    \item $\rho^{D_p}\leq j$,
    \item $\lVert \alpha_p\rVert_j=a, $
    \item $(\alpha_p)_j[C]=p$ and $(\alpha_p)_j(b)=\nu(p)$.
\end{enumerate}
Furthermore, we may suppose that there is $E\subseteq \omega$ so that $\pi[\llparenthesis p\rrparenthesis]=E$ for each $p\in \mathcal{A}$. Let $I=\{i\in\omega\,:\,E\cap P_i\not=\emptyset\}$. Note that $I$ is finite. Therefore, $$A=\omega\backslash (\bigcup\limits_{i\in I}P_i)\in \pi \mathcal{V}(\mathcal{F}).$$
In particular, $A\cap \pi[ \llparenthesis \alpha_p\,:\,p\in \mathcal{A}\rrparenthesis]\not=\emptyset$. Thus, there are $p_0,p_1\in \mathcal{A}$ so that $\{\alpha_{p_0},\alpha_{p_1}\}$ is captured at some level $l$ and $\pi ( \llbracket \alpha_{p_0}\alpha_{p_1}\rrbracket )\in A$. By arguing in the same way as in the previous theorem, we can conclude that both $\{\nu(p_0),\nu(p_1)\}$ and $\{p_0,p_1\}$ are captured. Even more, for $q=p_0\cup p_1$ we have that $$\llparenthesis q\rrparenthesis=\{\llbracket \alpha_{p_0},\alpha_{p_1}\rrbracket\}\cup \llparenthesis p_0\rrparenthesis \cup \llparenthesis p_1\rrparenthesis.$$
In virtue of this,  $\pi[\llparenthesis q\rrparenthesis]=\{\pi(\llbracket \alpha_{p_0},\alpha_{p_1}\rrbracket)\}\cup E$. As $\pi[\llbracket \alpha_{p_0},\alpha_{p_1}\rrbracket]\not\in A$, it easily follows that $q\in \mathbb{P}$. This finishes the proof.
\end{claimproof}

\end{proof}
\end{theorem}

%% file: chapters/A_deeper_analysis_of_construction_schemes.tex
\chapter{On the existence of construction and capturing schemes}\label{deeperanalysis}
The main goal of this chapter is to present a complete proof  that the capturing axiom $FCA(part)$ follows from  Jensen's $\Diamond$-principle. In order to do that, we will first analyse in more depth some structural properties regarding construction scIn fact, according tonals. Among other things, we will conclude that if the type is good, any construction scheme over $\omega_1$ can be reconstructed from an $\subseteq$-increasing chain of construction schemes over the countable limit ordinal and viseversa. That is, any increasing chain of construction scheme over the countable limit ordinals determines a construction scheme over $\omega_1$.  Furthermore, the first element of such chain is uniquely determined.  In other words, there is a unique construction scheme over $\omega$ (which we call $\mathcal{F}(\omega)$). In Section \ref{sectionpfforcing}, we will define the forcing $\mathbb{P}(\mathcal{F})$. We will use this forcing as a tool for extending a construction scheme from a countable limit ordinal to the next one. In Subsection \ref{subsectionschemeszfc}, we will use the forcing $\mathbb{P}(\mathcal{F})$ to show that construction schemes (over $\omega_1$) of any given good type exist without the assumption of any extra axioms. In Section \ref{sectiondiamond}, we will prove that the $\Diamond$-principle implies $FCA(part)$. Finally, in Section \ref{sectionpidschemes}, we will prove that $PID$ is incompatible with the existence of a $2$-capturing construction scheme.\\

For the rest of this chapter, we fix a type $\tau=\langle m_k,n_{k+1},r_{k+1}\rangle_{k\in\omega}.$
\begin{definition}[The restriction of a scheme] Let $\mathcal{F}$ be a construction scheme over $X$. For any $Y\subseteq X$, we define the \textit{restriction of $\mathcal{F}$ to $Y$} as $$\mathcal{F}|_Y=\{F\in \mathcal{F}\,:\,F\subseteq Y\}.$$ 
\end{definition}
\begin{rem}If $Y\in \mathcal{F}$, then $\mathcal{F}|_{Y}$ is a construction scheme over $Y$.
    
\end{rem}
\begin{rem}\label{remarklevelrestriction} If $Y\subseteq X$ and $k\in \omega$, then $(\mathcal{F}|_Y)_k=\{F\in \mathcal{F}_k\,|\, F\subseteq Y\}$.
\end{rem}

\begin{rem}\label{remarkXfiniteconstscheme}If $\mathcal{F}$ is a construction scheme over a finite set $X$, then $X\in \mathcal{F}$. This is due to the point (a) of Definition \ref{constructionschemedef}.
\end{rem}

\begin{rem}If $\mathcal{F}$ is a construction scheme over $X$, and $Y\subseteq X$ is such that $\mathcal{F}|_Y$ is a construction scheme over $Y$, then $\rho_{\mathcal{F}}|_{Y^2}=\rho_{\mathcal{F}|_Y}$. 
\end{rem}
In the two following propositions we characterize the finite sets of ordinals which admit construction schemes. 

\begin{proposition}\label{uniquefiniteschemes} If $X$ is a finite set of ordinals, then there is at most one construction scheme over $X$ of type $\tau$.
\begin{proof}The proof is carried by induction over $m=|X|$. If $m=1$, the result is immediate. For the inductive step, suppose that we have proved the proposition for each set of ordinals of cardinality at most $m$ and consider $X$ such that $|X|=m+1$. Let $\mathcal{F}$ and $\mathcal{F}'$ be construction schemes over $X$. According to the Remark \ref{remarkXfiniteconstscheme}, $X\in \mathcal{F}\cap \mathcal{F}'$. Thus, there are $k,k'\in\omega$ so that $X\in \mathcal{F}_k$ and $X\in \mathcal{F}'_{k'}$. By the point (b) of Definition \ref{constructionschemedef}, $m_k=|X|=m_{k'}$. In this way, $k=k'$. Furthermore, since $|X|=m+1>1$ then $k>0$. Therefore, due to the point (d) of Definition \ref{constructionschemedef}, we can take   $\langle X_i\rangle_{i<n_k}$ and $\langle X'_i\rangle_{i<n_k}$ to be the canonical decompositions of $X$ inside $\mathcal{F}$ and $\mathcal{F}'$ respectively. In the Remark \ref{remarkcanonicaldecomposition}, we explicitly described the canonical decomposition of a set in terms of the type. From this description it is easy to see that $X_i=X'_i$ for any $i<n_k$. By means of Proposition \ref{propnoname} we can deduce that $$\mathcal{F}=\{X\}\cup (\bigcup_{i<n_k} \mathcal{F}|_{X_i}),$$
$$\mathcal{F}'=\{X\}\cup (\bigcup_{i<n_k} \mathcal{F}'|_{X_i}).$$
But $\mathcal{F}|_{X_i}=\mathcal{F}'|_{X_i}$ for each $i$ due to the inductive hypotheses. Consequently, $\mathcal{F}=\mathcal{F}'$. 
\end{proof}
    
\end{proposition}

\begin{proposition}Let $k\in\omega$. If $X$ is a finite set of ordinals and $|X|=m_k$ then there is a unique construction scheme over $X$ of type $\tau$.
\begin{proof}The uniqueness part was showed in Proposition \ref{uniquefiniteschemes}. The \say{existence} part of proof will be carried by induction over $k$.\\\\
\underline{Base step}: If $k=0$, then $m_k=1$. It should be clear that $\mathcal{F}=\{X\}$ is the only construction scheme over $X$.\\\\
\underline{Inductive Step}: Suppose that the proposition is true for some $k\in \omega$ and each set of ordinals of size $m_k$. Let $X$ be a set of ordinals with $|X|=m_{k+1}$. Given $i<n_{k+1}$, let us define $a_i=r_{k+1}+ i\cdot(m_k-r_{k+1})$ and
$$F_i=X[r_{k+1}]\cup X[\,[a_i,\, a_{i+1})\,].$$

Note that each $F_i$ has cardinality $m_k$. Therefore, by the inductive hypotheses there is a unique construction scheme $\mathcal{F}_i$ over $F^i$ of type $\tau$. Let $$\mathcal{F}=\{X\}\cup \big(\bigcup\limits_{i<n_{k+1}}\mathcal{F}^i\big).$$
We claim that $\mathcal{F}$ is the construction scheme over $X$ that we are looking for. $\mathcal{F}$ satisfies condition (a) of Definition \ref{constructionschemedef} because $X\in \mathcal{F}$. It satisfies condition (b) of such definition because for each $l\in\omega$, $$\mathcal{F}_l=\begin{cases}\emptyset&\textit{if }l>k+1\\
\{X\}&\textit{if }l=k+1\\
\bigcup\limits_{i<n_{k+1}}\mathcal{F}^i_l&\textit{if }l\leq k   
\end{cases}$$
Let us fix $l\in \omega$. We proceed to prove that the conditions (c) and (d) of Definition \ref{constructionschemedef} applied to $l$ also hold.  Note that when $l>k+1$ there is nothing to verify. This is because, in this case,  $\mathcal{F}_l$ is empty. Therefore, we will assume that $l\leq k+1.$
\begin{claimproof}[Proof of $(c)$] If $l=k+1$, condition (c) trivially holds for $l$ because $\mathcal{F}_l=\{X\}$. Therefore, we can assume that $l\leq k$. In this case, let $F,E\in \mathcal{F}_l$. We need to show that $E\cap F\sqsubseteq F$. According to the formula for $\mathcal{F}_l$ 
written above, we know that there are $i,j<n_{k+1}$ for which $E\in \mathcal{F}^i_l$ and $F\in \mathcal{F}^j_l.$ If $i=j$, we are done since $\mathcal{F}^i$ is already a construction scheme. So suppose that $i\not=j$ and consider $\phi:F_i\longrightarrow F_j$ the increasing bijection. It is not hard to see that $\mathcal{F}^\phi=\{\,\phi[G]\,:\,G\in \mathcal{F}^i\,\}$ is a construction scheme over $F_j$ of type $\tau$. Therefore, $\mathcal{F}^j=\mathcal{F}^\phi$ by virtue of the inductive hypotheses. In particular, there is $G\in \mathcal{F}_l^i$ for which $F= \phi[G]$. According to the definitions of $F_i$ and $F_j$ it follows that $F\cap X[r_{k+1}]=G\cap X[r_{k+1}]$ and $E\cap F\subseteq X[r_{k+1}]$. In this way, $$E\cap F=E\cap F\cap X[R_{k+1}]=E\cap G\cap X[r_{k+1}].$$
Note that $E\cap G\sqsubseteq G$ because $E,G\in \mathcal{F}^i_l$. Furthermore, since $X[r_{k+1}]\sqsubseteq F_j$ and $F\subseteq F_j$, then $F\cap X[r_{k+1}]\sqsubseteq F$. Therefore,  $$E\cap F=E\cap G\cap X[r_{k+1}]\sqsubseteq G\cap X[r_{k+1}]=F\cap X[r_{k+1}]\sqsubseteq  F.$$
\end{claimproof}
\begin{claimproof}[Proof of $(d)$] Let $F\in \mathcal{F}_l$. If $l=k+1$ then $F=X$. In this case, it is easy to see that the sequence $\langle F_i\rangle_{i<n_{k+1}}$ satisfies all the required properties. On the other hand, if $l<k+1$ then $F\in \mathcal{F}^i$ for some $i<n_{k+1}$. Therefore, the existence of the canonical decomposition of $F$ is assured because $\mathcal{F}^i$ is already a construction scheme.  
\end{claimproof}
\end{proof}
\end{proposition}

\begin{definition}[Unique scheme (finite version)] If $k\in\omega$ and $X$ is a finite set of ordinals of size $m_k$, we call $\mathcal{F}(X)$ the unique construction scheme over $X$ of type $\tau.$
\end{definition}

\begin{rem}If $\mathcal{F}$ is a construction scheme over $X$ and $F\in \mathcal{F}$, then $\mathcal{F}|_F=\mathcal{F}(F).$
\end{rem}
\begin{rem}\label{remarkincreasingbijectionisofinite}If $X$ and $Y$ are finite sets of cardinality $m_k$ and $\phi:X\longrightarrow Y$ is the increasing bijection, then $\mathcal{F}(Y)=\{\,\phi[F]\,:F\in \mathcal{F}(X)\}$.
\end{rem}
Now, we will show that $\omega$ also admits a unique construction scheme (of type $\tau$). The uniqueness of the scheme  will be proved using the following result. 
\begin{lemma}\label{lemmarestrictionschemes}Suppose that $X$ and $Y$ are two sets of ordinals for which $X\subseteq Y$. Let $\mathcal{F}$ and $\mathcal{G}$ be construction schemes over $X$ and $Y$ respectively. If $\mathcal{F}\subseteq \mathcal{G}$ then $\mathcal{F}=\mathcal{G}|_X$.
\begin{proof}By definition we have that $\mathcal{F}\subseteq \mathcal{G}|_X$. We proceed to prove the remaining inclusion. For this, let $G\in \mathcal{G}|_X$, $k=\rho^G_\mathcal{G}$ and $\gamma=\max(G)$. By virtue of Proposition \ref{unionlevelscheme}, there is $F\in \mathcal{F}_k$ such that $\gamma\in F$. By the hypotheses, $F\in \mathcal{G}$. Furthermore, $F\in \mathcal{G}_k$ because $|F|=m_k$. In this way, $F\cap G\sqsubseteq G$ by means of the property (c) of Definition \ref{constructionschemedef}. Since $\gamma\in F\cap G$, this means that $F\cap G= G$. Finally, since $|F|=m_k=|G|$, we conclude that $F=G$. In other words, $G\in \mathcal{F}$. 
\end{proof}
\end{lemma}
\begin{rem}If $X$, $Y$, $\mathcal{F}$
 and $\mathcal{G}$ are as in the previous lemma, then $\rho_\mathcal{F}=\rho_\mathcal{G}|_{X^2}$.    
\end{rem}

\begin{lemma}\label{unionschemelemmadeeper}Let $\mathcal{B}$ be a family of sets of ordinals which is totally ordered with respect to $\subseteq$.  Suppose that $\langle \mathcal{F}^X\rangle_{X\in \mathcal{B}}$ is a sequence so that for any $X,Y\in \mathcal{B}$ with $X\subseteq Y$, the following conditions hold:
\begin{itemize}
    \item $\mathcal{F}^X$ is a construction scheme over $X$,
    \item $\mathcal{F}^Y|_X=\mathcal{F}^X$ ( or equivalently, $\mathcal{F}^X\subseteq \mathcal{F}^Y$).
\end{itemize}
Then $\mathcal{F}=\bigcup\limits_{X\in \mathcal{B}} \mathcal{F}^X$ is a construction scheme over $Z=\bigcup \mathcal{B}$.
\begin{proof}$\mathcal{F}$ satisfies the point (a) of Definition \ref{constructionschemedef} because $\mathcal{B}$ is totally ordered. In order to finish,  note that $\mathcal{F}_k=\bigcup\limits_{X\in \mathcal{B}}\mathcal{F}^X_k$ for any $k\in\omega$. The remaining points of Definition \ref{constructionschemedef} hold for $\mathcal{F}$ as a consequence of this observation.
\end{proof}
\end{lemma}
\begin{lemma}\label{lemmafiniteinterval}Let $X$ be a set of ordinals, $\mathcal{F}$ be a construction scheme over $X$ and $k<l\in\omega$. If $F\in \mathcal{F}_l$ and $\alpha\in F$, then there is (a unique) $E\in \mathcal{F}_k$ such that:\begin{enumerate}[label=$(\alph*)$]
    \item $E\subseteq F,$
    \item $\alpha\in E$,
    \item $E\backslash \alpha$ is an interval in $F.$
\end{enumerate}
\begin{proof}The uniqueness of $E$ holds because $E\cap \alpha=(\alpha)^-_k$ and $E\backslash \alpha$ is the unique interval of $F$ of size $m_k-\lVert \alpha\rVert_k$ starting at $\alpha$. The proof of the existence carried by induction over $l$ starting at $k+1$.\\\\
\underline{Base step}: Suppose that $l=k+1$. In the case where $\alpha\in R(F)$ let $E=F_0$. Otherwise, let $E=F_{\Xi_\alpha(l)}$. It is straightforward that $E$ satisfies all the required conditions in both cases.\\\\
\underline{Inductive step}: Suppose that we have proved the lemma for some $l>k$. We will now prove it for $l+1$. So let $F\in \mathcal{F}_{l+1}$ be such that $\alpha\in F$. First we deal with the case where $\alpha\in R(F)$. Note that $\alpha\in F_0$, so according to the inductive hypotheses, there is $E\in \mathcal{F}_k$ such that $E\subseteq F_0$, $\alpha\in E$ and $E\backslash \alpha$ is an interval in $F_0$. Since $F_0$ is already an interval in $F$, then $E\backslash \alpha$ is also an interval. Now we deal we the case where $\alpha\not \in R(F)$. In here, let $i=\Xi_\alpha(l)$. Then $\alpha\in F_i\backslash R(F)$. Again, by the inductive hypotheses, there is $E\in \mathcal{F}_k$ such that $E\subseteq F_i$, $\alpha\in E$ and $E\backslash \alpha$ is an interval in $F_i$. Observe that $F_i\backslash R(F)$ is an interval in $F$ and $E\backslash \alpha$ is contained in it. Therefore, $E\backslash \alpha$ must be an interval in $F$. This finishes the proof.
\end{proof}
\end{lemma}

\begin{corollary}\label{corollaryreductionunique} Let $\mathcal{F}$ be a construction scheme over $X$. Suppose that $\alpha\in X$ and $k\in \omega$. Then there is (a unique) $E\in \mathcal{F}_k$ such that $\alpha\in E$ and $E\backslash \alpha$ is an interval in $X$.
\begin {proof} Let $m= m_k-\lVert \alpha \rVert_k$. Note that $|\,[\alpha,\infty)_X|\geq m$. This is because there is at least one $F'\in \mathcal{F}$ with $\alpha\in F'$. For such $F'$ we have that $\alpha\cap F'=(\alpha)^-_k$. Hence, $m=|F\backslash \alpha|\leq |X|$.

Since $X$ is a well-ordered set and $|[\alpha,\infty)_X|\geq m$, we can take $P$ the unique interval in $[\alpha,\infty)_X$ of size $m$ with $\alpha\in P$. Let $l=\max(k+1,\rho^P)$. Then there is $F\in \mathcal{F}_l$ such that $P\subseteq F$. According to Lemma \ref{lemmafiniteinterval}, there is $E\in \mathcal{F}_k$ such that $E\subseteq F$, $\alpha\in F$ and $E\backslash \alpha$ is an interval in $F$. Note that both $E\backslash \alpha$ and $P $ are intervals in $F$ of size $m$ having $\alpha$ as their minimum. In this way, $E\backslash \alpha=P$. Therefore, $E\backslash \alpha$ is an interval in $X$. This finishes the proof.
\end{proof}
\end{corollary}

\begin{rem}\label{remarkintervalreduction}If we apply the previous corollary to the case where $X$ is an ordinal, say $\gamma$,  we get the following: For any $\alpha\in \gamma$ and $k\in\omega$, we have that $$(\alpha)^-_k\cup [\alpha,\alpha+(m_k-\lVert \alpha \rVert_k)\,)$$
is an element of $\mathcal{F}_k.$
\end{rem}

\begin{theorem}\label{uniqueschemeomegatheorem} There is a unique construction scheme over $\omega$ of type $\tau.$ Furthermore, $\{m_k\,:\,k\in\omega\}$ is contained in such scheme.

\begin{proof}\begin{claimproof}[Proof of existence]For any $k\in \omega$, let $\mathcal{F}(m_k)$ be the only construction scheme over $m_k$. Note that $m_k$ is the first element of the canonical decomposition of $m_{k+1}$ inside $\mathcal{F}(m_{k+1})$. This means that $m_k\in \mathcal{F}(m_{k+1})$. Therefore, $\mathcal{F}(m_k)=\mathcal{F}(m_{k+1})|_{m_k}$. In this way, by Lemma \ref{unionschemelemmadeeper}, the family $\mathcal{F}=\bigcup\limits_{k\in\omega}\mathcal{F}(m_k)$ is a construction scheme over $\omega=\bigcup\limits_{k\in\omega}m_k$.
\end{claimproof}
\begin{claimproof}[Proof of uniqueness] Suppose that $\mathcal{F}'$ is another construction scheme over $\omega$. Let $k\in\omega$. According to the Remark \ref{remarkintervalreduction}, the set $$(0)_k^-\cup [0, 0-(m_k-\lVert 0\rVert_k)\,)=m_k$$
Is an element of $\mathcal{F}'$. Thus, $\mathcal{F}(m_k)=\mathcal{F}'|_{m_k}\subseteq \mathcal{F}'$. Since this is true for any $k$, we have that $\mathcal{F}\subseteq \mathcal{F}'$. But then $\mathcal{F}=\mathcal{F}'$ due to Lemma \ref{lemmarestrictionschemes}.
\end{claimproof}
\end{proof}
\end{theorem}
\begin{definition}[Unique scheme over $\omega$] We will call $\mathcal{F}(\omega)$ the only construction scheme over $\omega$ of type $\tau$.  
\end{definition}

In order to build construction schemes over $\omega_1$, we will use the same technique as when building $\mathcal{F}(\omega)$. That is, we will define such schemes by taking unions of schemes over distinct sets of countable sets of ordinals. In the next two propositions, we show that there is essentially one way of choosing those countable sets.
\begin{proposition}\label{restrictionlimitordinalscheme} Let $\mathcal{F}$ be a construction scheme over an ordinal $\gamma$ of type $\tau$. If $\delta\leq \gamma$ is a limit ordinal, then $\mathcal{F}|_\delta$ is a construction scheme over $\delta$ of type $\tau$. In particular, $\mathcal{F}(\omega)\subseteq \mathcal{F}$ whenever $\gamma$ is infinite.
\begin{proof} It is straightforward that $\mathcal{F}|_\delta$ satisfies conditions (b), (c) and (d) of Definition \ref{constructionschemedef}. In order to prove that condition (a) also holds, take an arbitrary $A\in \text{FIN} (X)$ and let $\alpha=\max(A)$. We need to prove that there is $F\in \mathcal{F}|_\delta$ such that $A\subseteq F$. Indeed, according to the Remark \ref{remarkintervalreduction}, we know that $$F=(\alpha)_k^-\cup [\alpha,\alpha+(m_k-\lVert \alpha\rVert_k)\,)$$
is an element of $\mathcal{F}$ for $k=\rho^A$. Since $\delta$ is a limit ordinal it follows that $F\subseteq \delta$ (that is, $F\in \mathcal{F}|_\delta$). Furthermore, $A\subseteq (\alpha)_k$ so $A\subseteq F$. 
\end{proof}  
\end{proposition}
\begin{proposition}\label{restrictionlimitordinalscheme2}Let $\mathcal{F}$ be a construction scheme over an ordinal $\gamma$ of type $\tau$. If $X\subseteq \gamma$ is an infinite subset of $\gamma$ and $\mathcal{F}|_X$ is a construction scheme over $X$, then $X$ is an ordinal. Furthermore, if $\tau$ is a good type then $X$ is limit.
\begin{proof}In order to show that $X$ is an ordinal, it is enough to prove that $X$ is an initial segment of $\gamma$. For this, let $\alpha<\beta\in \gamma$  be such that $\beta\in X$. Consider $k=\rho(\alpha,\beta)$. Since $X$ is infinite, there is $F\in(\mathcal{F}|_X)_k$ such that $\beta\in F$. This is due to Proposition \ref{unionlevelscheme}. According to the Remark \ref{remarklevelrestriction}, $F\in \mathcal{F}_k$. Therefore, $\alpha\in F$. In particular $\alpha\in X$, so we are done.

Now, suppose that $\tau$ is a good type.  We will prove that $X$ is limit. Assume towards a contradiction that this is not the case. Let $\alpha=\max(X)$ and fix $k>\rho(0,\alpha)$ for which $r_k=0$. Observe that if $F\in \mathcal{F}_k$  is such that $\alpha\in F$, then $\alpha=\max(F)$. In particular,  $\alpha\in F_{n_k-1}$. Furthermore, since $k>\rho(0,\alpha)$ then $0$ also belongs to $F_{n_{k-1}}.$ Lastly, $F_0<F_{n_{k-1}}$ because $r_k=0$. This implies that every element of $F_0$ is smaller than $0$, which is a contradiction. Thus, the proof is over.
\end{proof}   
\end{proposition}

In summary, by means of the two previous propositions we know that if $\tau$ is a good type and $\mathcal{F}$  is a construction scheme over and ordinal $\gamma$ then $\gamma$ is either a $m_k$ for some $k\in \omega$, or $\gamma$ is limit. Furthermore, the set of all $X\in \mathscr{P}(\gamma)$ for which $\mathcal{F}|_X$ is a construction scheme can be fully described as $$\mathcal{F}\cup(\text{Lim}\cap (\gamma+1)).$$ 
In the case where $\gamma$ is infinite, the previous family has $\omega$ as a special member. This is because, according to the Theorem \ref{uniqueschemeomegatheorem}, $\mathcal{F}|_\omega$ is simply $\mathcal{F}(\omega)$. This scheme in $\omega$ stands out not only because of its uniqueness but because each $m_k$ belongs to it. Therefore, $m_k$ not only acts as the common cardinality of all elements in $\mathcal{F}_k$,  but it is the \say{most canonical} member of such family. Actually, this line of thought transfers to all natural numbers. Concretely, since the initial segments of closed sets are closed according to Corollary \ref{initialsegmenteclosedisclosed}, we have the following corollary.
\begin{corollary}Let $\mathcal{F}$ be a construction scheme over an infinite ordinal $\gamma$. Then $n\in \omega$ is a closed set with respect to $\rho$.  Furthermore, $n\in(m_{\rho^n-1}, m_{\rho^n}]$ whenever $n>0$.
\end{corollary}
The metric structure of closed sets is fully determined by their cardinality (see Lemma \ref{ballhomogeneitylemma} and Theorem \ref{closedsizetheorem}). In particular, for each $\beta\in \gamma$ and $l\in \omega$, the set $(\beta)_l$ is $\rho$-isomorphic to $\lVert\beta\rVert_l+1=[0,\lVert \beta\rVert_l]$ (recall $\lVert \beta\rVert_l=|(\beta)^-_l|$). It will be convenient to name and explicitly describe the increasing bijections between these two sets. 
\begin{definition}\label{functionphibetaldef}Let $\mathcal{F}$ be a construction scheme over a limit ordinal $\gamma$. Given $\beta\in \gamma$ and $l\in \omega$, we define $\phi^\beta_l:(\beta)_l\longrightarrow \lVert \beta\rVert_l+1$ as:
$$\phi^\beta_l(\alpha)=\lVert \alpha \rVert_l.$$
The inverse function of $\phi^\beta_l$ is given by:
$$(\phi^\beta_l)^{-1}(i)=(\beta)_l(i).$$    
\end{definition}

\begin{proposition}\label{kcardinalityisomorphismprop}Let $\mathcal{F}$ be a construction scheme  over an ordinal $\gamma$. Given $\alpha,\beta\in \gamma$ and $k\leq l\in \omega$, the following conditions hold:
\begin{enumerate}[label=$(\alph*)$ ]
    \item $\lVert \beta\rVert_k=\lVert\,\lVert\beta\rVert_l\,\rVert_k,$
    \item If $l\geq \rho(\alpha,\beta)$, then $\rho(\alpha,\beta)=\rho(\lVert\alpha\rVert_l,\lVert\beta\rVert_l)$
    \item If $l\geq\Delta(\alpha,\beta)$, then $\Delta(\alpha,\beta)=\Delta(\lVert \alpha\rVert_l,\lVert \beta\rVert_l).$
\end{enumerate}
\end{proposition}
We can rewrite the point (a) of the previous proposition in a more useful way.
\begin{lemma}\label{compositionphi}Let $\mathcal{F}$ be a construction scheme over an ordinal $\gamma$. Given $\beta\in \gamma$ and $k\leq l\in \omega$, we have that $$\phi^{\lVert \beta\rVert_l}_{k}\circ\phi^\beta_l=\phi^\beta_k.$$
     
\end{lemma}
    
\section{The forcing $\mathbb{P}(\mathcal{F})$ and schemes in ZFC}\label{sectionpfforcing}
In this section we will present an incarnation of the Cohen forcing which was first considered in \cite{schemenonseparablestructures}.  This forcing will serve the purpose of extending a construction scheme from a countable limit ordinal to the next one, and will play a key roll in the arguments of the next sections of this chapter. As an application, we will show in subsection \ref{subsectionschemeszfc} that there is a construction scheme over $\omega_1$ of any given type. All of this section is based in \cite{schemenonseparablestructures}.\\

For the rest of this section, let us fix a construction scheme $\mathcal{F}$ of type $\tau$ over a limit ordinal $\gamma$. 
\begin{definition}[The forcing $\mathbb{P}(\mathcal{F})$]\label{forcingPFdef} We define $\mathbb{P}(\mathcal{F})$ as the forcing consisting of the empty set and  of all $p\in \text{FIN}(\gamma+\omega)$ with the following properties:
\begin{enumerate}[label=(\Roman*)]
\item There is $k_p\in \omega$ such that $|p|=m_{k_p}$.
\item There is $F\in \mathcal{F}_{k_p}$ such that $p\cap \gamma\sqsubseteq F$.
\item $p\cap [\gamma,\gamma+\omega)$ is an initial segment of $[\gamma,\gamma+\omega)$.
\end{enumerate}
Whenever $p\cap\gamma\not=\emptyset$ (even if $p$ is not a condition of $\mathbb{P}(\mathcal{F})$), we let $\alpha_p=\max(p\cap \gamma)$. Additionally, for each $k\in\omega$ we let $\mathbb{P}_k(\mathcal{F})=\{\,p\in\mathbb{P}(\mathcal{F})\,:\,k_p=k\, \}.$  The order on $\mathbb{P}(\mathcal{F})$ is given by $$p\leq q\text{ if and only if }q\in \mathcal{F}(p)\textit{ or }q=\emptyset.$$
\end{definition}
Note that $\mathbb{P}(\mathcal{F})$ is always countable. Therefore, it is forcing equivalent to the Cohen forcing. 

\begin{center}
\begin{minipage}[b]{0.8\linewidth}
\centering
\includegraphics[width=12 cm, height=5.5cm]{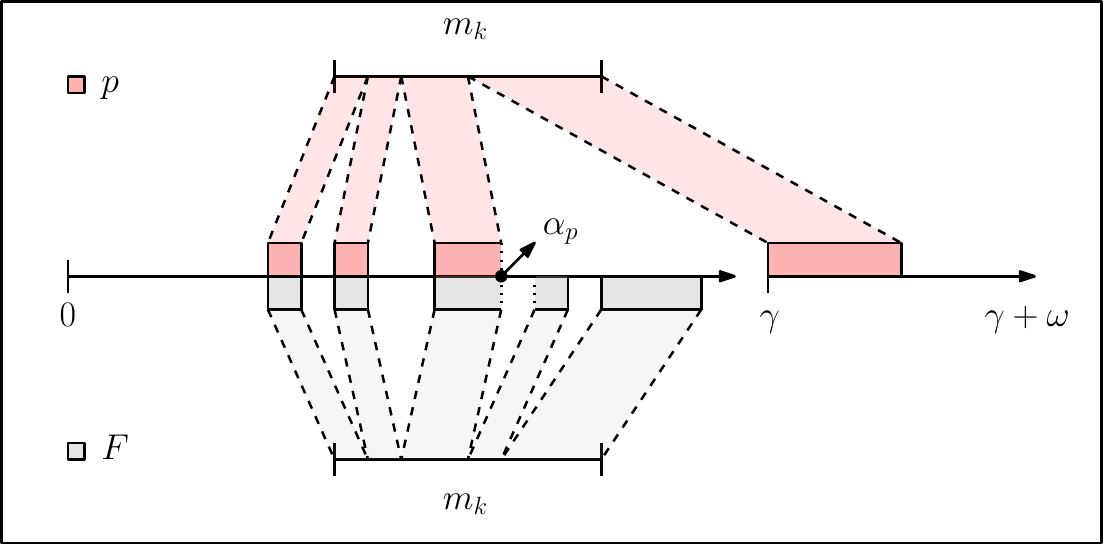}

\textit{\small In here, $p$ represents an element of $\mathbb{P}_k(\mathcal{F})$ and $F$ represents an element of $\mathcal{F}_k$ which testifies the condition (II) of Definition \ref{forcingPFdef} for $p$. }
\end{minipage}
\end{center}

\begin{rem}\label{remarkconditiontwoforcingequivalence} From Lemma \ref{closureschemelemma} it follows that the condition (II) of Definition \ref{forcingPFdef} is equivalent to: \begin{center}(II)$^*\:\:p\cap \gamma=(\alpha_p)_{k_p}$,
\end{center}
in the case where $p\cap \gamma\not=\emptyset.$
\end{rem}
In general, if $p\in \mathbb{P}(\mathcal{F})$, there are many elements $F\in \mathcal{F}_{k_p}$ testifying the condition (II) of Definition \ref{forcingPFdef} with respect to $p$. However, it turns out that there is a canonical one. Such $F$ will be useful in many of the arguments involving the forcing $\mathbb{P}(\mathcal{F})$. In the following definition, we explicitly describe it.
\begin{definition}[The reduction operation]If $p\in \text{FIN}(\gamma+\omega)$ and $\delta\leq \gamma$, we define the \textit{reduction of $p$ to $\delta$} as follows: $$red_\delta(p)=\begin{cases}(p\cap \delta)\cup[\max(p\cap \delta)+1,\max(p\cap \delta)+|p\backslash \delta|+1)&\textit{ if }p\cap\delta\not=\emptyset\\
|p|&\textit{ if }p\cap\delta=\emptyset
\end{cases}
$$
\end{definition}
\begin{rem}\label{remark1reductiondef}The reduction operation is closely related to the Corollary \ref{corollaryreductionunique}. Suppose that we take $\alpha\in \gamma$ and $k\in\omega$. According to the corollary, we know that there is a unique $E\in \mathcal{F}_k$ so that $\alpha\in E$ and $E\backslash \alpha$ is an interval in $\gamma$. It turns out that such $E$ can be described using the reduction operation. In order to do this, take an arbitrary $F\in \mathcal{F}_k$ with $\alpha\in F$. Then: \begin{align*}red_{\alpha+1}(F)&=F\cap(\alpha+1)\cup [\,\alpha+1,\,\alpha+1+|F\backslash (\alpha+1)|\,)\\
&=(\alpha)^-_k\cup [\,\alpha,\, \alpha +1+ (m_k-(\lVert \alpha\rVert_k+1))\,)\\
&=(\alpha)^-_k\cup [\,\alpha,\,\alpha+m_k-\lVert \alpha \rVert\,).
\end{align*}
Due to the Remark \ref{remarkintervalreduction}, it is straightforward that $E=red_{\alpha+1}(F)$.
\end{rem}
The reduction operation will be mainly used in the case where $p\in \mathbb{P}(\mathcal{F})$ and $\alpha=\gamma$. We now illustrate such situation for convenience of the  reader.

\begin{center}
\begin{minipage}[b]{0.8\linewidth}
\centering
\includegraphics[width=12 cm, height=5.5cm]{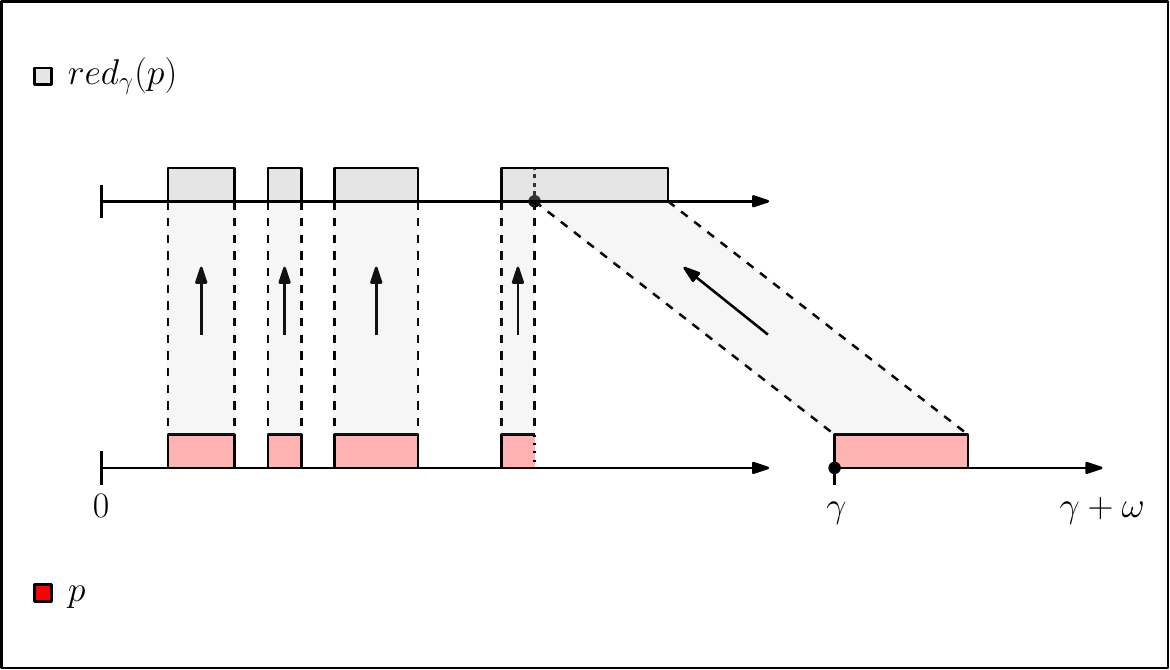}

\end{minipage}
\end{center}

By arguing in a similar manner as in the Remark \ref{remark1reductiondef}, we get the following lemma. The proof of it is left to the reader.
\begin{lemma}\label{equivalencereductioncondition} If $p\in \text{FIN}(\gamma+\omega)$ is such that $p\cap[\gamma,\gamma+\omega)$ is an initial segment of $[\gamma,\gamma+\omega)$, then the following conditions are equivalent:
\begin{enumerate}[label=$(\alph*)$]
\item $p\in \mathbb{P}(\mathcal{F})$,
\item $red_\gamma(p)\in \mathcal{F}.$
\end{enumerate}
\end{lemma}
Even though the previous lemma is useful, the main tool for defining and extending conditions in $\mathbb{P}(\mathcal{F})$ relies in the next definition.
\begin{definition}[The cut operation] Let $F\in \text{FIN}(\gamma)$ and $\alpha\in \gamma$. We define the \textit{cut of $F$ at $\alpha$} as follows: $$Cut_\alpha(F)=(F\cap \alpha)\cup [\gamma,\gamma+|F\backslash\alpha|)$$
\end{definition}
The following picture illustrates the previous definition for the particular case in which $\alpha\in F$.  This is the main situation in which the cut operation will be applied.

\begin{center}
\begin{minipage}[b]{0.8\linewidth}
\centering
\includegraphics[width=12cm, height=5.5cm]{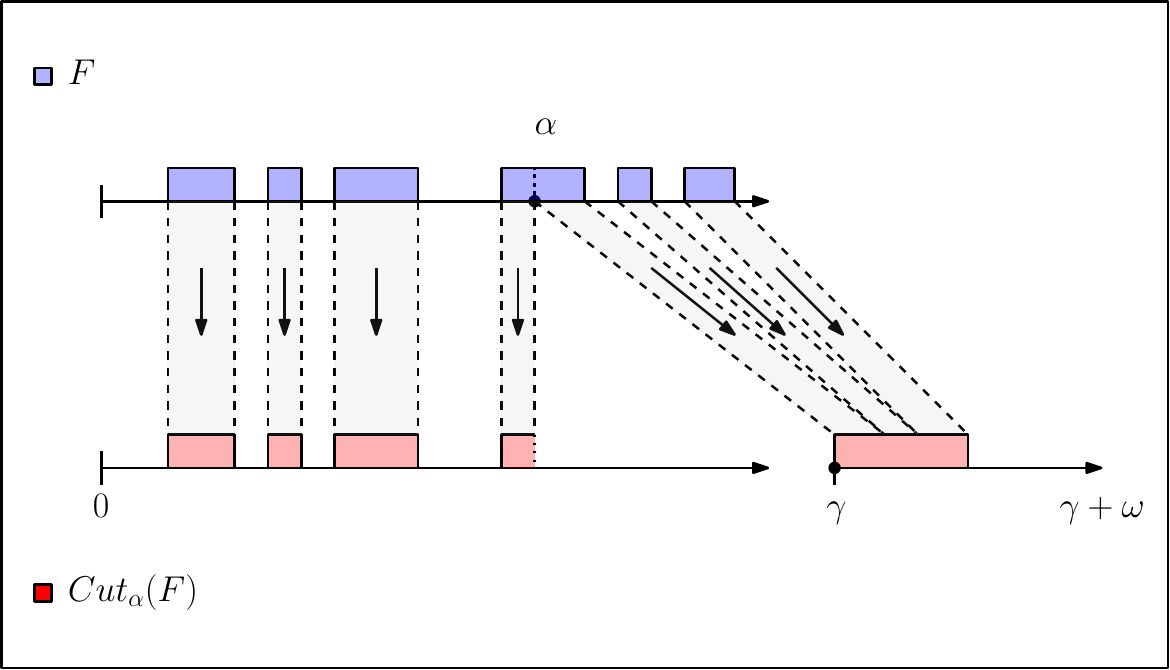}

\end{minipage}
\end{center}

\begin{rem}\label{remarkcutcannonicaldecomposition} Suppose that $F\in \mathcal{F}$ and $\alpha\in \gamma$. If $\phi:F\longrightarrow Cut_\alpha(F)$ is the increasing bijection, then $$\mathcal{F}(Cut_\alpha(F))=\{\, \phi[H]\,:\,H\in \mathcal{F}(F)\,\}$$
due to the Remark \ref{remarkincreasingbijectionisofinite}. In particular, this means  that $Cut_\alpha(F_i)$ is the $ith$ element of the canonical decomposition of $Cut_\alpha(F)$ for each $i<n_{\rho^F}$. In other words,   $Cut_\alpha(F_i)=(Cut_\alpha(F))_i.$ On the other hand, it is not necessarily true that $\phi[H]=Cut_\alpha(H)$ for a given $H\in \mathcal{F}(F)$. A particular case in which the previous equality holds is when $\alpha\in H$ and $H\backslash \alpha$ is an interval in $F$ (Thus, an initial segment of $F\backslash \alpha$).
    
\end{rem}
\begin{rem}\label{remarkcutreductioninverse}The reduction and cut operations are in some sense inverses of each other. That is, if $p\in \mathbb{P}(\mathcal{F})$ and $p\cap \gamma\not=\emptyset$, then 
$Cut_{\alpha_p+1}(red_\gamma(p))=p$. On the other hand, if $F\in \text{FIN}(\gamma)$ and $\alpha\in F$ then $red_\gamma(Cut_{\alpha+1}(F))=F$.
\end{rem}

\begin{lemma}\label{cutconditionlemma}Let $F\in \mathcal{F}$ and $\alpha\in \gamma$. Then $Cut_\alpha(F)\in \mathbb{P}(\mathcal{F})$. Furthermore, if $F,G\in \mathcal{F}$ are such that $F\subseteq G$ and $\alpha\in F$ then $Cut_\alpha(G)\leq Cut_\alpha(F)$.
\begin{proof}For the first part of the lemma, suppose that $F\in \mathcal{F}$ and $\alpha\in \gamma$. Let $k=\rho^F$. Then $m_k=|F|=|Cut_\alpha(F)|$. In this way, $Cut_\alpha(F)$ satisfies the condition (I) of Definition \ref{forcingPFdef}. Now, $Cut_\alpha(F)\cap \gamma= F\cap \alpha\sqsubseteq F$. Therefore, the condition (II) of the same definition  holds for $Cut_\alpha(F)$. Finally,  $Cut_\alpha(F)\cap [\,\gamma,\,\gamma+\omega\,)=[\,\gamma,\,\gamma+|F\backslash \alpha|\,)$. Thus, $Cut_\alpha(F)$ satisfies the condition (III) of Definition \ref{forcingPFdef}, which means that $Cut_\alpha(F)\in \mathbb{P}(\mathcal{F}).$

Next, we prove the second part of the lemma. Let $F,G\in \mathcal{F}$ be such that $F\subseteq G$ and let $\alpha\in F$. If $F=G$, the result is obvious. So we may assume that the inclusion between $F$ and $G$ is proper. It follows that $k<l$ where $k=\rho^F$ and $l=\rho^G$. Therefore, we are in the conditions of applying the lemma \ref{lemmafiniteinterval} to $\alpha, G$ and $k$. In this way, we get $E\in \mathcal{F}_k$ for which $\alpha\in E$, $E\subseteq G$ and $E\backslash \alpha$ is an interval in $G$. Note that $|E\backslash \alpha|=|F\backslash \alpha|$ and $F\cap \alpha=(\alpha)^-_k=E\cap \alpha$. Thus, $Cut_\alpha(F)=Cut_\alpha(E)$. Now, let $\phi:G\longrightarrow Cut_\alpha(G)$ be the increasing bijection. Due to the Remark \ref{remarkcutcannonicaldecomposition}, $$\mathcal{F}(Cut_\alpha(G))=\{\, \phi[H]\,:\,H\in \mathcal{F}(G)\,\}.$$
Finally, just note that $\phi[E]=Cut_\alpha(E)$ by virtue of the same remark.
\end{proof}
\end{lemma}
\begin{lemma}\label{lemmacut}Let $k\in\omega$, $p\in \mathbb{P}_k(\mathcal{F})$ (that is, $k=k_p$) and $\alpha\in \gamma$ be such that $(\alpha)^-_k=p\cap\gamma$. Then $Cut_\alpha(G)\leq p$ for each $G\in \mathcal{F}$ with $\alpha\in G$ and $\rho^G\geq k.$
\begin{proof}Let $G\in \mathcal{F}$ be as in the hypotheses. Consider $F\in \mathcal{F}_k$ for which $\alpha\in F$ and $F\subseteq G.$ According to the Lemma \ref{cutconditionlemma}, $Cut_\alpha(G)\leq Cut_\alpha(F)$. To finish, just note that $|Cut_\alpha(F)|=|p|$ and $Cut_\alpha(G)\cap \gamma=G\cap \alpha=(\alpha)^-_k=p\cap \gamma.$  In this way,  $Cut_\alpha(F)=p$ due to the condition (III) of Definition \ref{forcingPFdef}. 

\end{proof}
\end{lemma}
\subsection{Construction schemes in ZFC}\label{subsectionschemeszfc}
The next definition appeared for the first time in \cite{schemenonseparablestructures}. This hypothesis is sufficient to build construction schemes in a recursive manner.
\begin{definition}[The $IH_1$ property]Let $A\in \text{FIN}(\gamma)$,  $\alpha\in \gamma$ and $F\in \mathcal{F}$. We say that  $IH_1(\alpha,A,F)$ holds if:
\begin{enumerate}[label=$(\arabic*)$]
    \item $A\subseteq F_0,$
    \item $R(F)=F\cap \alpha$.
\end{enumerate}
Additionally, we say that $\mathcal{F}$ satisfies $IH_1$ if for all $A\in\text{FIN}(\gamma)$ and $\alpha\in \gamma$, there is $F\in \mathcal{F}$ for which $IH_1(\alpha,A,F)$ holds.
\end{definition}
\begin{proposition}\label{IH1Fomega}Suppose that $\tau$ is a good type. Then $\mathcal{F}(\omega)$ satisfies $IH_1$.
\begin{proof}Let $A\in \text{FIN}(\gamma)$ and $\alpha\in \omega$. Since $\tau$ is a good type, there is $k>\max(A)$ for which $r_{k+1}=\alpha$. Let $F=m_{k+1}$. Then $F\in \mathcal{F}(\omega)$, $F_0=m_k\supseteq A$ and $R(F)=r_{k+1}=\alpha$. This finishes the proof.
\end{proof}
\end{proposition}
The following lemma is easy. 
\begin{lemma}\label{unionschemelemmadeeperIH1}Suppose that $\gamma$ is a limit of limit ordinals and $\mathcal{F}|_\delta$ satisfies $IH_1$ for each limit $\delta<\gamma$. Then $\mathcal{F}$ also satisfies $IH_1$.
 \end{lemma}

\begin{definition} Given a filter $\mathcal{G}$ over $\mathbb{P}(\mathcal{F})$, we define  $\mathcal{F}^\mathcal{G}$ as $\bigcup_{p\in \mathcal{G}}\mathcal{F}(p).$ Finally, $\mathcal{F}^{Gen}$ denotes the name for $\mathcal{F}^\mathcal{G}$ where $G$ is a generic filter.
\end{definition}
\begin{lemma}\label{nonemptygammacondition}Suppose that $\tau$ is a good type and $\mathcal{F}$ satisfies $IH_1$. Then the set $\{\,p\in \mathbb{P}(\mathcal{F})\,:\,p\cap \gamma\not=\emptyset\,\}$ is
 open dense in $\mathbb{P}(\mathcal{F}).$
\begin{proof}Let $q\in \mathbb{P}(\mathcal{F})$. Without loss of generality we can assume that $q\not=\emptyset$. We need to find $p\leq q$ such that $p\cap \gamma\not=\emptyset$. If $q\cap \gamma\not=\emptyset$, we are done. So suppose that $q\cap \gamma=\emptyset.$ Then $q=[\gamma,\gamma+m_{k_q})$ due to the condition (III) of Definition \ref{forcingPFdef}. Since $\tau$ is a good type, then there is $k>k_q$ such that  $r_{k+1}=0$. Let $F=m_{k+1}$ and $\alpha=m_k$. Then $F\in \mathcal{F}(\omega)\subseteq \mathcal{F}$. By means of the \ref{cutconditionlemma}, $$m_k\cup [\gamma, \gamma +m_k+1-m_k)=Cut_\alpha(F)\in\mathbb{P}(\mathcal{F}).$$
Let $p=Cut_\alpha(F)$. Note that the second piece of the canonical decomposition of $p$, that is $p_1$, is equal to $[\gamma,\gamma+m_k)$. From this, it is easy to see that $[\gamma,\gamma+m_j)\in\mathcal{F}(p_1)$ for each $j\leq k$. In particular, $q\in \mathcal{F}(p_1)\subseteq \mathcal{F}(p)$. In other words, $p\leq q$. Since $p\cap \gamma=m_k\not=\emptyset$, the proof is over.
\end{proof}
\end{lemma}
\begin{lemma}\label{upslemmadups}Suppose that $\tau$ is a good type and $\mathcal{F}$ satisfies $IH_1$. Then the set $\{\,p\in \mathbb{P}(\mathcal{F})\,:\,p\cap [\gamma,\gamma+\omega)\not=\emptyset\,\}$ is open dense in $\mathbb{P}(\mathcal{F}).$
\begin{proof}Let $q\in \mathbb{P}(\mathcal{F})$. Without loss of generality we may assume that $q\not=\emptyset$. We need to find $p\leq q$ with $p\backslash \gamma\not=\emptyset$. If $q\backslash \gamma\not=\emptyset$, there is nothing to do. In this way, we may assume that $q\subseteq \gamma$. In this case, it follows that $q\in \mathcal{F}$. Let $\alpha=\alpha_p+1=\max(q)+1.$ Since $\mathcal{F}$ satisfies $IH_1$, there is $F\in \mathcal{F}$ with $q\cup\{\alpha\}\subseteq F_0$ and $\alpha\cap F=R(F)$. Note that $q$ is actually a subset of $R(F)$. Let $p=Cut_\alpha(F)$. According to the Lemma \ref{cutconditionlemma}, $p=Cut_\alpha(F)\leq Cut_\alpha(q)=q$. Furthermore, $p\backslash \gamma\not=\emptyset$. This finishes the proof.
\end{proof}
\end{lemma}

\begin{lemma}\label{densedomainlemma}Suppose that $\mathcal{F}$ is a construction scheme over $\gamma$ satisfying $IH_1$. For any $\alpha\in \gamma+\omega$, the set $$\mathcal{D}'_\alpha=\{p\in \mathbb{P}(\mathcal{F})\,:\,\alpha\in p\}$$
is open dense in $\mathbb{P}(\mathcal{F}).$
\begin{proof}Let $q\in \mathbb{P}(\mathcal{F})$. By Lemmas \ref{nonemptygammacondition} and \ref{upslemmadups}, we may assume that $q\cap \gamma\not=\emptyset$ and $q\backslash \gamma\not=\emptyset$. In this way, $\alpha_q$ is well defined and $\alpha_q+1\in F$ where $F=red_\gamma(q)$. In particular, $(\alpha_q+1)^-_{k_q}=q\cap \gamma$. Recall that our goal is to find $p\in \mathbb{P}(\mathcal{F})$ with $p\leq q$ and such that $\alpha\in p$. The proof of this fact is divided into two cases.\\\\
\underline{Case 1}: If $\alpha<\gamma$.
\begin{claimproof}[Proof of case]Since $\mathcal{F}$ satisfies $IH_1$, there is $G\in \mathcal{F}$ so that $\{\alpha\}\cup F\subseteq G_0$ and $(\alpha_q+1)\cap G=R(G)$. Let $k=\rho^G$ and $\beta=\min(G_1\backslash R(G))$. Note that $k>k_q$ because $F\subseteq G_0$.  Moreover, $(\beta)^-_{k-1}=R(F)=(\alpha_q+1)^-_{k-1}$.  Since $k_q\leq k$, it follows that $(\beta)^-_{k_q}=(\alpha_q+1)^-_{k_q}=q\cap \gamma$. Let $p=Cut_\beta(G)$. According to the Lemma \ref{lemmacut}, $p\leq q$. Moreover, $p\cap \gamma=G_0$ and this set contains $\alpha$. Hence, the proof of this case is done.
\end{claimproof}
\noindent\underline{Case 2}: If $\alpha\geq \gamma$.
\begin{claimproof}[Proof of case] Let $n=|\alpha\backslash \gamma|+1$. Again, using that $\mathcal{F}$ satisfies $IH_1$, we can get $G\in \mathcal{F}$ so that $F\cup [\alpha_q,\alpha_q+n]\subseteq G_0$ and $(\alpha_q+1)\cap G=R(G)$. By arguing in a similar way as we did in the first case, we may conclude that for $\beta=\min(G_1\backslash R(G))$, $p=Cut_\beta(G)$ is a condition in $\mathbb{P}(\mathcal{F})$ which is smaller than $q$. In order to finish, note that $[\alpha_q,\alpha_q+n]\subseteq G_0\backslash R(G)$. Therefore, $n\leq |G_0\backslash R(G)|=|G_1\backslash R(G)|$. According to the definition of $Cut_\beta(G)$, we get that $|p\cap [\gamma,\gamma+\omega)|\geq n$. Thus, $\alpha\in [\gamma,\gamma+n]\subseteq p$.
\end{claimproof}
\end{proof}
\end{lemma}

\begin{lemma}\label{lemmaFdenseinPF}Suppose that  $\tau$ is a good type and $\mathcal{F}$ satisfies $IH_1$. For any $F\in \mathcal{F}$, the set $$\mathcal{D}_F=\{\,p\in \mathbb{P}(\mathcal{F})\,:\,F\in \mathcal{F}(p)\,\}$$
is open dense in $\mathbb{P}(\mathcal{F}).$
\begin{proof}Let $q\in \mathcal{F}$. By applying Lemma \ref{densedomainlemma} to each element of $F$, we get that there is $p\in \mathbb{P}(\mathcal{F})$ so that $p\leq q$ and $F\subseteq p.$  According to the Lemma \ref{cutconditionlemma}, $Cut_{\alpha_p+1}(red_\gamma(p))\leq Cut_{\alpha_p+1}(F)$. That is, $Cut_{\alpha_p+1}(F)\in \mathcal{F}(Cut_{\alpha_p+1}(red_\gamma(p)))$. Note that $Cut_{\alpha_p+1}(F)=F$ because $F\subseteq \alpha_p+1$. Moreover, $Cut_{\alpha_p+1}(red_\gamma(p))=p$ by virtue of the Remark \ref{remarkcutreductioninverse}. This finishes the proof.
    
\end{proof}
\end{lemma}

\begin{lemma}\label{lemmaih1dense}Suppose that $\tau$ is a good type and $\mathcal{F}$ satisfies $IH_1$. For any $A\in \text{FIN}(\gamma+\omega)$ and each $\alpha\in \gamma+\omega$, the set $$\mathcal{E}_{\alpha,A}=\{\,p\in\mathbb{P}(\mathcal{F})\,:\,\alpha\in p, \, A\subseteq p_0,\textit{ and }\alpha\cap p=R(p)\,\}$$
is dense in $\mathbb{P}(\mathcal{F}).$
\begin{proof}Let $q\in \mathbb{P}$. By means of the previous lemmas, we can assume without loss of generality that $q\cap \gamma\not=\emptyset$, $q\backslash \gamma\not=\emptyset$ and $\{\alpha\}\cup A\subseteq q$. Let $\alpha'=red_\gamma(q)(|\alpha\cap q|)$. That is, $\alpha'$ is the element of $red_\gamma(q)$ which corresponds to $\alpha$ via the increasing bijection from $red_\gamma(q)$ to $q$. Since $\mathcal{F}$ satisfies $IH_1$, there is $F\in \mathcal{F}$ so that $red_\gamma(q)\subseteq F_0$ and $F\cap \alpha'= R(F).$ Let $p=Cut_{\alpha_q+1}(F)$.  We claim that $p$ is an element of $\mathcal{E}_{\alpha,A}$ which is smaller that $q$. Indeed, since $\alpha_p+1\in red_\gamma(q)$, then $(\alpha_q+1)^-_{k_q}=q\cap \gamma$. Furthermore, $\rho^F>k_q$ because $F_0$ contains the reduction of $q$ to $\gamma$. In this way, $p\leq q$ by virtue of the Lemma \ref{lemmacut}. Moreover, by virtue of the Remarks \ref{remarkcutcannonicaldecomposition} and \ref{remarkcutreductioninverse}, the following chain of inclusions hold: $$A\subseteq q=Cut_{\alpha_q+1}(red_\gamma(q))\subseteq Cut_{\alpha_q+1}(F_0)=p_0.$$
Now, let $\phi:F\longrightarrow Cut_{\alpha_q+1}(F)$ be the increasing bijection. Note that $red_\gamma(q)\backslash \alpha_q$ is an initial segment of $F$. According to the last part of the Remark \ref{remarkcutcannonicaldecomposition}, $\phi[red_\gamma(q)]=Cut_{\alpha_q+1}(red_\gamma(q))=q$. In particular, $\phi[\alpha']=\alpha.$  Since $\phi[R(F)]=R(p)$,  we have that $R(p)=\phi[\alpha'\cap F]=\alpha\cap p$.
\end{proof}
\end{lemma}

\begin{proposition}\label{schemeinZFCextension}Suppose that $\tau$ is a good type and $\mathcal{F}$ satisfies $IH_1$. Let $\mathcal{G}$ be a filter over $\mathbb{P}(\mathcal{F})$ intersecting $\mathcal{D}'_\alpha$, $\mathcal{D}_F$ and $\mathcal{E}_{\alpha,A}$ for all $\alpha\in \omega+\gamma$, $A\in \text{FIN}(\omega+\gamma)$ and $F\in \mathcal{F}$\footnote{Such $\mathcal{G}$ exists due to the Rasiowa-Sikorski Lemma}. Then $\mathcal{F}^\mathcal{G}$ is a construction scheme over $\gamma+\omega$ which contains $\mathcal{F}$ and satisfies $IH_1$.
\begin{proof} If we use that $\mathcal{G}$ intersects each $\mathcal{D}'_\alpha$, and argue in a similar fashion as we did in Lemma \ref{unionschemelemmadeeper}, it is easy to see that $\mathcal{F}^G$ is a construction scheme over $\gamma+\omega$. Furthermore, $\mathcal{F}\subseteq\mathcal{F}^\mathcal{G}$ because $\mathcal{G}$ intersects each $\mathcal{D}_\mathcal{F}$, and $\mathcal{F}^\mathcal{G}$ satisfies $IH_1$ because $\mathcal{G}$ intersects each $\mathcal{E}_{\alpha,A}.$
\end{proof}
\end{proposition}
Suppose that $\tau$ is a good type, By means of the results \ref{IH1Fomega}, \ref{schemeinZFCextension}, \ref{unionschemelemmadeeper}, and \ref{unionschemelemmadeeperIH1}, we can recursively construct, for each limit $\gamma<\omega_1$, a construction scheme $\mathcal{F}^\gamma$ which satisfies $IH_1$ and such that $\mathcal{F}^\delta\subseteq \mathcal{F}^\gamma$ whenever $\delta<\gamma$. If we use one more time the Lemmas \ref{unionschemelemmadeeper} and \ref{unionschemelemmadeeperIH1}, we can conclude that $$\mathcal{F}=\bigcup\limits_{\gamma\in Lim}\mathcal{F}^\gamma$$
is a construction scheme over $\omega_1$ which satisfies $IH_1$. Thus, we have proved Theorem \ref{theoremschemesinzfcstevo}.\\\\
\noindent
\textbf{Theorem \ref{theoremschemesinzfcstevo}.}\textit{ For any good type there is a construction scheme over $\omega_1$ of that type which satisfies $IH_1.$}

\section{Diamond principle and FCA(part)}\label{sectiondiamond}
The purpose of this section is to prove that Jensen's $\Diamond$-principle implies $FCA$. In order to do this, we will work in the same manner as we did in the Subsection \ref{subsectionschemeszfc}. That is, we want to find a suitable property $IH_2$ which is satisfied by $\mathcal{F}(\omega)$. Furthermore, given a construction scheme $\mathcal{F}$ over a limit ordinal $\gamma\leq \omega_1$, we want the two following things to happen:
\begin{enumerate}[label=$(\alph*)$]
    \item If $\gamma$ is a limit of limit ordinals and $\mathcal{F}|_\delta$ satisfies $IH_2$ for each limit $\delta<\gamma$, then $\mathcal{F}$ also satisfies $IH_2$.
    \item If  $\gamma$ is countable and $\mathcal{F}$ satisfies  $IH_1$ and $IH_2$ then there is a construction scheme $\mathcal{F}'$ over $\gamma+\omega$ containing $\mathcal{F}$ which satisfies $IH_1$ and $IH_2$.
\end{enumerate}
By doing this, we may conclude that there is a construction scheme over $\omega_1$ which satisfies $IH_2$. Finally, we want that: \begin{enumerate}[label=$(c)$]
   \item Any construction scheme over $\omega_1$ which satisfies $IH_2$ is fully capturing.
\end{enumerate}
While the points (a) and (c) are relatively easy to guarantee, it turns out that finding a property $IH_2$ which also satisfies the point (b) is a non-trivial problem. \\

For the rest of this section, we fix a countable limit ordinal $\gamma$ and a construction scheme $\mathcal{F}$ be over $\gamma$ of type $\tau.$ Moreover, we fix  a $\Diamond$-sequence $\langle D_\alpha\rangle_{\alpha\in \text{Lim}}$.

\begin{definition}[Block interval sequence]Let $j\in \omega$ and $A$ be a non-empty subset of $\omega_1$. We say that a sequence $\mathbb{I}=\langle \mathbb{I}(i)\rangle_{i<j}\subseteq \text{FIN}(A)$ is a \textit{block interval sequence over $A$} if:
\begin{itemize}
    \item For any $i<t$, $\mathbb{I}(i)$ is an interval in $A$.
    \item For all $i<i'<t$, $\mathbb{I}(i)<\mathbb{I}(i')$.
\end{itemize}
Given $\alpha\in \omega_1$, we let $Bl(\alpha,A)$ to be the set of all interval sequences over $A\backslash\alpha.$\footnote{Note that $Bl(\alpha,A)=Bl(0,A\backslash \alpha)$ for all $\alpha\in \omega_1$ and $A\in \text{FIN}(\omega_1).$} Given $\alpha\in X$, we denote $Bl(\,\lVert\alpha\rVert_k,\,m_k)$ simply as $Bl_k(\alpha)$ for each $k\in \omega$. 
\end{definition}
\begin{rem}Note that a block interval sequence $\mathbb{I}$ may be empty even though its elements are not. This case is not pathological, and will be highly relevant for the proof of Theorem \ref{ih2omega1fullycapturing}.
    
\end{rem}
\begin{definition}Suppose that $A$ and $B$ are non-empty subsets of $\omega_1$ and $\phi:A\longrightarrow B$ is an increasing function. Given $\mathbb{I}=\langle \mathbb{I}(i)\rangle_{i<j}\in Bl(0,A)$, we let $$\phi\circ\mathbb{I}=\langle\, \phi[\mathbb{I}(i)] \,\rangle_{i<j}.$$
Recall that if $B$ has size $m$ for some $m\in \omega$, then we identity the increasing bijection $\phi:m\longrightarrow B$ with $B$, and we denote the inverse of such function as $B^{-1}$. Following this convention, if $A=m$ and $\mathbb{I}$ is as before, $\phi\circ \mathbb{I}$ is denoted as $B\circ \mathbb{I}.$ Analogously, if $\mathbb{J}\in Bl(0,B)$,  we denote $\phi^{-1}\circ \mathbb{J}$  as $B^{-1}\circ \mathbb{J}$.
\end{definition}
\begin{rem}Suppose that $A,B,\phi$ and $\mathbb{I}$ are as in the previous definition. Since $\phi$ is increasing, then $\phi\circ \mathbb{I}(i)<\phi\circ\mathbb{I}(i')$ for all $i<i<j$. However, it is not necessarily true that  $\phi\circ\mathbb{I}(i)$ is an interval for a given $i.$ A particular case in which this occurs is when $im(\phi)$ is itself an interval in $B$.  In fact, it is easy to show that $im(\phi)$ is an interval in $B$ if and only if $$\{\phi\circ \mathbb{I}\,:\,\mathbb{I}\in Bl(0,A)\}\subseteq Bl(0,B).$$
Moreover, if $\phi$ is surjective then the equality between the previous sets will hold. Particularly, if  $\alpha\in \gamma$ and $k\in \omega$, then $$Bl(\alpha,F)=\{\,F\circ \mathbb{I}\,:\,\mathbb{I}\in Bl_k(\alpha)\,\},$$
$$Bl_k(\alpha)=\{F^{-1}\circ \mathbb{I}\,:\,\mathbb{I}\in Bl(\alpha,F)\,\}$$
for all $F\in \mathcal{F}_k$ with $\alpha\in F.$
    
\end{rem}

\begin{center}
\begin{minipage}[b]{0.8\linewidth}
\centering
\includegraphics[width=12 cm, height=5.5cm]{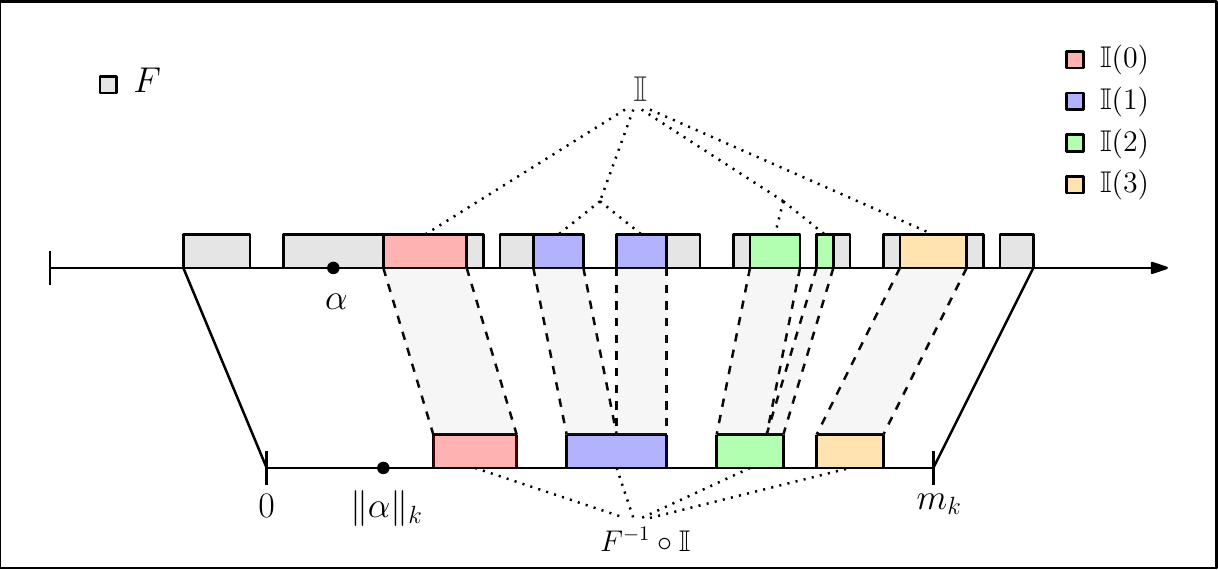}

\textit{\small In here, $F$ represents an element of $\mathcal{F}_k$ with $\alpha\in F$  and $\mathbb{I}$ is an element of $Bl(\alpha,F)$.}
\end{minipage}
\end{center}

\begin{definition}[Good pairs and good sets]Let $\alpha\in \omega_1$ and $A\in \text{FIN}(\omega_1)$. we say that a pair $(\mathbb{I},z)$ is \textit{$(\alpha,A)$-good} if the following conditions hold:
\begin{itemize}
\item $z\in A\backslash \alpha,$
\item $\mathbb{I}\in Bl(\alpha, A\cap (z+1))$. That is, $\mathbb{I}\in Bl(\alpha,A)$ and either $\mathbb{I}=\emptyset$ or $\max(\bigcup \mathbb{I})\leq z.$
\end{itemize}
A non-empty set of ordered pairs $T$ is said to be \textit{$(\alpha,A)$-Good} if every $(\mathbb{I},z)\in T$ is  $(\alpha,A)$-good. We denote the family of $(\alpha,A)$-Good sets as  $Good(\alpha,A)$. Given $\alpha\in \gamma$ and $k\in \omega$, we denote $Good(\lVert \alpha \rVert_k,m_k)$ as $Good_k(\alpha)$.

\end{definition}

\begin{center}
\begin{minipage}[b]{0.8\linewidth}
\centering
\includegraphics[width=10 cm, height=3cm]{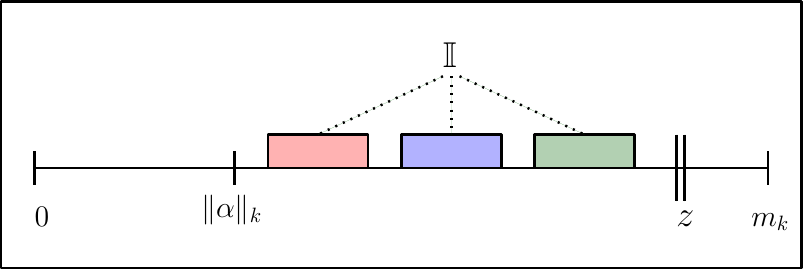}

\textit{\small In here, $(\mathbb{I},z)$ represents an $(\lVert \alpha\rVert_k,m_k)$-good pair.}
\end{minipage}
\end{center}

\begin{definition} Suppose that $A$ and $B$ are non-empty subsets of $\omega_1$ and $\phi:A\longrightarrow B$ is an increasing function. Given $\alpha\in \omega_1$ and $T\in Good(\alpha,A)$ we define $$\phi\bullet T=\{\,(\phi\circ \mathbb{I},\phi(z))\,:\,(\mathbb{I},z)\in T\,\}.$$
In the particular case where $|B|=m$ for some $m\in \omega$ and $A=m$, we denote $\phi \bullet T$ as $B\bullet T$. Analogously, if $T'\in Good(\alpha,B)$ then $\phi^{-1}\bullet T'$ is denoted as $B^{-1}\bullet T'.$
    
\end{definition}
\begin{rem}If $\alpha\in \gamma$ and $k\in \omega$ then
$$Good(\alpha,F)=\{F\bullet T\,:\,T\in Good_k(\alpha)\,\},$$
$$Good_k(\alpha)=\{F^{-1}\bullet T\,:\,T\in Good(\alpha,F)\,\}$$
for each $F\in \mathcal{F}_k$ with $\alpha\in F.$
\end{rem}
For the next the definition, it is convenient to recall the function $\phi^\beta_l$ that was presented in Definition \ref{functionphibetaldef}.
\begin{definition}Let $\beta\in \gamma$ and $k\leq l\in \omega$. We define $\pi^{k,l}_\beta:\lVert \alpha\rVert_k+1\longrightarrow \lVert \beta\rVert_l+1$ as: $$\pi^{k,l}_\beta=\phi^\beta_l\circ(\phi^\beta_k)^{-1}.$$
    
It is not hard to see that $\pi_\alpha^{k,l}$ is just the increasing bijection from $\lVert \alpha\rVert_k+1$ to $(\lVert \alpha\rVert_l)_k.$
 \end{definition}
\begin{definition}[The $\bigstar$ relation]\label{defstarrelation} Let $\beta, \xi\in \gamma$, $2\leq k<l\in \omega$ , and $T\in Good_k(\beta).$ We say that \textit{$T$ guesses $(\beta,\xi,k,l)$}, and write it as $T\bigstar (\beta,\xi,k,l)$, if there is $(\mathbb{I},z)\in T$ for which:
\begin{enumerate}[label=$(\alph*)$]
    \item $\lVert \beta \rVert_l\leq \lVert \xi\rVert_l$ and $\lVert \xi\rVert_k=z$,
    \item $(\xi)_k\circ \mathbb{I}\in Bl(0,(\xi)_l).$ Equivalently, if $\pi^{k,l}_\xi\circ \mathbb{I}\in Bl(0,m_l),$
    \item $(\xi)_l(\lVert\beta\rVert_l)\in (\xi)_k.$ Equivalently, $\lVert \beta\rVert_l\in (\lVert \xi\lVert_l)_k.$
\end{enumerate}
\end{definition}
\begin{center}
\begin{minipage}[b]{0.8\linewidth}
\centering
\includegraphics[width=12 cm, height=5cm]{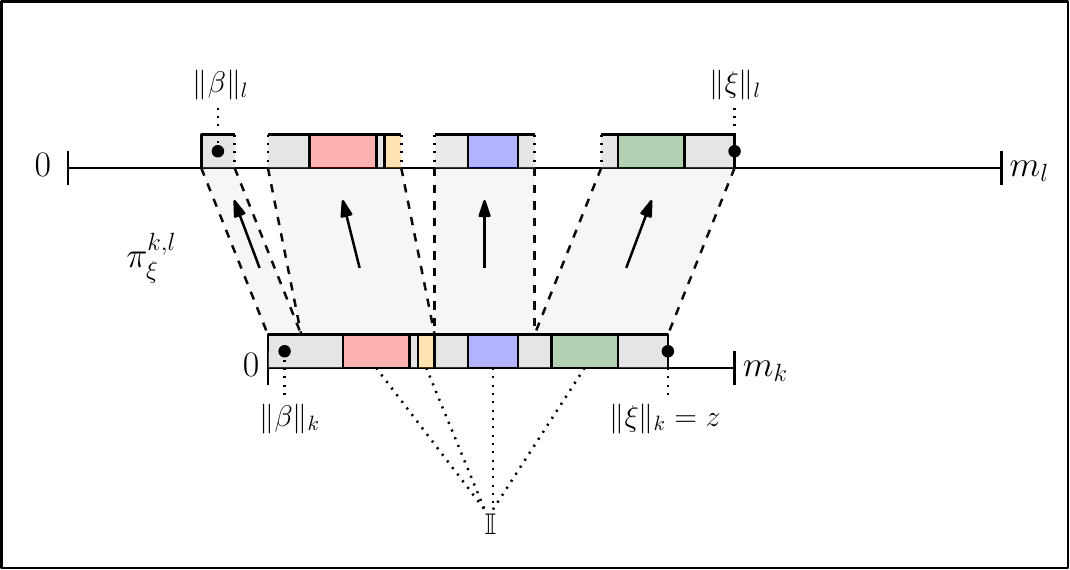}
\textit{\small{In here, the pair $(\mathbb{I},z)$ represent an element of $T$ which testifies that $T\star(\beta,\xi,k,l).$}}
\end{minipage}
\end{center}
\begin{rem}\label{bigstarremark} If $T\bigstar(\beta,\xi,k,l)$ and $\xi'$ is such that $\lVert \xi'\rVert_l=\lVert \xi\rVert_l$ then $T\bigstar (\beta,\xi',k,l).$
\end{rem}
\begin{rem}The previous definition will often be applied in cases where $\xi<\beta$.
    
\end{rem}
Recall that we fixed  $\langle D_\delta \rangle_{\delta\in \text{Lim}}$ a $\Diamond$-sequence at the start of this section.
\begin{definition}[The $\checkmark$ relation]\label{checkmarkdef} Let $\beta\in \gamma$, $k\leq l\in \omega$ and $\delta\in Lim$. Given $T\in Good_k(\beta)$ and $C\in \text{FIN}(\gamma)$, we say that \textit{$(C,\delta)$ accepts $(T,\beta,k,l)$}, and write it as $(C,\delta)\checkmark (T,\beta,k,l)$, if the following conditions hold:
\begin{enumerate}[label=$(\arabic*)$]
\item $C\subseteq D_\delta$ and $C$ is captured at level $l$,
\item $T\bigstar (\beta,\,C(0),\, k,\,l-1),$.  
\item $|(\beta)_{l-1}\cap \delta|=r_l$. Equivalently, if $F\in \mathcal{F}_l$ is such that $C\subseteq F$, then $(\beta)_{l-1}\cap \delta=R(F).$ 
\end{enumerate}
Whenever there is $C$ for which $(C,\delta)\checkmark(T,\beta,k,l)$, such $C$ has cardinality at least $1$ and at most $n_l$ due to the point (1).  Thus, we can define $j(\delta,T,\beta,k,l)$ as the maximum of such cardinalities. That is: $$j(\delta,T,\beta,k,l)=\max(\, j\leq n_l\,:\,\exists C\in [\gamma]^j\,(\,(C,\delta)\checkmark (T,\beta,k,l)\,)\,).$$
If there is no such $C$, then we define $j(\delta,T,\beta,k,l)$ as $0$.
\end{definition}
\begin{definition}[The $IH_2$ property]\label{ih2propertydef}Given $\delta\in\text{Lim}\cap \gamma$, we say that $IH_2(\delta,\mathcal{F})$ holds if one of the two following mutually excluding conditions occurs:
\begin{enumerate}[label=$(\Alph*)$]
\item There are infinitely many $l\in \omega$ for which there is $C\in \text{FIN}(D_\delta)$ which is fully captured at level $l$. In this case, we will say that $IH^A_2(\delta,\mathcal{F})$ holds.
\item For all $\beta\in [\delta,\gamma)$, $2\leq k\in \omega$  and $T\in Good_k(\beta)$ there are infinitely many $k\leq l\in \omega$ for which $j=j(\delta,T,\beta,k,l)<n_l$ and such that:
\begin{enumerate}[label=$(B.\arabic*)$]
    \item $|(\beta)_{l-1}\cap \delta|=r_l,$
    \item $\Xi_\beta(l)=j,$
    \item  If $j>0$ then there is $C\in [\delta]^j$ for which $(C,\delta)\checkmark(T,\beta,k,l)$ and $C\subseteq (\beta)_l.$
\end{enumerate}
In this case we will say that $IH^B_2(\delta,\mathcal{F})$ holds. 
\end{enumerate}  
Finally, we say that  $\mathcal{F}$ satisfies $IH_2$ if $IH_2(\delta,\mathcal{F})$ holds for any $\delta \in \text{Lim}\cap \gamma$.
\end{definition}
The two following results follow directly from the definition.
\begin{proposition}\label{IH2omega}$\mathcal{F}(\omega)$ satisfies $IH_2$.    
\end{proposition} 

\begin{lemma}\label{unionschemelemmadeeperIH2}
    Assume that $\gamma$ is a limit of limit ordinals and  for each limit $\gamma'<\gamma$, $\mathcal{F}|_{\gamma'}$ satisfies $IH_2$. Then $\mathcal{F}$ satisfies $IH_2$.
\end{lemma}

More important, the property $IH_2$ in fact implies full capturing when $\gamma=\omega_1$.
\begin{theorem}\label{ih2omega1fullycapturing}Suppose that $\gamma=\omega_1$ and $\mathcal{F}$ satisfies $IH_2$. Then $\mathcal{F}$ is fully capturing.
\begin{proof} We will prove the Theorem by appealing to the equivalence of fully capturing stated in Lemma \ref{equivalencecapturing}. Let $S\in [\omega_1]^{\omega_1}$. Since $\langle D_\delta \rangle_{\delta\in \text{Lim}}$ is a $\Diamond$-sequence, there is $\delta\in \text{Lim}$ so that:
\begin{enumerate}
    \item $S\cap \delta=\delta$.
    \item $(\delta,<,D_\delta,\mathcal{F}|_\delta)$ is an elementary submodel of $(\omega_1,<,S,\mathcal{F}).$
\end{enumerate}
Suppose towards a contradiction that there is no $C\in \text{FIN}(S)$ so that $C$ is fully captured. By elementary, the same is true for $D_\delta$. In other words, $IH_2^A(\delta,\mathcal{F})$ fails and consequently $IH_2^B(\delta,\mathcal{F})$ holds. Let us fix  $\beta\in S\backslash \delta$ and $k=2.$ We now define $\mathbb{I}=\emptyset$, $z=\lVert \beta\rVert_2$ and $T=\{ (\mathbb{I},z) \}\in Good_2(\beta)$. According to the point (B) in Definition \ref{ih2propertydef},  there is $2\leq l\in \omega$ for which $j=j(\delta,T,\beta,2,l)<n_l$ and such that:
\begin{itemize}
\item $|(\beta)_{l-1}\cap \delta|=r_l$,
    \item $\Xi_\beta(l)=j$,
    \item If $j>0$ then there is $C\in [\delta]^j$ for which $(C,\delta)\checkmark (T,\beta,2,l)$ and $C\subseteq (\beta)_l$.
\end{itemize}
We will arrive to the desired contradiction by finding  a finite set $C'\in \text{FIN}(\omega_1)$ for which $(C',\delta) \checkmark (T,\beta,2,l)$ and $|C'|=j+1$. The proof of this is divided into two cases:\\

\noindent
\underline{Case 1}: If $j=0.$
\begin{claimproof}[Proof of case] In this case, as $\Xi_\beta(l)=0$, it is straightforward that $\{\beta\}$ is captured at level $l$.  By elementarity, there is $\xi\in D_\delta$ so that $\lVert \xi\rVert_l=\lVert \beta\rVert_l$. In particular, $\{\xi\}$ is captured at level $l$. The proof of this case ends by proving the following claim.\\

\noindent
\underline{Claim 1}: $(\{\xi\},\delta)\checkmark (T,\beta,k,l)$.
\begin{claimproof}[Proof of claim] It is enough to show that $T\bigstar (\beta, \xi,2, l-1)$.  For this, note that the  point $(a)$ of Definition \ref{defstarrelation} is satisfied for the unique pair $(\mathbb{I},z)\in T$ because $\lVert \beta\rVert_l=\lVert \xi\rVert_l$. In particular, $\Delta(\beta,\xi)>l$ so $\lVert \xi\rVert_2=\lVert \beta\lVert_2= z$. The  point $(b)$ of Definition \ref{defstarrelation} holds because $\mathbb{I}=\emptyset$. Therefore, $\pi^{2,l}_\xi=\emptyset\in Bl(0,m_l)$. Finally, the point $(c)$ of Definition \ref{defstarrelation} is satisfied because $(\xi)_l(\lVert \beta\rVert_l)=(\xi)_l(\lVert \xi\rVert_l)=\xi$.
\end{claimproof}

\end{claimproof}

\noindent
\underline{Case 2}: If $j>0$.
\begin{claimproof}[Proof of case]

Let $C\in [\delta]^j$ be such that $(C,\delta)\checkmark (T,\beta,2,l)$ and $C\subseteq (\beta)_l$. By elementarity, there is $\xi\in D_\delta$ so that $\lVert \xi\rVert_l=\lVert \beta\rVert_l$ and $C\subseteq (\xi)_l$. Particularly, $\Xi_\xi(l)=j$ and $\rho(C(i),\xi)\leq l$ for each $i<j$. In fact, as $\Xi_{C(i)}(l)=i$ because $C$ is captured at level $l$,  we have that $\rho(C(i),\xi)=l$ due to Lemma \ref{lemmaxi}.
We will finish the proof of this case by showing the next claim.\\

\noindent
\underline{Claim 2}: $(C\cup\{\xi\},\delta)\checkmark (T,\beta,k,l)$.
\begin{claimproof}[Proof of claim] It suffices to prove that $C\cup \{\xi\}$ is captured at level $l. $ For this, we appeal to the Proposition \ref{deltarhoequalityprop}. By the previous observations, we only need to prove that $\Delta(C(i),\xi)=\rho(C(i),\xi)$ for each $i<j$. First note that, according to the point $(2)$ in Definition \ref{checkmarkdef}, $T\bigstar (\beta, C(0), 2,l-1)$. As $(\emptyset, \lVert \beta\rVert_2)$ is the only element of $T$, then $\lVert \beta\rVert_{l-1}\leq \lVert C(0)\rVert_{l-1}$, $\lVert C(0) \rVert_2=\lVert \beta\rVert_2$ and $\lVert \beta\rVert_{l-1}\in (\lVert C(0)\rVert_{l-1})_2$. In other words, $\lVert \beta\rVert_{l-1}$ is in the domain of $\phi^{\lVert C(0)\rVert_{l-1}}_2$. Moreover,  \begin{align*}\phi^{\lVert C(0)\rVert_{l-1}}_2(\lVert C(0)\rVert_{l-1})=\lVert\,\lVert C(0)\rVert_{l-1}\,\rVert_2 &= \lVert C(0)\rVert_2\\ &=\lVert \beta\rVert_2=\lVert\, \lVert \beta\rVert_{l-1}\, \rVert_2=\phi^{\lVert C(0)\rVert_{l-1}}_2(\lVert \beta\rVert_{l-1})
\end{align*}
  due to the part $(a)$ of Proposition \ref{kcardinalityisomorphismprop}. Since $\phi^{\lVert C(0)\rVert_l}_2$ is one to one, it follows that $\lVert C(0)\rVert_{l-1}=\lVert \beta\rVert_{l-1}$. Now, as $\lVert \beta\rVert_l=\lVert \xi\rVert_l$ and $\lVert C(0)\rVert_{l-1}=\lVert C(i)\rVert_{l-1}$, then $\lVert \xi\rVert_{l-1}=\lVert C(i)\rVert_{l-1}$. That is, $l\leq \Delta(C(i),\xi)\leq \rho(C(i),\xi)\leq l$. This finishes the proof.
\end{claimproof}
\end{claimproof}
\end{proof} 
\end{theorem}
Before starting the next subsection, we present the key concept needed for proving Theorem \ref{superprincipaltheorem}.

\begin{definition}[Transferring  Good sets] Let  $\beta\in \gamma$, $k\in \omega$ and $T\in Good_k(\beta)$. Given $k'>k$, we define $\mathcal{Q}_{k'}(\beta,T)$ as the set of all pairs $(F,\mathbb{I})\in \mathcal{F}_k(m_{k'})\times Bl_k(\beta)$ such that:
\begin{itemize}
\item There is $z\in m_k$ such that $(\mathbb{I},z)\in T$,
    \item $\lVert \beta\rVert_{k'}\in F$,
    \item $F\circ \mathbb{I}\in Bl_{k'}(\beta).$ In other words, $F[\mathbb{I}(i)]$ is an interval in $m_k'$ for every $i<\mathbb{I}.$
\end{itemize}
Given $\alpha\in (\beta)^-_{k'}($ and $(F,\mathbb{I})\in \mathcal{Q}_{k'}(\beta,T)$, let $\mathbb{I}_F=(\lVert \alpha\rVert_{k'},\lVert \beta,\rVert_{k'}]^\smallfrown (F\circ \mathbb{I})$\footnote{In principle, $\mathbb{I}_F$ also depends on $\alpha$. However, we will never use this concept with two different $\alpha$'s at the same time. For that reason, we omit any mention of that ordinal in the notation.}. That is, $\mathbb{I}_F=\langle \mathbb{I}_F(i)\rangle_{i<\mathbb{I}+1}$ is given by: $$\mathbb{I}_F(i)=\begin{cases}(\lVert \alpha\rVert_{k'},\lVert \beta\rVert_{k'}]&\textit{ if }i=0\\
F[\mathbb{I}(i-1)]&\textit{ if }i>0
\end{cases}$$ We define \textit{the transferring of $T$  relative to $k',\alpha$ and $\beta$} as:
$$\text{Transf}_{k'}(T,\alpha,\beta)=\{\,(\mathbb{I}_F,\,F(z))\,:\,(\mathbb{I},z)\in T\textit{ and }(F,\mathbb{I})\in \mathcal{Q}_{k'}(\beta,T)\,\}$$
\end{definition}

\begin{center}
\begin{minipage}[b]{1\linewidth}
\centering
\includegraphics[width=15cm, height=3.5cm]{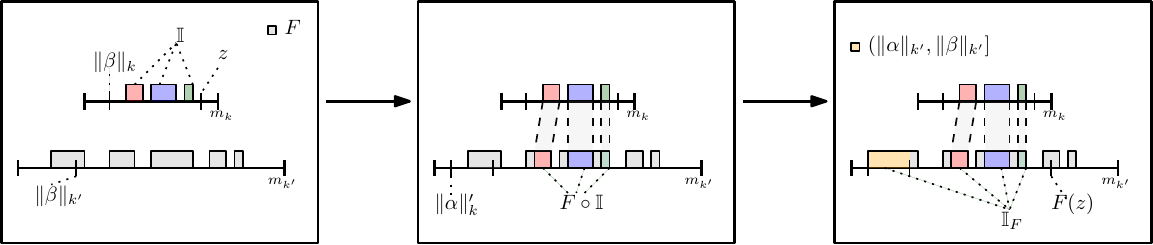}
\end{minipage}
\textit{\small In the first picture, we considered a a pair $(\mathbb{I},z)\in T$ and $F\in \mathcal{F}_k(m_{k'})$ with $\lVert \beta\rVert_{k'}\in F$. In the second picture, we verified that $F\circ \mathbb{I}$ is an element of  $Bl_{k'}(\beta)$. As this is the case, we can conclude that $(\mathbb{I}_F,F(z))\in \text{Transf}_{k'}(T,\alpha,\beta)$.}
\end{center}

\begin{rem}Note that if $T\in Good_k(\beta)$, $\alpha\in (\beta)^-_k$ and $k'>k$ then $Transf_{k'}(T,\alpha,\beta)\in Good_{k'}(\alpha).$  
\end{rem}

\begin{lemma}\label{transfequiv}Let $\beta\in \gamma$, $k\in \omega$, $T\in Good_k(\beta)$, $k'>k$ and $\alpha\in (\beta)^-_{k'}$. Suppose that $l>k'$ and $\xi\in \gamma$ are such that:\begin{enumerate}[label=$(\arabic*)$] \item $\lVert \beta\lVert_l\leq \lVert \xi\rVert_l,$ 
\item $(\xi)_l(\lVert \beta\rVert_l)\in (\xi)_{k'}.$ In other words, $\lVert \beta\rVert_l\in (\lVert \xi\rVert_l)_{k'}$
\item $(\beta)_{k'}[\,(\lVert\alpha\rVert_{k'},\lVert\beta\rVert_{k'}]\,]=(\beta)_{k'}\backslash (\alpha+1)$ is an interval in $(\beta)_l.$
\end{enumerate}
Then the following statements are equivalent:
\begin{enumerate}[label=$(\alph*)$]
    \item $T\bigstar (\beta,\xi,k,l).$
    \item $\text{Transf}_{k'}(T,\alpha,\beta)\bigstar (\alpha,\xi,k',l).$
\end{enumerate}
\begin{proof}Let $\mathcal{Q}=\mathcal{Q}_{k'}(\beta,T)$ and $T'=\text{Transf}_{k'}(T,\alpha,\beta)$.
\begin{claimproof}[$(a)\Rightarrow (b)$] Let $(\mathbb{I},z)\in T$ so that $\lVert\beta\rVert_k\leq \lVert \xi \rVert_k=z$, $\pi^{k,l}_\xi\circ \mathbb{I}\in Bl(0,m_l)$ and $\lVert \beta\rVert_l\in (\lVert\xi\rVert_l)_k$. Now, let us take $F\in \mathcal{F}_{k}(m_{k'})$ having $\lVert \xi \rVert_{k'}$ as an element. Before we continue, we will highlight some useful facts:\begin{itemize}
    \item  $\lVert \xi\rVert_{k'}\in F$ so $F\cap (\xi+1)=(\,\lVert \xi\rVert_{k'})_k$. Consequently, $F(a)=(\phi_k^{\lVert \xi\rVert_{k'}})^{-1}(a)$ for each $a\leq \lVert \xi\rVert_k$.
    \item $\lVert\xi\rVert_k=\phi^{(\lVert{\xi}\rVert_{k'})}_k\circ \phi^{\xi}_k(\xi)$ according to Lemma \ref{compositionphi}. Therefore, $F(\lVert \xi\rVert_k)=\lVert \xi\rVert_{k'}$ by virtue of the previous point.
    \item Again, by Lemma \ref{compositionphi}: $$\phi^{\lVert \xi\rVert_l}_{k'}\circ \pi^{k,l}_\xi=\phi^{\lVert \xi\rVert_l}_{k'}\circ  \phi^\xi_l\circ(\phi^\xi_k)^{-1}=\phi^\xi_{k'}\circ (\phi^\xi_k)^{-1}=(\phi^{\lVert \xi\rVert_{k'}}_k)^{-1}. $$
    
\end{itemize}
\underline{Claim 1}: $(F,\mathbb{I})\in \mathcal{Q}$.
\begin{claimproof}[Proof of claim] First we will show that $\lVert \beta\rVert_{k'}\in F$. Indeed, since $k\leq k'<l$ and $\lVert \beta \rVert_l\in (\lVert \xi\rVert_l)_k$ then, by Proposition \ref{kcardinalityisomorphismprop} and Lemma \ref{compositionphi}, we can conclude that \begin{align*}\lVert \beta \rVert_{k'}=\lVert\,\lVert \beta\rVert_l\,\rVert_{k'}=\phi^{\lVert \xi\rVert_l}_{k'}(\lVert \beta \rVert_l)&\in \phi^{\lVert \xi\rVert_l}_{k'}[\,(\,\lVert \xi\rVert_l\,)_k\,]\\&=(\,\phi^{\lVert \xi\rVert_l}_{k'}(\lVert \xi\rVert_l)\,)_k\\&=(\,\lVert\,\lVert \xi\rVert_l\,\rVert_{k'}\,)_k=(\lVert \xi\rVert_{k'})_k\subseteq F.
\end{align*}
Now, let $i<|\mathbb{I}|$. The claim will be over once we prove that $F[\mathbb{I}(i)]$ is an interval in $m_{k'}$. As $\mathbb{I}(i)\leq z=\lVert\xi\rVert_k$ then $F[\mathbb{I}(i)]=(\phi^{\lVert \xi\rVert_{k'}}_k)^{-1}[\mathbb{I}(i)]=\phi^{\lVert \xi\rVert_l}_{k'}\circ\pi^{k,l}_\xi[\mathbb{I}(i)].$ By virtue of the hypotheses, $\pi^{k,l}_\xi[\mathbb{I}]$ is an interval in $m_l$ which is contained in $(\,\lVert \xi\rVert_l\,)_k$. Hence, such interval is also included in $(\,\lVert \xi\rVert_l\,)_{k'}$. Since, $\phi^{\lVert \xi\rVert_l}_{k'}$ is the increasing bijection from $(\,\lVert \xi\rVert_l\,)_{k'}$ to $\lVert \xi\rVert_{k'}$, it follows that $\phi^{\lVert \xi\rVert_l}_{k'}[\pi^{k,l}_\xi[\mathbb{I}]]$ is an interval in $\lVert \xi\rVert_{k'}+1$. From this it is straightforward that $F[\mathbb{I}(i)]$ is an interval in $m_{k'}$. This finishes the claim.
\end{claimproof}
Since $(\mathbb{I},z)\in T$ and $(F,\mathbb{I})\in \mathcal{Q}$ then $(\mathbb{I}_F,F(z))\in T'$. We will show that such pair testifies that $T'\bigstar (\alpha,\xi,k',l)$. This will be achieved by proving the three following claims.\\

\noindent
\underline{Claim 2}: $\lVert \alpha\rVert_{l}\leq \rVert \xi\rVert_{l}$ and $\lVert \xi\rVert_{k'}=F(z).$
\begin{claimproof}[Proof of claim] We know that $z=\lVert \xi\rVert_k$, and it has already been proved that $F(\lVert \xi\rVert_k)=\lVert \xi\rVert_{k'}$. As $\alpha\in (\beta)_{k'}$ then $\lVert \alpha\rVert_l\leq \lVert \beta\rVert_l$. But we know by the hypotheses that $\lVert \beta\rVert_l\leq \lVert \xi\rVert_l$. So we are done 
    
\end{claimproof}
\noindent 
\underline{Claim 3}: $\pi^{k',l}_\xi\circ \mathbb{I}_F\in Bl(0,m_l).$
\begin{claimproof}[Proof of claim] Let $i<|\mathbb{I}_F|$.  We want to show that $\pi^{k',l}_\xi[\mathbb{I}_F(i)]$ is an interval in $m_l$. First we deal with the case where $i=0$. Here,  $\mathbb{I}_F(i)=(\lVert\alpha\rVert_{k'},\lVert \beta\rVert_{k'}]$. We will first argue that, in this case, $\pi^{k',l}_\xi[\mathbb{I}_F(i)]=\pi^{k',l}_\beta[\mathbb{I}_F(i)]$.  Indeed, since $\lVert \beta\rVert_l\in (\lVert \xi\rVert_l)_{k'}$ then $\phi^{\lVert \xi\rVert_l}_{k'}|_{\lVert \beta\rVert_l+1}= \phi^{\lVert\beta\rVert_l}_{k'}.$ In particular, this means that $(\phi^{\lVert \xi\rVert_l}_{k'})^{-1}[a]=(\phi^{\lVert \beta\rVert_l}_{k'})^{-1}[a]$ for each $a\leq \lVert \beta\rVert_{k'}$. In this way, by using Lemma \ref{compositionphi}, we conclude that \begin{align*}\pi^{k',l}_\xi[\mathbb{I}_F(i)]&=\phi^\xi_l\circ (\phi^\xi_{k'})^{-1}[\mathbb{I}_F(i)]
\\&=\phi^\xi_l\circ ((\phi^{\xi}_l)^{-1}\circ (\phi^{\lVert \xi\rVert_l}_{k'})^{-1})[\mathbb{I}_F(i)]
\\&=(\phi^{\lVert \beta\rVert_l}_{k'})^{-1}[\mathbb{I}_F(i)]=(\phi^{\lVert \beta\rVert_l}_{k'})^{-1}[\mathbb{I}_F(i)]\\&=(\phi^\beta_l\circ ((\phi^{\beta}_l)^{-1})\circ (\phi^{\lVert \beta\rVert_l}_{k'})^{-1}[\mathbb{I}_F(i)]\\&=\phi^\beta_l\circ(\phi^\beta_{k'})^{-1}[\mathbb{I}_F(i)]=\pi^{k',l}_\beta[\mathbb{I}_F(i)]
\end{align*}
Observe that $(\phi^{\beta}_{k'})^{-1}[\mathbb{I}(i)]=(\beta)_{k'}[\,(\lVert \alpha\rVert_{k'},\lVert \beta\rVert_{k'}]\,]$ is an interval in $(\beta)_l$, according to the point (3) in the hypotheses of this lemma. Since $\phi^{\beta}_l$ is the increasing bijection from $(\beta)_l$ to $\lVert \beta\rVert_l+1\sqsubseteq m_l$, then $\pi^{k',l}_\beta[\mathbb{I}(i)]$ is an interval in $m_l$. This finishes the first case. Now we deal with the case where $i>0$. Here, $\mathbb{I}_F(i)=F[\mathbb{I}(i-1)]$. Thus,\begin{align*}
    \pi^{k',l}_\xi[\mathbb{I}_F(i)]&=\phi^{\xi}_l\circ (\phi^{\xi}_{k'})^{-1}[F(i)]\\&=\phi^{\xi}_l\circ (\phi^{\xi}_{k'})^{-1}[(\phi^{\lVert \xi\rVert_{k'}}_k)^{-1}[\mathbb{I}(i-1)]]\\&=\phi^\xi_l\circ((\phi^\xi_{k'})^{-1}\circ (\phi^{\lVert \xi\rVert_{k'}}_k)^{-1})[\mathbb{I}(i-1)]\\&=\phi^\xi_l\circ(\phi^\xi_k)^{-1}[\mathbb{I}(i-1)]=\pi_\xi^{k,l}[\mathbb{I}(i-1)]\end{align*} 
 Due to the hypotheses, $\pi^{k,l}_\xi[\mathbb{I}(i-1)]$ is already an interval in $m_l$. Consequently, this case is over. 
\end{claimproof}
\noindent
\underline{Claim 4}: $\lVert \alpha\rVert_l\in (\lVert \xi\rVert_l)_{k'}.$
\begin{claimproof}[Proof of claim] Recall that $\alpha\in (\beta)_k$. Therefore, $\lVert\alpha\rVert_l\in (\lVert \beta\rVert_l)_k$ due to the point (b) of Proposition \ref{kcardinalityisomorphismprop}. As $\lVert \beta\rVert_l\in (\lVert \xi \rVert_l)_{k'}$, we are done.
\end{claimproof}
\end{claimproof}
\begin{claimproof}[$(b)\Rightarrow (a)$] Let $(\mathbb{I},z)\in T$ and $F\in \mathcal{F}_k(m_{k'})$ with $(F,\mathbb{I})\in \mathcal{Q}$ be such that  $\lVert \alpha\rVert_{k'}\leq \lVert \xi \rVert_{k'}=F(z)$, $\pi^{k',l}_\xi\circ \mathbb{I}_F\in Bl(0,m_l)$ and $\lVert \alpha\rVert_l\in (\lVert \xi\rVert_l)_{k'}$. We will show in the next three claims that $(\mathbb{I},z)$ testifies that $T\bigstar (\beta,\xi,k',l)$.\\

\noindent
\underline{Claim 5}: $\lVert \beta \rVert_l\leq \lVert \xi\rVert_l$ and $\lVert \xi\rVert_k=z.$
\begin{claimproof}[Proof of claim]
Due to the hypotheses of the lemma, we already know that $\lVert \beta\rVert_l\leq \lVert \xi\rVert_l$. Therefore it suffices to prove that $\lVert \xi\rVert_k=z$. Indeed, as $F\in \mathcal{F}_k(m_{k'})$ then $\lVert F(z)\rVert_k=z$. Moreover, we know that $F(z)=\lVert \xi\rVert_{k'}$. Thus, $z=\lVert F(z)\rVert_k=\lVert \,\lVert \xi\rVert_{k'}\,\rVert_k=\lVert \xi\rVert_k$ by virtue of the point (a) in Proposition \ref{kcardinalityisomorphismprop}.
    
\end{claimproof}
\noindent
\underline{Claim 6}: $\pi^{k,l}_\xi\circ \mathbb{I}\in Bl(0,m_l).$
\begin{claimproof}[Proof of claim] Let $i<|\mathbb{I}|$. We need to prove that $\pi^{k,l}_\xi[\mathbb{I}(i)]$ is an interval in $m_l$. First note that $\lVert \xi\rVert_{k'}\in F$, so as in the proof of the other direction of the lemma,  $F\cap (\xi+1)=(\lVert \xi\rVert_{k'})_k$. In this way, $F(a)=(\phi^{\lVert \xi\rVert_{k'}}_k)^{-1}(a)$ for each $a\leq \lVert \xi \rVert_k.$ From this fact, we get the following chain of equalities:
\begin{align*}\pi^{k,l}_\xi[\mathbb{I}(i)]=\phi^\xi_l\circ(\phi^\xi_k)^{-1}[\mathbb{I}(i)]&=\phi^\xi_l\circ ( (\phi^\xi_{k'})^{-1}\circ(\phi^{\lVert \xi\rVert_{k'}}_k)^{-1})[\mathbb{I}(i)]\\&=\phi^\xi_l\circ \phi^\xi_{k'}[F[\mathbb{I}(i)]]\\&=\pi^{k',l}_\xi[\mathbb{I}_F(i+1)]
\end{align*}
As $\pi^{k',l}_\xi\circ \mathbb{I}\in Bl(0,m_l)$, we are done.
\end{claimproof}
\noindent
\underline{Claim 7}: $\lVert \beta\rVert_l\in (\lVert \xi\rVert_l)_k$.
\begin{claimproof}[Proof of claim] On one hand, $\lVert \beta\rVert_{k'}\in F$ since $(F,\mathbb{I})\in \mathcal{Q}$. On the other hand, we already know that $\lVert \xi \rVert_k\in F$. Therefore, $\rho(\lVert \beta\rVert_{k'},\lVert \xi\rVert_{k'})\leq k$. Now, according to the hypotheses $\rho(\lVert \xi\rVert_l,\lVert \beta\rVert_l)\leq {k'}$. Thus, $$\rho(\lVert \beta\rVert_l,\lVert \xi\rVert_l)=\rho(\lVert\,\lVert \beta\rVert_l\,\rVert_{k'},\lVert \,\lVert \xi\rVert_l\,\rVert_{k'})=\rho(\lVert \beta\rVert_{k'},\lvert \xi\rVert_{k'})$$
by virtue of the points (a) and (b) of Proposition \ref{kcardinalityisomorphismprop}. In this way, $\rho(\lVert \beta\rVert_l,\lVert \xi\rVert_l)\leq k$.
    
\end{claimproof} 
\end{claimproof}
\end{proof}   
\end{lemma}

\begin{corollary}\label{trasnfcoro}Let $\mathcal{F}$ be a construction scheme over an ordinal $\gamma$. Also, let $\beta\in \gamma$, $T\in Good_k(\beta)$, $k'>k$ and  $\alpha\in (\beta)^-_{k'}$. Suppose that $(\beta)_{k'}\backslash (\alpha+1)$ is an interval in $(\beta)_{l-1}$. Furthermore, assume that $l>k'$ and $\delta\in \text{Lim}\cap \alpha$ are such that $|(\alpha)_{l-1}\cap \delta|=r_l$. If $C\in \text{FIN}(\gamma)$ is such that $(C,\delta)\checkmark (k,l,\beta,\delta,T)$ then $(C,\delta) \checkmark(k',l,\beta,\delta,T')$ where $T'=\text{Transf}_{k'}(\alpha,\beta,T)$. In particular, $j(\delta,T,\beta,k,l)\leq j(\delta,T',\alpha,k',l)$.
\begin{proof} Let $C\in \text{FIN}(\gamma)$ be such that $(C,\delta)\checkmark (k,l,\beta, \delta ,T')$.   Since $\alpha\in (\beta)_{k'}$ and $l>k'$ then $(\beta)_l\cap(\alpha+1)=(\alpha)_l$. By the hypotheses, $\alpha\geq \delta$. This means that $(C(0))_l[r_l]=(\beta)_l\cap \delta=(\alpha)_l\cap \delta$. The only thing left to prove is that $T'\bigstar (\alpha,C(0), k',l-1)$. For this, it suffices to show that the hypotheses of Lemma \ref{transfequiv} hold for $\xi=C(0)$. Indeed, since $T\bigstar(\beta,\xi,k,l-1)$ then $\lVert \beta\rVert_{l-1}\in (\lVert \xi\rVert_{l-1})_k$. This implies both that $\lVert \beta \rVert_{l-1}\leq \lVert \xi\rVert_{l-1}$ and $\lVert \beta\rVert_{l-1}\in (\lVert \xi\rVert_{l-1})_{k'}$. According to the assumptions $(\beta)_{k'}[\,(\lVert\alpha\rVert_{k'},\,\lVert \beta\rVert_{k'}]\,]$ is an interval in $(\beta)_{l-1}$. Thus, we are done.
    
\end{proof}
    
\end{corollary}
\subsection{Forcing the property $IH_2$}

In this subsection we aim  to prove that if $\gamma$ is a countable ordinal and $\mathcal{F}$ satisfies $IH_1$ and $IH_2$, then we can extend $\mathcal{F}$ to a construction scheme over $\gamma+\omega$ which also satisfies $IH_1$ and $IH_2$. For this, we will make use of the forcing $\mathbb{P}(\mathcal{F})$. In Section \ref{sectionpfforcing} we already proved that there is a countable family of dense sets over $\mathbb{P}(\mathcal{F})$ so that if $G$ is a filter intersecting each of those sets, then $\mathcal{F}^\mathcal{G}$ is construction scheme over $\gamma+\omega$ which also satisfies $IH_1$. Here we want to show the same result but for the Property $IH_2$. However, due to the large amount of variables appearing in the Definition \ref{ih2propertydef}, the explicit formula for each of those dense sets is quite messy. For that reason, we will instead show that $\mathbb{P}(\mathcal{F})$ forces $\mathcal{F}^{Gen}$ to satisfy $IH_2$. Before doing this, it is worth  pointing out some last remarks regarding the forcing $\mathbb{P}(\mathcal{F}).$
\begin{rem}By virtue of the the Lemmas \ref{densedomainlemma},and \ref{lemmaih1dense}, we have that $$\mathbb{P}(\mathcal{F})\Vdash\text{\say{ $\mathcal{F}^{Gen}$ is construction scheme over $\gamma+\omega$ which satisfies $IH_1$}}.$$
\end{rem}
\begin{rem} $\mathbb{P}(\mathcal{F})\Vdash\text{\say{ $\mathcal{F}\subseteq \mathcal{F}^{Gen}$}}$ due to the Lemma \ref{lemmaFdenseinPF}. This means that $\rho_{F^{Gen}}|_{\gamma^2}$ is forced to be equal to $\rho_\mathcal{F}$. Therefore, there is no risk of confusion by referring to both of these two ordinal metrics simply as $\rho$.
\end{rem}
\begin{rem}If $p\in \mathbb{P}_k(\mathcal{F})$ for some $k\in\omega$, then $p\vDash\text{\say{ $p\in \mathcal{F}^{Gen}_k$ }}.$ In particular, $\mathcal{F}(p)$ is forced to be equal to $\mathcal{F}^{Gen}|_p$. In other words, $\rho|_{p^2}=\rho_{\mathcal{F}(p)}$. So, again, there is no risk of confusion on referring to $\rho_{\mathcal{F}(p)}$, simply as $\rho$.    
\end{rem}
\begin{rem}\label{remarkdecidingkclosure}Suppose that $p\in \mathbb{P}_l(\mathcal{F})$, $\beta\in p$ and $k\leq l$. Then for each $H\in \mathcal{F}_j(p)$ it happens that $p\Vdash\text{\say{ $H\cap (\beta+1)=(\beta)_k $ }}. $\footnote{Here, by $(\beta)_k$ we mean $k$-closure of $\beta$ calculated inside $\gamma+\omega$ with respect to $\rho_{\mathcal{F}^{Gen}}$.} That is, $p$ decides the value of $(\beta)_k$. Furthermore, this value can be calculated using $\mathcal{F}(p)$. From this, it follows that $p$ also knows who is $\lVert \beta\rVert_k$, $\Xi_\beta(k)$, $Bl_k(\beta)$ ($=Bl(\lVert \beta\rVert_k,m_k)$) and $Good_k(\beta)$. Therefore, there should not be any confusion while working explicitly with this sets instead of with their names, of course, as long as the conditions at the beginning of this remark hold.
\end{rem}
\begin{lemma}\label{lemmaforcingroot} Let $q\in \mathbb{P}(\mathcal{F})$. Then  there is $p\leq q$ so that $|p\cap \gamma|=r_{k_p+1}$.
\begin{proof} From Lemma \ref{densedomainlemma}, it is easy to get $q'\leq q$ so that $k_{q'}>k_q$.  Now, according to Lemma \ref{lemmaih1dense}, there is $p'\in \mathbb{P}(\mathcal{F})$ so that $p'\leq q'$ and $p'\cap \gamma=R(p')$. Put $p=p'_0$. Then $p\in \mathbb{P}(\mathcal{F})$ and $p\leq q$. Furthermore, $|p\cap \gamma|=|R(p')|=r_{k_{p'}}=r_{k_p+1}.$
    
\end{proof}
    
\end{lemma}
\begin{theorem}\label{superprincipaltheorem}Suppose  that $\gamma$ is countable and $\mathcal{F}$ satisfies $IH_1$ and $IH_2$. Then $$\mathbb{P}(\mathcal{F})\Vdash \text{\say{ $\mathcal{F}^{Gen}$ satisfies $IH_2$ }}.$$
\begin{proof}According to the Definition \ref{ih2propertydef}, our goal is to prove that $$\mathbb{P}(\mathcal{F})\Vdash\text{\say{ $IH_2(\delta,\mathcal{F}^{Gen}) $ holds }}$$ for any $\delta\in \text{Lim}\cap (\gamma+\omega)$. If $\mathcal{F}$ satisfies that there are infinitely many $l\in \omega$ for which there is $C\in \text{FIN}(D_\delta)$ which is fully captured at level $l$,  it is easy to see that  $\mathbb{P}(\mathcal{F})\vDash\text{\say{ $IH^A_2(\gamma,\mathcal{F}^{Gen})$ holds }}$. This is because $\mathbb{P}(\mathcal{F})\Vdash \text{\say{  $\mathcal{F}\subseteq \mathcal{F}^{Gen}$ }}$. Therefore,  the only interesting case happens when:
\begin{center}With respect to $\mathcal{F}$, there  is $k''\in \omega$ such that for each $l\geq k''$ there is no $C\in \text{FIN}(D_\gamma)$ which is fully captured at level $l$.
\end{center}
From now on, let us assume that this case holds. \\

\noindent
\underline{Claim}: $\mathbb{P}(\mathcal{F})\Vdash \text{\say{ $IH^B_2(\delta,\mathcal{F}^{Gen})$ holds }}.$
\begin{claimproof}[Proof of claim]Let $q\in \mathbb{P}(\mathcal{F})$, $\beta\in [\delta,\gamma+\omega)$, $2\leq k\in \omega$ and a name $\dot{T}$ for an element of $Good_k(\beta)$. Because of Lemma \ref{densedomainlemma} , we may assume without loss of generality that $\beta\in q$ and $q\in \mathbb{P}_{k'}(\mathcal{F})$ for some $k'>k,k''$. According to the Remark \ref{remarkdecidingkclosure}, $Good_k(\beta)$ is fully determined by $q$. In particular, this means that  there is $T$ such that $q\Vdash\text{\say{ $T=\dot{T}$ }}$.

We need to show that there is $l>{k'}$ and $p\leq q$ with the following properties:
\begin{itemize}
    \item $p\Vdash\text{\say{ $|(\beta)_{l-1}\cap \delta|=r_l$ }}.$
    \item There is $j\in \omega$ so that $p\Vdash \text{\say{ $j=j(\delta,T,\beta,k,l)$ }}$.
    \item $p\Vdash\text{\say{ $\Xi_\beta(l)=j$ }}$.
    \item If $j>0$ then there is $C\in [\gamma+\omega]^j$ so that $p\Vdash \text{\say{ $(C,\delta)\checkmark (T,\beta,k,l)$ and $C\subseteq (\beta)_l$}}$.
\end{itemize}

We will divide the rest of the proof into two cases.\\

\noindent
\underline{Case 1}: $\delta=\gamma$.
\begin{claimproof}[Proof of case] By virtue of the Lemma \ref{lemmaforcingroot}, there is $q'\in \mathbb{P}(\mathcal{F})$ so that $q'\leq q$ and $|q'\cap \gamma|=r_{k_q+1}$. Let $l=k_q+1$. Note that the objects needed to calculate $j(\gamma,T,\beta,k,l)$ are $\mathcal{F}^{Gen}|_{\gamma}$, $T$, $(\beta)_{l-1}$, $k $ an $l$. All of this objects are already decided by $q'$. Therefore, there is $j\in \omega$ so that $q'\Vdash\text{\say{ $j=j(\gamma,T,\beta,k,l)$}} $.\\

\noindent
\underline{Subcase 1.1}: If $j=0$.
\begin{claimproof}[Proof of subcase] Consider $F\in \mathcal{F}_l$ so that $q'\cap \gamma\sqsubseteq F$. Then $R(F)=q'\cap \gamma$. Let $\alpha=\min(F_0\backslash R(F))$ and $p=Cut_\alpha(F)$. It is straightforward that $R(p)=q'\cap\gamma$ and $q'=p_0$. This means that $p\leq q$. Moreover, as $\beta\in q\backslash \gamma$ then $p\Vdash\text{\say{$|(\beta)_{l-1}\cap \delta|=r_l$ }}$ and  $p\Vdash\text{\say{ $\Xi_\beta(l)=0=j$ }}$. In this case, there is nothing more to prove.
\end{claimproof}
\noindent
\underline{Subcase 1.2}: If $j>0$.
\begin{claimproof}[Proof of subcase] By definition of $j$, there is $C\in [\gamma+\omega]^{j}$ so that $(C,\gamma)\checkmark (T,\beta,k,l)$. According the point (1) of Definition \ref{checkmarkdef}, $C\subseteq D_\gamma$ and $C$ is captured at level $l$. Let $F\in \mathcal{F}_l$ be so that $C\subseteq F$. Then, by the point (3) of Definition \ref{checkmarkdef}, it follows that $R(F)=(\beta)_{l-1}\cap\gamma=q'\cap \gamma$. In this case, let $\alpha=\min(F_j\backslash R(F))$ and $p=Cut_\alpha(F)$. As in the previous case,  we have  $p_j=q'$ and $R(p)=q'\cap \gamma$. Thus, $q'\leq p_j$. Again, as $\beta\in q\backslash \gamma$ then $p\Vdash\text{\say{ $|(\beta)_{l-1}\cap \delta|=r_l$ }}$ and $p\Vdash \text{\say{ $\Xi_\beta(l)=j$ }}$. In order to finish, note that $p\cap \gamma=\bigcup\limits_{i<j} F_i$. In particular, $C\subseteq p\cap \gamma\subseteq (\beta+1)\cap p=(\beta)_l$. That is, $$p\Vdash \text{\say{ $(C,\gamma)\checkmark (T,\beta,k,l)$ and $C\subseteq (\beta)_l$ }}.$$

\end{claimproof}
\end{claimproof}
\noindent
\underline{Case 2}: If $\delta<\gamma$.
\begin{claimproof}[Proof of case] If $\beta<\gamma$ there is nothing to do. This is because $\mathcal{F}$ is forced to be contained in $\mathcal{F}^{Gen}$ and $\mathcal{F}$ already satisfies $IH_2$. Therefore, we may assume that $\beta\geq \gamma$. Appealing to Lemma \ref{densedomainlemma}, we may assume without loss of generality that $q\cap [\delta,\gamma)\not=\emptyset$. Let $\alpha=\alpha_q=\max(q\cap \gamma)\geq \delta$. The plan is apply the property $IH_2$ to $\alpha$ (inside $\mathcal{F}$). For this purpose, let $\mathcal{Q}=\mathcal{Q}_{k'}(\beta,T)$ and $T'= \text{Transf}_{k'}(T,\alpha,\beta)$. 

By our assumptions $IH^A_2(\delta,\mathcal{F})$ does not hold. Therefore, $IH^B_2(\delta,\mathcal{F})$ holds. This means that there is $l\geq k'$ for which $j'=j(\delta,T',\alpha,k',l)<n_l$ and such that:
\begin{enumerate}[label=(\Roman*)]
    \item $|(\alpha)_{l-1}\cap \delta|=r_l,$
    \item $\Xi_\alpha(l)=j,$
    \item If $j>0$ then there is $C\in [\gamma]^j$ for which $(C,\delta)\checkmark (T',\alpha,k',l)$ and $C\subseteq (\alpha)_l$.
\end{enumerate}
Before we continue, we remark that:
\begin{itemize}
    \item $q\Vdash \text{\say{ $\alpha\in\beta\cap q=(\beta)^-_{k'}$ }}$.
    \item $q\Vdash \text{\say{ $(\beta)_{k'}\backslash (\alpha+1)=[\gamma,\beta+1)$ }}$. Since $[\gamma,\beta+1)$ is already an interval in $\gamma+\omega$, it follows that $q\Vdash\text{\say{ $(\beta)_{k'}$ is an interval in $(\beta)_{l-1}$}}$.
    \item $q\Vdash \text{\say{ $(\beta)^-_{l-1}\cap \alpha=(\alpha)_{l-1}$ }}$ because $l>k'$. In particular, this means that $q\Vdash \text{\say{ $|(\beta)_{l-1}\cap \delta|=|(\alpha)_{l-1}\cap \delta|=r_l$ }}.$
    \item  Since $j'\geq 0$ and $\Xi_\alpha(l)=j'$ then $q\Vdash \text{\say{ $\Xi_\beta(l)=j'$ }}$ by means of the Lemma \ref{lemmaxi}.
    
\end{itemize}
Due to the first two points, the hypotheses of Corollary \ref{trasnfcoro} are fulfilled. In particular, this means that $$q\Vdash \text{\say{$j(\delta,T,\beta,k,l)\leq j'$ }}.$$

We will divide the rest of the proof into two subcases.\\

\noindent
\underline{Subcase 2.1}: If $j'=0$.
\begin{claimproof}[Proof of subclaim]

As $j(\delta,T,\beta,k,l)$ is always non-negative then $$q\Vdash \text{\say{$j(\delta,T,\beta,k,l)=0$ }}.$$ Because of this and by the last two points of the previous remark, $p=q$ forces everything that we want.
\end{claimproof}
\noindent
\underline{Subcase 2.2}: If $j'>0$.
\begin{claimproof}[Proof of subcase]
Let $F\in \mathcal{F}_l$ be such that $\alpha\in F$. Since $j>0$, there is $C\in [\gamma]^j$ for which $(C,\delta)\checkmark (T',\alpha,k',l)$ and $C\subseteq (\alpha)_l$. Note that $q\vDash\text{\say{ $(\alpha)_l\subseteq (\beta)_l$}}$. In this way, $q\Vdash \text{\say{ $C\subseteq (\beta)_l$ }}$. Our next task will be to extend $q$ to a condition forcing that $T\bigstar(\beta, C(0),k,l-1)$. By appealing to Lemma \ref{transfequiv}, it suffices to find $p\leq q$ which forces that $\lVert \beta\rVert_{l-1}\leq \lVert C(0)\rVert_{l-1}$ and $\lVert \beta\rVert_{l-1}\in (\lVert \xi\rVert_{l-1})_{k'}$.

As $(C,\delta)\checkmark (T',\alpha,k',l)$ then $T'\bigstar (\alpha,C(0),k',l-1)$. Let $\xi\in F_j$ be such that $\lVert \xi\rVert_{l-1}=\lVert C(0)\rVert_{l-1}$. By Remark \ref{bigstarremark}, $T'\bigstar (\alpha,\xi,k',l-1)$. Thus, there are $(\mathbb{I},z)\in T$ and $(G,\mathcal{I})\in \mathcal{Q} $ so that $(\mathbb{I}_G,G(z))$ testifies the previous guessing relation. That is:
\begin{enumerate}[label=$(\alph*)$]
    \item $\lVert \alpha\rVert_{l-1}\leq \lVert \xi\rVert_{l-1}$ and $\lVert \xi\rVert_{k'}=G(z)$. Since both $\xi$ and $\alpha$ are in $F_j$ then $\rho(\alpha,\xi)\leq l-1$. Thus, this point implies that $\alpha\leq \xi$.
    \item $(\xi)_{k'}\circ \mathbb{I}_G\in Bl(0,(\xi)_{l-1}).$
    \item $(\xi)_{l-1}(\lVert \alpha\rVert_{l-1})\in (\xi)_{k'}.$ As $\alpha,\xi\in F_j$ then $(\xi)_{l-1}(\lVert \alpha\rVert_{l-1})=F_j(\lVert \alpha\rVert_{l-1})=\alpha$. Therefore, this point implies that $\alpha\in (\xi)_{k'}$.
\end{enumerate}
According to the point (b), $(\xi)_{k'}[\mathbb{I}_G(0)]=(\xi)_{k'}[\,(\lVert \alpha\rVert_{k'},\lVert \beta\rVert_{k'}]]$  is an interval in $(\xi)_{l-1}$. Moreover, by the point (a), this interval is contained in $F_j\backslash R(F)$. Therefore, it is also an interval in $F$. Let $\alpha'=(\xi)_{k'}[\lVert \alpha\rVert_{k'}+1]$. In other words, $\alpha'$ is the successor of $\alpha$ inside $(\xi)_{k'}$. Because of this, it follows that $(\alpha')^-_{k'}=(\alpha)_{k'}=q\cap \gamma$. In this way, $p=Cut_{\alpha'}(F)\leq q$ due to the Lemma \ref{lemmacut}. \\

\noindent
\underline{Subclaim 1}: $p\Vdash\text{\say{ $\lVert \beta\rVert_{l-1}\leq \lVert C(0)\rVert_{l-1}$ and $\lVert \beta\rVert_{l-1}\in (\lVert C(0)\rVert_{l-1})_{k'}$ }}$.
\begin{claimproof}[Proof of subclaim] 
Let $\phi: F\longrightarrow p$ be the increasing bijection. By the definition of the Cut function, it follows that $\phi(\alpha')=\gamma$. Now, since $(\xi)_{k'}[\mathbb{I}_F(0)]$ is an interval in $F$, then $\phi[\,(\xi)_{k'}[\mathbb{I}_G(0)]\,]=(\phi(\xi))_{k'}[\mathbb{I}_G(0)]$ is an interval in $p$. Furthermore, the first point of this interval is $\gamma$. In this way, as $p\backslash \gamma$ is an initial segment of $[\gamma,\gamma+\omega)$, we conclude that $\phi((\xi)_k'(\lVert \beta\rVert_{k'}))$ (the last point of the interval $\phi[(\xi)_{k'}[\,\mathbb{I}_G(0)]\,]$ ) is equal to $$\gamma+(\lVert \beta\rVert_{k'}-\lVert \alpha\rVert_{k'})=\beta.\footnote{All the definitions related to $IH_2$ were originally formulated to achieve this single step of the proof.} $$

As direct consequence of this equality, p forces that $\beta\leq\phi(\xi)$ and $\rho(\beta,\phi(\xi))\leq {k'}$. It then follows that $\lVert \beta\rVert_{l-1}\leq \lVert \phi(\xi)\rVert_{l-1}$  and $\lVert \beta\rVert_{l-1}\in (\lVert\phi( \xi)\rVert_{l-1})_{k'}$. Since $\lVert \phi(\xi)\rVert_{l-1}=\lVert \xi\rVert_{l-1}=\lVert C(0)\rVert_{l-1}$, we are done.
\end{claimproof}
By virtue of the previous subclaim,  $p\Vdash\text{\say{ $ T\bigstar (\beta,C(0),k,l-1)$ }}$. It then follows that $p\Vdash\text{\say{ $(C,\delta)\checkmark (\delta,T,\beta,k,l)$ }}$. As $|C|=j'$, we get that $$p\Vdash \text{\say{ $j'\leq j(\delta,T,\beta,k,l)$ }}.$$ Hence, these two numbers are equal. In conclusion, we have shown that there are $p\leq q$ and  $j\in \omega$ (namely, $j'$), so that:
\begin{itemize}
    \item $p\Vdash\text{\say{ $|(\beta)_{l-1}\cap \delta|=r_l$ }}$ and $p\Vdash\text{\say{ $\Xi_\beta(l)=j$ }}$.\footnote{ This was already forced by $q$.}

    \item There is $C\in [\gamma+\omega]^j$ so that $p\Vdash \text{\say{ $(C,\delta)\checkmark (T,\beta,k,l)$ and $C\subseteq (\beta)_l$}}$. 
\end{itemize}
Thus, the proof of this subcase is over.
\end{claimproof}
\end{claimproof}
\end{claimproof}
\end{proof}
\end{theorem}

It is not hard to see that the property forced in the previous theorem can be coded by countable many dense subsets of $\mathbb{P}(\mathcal{F})$. In this way, we have the following corollary.
\begin{corollary}\label{lastcorollary}Suppose that $\gamma$ is a limit ordinal and $\mathcal{F}$ satisfies $IH_1$ and $IH_2$. There is a countable family of dense sets in $\mathbb{P}(\mathcal{F})$ so that if $\mathcal{G}$ is a filter intersecting all of them, then $\mathcal{F}^\mathcal{G}$ is a construction scheme over $\gamma+\omega$ which contain $\mathcal{F}$ and satisfies both $IH_1$ and $IH_2$.
\end{corollary}

Now we have proved all the necessary results to show that there is a fully capturing construction scheme over $\omega_1$ of type $\tau$. Namely, Propositions \ref{IH1Fomega} and \ref{IH2omega} , Lemmas \ref{unionschemelemmadeeperIH1} and \ref{unionschemelemmadeeperIH2}, and Corollary \ref{lastcorollary}. From them we conclude that there is a sequence $\langle \mathcal{F}^\delta\rangle_{\delta\in \text{Lim}}$ so that for all $\delta<\gamma\in \text{Lim}$, the following properties hold:
\begin{itemize}
    \item $\mathcal{F}_\delta$ is a construction scheme over $\delta$ satisfying $IH_1$ and $IH_2$,
    \item $\mathcal{F}_\delta\subseteq \mathcal{F}_\gamma$.
\end{itemize}
It follows that $\mathcal{F}=\bigcup\limits_{\delta\in \text{Lim}}\mathcal{F}^\delta$ is a fully capturing construction scheme. This is due to the Theorem \ref{ih2omega1fullycapturing}. 
\begin{theorem}The $\Diamond$-principle implies $FCA$. 
\end{theorem}
Now fix $\mathcal{P}$ an aribtrary partition of $\omega$ compatible with $\tau$. Consider the following variations of the the properties $IH_1$ and $IH_2$.

\begin{definition}[The $IH_1(\mathcal{P})$ property] We say that $\mathcal{F}$ satisfies $IH_1(\mathcal{P})$ if for all $A\in\text{FIN}(\gamma)$, $P\in \mathcal{P}$ and $\alpha\in \gamma$, there is $F\in \mathcal{F}$ with $\rho^F\in P$ such that $IH_1(\alpha,A,F)$ holds.
\end{definition}
\begin{definition}[The $IH_2(\mathcal{P})$ property]Given $\delta\in\text{Lim}\cap \gamma$ and $P\in \mathcal{P}$, we say that $IH_2(\delta,\mathcal{F})$ holds if one of the two following mutually excluding conditions occurs:
\begin{enumerate}[label=$(\Alph*)$]
\item There are infinitely many $l\in P$ for which there is $C\in \text{FIN}(D_\delta)$ which is fully captured at level $l$. In this case, we will say that $IH^A_2(P,\delta,\mathcal{F})$ holds.
\item For all $\beta\in [\delta,\gamma)$, $2\leq k\in \omega$  and $T\in Good_k(\beta)$ there are infinitely many $k\leq l\in P$ for which $j=j(\delta,T,\beta,k,l)<n_l$ and such that:
\begin{enumerate}[label=$(B.\arabic*)$]
    \item $|(\beta)_{l-1}\cap \delta|=r_l,$
    \item $\Xi_\beta(l)=j,$
    \item  If $j>0$ then there is $C\in [\delta]^j$ for which $(C,\delta)\checkmark(T,\beta,k,l)$ and $C\subseteq (\beta)_l.$
\end{enumerate}
In this case we will say that $IH^B_2(P,\delta,\mathcal{F})$ holds. 
\end{enumerate}  
Finally, we say that  $\mathcal{F}$ satisfies $IH_2(\mathcal{P})$ if $IH_2(\delta,\mathcal{F})$ holds for any $\delta \in \text{Lim}\cap \gamma$ and each $P\in \mathcal{P}$.
\end{definition}
All the results stated in this section are true when $IH_1$ is changed by $IH_1(\mathcal{P})$ and $IH_2$ is changed by $IH_2(\mathcal{P})$. Moreover, the proofs are completely similar. In this way, we deduce theorem \ref{diamondfcaparttheorem}.

\begin{theorem}\label{diamondfcapartheoremproooof} The $\Diamond$-principle implies $\text{FCA}(part)$.
\end{theorem}

\section{P-ideal Dichotomy and capturing axioms}\label{sectionpidschemes}
Throughout this work, we have seen a handful of applications of construction schemes. Many of the objects that we have constructed using capturing construction schemes already contradict to some degree the $P$-ideal dichotomy. For example, since the Suslin Hypothesis follows from $PID$ (see \cite{partitionpropertiesch} and \cite{dichotomypidfirst}), this axiom is incompatible with $FCA$. In addition, since $PID$ implies that any gap is indestructible, $CA_3$ also contradicts $PID$. However, note that none of the applications that we have of $CA_2$ is enough to contradict the $P$-ideal dichotomy. For now, the best we have is that $PID+\min(\mathfrak{b},\,cof(\mathcal{F}_\sigma)\,)>\omega_1$ is incompatible with the existence of a $2$-capturing construction scheme. This is because such assertion is equivalent to the nonexistence of a sixth Tukey type (see \cite{CombinatorialDichotomiesandCardinalInvariants}). The purpose of this section is to prove that $PID$ is in fact incompatible with $CA_2$.

For the remainder of this section, we will assume that $PID$ holds and that $\mathcal{F}$ is a construction scheme. Our plan is to define a $P$-ideal $\mathcal{I}$ from $\mathcal{F}$ that contradicts $PID$ whenever $\mathcal{F}$ is $2$-capturing. This ideal $\mathcal{I}$ will be defined as family of countable sets of the orthogonal ideal of some other ideal $\mathcal{H}$. When applying the second alternative of $PID$ to the ideal $\mathcal{I}$, we will need to calculate who is $\mathcal{I}^\perp$. In general, this may be a difficult task. However, the following concept will facilitate these calculations.

\begin{definition}Let $\mathcal{H}$ be an ideal of countable sets in $\omega_1$. We will say that $\mathcal{H}$ is \textit{Frechet in $[\omega_1]^{\omega}$} if for any $A\in \mathcal{H}^+$, there is an infinite $B\subseteq A$ such that $B\in \mathcal{H}^\perp$.
    
\end{definition}

\begin{proposition}Let $\mathcal{H}$ be an ideal of countable sets in $\omega_1$ and assume that $\mathcal{H}$ is Frechet in $[\omega_1]^{\omega_1}$. Then $(\mathcal{H}^\perp\cap [\omega_1]^{\leq\omega})^\perp\cap [\omega_1]^{\omega}=\mathcal{H}\cap [\omega_1]^{\omega}$.
\begin{proof}We will only show that the  inclusion from right to left holds. For this purpose, let $A\in (\mathcal{H}^\perp\cap[\omega_1]^{\leq \omega})^\perp\cap[\omega_1]^{\omega}$. Assume towards a contradiciton that $A\not\in \mathcal{H}$. Then $A\in \mathcal{H}^+\cap [\omega_1]^\omega$. Since $\mathcal{H}$ is Frechet, there is an infinite $B\subseteq A$ such that $B\in \mathcal{H}^\perp$. In particular, $A\cap B$ is infinite, and this contradicts the fact that $A\in (\mathcal{H}^\perp)^\perp$. Thus, the proof is over.
    
\end{proof}
\end{proposition}
We now define the ideal $\mathcal{I}$ which will be used in the proof of Theorem \ref{fnotcapturingpidtheorem}.

\begin{definition} Given $n\in\omega$ and $\alpha\in \omega_1$, we define $$H_n(\alpha)=\{\xi<\alpha\,:\,\forall m>n\,(\,\Xi_\alpha(m)=-1\textit{ or }\Xi_\xi(m)\leq \Xi_\alpha(m)\,)\}.$$
Additionally, we define the ideal generated by the family $\{H_n(\alpha)\,:\,n\in\omega\textit{ and }\alpha\in\omega_1\}$ as $\mathcal{H}$. Lastly, we define $\mathcal{I}=\mathcal{H}^\perp\cap [\omega_1]^{\leq \omega}.$
\end{definition}
\begin{rem}$[\omega_1]^{<\omega}$ is included in both $\mathcal{H}$ and $\mathcal{I}.$
    
\end{rem}
The following lemma is easy.
\begin{lemma}\label{lemmaHmpideal1}Let $\alpha,\beta,\gamma<\omega_1$ and $m,n\in \omega$. Then:
\begin{enumerate}[label=$(\arabic*)$]
    \item If $m<n$, then $H_m(\alpha)\subseteq H_n(\alpha).$
    \item If $\beta<\gamma$ and $n>\rho(\beta,\gamma)$, then $H_n(\gamma)\cap \beta\subseteq H_n(\beta).$
\end{enumerate}
\begin{proof}The point (1) follows directly from the definition. In order to prove (2) we will show that $\omega_1\backslash H_n(\beta)\subseteq \omega_1 \backslash (H_n(\gamma)\cap \beta)$. Let $\xi \in \omega_1\backslash H_n(\beta) $. If $\xi\geq \beta$ we are done, so let us assume that $\xi<\beta.$ Then there is $m>n$ for which $\Xi_\xi(m)>\Xi_\beta(m)\geq 0$. Note that $m>\rho(\beta,\gamma)$. Thus, by the point (c) of Lemma \ref{lemmaxi} we conclude that $\Xi_\beta(m)=\Xi_\gamma(m)$. In this way, $\xi\not\in H_n(\gamma)$. This finishes the proof.
    
\end{proof}
    
\end{lemma}
\begin{rem}By the point (1) of the previous lemma we have that for any $A\in \mathcal{H}$ there are $n\in\omega$ and $\alpha_1,\dots,\alpha_n\in \omega_1$ for which $A\subseteq \bigcup\limits_{i<n}H_n(\alpha_i).$
    
\end{rem}

\begin{definition}Given $\beta\in \omega_1$ we define $\mathcal{H}|_\beta$ as the ideal generated by the family $\{H_n(\alpha)\,:\,\alpha\leq \beta \textit{ and }n\in\omega\}.$
    
\end{definition}
\begin{lemma}\label{lemmaHmpideal2}Let $\beta<\omega_1$ and $A\subseteq \beta$.
\begin{enumerate}[label=$(\alph*)$]
    \item If $A\in (\mathcal{H}|_\beta)^{\perp}$, then $A\in \mathcal{H}^\perp.$ Furthermore, since $A$ is countable then $A\in \mathcal{I}.$
    \item If $A\in (\mathcal{H}|_\beta)^{+}$, then $A\in \mathcal{H}^+.$
\end{enumerate}
\begin{proof}\begin{claimproof}[Proof of $(a)$] It is enough to show that $A\cap H_n(\gamma)=^*\emptyset$ for any $\beta\leq\gamma\in \omega_1$ and each $n\in\omega$. Furthermore, we may assume that $n>\rho(\beta,\gamma)$ due to the point (1) of Lemma \ref{lemmaHmpideal1}. Indeed, by virtue of the point (2) of the same lemma, we have that $$A\cap H_n(\gamma)=A\cap \beta\cap H_n(\gamma)\subseteq A\cap H_n(\beta)=^*\emptyset.$$
\end{claimproof}
\begin{claimproof}[Proof of $(b)$]Suppose towards a contradiction that $A\notin \mathcal{H}^+$. Then there are $n\in\omega$ and $\alpha_1,\dots,\alpha_n\in \omega_1$ for which $A\subseteq \bigcup\limits_{i<n}H_n(\alpha_i).$ Without any loss of generality we may assume that there is $j<n$ so that, for any $i<n$, $\alpha_i\leq \beta $ if and only if $i<j$. Furthermore, we may assume that $n>\rho(\beta,\alpha_i)$ for any $i\geq j$. From this, it follows that $$A\subseteq\big( \bigcup\limits_{i<j}H_n(\alpha_i)\big)\cup \big(\bigcup\limits_{j\leq i<n}H_n(\alpha_i)\cap \beta\big)\subseteq \big(\bigcup\limits_{i<j}H_n(\alpha_i)\big)\cup \big(\bigcup\limits_{j\leq i<n}H_n(\beta)\big). $$
In this way, $A\in \mathcal{H}|_\beta$ which is a contradiction.
    
\end{claimproof}

\end{proof}
    
\end{lemma}

\begin{proposition}\label{propIfrechet}$\mathcal{H}$ is Frechet in $[\omega_1]^\omega.$
\begin{proof}Let $A\in \mathcal{H}^+\cap [\omega_1]^\omega$ and $\beta\in\omega_1$ be such that $A\subseteq \beta$. Note that $A\in (\mathcal{H}|_\beta)^+$. Since $\mathcal{H}|_\beta$ is a countably generated ideal, there is $B\in [A]^\omega$ such that $B\in (\mathcal{H}|_\beta)^\perp$. According to the point $(a)$ of Lemma \ref{lemmaHmpideal2}, $A\in \mathcal{I}$. Thus, the proof is over.
    
\end{proof}
\end{proposition}

\begin{proposition}\label{propIpideal}
$\mathcal{I}$ is a $P$-ideal.
\begin{proof}Let $\langle A_n\rangle_{n\in\omega}$ be an increasing sequence of elements of $\mathcal{I}$. Then there is $\beta\in\omega_1$ for which $\bigcup\limits_{n\in\omega}A_n\subseteq \beta$. Let us enumerate $\beta+1$ as $\langle \beta_n\rangle_{n\in\omega}$. We then define $$A=\bigcup\limits_{n\in\omega}\big(A_n\backslash \big(\bigcup\limits_{i\leq n}H_n(\beta_i)\big)\big).$$
By the Lemma \ref{lemmaHmpideal1}, it follows that $A$ is a pseudo-union of $\langle A_n\rangle_{n\in\omega}$ and $A\in (\mathcal{H}|_\beta)^\perp$. In particular, this means that $A\in \mathcal{I}$ due to the point (a) of Lemma \ref{lemmaHmpideal2}
    
\end{proof}
\end{proposition}
\begin{theorem}\label{fnotcapturingpidtheorem}$\mathcal{F}$ is not $2$-capturing.
\begin{proof} Suppose towards a contradiction that $\mathcal{F}$ is $2$-capturing. According to the Proposition \ref{propIpideal}, $\mathcal{I}$ is a $P$-ideal. Since we are assuming $PID$, there are two possibilities for $\mathcal{I}$. We will finish the proof by showing that both of them lead to a contradiciton.\\\\
\noindent
\underline{Case 1}: There is an uncountable $Y\subseteq \omega_1$ such that $[Y]^{\omega}\subseteq \mathcal{I}.$
\begin{claimproof}[Non satisfaction of case]Consider the coloring $c:[Y]^2\longrightarrow 2$ given by:
$$c(\{\alpha,\beta\})=0\text{ if and only if }\{\alpha,\beta\}\text{ is captured}.$$
Using the partition relation $\omega_1\longrightarrow (\omega_1,\omega+1)^2_2$ and the fact that $\mathcal{F}$ is $2$-capturing, we can conclude that there is an increasing sequence $\langle \alpha_i\rangle_{i<\omega+1}\subseteq Y$ such that $\{\alpha_i,\alpha_j\}$ for any two distinct $i,j< \omega+1$. In particular, $\{\alpha_i,\alpha_\omega\}$ is captured for any $i<\omega$. By lemma \ref{lemmaxi}, we get that $A=\{\alpha_i\,:\,i<\omega\}\subseteq H_0(\alpha_\omega)$. In other words, $A\in \mathcal{I}$. This is a contradiction, so this case can not occur.
    
\end{claimproof}
\noindent
\underline{Case 2}: There $Y$ an uncountable subset of $\omega_1$ such that $[Y]^{\omega}\subseteq \mathcal{I}^\perp.$
\begin{claimproof}[Non satisfaction of case]
According to Proposition \ref{propIfrechet}, $\mathcal{I}$ is Frechet over $[\omega_1]^\omega$. From this it follows that  $[Y]^{\omega}\subseteq \mathcal{H}$. Let $$P=\{l\in\omega\,:\,\exists \alpha\in Y\,(\,\Xi_\alpha(l)>0\,)\}.$$ By refining $Y$,  we may assume that for any $l\in P$ there are uncountable many $\alpha\in Y$ for which $\Xi_\alpha(l)>0$. In this way, we can recursively construct sequences $\langle A_\beta\rangle_{\beta\in\omega_1}\subseteq [Y]^\omega$, $\langle B_\beta\rangle_{\beta\in \omega_1}\subseteq \text{FIN}(\omega_1)$ and $\langle n_\beta\rangle_{\beta\in\omega_1}$ so that the following properties hold for any two $\beta<\gamma\in \omega_1$:
\begin{enumerate}[label=$(\arabic*)$]
    \item $A_\beta<B_\beta<A_\gamma$. In particular $B_\beta\cap B_\gamma=\emptyset.$
    \item For all $l\in P$ there is $\alpha\in A_\beta$ for which $\Xi_\alpha(l)>0$.
    \item $A_\beta\subseteq \bigcup\limits_{\alpha\in B_\beta}H_{n_\beta}(\alpha).$
    \item $B_\beta\cap Y\not=\emptyset.$ 
\end{enumerate}
By refining such sequence, we can suppose that there is $n\in\omega$ such that $n_\beta=n$ for any $\beta\in\omega_1$. Since $\mathcal{F}$ is $2$-capturing, there are $\beta<\gamma\in\omega_1$ such that $\{B_\beta,B_\gamma\}$ is captured at some level $l>n$. Now, $B_\beta$ and $B_\gamma$ are disjoint. Thus, $\Xi_\alpha(l)=1$ for any $\alpha\in B_\gamma$. Therefore, $l\in P$ due to the point (4). By virtue of the point (2), there must be $\xi\in A_\beta$ such that $\Xi_\xi(l)>0$. On the other hand, $\Xi_\alpha(l)=0$ for all $\alpha\in B_\beta.$ Consequently, $\xi \not\in H_n(\alpha)$ for any such $\alpha$. This contradicts the point (3), so this case can not occur.
    
\end{claimproof}
    
\end{proof}
\end{theorem}
\begin{corollary}(Under $PID$)\label{coropidnotca2} There are no $2$-capturing construction schemes.
    
\end{corollary}

%% file: chapters/problems.tex
\chapter{Open problems}\label{problemschapter}

\section{Problems about gaps}

One of the most important concepts introduced in Section \ref{sectiongaps} was the Luzin representation of a partial order and what it means for a Luzin family to code a partial order (see Definition \ref{Luzinrepdef}) Indeed, in order to prove that there is a Luzin-Jones family $\mathcal{A}$ in $ZFC$, we ended up proving that such family $\mathcal{A}$we can differentiatet is natural to ask if the same is true for any Luzin-Jones family.

\begin{problem}\label{problemomega1code}Is there a Luzin-Jones family which doesn't code $\omega_1$?
\end{problem}

A partial answer to this question is \say{no}. If $\mathfrak{b}>\omega_1$ then any Luzin family is Jones and it codes $\omega_1$. Thus, a related problem would be whether the existence of a Luzin-Jones family which doesn't code $\omega_1$ follows from the equality $\mathfrak{b}=\omega_1$. Although such a family might exist, it is still possible that there is a partial order $X$ which characterizes the Jones property in terms of coding: Formally:

\begin{problem}Is there a partial order $X$ so that for any Luzin family, coding $X$ and being Jones are equivalent statements?
\end{problem}

Given an almost disjoint family $\mathcal{A}$, let us define the \textit{spectrum of $\mathcal{A}$}, namely $Spec(\mathcal{A})$, as the class of all partial orders $(X,<)$ such that $\mathcal{A}$ codes $X$ (Definition  \ref{Luzinrepdef}). The class $Spec(\mathcal{A})$ can be regarded as a measurement of the \say{malleability} of $\mathcal{A}$. From the point of view of the author, it is an interesting problem to know how much we can differentiate an almost disjoint family from another one by analysing their spectrums.

\begin{problem}Which are the possible values of $Spec(\mathcal{A})$ for a given Luzin family $\mathcal{A}$? 
\end{problem}

Let $\mathcal{X}_{\omega_1}$ denote the class of $\omega_1$-like orders. In Theorem \ref{luzinreptheorem}, we have seen that the Luzin-Jones family $\mathcal{A}$ that we constructed codes each element of $\mathcal{X}_{\omega_1}$. Following with Problem \ref{problemomega1code}, we may ask how strong is the property of coding each $\omega_1$-like order in relation to the property of being Luzin-Jones. We may ask also ask how much of the spectrum do we need to know in order to assure that it contains $\mathcal{X}_{\omega_1}$. That is:

\begin{problem}Is there a finite family $\mathcal{Y}\subseteq \mathcal{X}_{\omega_1}$ so that, for any Luzin-Jones family $\mathcal{A}$, if $\mathcal{Y}\subseteq Spec(\mathcal{A})$, then $\mathcal{X}_{\omega_1}\subseteq Spec(\mathcal{A})$?
\end{problem} 

Now we move to problems about gaps in a more general context. It is a theorem of Rothberger that $\mathfrak{b}$ is the minimal cardinal $\kappa$ for which $(\omega,\kappa)$-gaps exist. In fact, when we restrict to cardinals, the only gaps that we can construct in $ZFC$ alone are the ones of type $(\omega_1,\omega_1)$, $(\omega,\mathfrak{b})$ and $(\mathfrak{b},\omega)$. Indeed, if we assume $PFA$ then these are the only posibilities (see \cite{StevoIlias}). However, as we have shown in Theorem \ref{allgapsthm}, the situation becomes much more complex when we consider $(X,Y)$-gaps where $X$ and $Y$ have cofinality $\omega_1$. We wonder if the same level of complexity carries to $(X,Y)$-gaps when $|X|=\omega$ and the cofinality of $Y$ is $\mathfrak{b}$. Before stating our problem, let us consider the natural generalisation of orders similar to $\omega_1$.  
\begin{definition}
Let $\kappa$ be a cardinal and $(X,<)$ be a partial order. Let us call $(X,<)$ a \textit{$\kappa$-like order} whenever $X$ is well-founded, $|X|=\kappa$, and $|(-\infty,x)|<\kappa$ for each $x\in X$. 
\end{definition}
\noindent
We then ask:

\begin{problem}Is there is an $(\omega,X)$-gap for any $\mathfrak{b}$-like order $X$? Does $PFA$ decides this question?    
\end{problem}
\noindent
It can been seen that the answer to this question is \say{yes} under $CH.$\\

The study of the Gap Cohomology group of $*$-lower semi-lattices was originally motivated due to the work of Talayco regarding towers. We have already mentioned that Todor\v{c}evi\'c showed that $|\mathcal{G}(\mathcal{T})|=2^{\omega_1}$ for any $\omega_1$-tower $\mathcal{T}$\footnote{Recall that our towers do not need to be maximal.}. Thus, it is natural to ask whether the same result holds for any $*$-lower semi-lattice. A partial answer to this question lies in almost disjoint families known as strong-$Q$-sequences (also known as uniformizable $AD$-systems). For more about these objects see \cite{Strongqdavid},  \cite{properqsequenceshelah} and \cite{Strongqjuris}.

\begin{definition}[Strong-Q-sequence] Let $\mathcal{A}$ be an almost disjoint family. We say that $\mathcal{A}$ is a \textit{strong-$Q$-sequence} if for each sequence $\langle f_A\rangle_{A\in \mathcal{A}}$ with $f:A\longrightarrow 2$ for each $A\in \mathcal{A}$, there is $f:\omega\longrightarrow 2$ such that: $$f|_A=^*f_A$$ for each $A\in \mathcal{A}.$
\end{definition}

It is easy so see that the existence of a strong-$Q$-sequence of size $\omega_1$ (in $\mathscr{P}(\omega)/\text{FIN}$) is equivalent to the existence of a $*$-Lower semi-lattice $\mathcal{T}$ isomorphic to $[\omega_1]^{<\omega}$ for which the gap cohomology group of $\mathcal{T}$, namely $\mathcal{G}(\mathcal{T})$, is trivial. Thus, Todor\v{c}evi\'c result about the size of the gap cohomology group of towers cannot be extended further without assuming additional axioms. The following problems arise:

\begin{problem}Is there (in $ZFC$) an $\omega_1$-like $*$-lower semi-lattice $\mathcal{T}$ such that $|\mathcal{G}(\mathcal{T})|=\omega$?
    
\end{problem}

\begin{problem}Is it consistent that  $\mathcal{G}(\mathcal{T})$ is non-trivial for any $\omega_1$-like $*$-lower semi-lattice $\mathcal{T}$?
    
\end{problem}
\begin{problem}Let $X$ be an $\omega_1$-like lower semi-lattice and suppose that $\mathcal{G}(\mathcal{T})$ is non-trivial for any $*$-lower semi-lattice $\mathcal{T}$ isomorphic to $X$. Is it true that $|\mathcal{G}(\mathcal{T})|=2^{\omega_1}$ for any such $\mathcal{T}$?
\end{problem}

As a direct corollary of Theorem \ref{independentcoherentschemefca}, we have that $FCA$ implies the existence of a family of size $\mathfrak{c}$ of independent gaps. Note that the main result of Yorioka in \cite{yoriokadestructiblegaps} under $\Diamond$ is stronger than ours. Thus, we ask:

\begin{problem}Does $FCA$ imply the existence of a family of size $2^{\omega_1}$ of independent gaps? What about $FCA(part)$? 
\end{problem} 

\section{Problems about trees and lines}

An uncountable linear order $(X,<)$ is said to be \textit{minimal} if for each uncountable $Y\subseteq X$, there is a strictly increasing function from $X$ to $Y$. It is a well-known theorem that $\mathfrak{m}>\omega_1$ implies that every Countryman line is minimal (see \cite{decomposingsquares}, \cite{justinomega1minimal}, \cite{ordertypesuncountablebaumgartner} and \cite{mapsontrees}). Is this statement true when we use any of the parameterized Martin's axioms $\mathfrak{m}^n_\mathcal{F}$ instead of $\mathfrak{m}$? 
\begin{problem}\label{problemmfallminimal}Let $2\leq n\in \omega$ and $\mathcal{F}$ be an $n$-capturing construction scheme. Suppose that $\mathfrak{m}^n_\mathcal{F}>\omega_1$. Is it true that any Countryman line is minimal?
    
\end{problem}

In this thesis, we showed that given a construction scheme $\mathcal{F}$, there is a natural way to define an order $<_\mathcal{F}$ over $\omega_1$ such that $(\omega_1,<_\mathcal{F})$ is a Countryman line (see Definition \ref{definitioncountrymanscheme}). So even if Problem \ref{problemmfallminimal} has a negative answer, it is interesting to know whether $(\omega_1,<_\mathcal{F})$ is minimal or not.
\begin{problem}Let $2\leq n\in \omega$ and $\mathcal{F}$ be an $n$-capturing construction scheme. Suppose that $\mathfrak{m}^n_\mathcal{F}>\omega_1$. Is it true that $(\omega_1,<_\mathcal{F})$ is a minimal Countryman line?
\end{problem}
One of the main applications of the cardinal $\mathfrak{m}_\mathcal{F}$ was to show that under $\mathfrak{m}_\mathcal{F}>\omega_1$, there are no Suslin trees. Baumgartner proved something stronger by using the cardinal $\mathfrak{m}$. He proved that under $\mathfrak{m}>\omega_1$, every Aronszajn tree is special. Thus, a natural question is the following:

\begin{problem}Let $\mathcal{F}$ be a $2$-capturing construction scheme. Does the inequality $\mathfrak{m}_\mathcal{F}>\omega_1$ implies that every Aronszajn tree is special?
\end{problem}
In some sense, the non-existence of Suslin trees under $\mathfrak{m}_\mathcal{F}>\omega_1$ says that the property of being $2$-capturing is not strong enough to imply the existence of a Suslin tree. On the other hand $FCA$ does imply the existence of a Suslin tree. So how much capturing is it really needed for this result? 

\begin{problem}Let $3\leq n\in\omega$ and $\mathcal{F}$ be a $3$-capturing construction scheme. Does the inequality $\mathfrak{m}^n_\mathcal{F}>\omega_1$ implies Suslin's hypothesis?
\end{problem}
\noindent
The same question applies to entangled sets.

\begin{problem}Does the existence of $2$-entangled sets follow from either $CA_n$ or $CA?$ In other words, which is the weakest capturing axiom implying the existence of $2$-entangled sets?
\end{problem}

\section{More problems regarding applications}
In Theorem \ref{policromaticscheme} we proved that, from $CA_3$, we can construct a coloring $c:[\omega_1]^2\longrightarrow 2$ which is $2$-bounded, has no uncountable injective sets, but is $ccc$ destructible with respect to this last property. 
\begin{problem} Does the existence of a coloring such as the one in Theorem \ref{policromaticscheme} follow from $CA_2$? Given a construction scheme $\mathcal{F}$, does $\mathfrak{m}_\mathcal{F}>\omega_1$ imply the non-existence of such a coloring?
\end{problem}
\noindent
Regarding Tukey types, we propose the following question.
\begin{problem}How many distinct Tukey types on $\omega_1$ exist under $CA_2$?
\end{problem}
\noindent
Particularly, we can ask:
\begin{problem} Suppose that $\mathcal{F}$ and $\mathcal{F}'$ are distinct  $2$-capturing construction schemes. Find sufficient conditions under which $(\mathcal{B}_\mathcal{F},\leq )$ is incomparable with $(\mathcal{B}_{\mathcal{F}'},\leq)$.
\end{problem}

Under Martin's Axiom, a great variety of ultrafilters can be constructed. Is the same true for $\mathfrak{m}_\mathcal{F}$?

\begin{problem}Let  $2\leq n\in\omega$ and $\mathcal{F}$ be an $n$-capturing construction scheme. Does $\mathfrak{m}^n_\mathcal{F}>\omega_1$ imply the existence of union ultrafilters?
    
\end{problem}

Let $S\subseteq \omega_1$ be stationary. Recall that a \textit{ladder system} on $S$ is a sequence $\langle L_\alpha\rangle_{\alpha\in S}$ such that $L_\alpha$ is an unbounded subset of $\alpha$ of order type $\omega$ for each $\alpha\in S.$
\begin{definition}[Uniformization of ladders] A ladder system $\langle L_\alpha\rangle_{\alpha\in S}$ is said to be \textit{uniformizable} if of each sequence $\langle f_\alpha\rangle_{\alpha\in S}$ of functions $f_\alpha:L_\alpha\longrightarrow \omega$, there is $F: \omega_1\longrightarrow \omega$ such that: $$f_\alpha=^*f$$
for any $\alpha\in S$.
\end{definition}
It is known that $\mathfrak{m}>\omega_1$ implies that every ladder system is uniformizable. On the other hand, $2^\omega<2^{\omega_1}$ implies that no ladder system is uniformizable (see \cite{devlinshelahweakdiamond}). Throughtout the course of the investigation of construction schemes, we tried to construct ladder systems with interesting properties. In particular, we tried to construct a ladder system which is no uniformizable by using $FCA$. Nevertheless we didn't succeed in that task. So:

\begin{problem}Does the existence of a non-uniformizable ladder system follow from $FCA$? What about $FCA(part)$?
\end{problem}
\begin{problem}Let $2\leq n\in\omega$ and $\mathcal{F}$ be an $n$-capturing construction scheme. Does $\mathfrak{m}^n_\mathcal{F}>\omega_1$ imply that every ladder system is uniformizable?
\end{problem}
To learn more about uniformizability of ladder systems see \cite{pauluniformization} and \cite{Cesaruniformization}.

\section{Problems about construction schemes}
 Recently, in \cite{clubadassaf}, Assaf Rinot and Roy Shalev introduced the principle $\clubsuit_{AD}$. This principle follows both from $\clubsuit$ and from the existence of a Suslin Tree. Thus, it also follows from $FCA$. In this way, we may ask:

\begin{problem}Is any instance of $\clubsuit_{AD}$ a consequence of $CA_n$ for some $2\leq n\in\omega?$ What about $CA_n(part)?$
\end{problem}

One of the most famous problems outside set theory, in which set theory played an important role, was the solution of Naimark's problem regarding $C^*$-algebras. In \cite{akemannweaver}, Charles Akemann and Nik Weaver constructed a counterexample to Naimark's problem using $\Diamond$-principle. In \cite{Diamondcohen}, Daniel Calder\'on and Ilijas Farah isolated the guessing principle  $\Diamond^{Cohen}$, which follows from $\Diamond$. There, they constructed a counterexample to Naimark's problem using  $\Diamond^{Cohen}$.
\begin{problem} Is $CA_n$, $CA$ or $FCA$ a consequence of $\Diamond^{Cohen}$? What about $CA_n(part)$, $CA(part)$ and $FCA(part)$?
\end{problem}
\begin{problem}Does the existence of a counterexample to Naimark's problem follow from $FCA(part)$?
    
\end{problem}

In \cite{ParametrizedDiamonds}, Mirna D\v{z}amonja, Michael Hru\v{s}\'ak and Justin Moore introduced parametrized diamonds, which are weakenings of $\Diamond$ that may hold in models where $CH$ fails. Both capturing construction schemes and parametrized diamonds are useful tool for building interesting objects of size $\omega_1$. We ask:

\begin{problem}Is there any relation between $FCA$ (respectively $CA_n$) and any parametrized $\Diamond$-principle?
\end{problem}
It would be specially interesting to know if there is any realtionship between $FCA$ and $\Diamond(non(\mathcal{M})).$\\

We already know that the forcing $\mathbb{C}_\kappa$ forces $FCA(part)$ whenever $\kappa$ is an uncountable ordinal. But what happens when we only add one Cohen real?
\begin{problem}Does the Cohen forcing $\mathbb{C}$ force $FCA$? What about $FCA(part)$?
\end{problem}

By far, the most interesting problem from the point of view of the author is the following:
\begin{problem}Extend the theory of capturing construction schemes to higher cardinals.
    
\end{problem}